\documentclass{amsart}
\usepackage{amsthm,amsmath,amssymb,amscd,graphics,psfig,epic,eepic}

{\end{list}}
{
   \newtheorem{theorem}{Theorem}[subsection]                     
   \newtheorem{proposition}[theorem]{Proposition}     
   \newtheorem{lemma}[theorem]{Lemma}

   \newtheorem{conjecture}[theorem]{Conjecture}

}
{\theoremstyle{definition}
   
   \newtheorem{example}[theorem]{Example}
   \newtheorem{definition}[theorem]{Definition}
}
{\theoremstyle{remark}
   \newtheorem*{remark}{Remark}
   \newtheorem*{remarks}{Remarks}
}

\newcommand{\ssigma}{{{\sigma_s}}}
\newcommand{\RR}{{\mathbb{R}}}

\newcommand{\QQ}{{\mathbb{Q}}}

\newcommand{\ZZ}{{\mathbb{Z}}}

\newcommand{\cI}{{\mathcal I}}
\newcommand{\cJ}{{\mathcal J}}

\newcommand{\cO}{{\mathcal O}}

\newcommand{\cU}{{\mathcal U}}

\newcommand{\prim}{{\operatorname{prim}}}
\newcommand{\supp}{{\operatorname{supp}}}
\newcommand{\reg}{{\operatorname{reg}}}
\newcommand{\can}{{\operatorname{\rm can}}}
\newcommand{\val}{{\operatorname{val}}}
\newcommand{\ce}{{\operatorname{Z}}}
\newcommand{\card}{{\operatorname{card}}}

\newcommand{\bl}{{\operatorname{bl}}}

\newcommand{\stab}{{\operatorname{stab}}}
\newcommand{\strat}{{\operatorname{strat}}}
\newcommand{\Stab}{{\operatorname{Stab}}}
\newcommand{\Ver}{{\operatorname{Vert}}}
\newcommand{\ort}{{\operatorname{ord}}}
\newcommand{\Inv}{{\operatorname{Inv}}}
\newcommand{\Tan}{{\operatorname{Tan}}}
\newcommand{\inte}{{\operatorname{int}}}
\newcommand{\lin}{\operatorname{lin}}
\newcommand{\sing}{\operatorname{sing}}
\newcommand{\Cl}{\operatorname{Cl}}
\newcommand{\Ga}{\operatorname{Gal}}
\newcommand{\Gl}{\operatorname{Gl}}
\newcommand{\Spec}{\operatorname{Spec}}
\newcommand{\Proj}{\operatorname{Proj}}

\newcommand{\Sing}{\operatorname{Sing}}
\newcommand{\Star}{\operatorname{Star}}

\newcommand{\Pic}{{\operatorname{Pic}}}

\newcommand{\Aut}{{\operatorname{Aut}}}

\newcommand{\semic}{{\operatorname{semic}}}

\newcommand{\ord}{{\operatorname{ord}}}

\newcommand{\Mid}{{\operatorname{Mid}}}
\newcommand{\Ctr}{{\operatorname{Ctr}}}
\newcommand{\Card}{{\operatorname{Card}}}
\newcommand{\inv}{{\operatorname{inv}}}
\newcommand{\ini}{{\operatorname{in}}}
\newcommand{\id}{{\operatorname{id}}}

\newcommand{\Dep}{{\operatorname{Dep}}}
\newcommand{\para}{{\operatorname{par}}}
\newcommand{\Ind}{{\operatorname{Ind}}}

\newcommand{\cW}{\cite{Wlodarczyk2}}

\begin{document}


\title[Toroidal varieties and the weak Factorization Theorem]{Toroidal varieties and the weak Factorization Theorem}
\author{Jaros{\l}aw W{\l}odarczyk}
\thanks{The author was supported in part by  Polish KBN
grant 2 P03 A 005 16 and NSF grant DMS-0100598.}
\address{Department of Mathematics\\Purdue University\\West
Lafayette, IN-47907\\USA}

\email{wlodar@math.purdue.edu, jwlodar@mimuw.edu.pl}
\date{\today}
\begin{abstract}
We develop the theory of stratified toroidal varieties, which
gives, together with the theory of birational cobordisms \cW, 
 a proof of the weak
factorization conjecture for birational maps in characteristic zero: 
a birational map between complete
nonsingular varieties over an algebraically closed field $K$ of characteristic
zero is a composite of blow-ups and blow-downs with smooth centers.  
\end{abstract}
\maketitle

\tableofcontents
\addtocounter{section}{-1}

\section{Introduction}

The main goal of the present paper is two-fold. First we extend
the theory of toroidal embeddings introduced by Kempf,
Knudsen, Mumford, and Saint-Donat to the class of toroidal
varieties with stratifications.  Second we give a proof of the following weak factorization conjecture as an application and illustration of the theory.
\begin{conjecture} {\bf The Weak Factorization conjecture}
 \begin{enumerate}
\item  Lef $f:X{{-}{\to}} Y$ be a 
birational map  of
smooth complete varieties 
over an algebraically closed field of characteristic zero, which is an
isomorphism over an open set $U$. Then $f$ can be factored
 as
$$X=X_0\buildrel f_0 \over {{-}{\to}}  X_1
\buildrel f_1 \over {{-}{\to}} \ldots \buildrel f_{n-1} \over
{{-}{\to}} X_n=Y  ,$$
where each $X_i$ is a smooth complete variety and $f_i$ is a blow-up
or blow-down at a smooth center which is an isomorphism
over $U$.
\item Moreover, if $X\setminus U$ and $Y\setminus U$ are divisors
with normal crossings, then each $D_i:=X_i\setminus U$ 
is a divisor with normal crossings and $f_i$ is a blow-up
or blow-down at a smooth center which has normal
crossings with components of $D_i$.
\end{enumerate}
\end{conjecture}

This theorem extends a theorem of Zariski, which states that
any birational map between two smooth complete
surfaces can be factored into a succession of blow-ups at
points followed by a succession of blow-downs at points.  
A stronger version of the above theorem, called the strong factorization 
conjecture remains open. 

\begin{conjecture} {\bf Strong factorization conjecture}.
Any birational map $f:X\mathrel{{-}\,{\to}}Y$ of  smooth
complete varieties can be factored into a succession of blow-ups
at
smooth centers followed by a succession of 
blow-downs at
smooth centers.
\end{conjecture}

One can find the statements of both conjectures in many papers.  Hironaka \cite{Hironaka1} formulated the strong
factorization conjecture. The weak factorization
problem was stated by Miyake and Oda \cite{Miyake-Oda}. The toric versions of the
strong and weak factorizations were also conjectured by
Miyake and Oda \cite{Miyake-Oda} and are called the strong and weak Oda
conjectures. The $3$-dimensional toric version of the weak form was
established by Danilov \cite{Danilov2} (see also Ewald \cite{Ewald}). 
The weak  toric conjecture 
in arbitrary dimensions was proved in   \cite{Wlodarczyk1} and later
independently by Morelli \cite{Morelli1},  who proposed  a proof of the
strong factorization conjecture 
(see also Morelli\cite{Morelli2}).
 Morelli's
proof of the weak Oda conjecture was completed, revised and
generalized to the toroidal case by Abramovich, Matsuki and Rashid
in \cite{Abramovich-Matsuki-Rashid}.  A gap in Morelli's proof of the strong Oda conjecture,
 which was
not noticed in \cite{Abramovich-Matsuki-Rashid}, was later
found by K. Karu.

The local version of the strong factorization problem was posed
by Abhyankar in dimension 2 and by Christensen  in general, 
who has solved it for 3-dimensional toric varieties \cite{Christensen}.
The local version of the weak factorization problem (in
characteristic 0) was 
solved by Cutkosky in \cite{Cutkosky1}, who also 
showed that Oda's strong conjecture implies the local
version of the strong conjecture for proper\ birational
morphisms \cite{Cutkosky2} and proved the local strong factorization conjecture
in dimension 3 (\cite{Cutkosky2}) via Christensen's theorem.

Shortly after the present proof was found, another proof 
of the weak factorization conjecture was conceived  
 by
Abramovich, Karu, Matsuki and the author \cite{AKMW}. Unlike the proof in \cite{AKMW} the present proof does not refer to the weak factorization theorem  
 for toric varieties.

Another 
application of the theory of stratified toroidal varieties is
given in Section \ref{se: resolution}, where we show 
the existence of a 
 resolution of
singularities of toroidal varieties in arbitrary characteristic by 
blowing up  ideals determined by
valuations (see Theorem \ref{th: resolution}).  

All schemes and varieties in the present paper are
considered over an algebraically closed field $K$.
The assumption of characteristic 0 is needed  for the
results in Sections \ref{se: construction} and \ref{se: factorization} only, where we use Hironaka's canonical 
resolution of singularities and canonical principalization for the second part of the theorem (see Hironaka \cite{Hironaka1},
Villamayor \cite{Villamayor} and Bierstone-Milman \cite{Bierstone-Milman}).

\section{Main ideas}\label{se: main ideas}

\subsection{Toroidal embeddings} \label{se: toroidal embeddings} 
The theory of toroidal embeddings was introduced and developed
by Kempf, Knudsen, Mumford and Saint-Donat in \cite{KKMS}. 
Although it was originally conceived for the purpose of
compactifying symmetric spaces and certain moduli spaces it has been successfully applied to 
 various problems concerning resolution of singularities
and factorization of morphisms (see
\cite{KKMS}, \cite{Abramovich-de-Jong}, \cite{Abramovich-Karu}). It is an 
efficient tool in many situations, turning  complicated or
even "hopeless" algebro-geometric problems into relatively "easy"
combinatorial problems. Its significance lies in
a simple combinatorial description of
the varieties considered carrying rich  and precise information about
singularities and  stratifications.

By {\it a
toroidal embedding} (originally  a
toroidal embedding without self-intersections) we mean a variety $X$ with an open subset
$U$ with the following property: for any $x\in X$ there is an open
neighborhood $U_x$ and an \'etale morphism
$\phi : U_x\to X_{\sigma_x}$ into a toric variety $X_{\sigma_x}$
containing a torus $T$ such that $U_x\cap U=\phi^{-1}(T)$.
Such a morphism is called a {\it chart}.

Let $D_i, i\in I$, be the  irreducible
components of the complement divisor $D=X\setminus U$. We
define a stratification $S$ on $X$
with strata $s\in S$ such that either $\overline{s}$ is an
irreducible component of $\bigcap_{i\in
J} D_i$ for some $J\subset I$ or $\overline{s}=X$. The strata are naturally ordered
by the following relation: $s\leq s'$ iff the closure $\overline{s}$
contains $s'$. 
In particular
$s=\overline{s}\setminus\bigcup_{s'>s}\overline{s'}$. Any
chart $\phi : U_x\to X_{\sigma_x}$ preserves strata: as a
matter of fact, all strata $s\subset U_x$ for $s\in S$ are
preimages of $T$-orbits on $X_{\sigma_s}$. In fact a toroidal
embedding can be defined as a stratified variety such that
for any point $x\in X$ there exists a neighborhood $U_x$
of $x$ and an \'etale morphism $\phi : U_x\to X_{\sigma_x}$
preserving strata.

The cone ${\sigma_x}$ in the above definition can be
described using Cartier divisors. Let $x$ belong to a
stratum $s\in S$. 
  We associate with  $s\in S$
the following data (\cite{KKMS}):

${\bf M}^s$: the group of Cartier divisors in $
\Star(s,S):=\bigcup_{s'\leq s}{s'}$, supported in $\Star(s,S)\setminus U$,

${\bf N}_s:={\rm Hom}({\bf M}_s,\ZZ)$,

${\bf M}_+^s\subset {\bf M}^s$: effective Cartier divisors,

$\sigma_s\subset {\bf N}_\RR^s$: the dual of ${\bf M}_+^s$. 

This gives a correspondence
between strata $s\in S$ and cones.  If $s'\leq s$ then
$\sigma_{s'}$ is a face of $\sigma_{s}$. 
By glueing cones $\sigma_s$ along their subfaces we
construct a {\it
conical complex} $\Sigma$. In contrast to in the standard theory of
complexes two faces might share a few common faces, not
just one. 
For any stratum $s\in S$ the integral vectors in the relative interior of 
$\sigma_s$ define monomial valuations on $X_{\sigma_s}$.
These valuations induce, via
the morphisms $\phi$ monomial 
valuations centered at $\overline{s}$ (independent of charts). 

The natural class of
morphisms of toroidal embeddings are called
{\it toroidal morphisms} and are exactly those which are locally
determined by charts $\phi$ and toric morphisms: 

 A birational  morphism of toroidal
embeddings $f: (Y,U) \to (X, U)$ is {\it  toroidal} if
for any $x\in s\subset X$ there exists a chart $x\in U_x\to
X_{\Delta^\sigma}$  and a subdivision
$\Delta^\sigma$ of $\sigma_s$ and a 
 fiber square   of morphisms of stratified varieties 

\[\begin{array}{rcccccccc}
&& & (U_x,U_x\cap U) &  \longrightarrow & (X_{\sigma_s},T)&&&\\

&&&\uparrow f & & \uparrow &&& \\

U_x \times_{X_{\sigma_s}}
X_{\Delta^\sigma} &&\simeq& (f^{-1}(U_x),f^{-1}(U_x)\cap U)
 & \rightarrow &
(X_{\Delta^\sigma} , T)&&&

\end{array}\]

Note that
toroidal morphisms are well defined and do not depend upon
the  choice of the charts. 
This fact can be described nicely using
the following Hironaka condition:

For any  points $ x, y$ which are in the same stratum    every
isomorphism $ \alpha : \widehat{X}_x\to 
\widehat{X}_y $ preserving stratification 
  can be lifted
to an isomorphism $\alpha':Y\times_X\widehat{X}_x
\to Y\times_X\widehat{X}_y$ preserving stratification.

Birational toroidal morphisms with fixed target are in bijective correspondence with subdivisions of
the associated complex.

\subsection{Stratified toroidal varieties}
\label{sse: stratified toroidal}
A more general class of objects occurring naturally in
applications is the class of stratified toroidal varieties
defined as follows.

A {\it stratified toric variety} is a toric variety with an 
invariant equisingular stratification. The important
difference between toric varieties and 
stratified toric varieties is that the latter come with a
stratification which may be coarser than the one given by
orbits.   As a consequence, the
combinatorial object associated to a stratified toric
variety, called a {\it semifan}, consists of those faces of
the fan of the toric variety 
corresponding to strata. The faces  which do not
correspond to strata  are
ignored (see Definition \ref{de: embedded}). 

By a {\it stratified toroidal variety} we mean a stratified
variety $(X,S)$ such that for any $x\in s$, $s\in S$ there
is an \'etale map called a {\it chart} 
$\phi_x:U_x\to X_{\sigma}$ 
from an open neighborhood to a stratified toric variety such that all strata in $U_x$ are preimages of
strata in $X_{\sigma}$ (see Definition 
\ref{de: stratified toroidal}). A toroidal embedding is a particular
example of a stratified toroidal variety if we consider the
stratification described in section \ref{se: toroidal embeddings}  
(see also Definition \ref{de: toroidal embeddings}).

Each chart associates with a
stratum $s$ the {\it semicone} $\sigma$ i.e the  semifan consisting of
faces of $\sigma$
corresponding to strata in ${\rm Star(s,S)}$ for the
corresponding stratum $s$. 

If one chart associates to a stratum a semicone $\sigma_1$, and
another associates $\sigma_2$ then  Demushkin's  theorem   says that
$\sigma_1\simeq\sigma_2$.

This gives us the
correspondence between  strata and cones with semifans. 
Consequently, if $s\leq s'$, then 
$\sigma$ is a face of $\sigma'$ and the semicone
$\sigma$ is a subset of the semicone
$\sigma'$.  By glueing
the semicones $\sigma$ and $\sigma'$ along the common faces
$\sigma''$
 we obtain the so called {\it semicomplex} $\Sigma$
associated to the stratified toroidal variety. The
semicomplex is determined uniquely up to isomorphism.
Its  faces of
the semicomplex are indexed by the stratification. If a
stratified toroidal variety is a toroidal embedding then
the associated semicomplex is the usual conical complex. 

 Note that in contrast in the case of toroidal
embeddings only some integral vectors $v$ in  
faces of this semicomplex determine unique valuations on $X$ which
do not depend upon a chart. These vectors are called {\it stable}
and the valuations they induce are also called {\it stable}
(see Definitions \ref{de: stable valuations}, \ref{de: stable vectors}). The stable vectors  form a subset in 
the support of the semicomplex,
which we call the {\it stable support} (see Definition \ref{de:
stable support}). 

In particular, let $X$ be a smooth surface and $S$ be its
stratification consisting of a point $p\in X$ and its complement. The associated
semicomplex consists of the regular cone ${\bf Q}_{\geq 0}\cdot
e_1+{\bf Q}_{\geq 0}\cdot e_2$ corresponding to the point $p$ and of the
vertex corresponding to a big stratum. The only stable
valuation is the valuation of the point $p$  which corresponds to the
unique (up to proportionality) stable vector $e_1+e_2$.

The notion of stable support plays a key role in the theory of stratified
toroidal varieties. The stable support has many nice properties. 
It intersects  the relative interiors of all faces of
the semicomplex and this intersection is polyhedral. The relative interiors of faces which are not
in the semicomplex  do  not intersect the stable support.
Consequently, all stable valuations
are centered at the closures of strata. The part of the stable support
which is contained in a face of the semicomplex is
convex and is a union of relative interiors of a finite
number of polyhedral cones.  Also some ''small'' vectors
like those which are in the parallelogram determined by ray
generators belong to the stable support. If we consider the 
canonical resolution of toric singularities then all divisors
occurring in this resolution determine stable valuations.
In section \ref{se: comparison}, where we discuss another
approach to the weak factorization conjecture in \cite{AKMW}, we also show a method of
computing the  stable support which is based on the idea of
torification of Abramovich and  de Jong \cite{Abramovich-de-Jong}. 

The stability condition for vectors in a face $\sigma$ of
the semicomplex $\Sigma$ can be expressed in
terms of invariance under the action of the group $G^{\sigma}$ of
automorphisms  of a scheme $\widehat{X}_\sigma$ 
associated to this face (see Definition \ref{de: stable
valuations})  (see
section \ref{se: -semicomplexes}).  
The vector $v$ is stable if the valuation determined by this
vector on $\widetilde{X}_\sigma$ is $G^{\sigma}$-invariant. 
(Definition \ref{de: stable vectors}).

Birational toroidal mophisms or simply {\it toroidal
morphisms}  between stratified toroidal varieties
are those induced locally by toric morphisms via toric
charts and satisfying the Hironaka condition (see Definition 
\ref{de: canonical stratification}).
The Hironaka condition implies that toroidal morphisms, which
are a priori induced by toric charts, in fact do not depend upon the charts
and depend eventually only on the subdivisions of the
semicomplex associated to the original stratified toroidal variety. 
                
This leads to a  $1-1$ correspondence between toroidal
morphisms and certain subdivisions, 
called {\it canonical},  of the semicomplex associated to the
original variety  
(see Theorems \ref{th: modifications} and \ref{th: correspondence}).

The canonical subdivisions are
 those  satisfying the following 
condition: all "new" rays of the subdivision are in the stable
support  (see Proposition
\ref{pr: simple}). It is not difficult to see the necessity
of the condition  since the "new" rays correspond to
invariant divisors which determine invariant valuations.
 It is more complicated to show the sufficiency of the
condition.

Let $(Y,R)\to (X,S)$ be a toroidal morphism of stratified
toroidal varieties associated to the canonical subdivision
$\Delta$ of the semicomplex $\Sigma$ associated to
$(X,S)$. Then the semicomplex $\Sigma_R$ associated to
$(Y,R)$ consists of all faces in $\Delta$ whose relative
interiors intersect the stable support of $\Sigma$ (see
Theorem \ref{th: correspondence}).

Toroidal morphisms between stratified toroidal varieties
generalize those defined for toroidal embeddings. In that
case, all subdivisions are canonical and we get the
correspondence between toroidal embeddings and subdivisions
of associated semicomplexes mentioned in Section 
\ref{se: toroidal embeddings} above. 

 In our considerations we
require some condition, called orientability, of compatibility of charts on
stratified toroidal varieties. This condition is 
analogous to the orientability of charts on an oriented
differentiable manifold. It says that two \'etale charts 
$\phi: U\to X_\ssigma$ and $\phi': U\to X_\ssigma$
define an automorphism
$\widehat{\phi}'\circ\widehat{\phi}^{-1}$ in 
the identity component of the group of  automorphisms of 
$\Spec(\widehat{\cO}_{X_{\sigma}})$
 (see Definition \ref{de:
orientation}). In the differential setting, the group of
local automorphisms consists of two components. In the
algebraic situation the number of components depends upon 
the singularities. It is finite for toric singularities.
 Smooth toroidal varieties are oriented, since
the group of automorphisms of the completion of a  regular local
ring is connected.

Let $X$ be a toroidal variety with isolated
singularity $x_1\cdot x_2=x_3\cdot x_4$  or equivalently
after some change of coordinates $y_1^2+y_1^2+
y_3^2+y_4^2=0$. In this case the group of automorphisms of
the 
formal completion of $X$ at the singular point consists of two 
components. Consequently, 
there are two orientations on $X$ (for more details see Example
\ref{ex: isolated}).

Toroidal embeddings are oriented stratified toroidal
varieties (see Example \ref{ex: toroidal embeddings}). 
Stratified toroidal varieties are orientable 
(see Section \ref{se: existence}).

\subsection{Toroidal varieties with torus action}
Let $\Gamma$ denote any algebraic subgroup of an algebraic
torus $T:=K^*\times\ldots\times K^*$.
If $\Gamma$ acts on a stratified toric variety via a group
homomorphism $\Gamma \to T$, then we require 
 additionally,  that the stabilizers at closed
points in a fixed  stratum $s$ are all the same and equal to a group
$\Gamma_s$ associated to $s$, and there is a
$\Gamma_s$-equivariant isomorphism  of the local rings at these points.

A {\it $\Gamma$-stratified toroidal variety}
is a stratified toroidal variety with a $\Gamma$-action which
is locally $\Gamma$-equivariantly \'etale isomorphic to a 
$\Gamma$-stratified toric variety. We associate with it a
$\Gamma$-semicomplex which is the semicomplex
$\Sigma$
with
associated groups $\Gamma_\sigma\subset \Gamma$ acting on $X_{\sigma}$.

$\Gamma$-semicomplexes determine local projections 
$\pi:\sigma\to\sigma^\Gamma$ defined by good 
quotients $X_{\sigma_s}\to
X_{\sigma_s^\Gamma}:=X_{\sigma}/\Gamma_\sigma$. These local
projections are coherent in the sense that they commute
with face restrictions and subdivisions.

All the definitions in this category are adapted from the
theory of stratified toroidal varieties. We additionally require that
charts, morphisms and automorphisms occurring in the definitions
are $\Gamma_\sigma$ or $\Gamma$-equivariant.

\subsection{Birational cobordisms, Morse theory and
polyhedral cobordisms of Morelli}  The theory of birational
cobordisms, which 
was developed in \cW,
was inspired by Morelli's proof of the weak
factorization theorem for toric  varieties,
where the notion of combinatorial cobordism was introduced (see
 \cite{Morelli1}).  

By a {\it birational cobordism} between birational
varieties $X$ and $X'$ we understand a variety $B$ 
with the action of the multiplicative group $K^*$ such that 
the "lower boundary" $B_-$ of $B$ (resp. the "upper boundary" $B_+$)
is an open subset which consists of orbits with no limit at
$0$ (resp. $\infty$).
The cobordant varieties $X$ and $X'$ are isomorphic 
to the spaces of orbits (geometric quotients) $B_-/K^*$ and
$B_+/K^*$ respectively. In  differential cobordism the
action of the $1$-parameter group of diffeomorphisms $G\simeq
({\bf R},+)\simeq ({\bf R^*},\cdot)$, defined by
the gradient field of a Morse function  gives us an
analogous interpretation of Morse theory. The bottom and top
boundaries determined by a Morse function are isomorphic to
the spaces of all orbits with no limit at $-\infty$ and
$+\infty$ respectively (or $0$ and $+\infty$ 
in multiplicative notation). The critical points of the
Morse function are the fixed points of the action. ''Passing
through'' these points induces  simple birational
transformations analogous to spherical transformations in
differential geometry. If $B$ is a smooth cobordism, then
these birational transformations locally replace one weighted
projective space with a complementary weighted projective
space. These birational transformations can be nicely  described 
using the language of toric varieties and associated fans.
The semiinvariant parameters in the neighbouhood of a fixed
point provide a chart which is a locally analytic
isomorphism of the tangent
spaces with the induced linear action. Hence locally, the
cobordism is isomorphic to an affine space with a  linear
hyperbolic action. Locally we can identify $B$ with an
affine space ${\bf A}^n$, which  is a toric variety 
corresponding to a regular cone $\sigma$ in the vector space
$N^{\bf Q}$. The action of $K^*$ determines a $1$-parameter
subgroup, hence an integral vector $v\in N^{\bf Q}$. The sets $B_-$
and $B_+$
correspond to the sets $\sigma_+$ and $\sigma_-$ of all 
faces of $\sigma$ visible from above or
below (with respect to the direction of $v$), hence to the
''upper'' and ''lower'' boundaries of $\sigma$. The vector $v$
defines a projection $\pi: N^{\bf Q}\to N^{\bf Q}/{\bf Q}v$. The quotient spaces
$B_-/K^*$, $B_+/K^*$ correspond to two ''cobordant'' (in the
sense of Morelli) subdivisions
$\pi(\sigma_+)$ and $\pi(\sigma_-)$ of the cone $\pi(\sigma)$.
The problem lies in the fact that the fans $\pi(\sigma_+)$, 
$\pi(\sigma_-)$ are singular (which means that they are not spanned by a
part of an integral basis). Consequently,   the corresponding
birational transformation is a composition of weighted
blow-ups and blow-downs between singular varieties. 
The process of resolution of
singularities of the quotient spaces is called
$\pi$-desingularization and can be achieved locally by
a combinatorial algorithm. The difficulty is now in patching
these local $\pi$-desingularizations. This problem can be
solved immediately  once the theory of stratified toroidal
varieties (with $K^*$-action) is established. 

\subsection{Sketch of  proof}
Let $B=B(X,X')$ be a smooth cobordism between smooth
varieties $X$ and $X'$. Let $S$ be the stratification determined by
the isotropy groups. Then $(B,S)$ is 
a stratified toroidal variety with $K^*$-action. We can
associate with it a $K^*$-semicomplex $\Sigma$. 
The $K^*$-semicomplex $\Sigma$ determines for any
face $\ssigma\in\Sigma$ a projection
$\pi_\sigma: \sigma\to \sigma^\Gamma$ determined by the good 
quotients
$X_{\sigma}\to X_{\sigma}/\Gamma_\sigma=X_{\sigma}^\Gamma$.   
We   apply a combinatorial algorithm - the so
called $\pi$-desingularization Lemma of Morelli -  to
the semicomplex associated to the cobordism. The algorithm
consists of starf subdivisions at stable vectors. Such
subdivisions are canonical. As a result we obtain a 
semicomplex with the property that the projections of all
simplices are either regular (nonsingular) cones or cones
of smaller dimension (than the dimension of the projected cones)
(see Lemma \ref{le: pi-lemma} for details).
The semicomplex corresponds to a $\pi$-regular toroidal
cobordism, all of  whose   
open affine fixed point free subsets have smooth geometric
quotients. The existence of such a cobordism
easily implies the weak factorization theorem. (The
blow-ups, blow-downs and flips induced by elementary
cobordisms are regular  (smooth)). Since each flip is a
composition of a blow-up and a blow-down at a smooth center we
come to a factorization of the map $X{-}{\to}X'$ into
blow-ups and blow-downs at smooth centers (see Proposition
\ref{pr: regular factorization}).  
\section{Preliminaries}
\subsection{ Basic notation and terminology}

Let $N\simeq
{\bf Z}^k$ be a lattice contained in 
the vector space $N^{\bf Q}:=N\otimes {\bf
Q}\supset N$.  
 By a {\it cone} in this paper 
we  mean a convex set $\sigma = {\bf Q}_{\geq 0}\cdot
v_1+\ldots+{\bf Q}_{\geq 0}\cdot v_k\subset N^{\bf Q}$. 
By abuse of language we shall speak of a {\it cone $\sigma$ in a
lattice $N$}. To avoid confusion we 
shall sometimes write $(\sigma, N)$ for the cone $\sigma$
in $N$. For a cone $\sigma$ in $N$ denote by $N_\sigma:=N\cap\lin(\sigma)$ the
sublattice generated by $\sigma$. Then by
$$\underline{\sigma}:=(\sigma,N_\sigma)$$ 
we denote the corresponding cone in $N_\sigma$. (Sometimes the cones $\sigma$ and 
$\underline{\sigma}$ will be identified).

We call a vector $v\in N$ {\it primitive} if it generates
sublattice ${\bf Q}_{\geq 0}v_i\cap N$. 
Each ray $\rho\in N^{\bf Q}$ contains a unique primitive
ray $\prim(\rho)$. 
 If $v_1,\ldots,v_k$ form a
minimal set of primitive vectors generating $\sigma\subset N^{\bf Q}$
 in the above sense, then we write $$\sigma =\langle
v_1,\ldots,v_k\rangle .$$  If
$\sigma$ contains no line we call it {\it strictly convex}. All
cones considered in this paper are strictly convex. For any
$\sigma$ denote by ${\rm lin}(\sigma)$ the linear span of $\sigma$. For any
cones $\sigma_1$ and $\sigma_2$ in $N$ we write $$\sigma
=\sigma_1+ \sigma_2$$ if $\sigma=\{v_1+v_2\mid  v_1\in
\sigma_1 ,v_2\in
\sigma_2\}$, and 
$$\sigma =\sigma_1\oplus \sigma_2$$ if
${\rm lin}(\sigma_1)\cap {\rm lin}(\sigma_2)=\{0\}$ and for any $v\in
\sigma\cap N$ there exist $v_1\in \sigma_1\cap N$ and 
$v_2\in \sigma_2\cap N$ such that $v=v_1+v_2$.
For any
cones $\sigma_1$ in $N_1$ and $\sigma_2$ in $N_2$ we
define the cone $\sigma_1\times\sigma_2$ in $N_1\times N_2$
to be
$$\sigma_1\times\sigma_2:=\{(v_1,v_2)\mid v_i\in\sigma_i\quad \mbox{for} \quad i=1,2 \}.$$
 We say that a cone
$\sigma$ in $N$ is {\it regular} if there exist vectors
$e_1,..,e_k\in N$ such that
$$\sigma =\langle e_1\rangle \oplus \ldots\oplus \langle e_k\rangle. $$  
A cone $\sigma$ is {\it simplicial} if 
$\sigma =\langle v_1,\ldots,v_k\rangle $ is generated by linearly
independent vectors. 
We call $\sigma$
{\it indecomposable} if it cannot be 
represented as $\sigma =\sigma'\oplus \langle e \rangle $ for
some nonzero vector $e\in N$.
\begin{lemma}
Any cone $\sigma$ in $N$ is uniquely
represented  as $$\sigma ={\rm sing}(\sigma )\oplus \langle e_1,\ldots,e_k\rangle,
 $$
 where $\langle e_1,\ldots,e_k\rangle $ is a regular cone and
${\rm sing}(\sigma )$ is the maximal indecomposable face of
$\sigma$.\qed  
\end{lemma}

By $\inte(\sigma)$ we denote the relative interior of  $\sigma$.

For any simplicial cone
$\sigma =\langle v_1,\ldots,v_k\rangle $ in $N$ set $${\rm par}(\sigma ):=\{ v\in
\sigma\cap N_{\sigma}\mid  v=\alpha_1v_1+\ldots+\alpha_kv_k,
\mbox{where}\,\, 0\leq\alpha_i< 1\},$$  
$$\overline{{\rm par}(\sigma )}:=\{ v\in
\sigma\cap N_{\sigma}\mid  v=\alpha_1v_1+\ldots+\alpha_kv_k,
\mbox{where}\,\, 0\leq\alpha_i\leq 1\}.$$
For any simplicial cone $\sigma =\langle v_1,\ldots,v_k\rangle $ in $N$, by
$\det(\sigma) $ we mean $\det(v_1,\ldots,v_k)$, where all
vectors are considered in some
 basis of $N\cap {\rm lin}(\sigma)$ (a change
 of basis can only change the sign of the determinant).

\begin{definition} \label{de: minimal}
A {\it minimal
generator} of a cone $\sigma$ is a vector not contained in a
one dimensional face of $\sigma$ and which cannot be
represented as the sum of two nonzero integral vectors in
$\sigma$. A {\it minimal
internal vector} of a cone $\sigma$ is a vector in ${\rm
int}(\sigma)$ which cannot be
represented as the sum of two nonzero integral vectors in
$\sigma$  such that at least one of them belongs to ${\rm
int}(\sigma)$. 
\end{definition}
Immediately from the definition we get

\begin{lemma} \label{le: minimal0} Any minimal generator $v$ of $\sigma$ 
is a minimal internal vector of the face
$\sigma_v$ of $\sigma$, containing $v$ in its relative interior.\qed
\end{lemma}

\begin{lemma} For any simplicial $\sigma$ each vector
from ${\rm par}(\sigma)$ can be represented as a nonnegative
integral combination of minimal generators.
\end{lemma}
 \qed

\subsection{Toric varieties}

\begin{definition} (see \cite{Danilov1},
\cite{Oda}).\label{de: fan} By a {\it fan}
$\Sigma $ in
$N$ we mean a finite collection of finitely 
generated strictly convex cones $\sigma$ in $N$ such 
that 

$\bullet$ any face of a cone in $\Sigma $ belongs to $\Sigma$,

$\bullet$ any two cones of $\Sigma $ intersect in a common face. 

By the {\it support} of the fan we mean the union of all
its faces, 
$|\Sigma|=\bigcup_{\sigma\in \Sigma}\sigma$.

If $\sigma$ is a face of $\sigma'$ we shall write $\sigma\preceq\sigma'$.

If $\sigma\preceq\sigma'$ but  $\sigma\neq\sigma'$ 
we shall write $\sigma\prec\sigma'$.

For any set of cones $\Sigma$ in $N$ by $\overline{\Sigma}$
we denote the set $\{\tau \mid \tau\prec \tau'\quad \mbox{for some}\quad\tau'\in\Sigma\}$
\end{definition}

\begin{definition}\label{de: star} Let $\Sigma$ be a fan
and $\tau \in \Sigma$. The {\it star} of the
cone $\tau$ and the {\it closed star} of $\Sigma$ are
defined as follows:
 $${\rm Star}(\tau ,\Sigma):=\{\sigma \in \Sigma\mid 
\tau\preceq \sigma\},$$ 
$$\overline{{\rm Star}}(\tau ,\Sigma):=\{\sigma \in
\Sigma\mid  \sigma'\preceq \sigma \, \mbox{for some} \, \sigma'\in
{\rm Star}(\tau ,\Sigma)\},$$ 
\end{definition} 

\begin{definition}\label{de: product} Let $\Sigma_i$ be a
fan in $N_i$ for $i=1,2$. Then the {\it product} of
$\Sigma_1$ and $\Sigma_2$ is a fan $\Sigma_1\times\Sigma_2$
in $N_1^{\bf Q}\times N_2^{\bf Q}$ defined as follows:
$$\Sigma_1\times\Sigma_2:=\{\sigma_1\times 
\sigma_2 \mid \sigma_1\in\Sigma_1, \sigma_2\in\Sigma_2\}.$$
\end{definition}
To a fan $\Sigma $ there is associated a toric variety
$X_{\Sigma}\supset {\bf }T$, i.e. a normal variety on which a torus
$T$ acts effectively with an open dense orbit 
(see \cite{KKMS}, \cite{Danilov2}, \cite{Oda}, \cite{Fulton}). To 
each cone $\sigma\in \Sigma$
corresponds an  open affine invariant subset
$X_{\sigma}$ and its unique closed orbit $O_{\sigma}$. The
orbits in the
closure of the  orbit $O_\sigma$ correspond to the cones of 
${\rm Star}(\sigma ,\Sigma)$.

Denote by $$M:={\rm Hom}_{alg.gr.}(T,K^*)$$ \noindent the lattice of
group homomorphisms to $K^*$, i.e.  characters of $T$. Then
the dual lattice $N=
{\rm Hom}_{alg.gr.}(K^*,T)$ can be 
identified with the lattice of $1$-parameter subgroups of $T$. Then
vector space $M^{\bf Q}:=M\otimes{\bf Q}$ is dual to $N^{\bf Q}:=
N\otimes{\bf Q}$. Let $(v,w)$ denote the relevant pairing for $v\in N,w\in M$ 

For any $\sigma\subset N^{\bf Q}$ we denote by
$$\sigma^\vee:=\{m\in M \mid (v,m) \geq 0 \,\, {\rm
for\,\, 
any} \,\,\,    m\in \sigma\}$$ 
\noindent the set of integral vectors of the dual cone to $\sigma$. Then the
ring of regular functions $K[X_\sigma]$ is $K[\sigma^\vee]$.  

Each vector $v\in N$ defines a linear function on $M$ which
determines a 
valuation ${\rm val}(v)$ on
$X_{\Sigma}$.

 For any regular function $f=\sum_{w\in M} a_wx^w\in K[T]$ set
$$\val(v)(f):=\min\{(v,w)\mid a_w\neq 0\}.$$ 
Thus $N$
can be perceived as the lattice of all $T$-invariant
integral valuations of the function field of $X_{\Sigma}$.

For any $\sigma\subset N^{\bf Q}$ set
$$\sigma^\perp:=\{m\in M \mid (v,m) = 0 \,\, {\rm
for\,\, 
any} \,\,\,    m\in \sigma\}.$$ 
\noindent The latter set represents all  invertible characters
on $X_\sigma$. All noninvertible characters are $0$ on $O_\sigma$.
The ring of regular functions on $O_\sigma\subset X_\sigma$ can be
written as $K[O_\sigma]=K[\sigma^\perp]\subset
K[\sigma^\vee]=K[\sigma^\perp]\otimes K[\underline{\sigma}]$. 
Thus
$K[X_\sigma]=K[O_\sigma][\underline{\sigma}^\vee]$

Let $T_\sigma\subset T$ be the subtorus corresponding the
sublattice $N_\sigma:=\lin(\sigma)\cap N$ of $N$. Then by
definition, $T_\sigma$ acts trivially on
$K[O_\sigma]=K[\sigma^\perp]=K[X_\sigma]^{T_\sigma}$. Thus
$T_\sigma=\{t\in T\mid tx=x, x\in O_\sigma\}$.

This leads us to the lemma
\begin{lemma} \label{le: und} Any toric variety
$X_\sigma$ is isomorphic to $X_{\underline{\sigma}}\times
O_\sigma$, where $O_\sigma\simeq T/T_\sigma$. 
\qed
\end{lemma}
For any $\sigma\in\Sigma$ 
the closure $\overline{O}_\sigma$ of the orbit 
$O_\sigma\subset X_\Sigma$ is a toric variety 
 with the
big torus $T/T_\sigma$. Let $\pi: N\to N/N_\sigma$ be the natural projection.
Then $\overline{O}_\sigma$ corresponds to the fan $\Sigma':=\{
\frac{\tau+\lin(\sigma)}{\lin(\sigma)}\mid \tau \in {\rm Star}(\tau ,\Sigma)\}$
in $N/N_\sigma$. 
\bigskip
\subsection {Morphisms of toric varieties}

\begin{definition}(see \cite{KKMS}, \cite{Oda},
\cite{Danilov2}, \cite{Fulton}). A
{\it birational toric morphism} or simply a {\it toric morphism} of toric
varieties $X_\Sigma \to X_{\Sigma'}$ is a 
morphism identical on $T\subset X_\Sigma, X_{\Sigma'}$ 
\end{definition}

\begin{definition} (see \cite{KKMS}, \cite{Oda},
\cite{Danilov2}, \cite{Fulton}). 
A {\it subdivision} of a fan
$\Sigma$ is a fan $\Delta$ such that $|\Delta|=|\Sigma|$
and any cone $\sigma\in
\Sigma $ is a union of cones $\delta\in
\Delta$. 
\end{definition}

\begin{definition}\label{de: star subdivision} Let $\Sigma$ be a fan and
$\varrho$ be a ray passing in the
relative interior of $\tau\in\Sigma$. Then the {\it star
subdivision} $\varrho\cdot\Sigma$ of $\Sigma$ with respect to
$\varrho$ is defined to be

$$\varrho\cdot\Sigma=(\Sigma\setminus {\rm Star}(\tau ,\Sigma) )\cup
\{\varrho+\sigma\mid   \sigma\in \overline{\rm Star}(\tau
,\Sigma)\setminus {\rm Star}(\tau
,\Sigma)\}.$$ If $\Sigma$ is regular, i.e. all its cones are
regular, $\tau=\langle v_1,\ldots,v_l\rangle $ and 
$\varrho=\langle v_1+\ldots+v_l\rangle $
then we call the star
subdivision $\varrho\cdot\Sigma$ {\it regular}. 
\end{definition}

\begin{proposition} (see \cite{KKMS}, \cite{Danilov2},
\cite{Oda}, \cite{Fulton}). Let
$X_\Sigma$ be a toric variety. There is a 1-1 correspondence
between subdivisions of the fan $\Sigma$ and proper toric
morphisms $X_{\Sigma'} \to X_{\Sigma}$.\qed
\end{proposition}

\begin{remark} Regular star subdivisions from
\ref{de: star subdivision} correspond to blow-ups of smooth varieties
at closures of orbits (\cite{Oda}, \cite{Fulton}). Arbitrary
star subdivisions correspond to blow-ups of some ideals
associated to valuations  (see Lemma \ref{le: blow-up valuation}).
\end{remark}

\subsection{Toric varieties with $\Gamma$-action }

In further considerations let $\Gamma$ denote any algebraic
subgroup of an algebraic torus  $T=K^*\times\ldots\times K^*$.
By $X/\Gamma$ and $X//\Gamma$ we denote respectively geometric and  good
quotients. If $\Gamma$ acts on an algebraic variety $X$ and
$x\in X$ is a closed point, then  $\Gamma_x$ denotes  the
isotropy group of $x$.

\begin{definition} \label{de: G-indecomposable} Let 
$\Gamma$ act on $X_{\sigma}\supset T$ via
a group homomomorphism  $\Gamma\to T$. 
Set $$\Gamma_\sigma:=\{g\in\Gamma\mid g(x)=x \quad \mbox{
for\,\, any} \quad  x \in O_\sigma\}
.$$ We shall write $$\sigma=\sigma'\oplus^\Gamma\langle
e_1,\ldots,e_k \rangle$$ \noindent if 
$\sigma=\sigma'\oplus\langle
e_1,\ldots,e_k \rangle$, $\langle
e_1,\ldots,e_k \rangle$ is regular and $\Gamma_{\sigma}=\Gamma_{\sigma'}$.  
We say that $\sigma$ is 
{\it $\Gamma$-indecomposable} if it cannot be  
represented as $\sigma=\sigma'\oplus^\Gamma\langle
e_1,\ldots,e_k \rangle$ (or equivalently $X_{\sigma}$ is not
of the form $X_{\sigma'}\times {\bf A}^k$ with $\Gamma_\sigma$ acting
trivially on ${\bf A}^k$). By 
$${\rm sing}^{\Gamma}(\sigma)$$ we mean the maximal 
$\Gamma$-indecomposable face of $\sigma$.
\end{definition}
If $\Gamma$ acts on a toric vatiety via a
group homomomorphism $\Gamma\to T$ we shall speak of a
{\it toric action} of $\Gamma$.
\subsection{Demushkin's Theorem}
For any algebraic variety $X$ and its (in general
nonclosed) point $x\in X$ we
denote by
$\widehat{\cO}_{X,x}$  the completion of the local
ring ${\cO}_{X,x}$ at the maximal ideal of $x$. We also set
$$\widehat{X}_x:={\rm Spec}(\widehat{\cO}_{X,x}).$$
For the affine toric variety $X_\sigma$ define 
$$\widehat{X}_\sigma:= {\rm Spec}(\widehat{\cO}_{X_\sigma,O_\sigma}).$$

 We shall use the following Theorem of
Demushkin (\cite{Demushkin}):

\begin{theorem}\label{th: Dem0} Let $\sigma$ and
$\tau$ be
two cones of maximal dimension in isomorphic lattices
$N_\sigma\simeq N_\tau$.

Then the following conditions are equivalent: 
\begin{enumerate}

\item $\sigma \simeq \tau$.

\item $\widehat{X}_{\sigma} \simeq \widehat{X}_{\tau}$. 
\end{enumerate}
\end{theorem}

\noindent{\bf Proof.} For the proof see \cite{Demushkin} or
the proof of
\ref{le: Dem2}. \qed

The above theorem can be formulated for affine toric
varieties with $\Gamma$-action.

\begin{theorem}\label{th: Dem00} Let $\sigma$ and
$\tau$ be
two cones of maximal dimension in isomorphic lattices
$N_\sigma\simeq N_\tau$.
Let $\Gamma$ act on $X_{\sigma}\supset{
T}_\sigma$ and on  $X_{\tau}\supset{
T}_\tau$ via group homomorpisms $\Gamma\to T_\sigma$ and
$\Gamma\to T_\tau$. 
 
Then the following conditions are equivalent:
\begin{enumerate}

\item There exists an isomorphism of cones $\sigma \simeq
\tau$ inducing a $\Gamma$-equivariant isomorphism of toric
varieties.

\item There exists a $\Gamma$-equivariant isomorphism 
$\widehat{X}_{\sigma} \simeq \widehat{X}_{\tau}$. 
\end{enumerate}
\end{theorem}
\noindent{\bf Proof.} For the proof see the proof of
\ref{le: Dem2} or
\cite{Demushkin} . \qed

\begin{definition} Let $\Gamma$ act on $X_{\sigma}\supset{
T}_\sigma$ and on  $X_{\tau}\supset{
T}_\tau$ via  group homomorpisms $\Gamma\to T_\sigma$ and
$\Gamma\to T_\tau$. 
 We say that cones $\sigma$ and
$ \tau$ are {\it $\Gamma$-isomorphic} if there exists
an isomorphism of cones $\sigma \simeq
\tau$ inducing a $\Gamma$-equivariant isomorphism of toric
varieties.
\end{definition}

The above theorem can be formulated as follows:

\begin{theorem} \label{th: Dem} 
Let $\sigma$ and
$\tau$ be
two cones  in isomorphic lattices
$N_\sigma\simeq N_\tau$.
Let $\Gamma$ act on $X_{\sigma}\supset{
T}_\sigma$ and on  $X_{\tau}\supset{
T}_\tau$ via group homomorpisms $\Gamma\to T_\sigma$ and
$\Gamma\to T_\tau$. Assume that $\Gamma_{\sigma}=\Gamma_{\tau}$.

Then the following conditions are equivalent:
\begin{enumerate}
\item  ${\rm sing}^\Gamma(\sigma)$ and ${\rm
sing}^\Gamma(\tau)$ are $\Gamma_\sigma$-isomorphic. 

\item
 For any closed points $x_\sigma\in O_{\sigma}$ and $x_\tau\in
O_{\tau}$, there exists a $\Gamma_\sigma$-equivariant 
isomorphism of the local rings ${\cO}_{X_{\sigma},x_\sigma}$
and ${\cO}_{X_{\tau},x_\tau}$. 

\item
 For any closed points $x_\sigma\in O_{\sigma}$ and $x_\tau\in
O_{\tau}$, there exists a $\Gamma_\sigma$-equivariant 
isomorphism of the completions of the local rings $\widehat
{\cO}_{X_{\sigma},x_\sigma}$ and $\widehat{\cO}_{X_{\tau},x_\tau}$ .

\end{enumerate}
\end{theorem}

\noindent{\bf Proof.} $(1)\Rightarrow(2)$ Write $\sigma={\rm
sing}^\Gamma(\sigma)\oplus^\Gamma {\rm r}(\sigma)$ and $\tau={\rm
sing}^\Gamma(\tau)\oplus^\Gamma {\rm r}(\tau)$, where
${\rm r}(\sigma)$, ${\rm r}(\tau)$ denote regular cones. By Lemma
\ref{le: und},  
$X_\sigma=X_{\underline{\sigma}}\times O_\sigma=X_{\underline{\rm
sing}^\Gamma(\sigma)}\times X_{\underline{{\rm r}(\sigma)}}\times O_\sigma$,
where $\Gamma_\sigma$ acts trivially on $X_{\underline{{\rm r}(\sigma)}}\times O_\sigma$.
Thus there exists a $\Gamma_\sigma$-equivariant isomorphism 
${\cO}_{X_{\sigma},x_\sigma}\simeq{\cO}_{X_{\underline{\rm
sing}^\Gamma(\sigma)}\times X_{{\rm r'}(\sigma)}}
$, where ${\rm r'}(\sigma)$ is a regular cone of dimension
$\dim(N_\sigma)-\dim({\sing^\Gamma(\sigma)})$ and
$\Gamma_\sigma$ acts trivially on $X_{{\rm r'}(\sigma)}$.
Analogously ${\cO}_{X_{\tau},x_\tau}\simeq{\cO}_{X_{\underline{\rm
sing}^\Gamma(\tau)}\times X_{{\rm r'}(\tau)}}\simeq {\cO}_{X_{\underline{\rm
sing}^\Gamma(\sigma)}\times X_{{\rm r'}(\sigma)}}={\cO}_{X_{\sigma},x_\sigma}
$.

The implication $(2)\Rightarrow(3)$
is trivial. 
For the proof of $(3)\Rightarrow (1)$ we find a regular
cone ${\rm re}(\sigma)$ in $N_\sigma$ such that
$\sigma':= \sigma\oplus {\rm re}(\sigma)$
is of maximal dimension in $N_\sigma$. Since $\Gamma_\sigma$
acts trivially on $O_\sigma$ and on $O_\tau$ we have 
$\sigma':= \sigma\oplus^{\Gamma_\sigma} {\rm re}(\sigma)$.
Analogously $\tau':= \tau\oplus^{\Gamma_\tau}{\rm  re}(\tau)$.
By Theorem \ref{th: Dem00}, $\sigma'$ and $\tau'$ are
$\Gamma_\sigma$-isomorphic.  
Consequently, 
${\rm sing}^\Gamma({\sigma})=
{\rm sing}^{\Gamma_\sigma}({\sigma'})\simeq {\rm
sing}^{\Gamma_\tau} ({\tau'})={\rm
sing}^\Gamma({\tau})$. \qed

Theorem \ref{th: Dem} allows us to assign a singularity type to
any closed point $x$ of a toric variety $X$ with a
toric action of the group $\Gamma$:

\begin{definition} \label{de: singularity type}
By the {\it singularity type} of
a point $x$ of a toric variety $X$ we mean the function 
$$\sing(x):=\underline{\sing(\sigma_x)},$$
\noindent 
where $\sigma_x$ is a cone of maximal dimension such that 
$\widehat{X}_x\simeq\widehat{X}_{\sigma_x}$.

By the {\it singularity type} of
a point $x$ of a toric variety $X$ with a
toric action of group $\Gamma$ we mean
$$\sing^\Gamma(x):=(\Gamma_x,\underline{\sing^{\Gamma_x}(\sigma_x)}),$$
\noindent 
where $\sigma_x$ is a cone of maximal dimension with toric action of  
$\Gamma_x$  on $\widehat{X}_{\sigma_x}$ and such that there exists a 
$\Gamma_x$-equivariant
isomorphism
$\widehat{X}_x\simeq\widehat{X}_{\sigma_x}$.
\end{definition}

\section{Stratified toric varieties and semifans}

\subsection{Definition of a stratified toric variety}

In this section, we give a combinatorial description of stratified
toric varieties in terms of so called (embedded) semifans. 

\begin{definition}\label{de: stratification}
Let $X$ be a noetherian scheme $X$ over ${\rm Spec}(K)$. A 
{\it stratification\/} of $X$ is a decomposition of $X$ into a
finite collection $S$ of pairwise disjoint locally closed irreducible
smooth subschemes $s \subset X$, called strata, with the following
property: For every $s \in S$, the closure $\overline{s} \subset X$
is a union of strata. 
\end{definition}

\smallskip
\begin{definition}
Let $S$ and $S'$ be two stratifications of $X$.
We say that $S$ is {\it finer} than $S'$ if any stratum
in $S'$ is a union of strata in $S$. In this case we shall
also call $S'$  {\it coarser} than $S$.
\end{definition}

\noindent \begin{definition}\label{de: stratified toric} \enspace
Let $X$ be a toric variety with big torus $T \subset X$. 
A {\it toric stratification\/} of $X$ is a stratification $S$ of 
$X$ consisting of $T$-invariant strata $s$ such that for any two
closed points $x,x' \in s$ their local rings ${\cO}_{X,x}$ and
${\cO}_{X,x'}$ are isomorphic.
\end{definition}

\begin{definition}\label{de: stratified toric2}
Let $X$ be a toric variety with a toric action of $\Gamma$.
 We say that a toric stratification $S$ of $ X$ is {\it
compatible with the action of $\Gamma$} if
\begin{enumerate}

\item All points in the  same stratum $s$ have the
same isotropy group $\Gamma_s$.

\item For any two closed points  $x, y$ from one
stratum $s$ there exists a $\Gamma_s$-equivariant
isomorphism $\alpha:\widehat{X}_x\to\widehat{X}_y$,
preserving all strata.
\end{enumerate}
\end{definition}
\smallskip

If $S$ is a toric stratification of $X$, then we shall also speak of a
{\it stratified toric variety} $(X,S)$. If $X$ is a toric
variety with a toric stratification , compatible with
$\Gamma$, then we shall also speak of a
{\it $\Gamma$-stratified toric variety} $(X,S)$.

The combinatorial objects we shall use
in this context are the following:

\begin{definition} \label{de: embedded} \enspace
An {\it embedded semifan\/} is a subset $\Omega \subset \Sigma$ of a
fan $\Sigma$ in a lattice $N$ such that for every $\sigma
\in \Sigma$ there is an $\omega(\sigma) \in \Omega$ satisfying
\begin{enumerate}

\item $\omega(\sigma) \preceq \sigma$ and any other 
$\omega \in \Omega$ with $\omega \preceq \sigma$ is a face of
$\omega(\sigma)$,

\item $\sigma = \omega(\sigma) \oplus {\rm r}(\sigma)$ for some
regular cone ${\rm r}(\sigma) \in \Sigma$.
\end{enumerate}

A {\it semifan\/} in a lattice $N$ is a set $\Omega$ of cones in $N$
such that the set $\Sigma$ of all faces of the cones of $\Omega$ is a
fan in $N$ and $\Omega \subset \Sigma$ is an embedded semifan.
\end{definition}

\smallskip

\begin{definition} \label{de: G-semifan} An embedded
semifan $\Omega\subset\Sigma$ with a toric action of
$\Gamma$ on $X_\Sigma$ will be called an
{\it embedded $\Gamma$-semifan} if for every
$\sigma\in\Sigma$, $\sigma=\omega(\sigma)\oplus^\Gamma {\rm r}(\sigma)$.

\end{definition}

\begin{remark}
A semifan or an embedded semifan can be viewed as a $\Gamma$-semifan
or an embedded $\Gamma$-semifan with trivial group $\Gamma$. 
\end{remark}

Some examples are discussed at the end of this section. The main
statement of this section says that stratified toric varieties are
described by embedded semifans:

\smallskip

\begin{proposition} \label{le: semifans correspondence} \enspace
 Let $\Sigma$ be a fan in a lattice $N$, and let $X$ denote the
associated toric variety with a toric action of $\Gamma$. 
There is a canonical 1-1 correspondence
between the toric stratifications of $X$ compatible withn the  
action of $\Gamma$ and the embedded
 $\Gamma$-semifans $\Omega \subset \Sigma$:

\begin{enumerate}
\item If $S$ is a toric stratification of $X$, compatible
with the action of $\Gamma$, then the
corresponding embedded $\Gamma$-semifan $\Omega \subset \Sigma$ consists of
all those cones $\omega \in \Sigma$ that describe the big orbit of
some stratum $s \in S$.  

\item If $\Omega \subset \Sigma$ is an embedded $\Gamma$-semifan, then
the strata of the associated toric stratification $S$ of $X$ arise
from the cones of $\Omega$ via 
$$ \omega \mapsto {\rm strat}(\omega) 
:= \bigcup_{\omega(\sigma) = \omega} O_\sigma. $$
\end{enumerate}

\end{proposition}

\smallskip
\noindent{\bf Proof.} 
$(1)\Rightarrow (2)$ Since strata of $S$ are $T$-invariant and
disjoint, each
orbit $O_\tau$ belongs to a
unique stratum $s$. Let $\omega \in \Omega$ describe the
big open orbit of $s$. Then $O_\tau$ is contained in the
closure of $O_{\omega}$. Hence $\omega$ is a face of $\tau$. 
By Definition {de: embedded}, $\Gamma_\tau=\Gamma_\omega$
and there exists
a $\Gamma_\tau$-equivariant isomorphism of the  local
rings   
${\cO}_{X_\tau,x}$ and ${\cO}_{X_\omega,y}$ of two 
points $x\in O_{\tau}$ and $
y\in O_{\omega} $.
By Theorem \ref{th: Dem} and since 
$\sing^\Gamma(\tau)\supset\sing^\Gamma(\omega)$, we infer that
$\sing^\Gamma(\tau)=\sing^\Gamma(\omega)$. Hence
$\omega=\sing^\Gamma(\omega)\oplus^\Gamma {\rm r}(\omega)$ and 
$\tau=\sing^\Gamma(\tau)\oplus^\Gamma {\rm r}(\tau)=
\omega\oplus^\Gamma {\rm r'}(\tau)$,
 where ${\rm r}(\tau)={\rm r}(\omega)\oplus {\rm r'}(\tau)$.
  
$(2)\Rightarrow (1)$. By definition
$\{\sigma\in\Sigma\mid \omega(\sigma)=\omega)\}=
{\rm Star}(\omega,{\Sigma})\setminus \bigcup_{\omega\prec
\omega' \in \Omega}
{\rm Star}(\omega', {\Sigma})$.
Hence all the defined subsets
$\strat(\omega)$ are
locally closed. The closure of each subset
$\strat(\omega)$ corresponds to ${\rm
Star}(\omega,\Sigma)$ and hence it is a union of
the sets  $\strat(\omega')$, where $\omega\prec\omega'$. 
Since $\tau=\omega\oplus^\Gamma {\rm r}(\tau)$ we have $\Gamma_\tau=\Gamma_\omega$.
Each subset $\strat(\omega)$ is a toric
variety with a fan
$\Sigma':=\{\frac{\tau+\lin(\omega)}{\lin(\omega)}|\omega(\tau)=\omega\}$
in $(N^Q)':=N^Q/\lin(\omega)$.
Since
$\tau=\omega\oplus^\Gamma {\rm r}(\tau)$ the cone
$\frac{\tau+\lin(\omega)}{\lin(\omega)}$ is isomorphic to
the regular cone ${\rm r}(\tau)$ in $(N^Q)'$.
 Thus the strata
$\strat(\omega)$ is  smooth. 
The local rings of closed points of the the strata
$\strat(\omega)$ have the same isotropy group $\Gamma_\omega$.
Moreover if $O_\tau\subset\strat(\omega)$ then
$\sing^\Gamma(\tau)=\sing^\Gamma(\omega)$. 
By Theorem \ref{th: Dem}, we conclude  that there exists 
a $\Gamma$-equivariant isomorphism of the local rings of any
 two points $x\in O_{\tau}$ and $
y\in O_{\omega} $.   
\qed

As a corollary from the above we obtain the following lemmas:

\begin{lemma} \label{le: strata}
${\rm strat}(\omega)=\overline{{\rm
strat}(\omega)}\setminus\bigcup_ {\omega'\prec\omega}
\overline{ {\rm strat}(\omega')}$. \qed
\end{lemma}

\begin{lemma} \label{le: strata2} Let $S_1$ and $S_2$ be two
toric stratifications on a toric variety $X$ corresponding to
two embedded semifans $\Omega_1,\Omega_2\subset\Sigma$.
Then $S_1$ is coarser than $S_2$ iff $\Omega_1\subset \Omega_2$.
\end{lemma}
\noindent{\bf Proof} $(\Rightarrow)$ Assume that $S_1$ is coarser
than $S_2$. Let $\omega\in\Omega_1$ be the cone
corresponding to  a stratum $s_1\in S_1$.  The closure $\overline{s_1}$
 is a union of strata in $S_2$. 
There is a generic stratum $s_2\in S_2$ contained in
$\overline{s_1}$. Then $\overline{s_1}=\overline{s_2}$ and
$\omega\in\Omega_2$ corresponds to the stratum $s_2$.

$(\Leftarrow)$ Now assume that $\Omega_1\subset \Omega_2$. 
Then the closure of any stratum $s_1$ equals 
$\overline{s_1}=\overline{\omega_1}$, where
$\omega_1 \in \Omega_1$ is the cone
corresponding  to $s_1$. But then
$\overline{\omega_1}=\overline{s_2}$, where $s_2\in S_2$ is the stratum
corresponding to $\omega_1\in\Omega_2$. Thus
$\overline{s_1}=\overline{s_2}$ and it suffices to apply
Lemma \ref{le: strata}. \qed

\begin{remark}\enspace
Let $(X,S)$ be the stratified toric variety arising from a semifan
$\Omega \subset \Sigma$. The above correspondence between the cones of
$\Omega$ and the strata of $S$ is order reversing in the sense that 
$\omega \preceq \omega'$ iff ${\rm strat}(\omega')$ is contained in
the closure of ${\rm strat}(\omega)$. 
\end{remark}

\smallskip

\noindent\begin{lemma}\label{le: stratified toric}
Let $\Omega\subset\Sigma$ be an embedded $\Gamma$-semifan.
Then all  $\Gamma$-indecomposable faces
 of $\Sigma$ are in $\Omega$. 
\end{lemma}
\noindent{\bf Proof.}
For a $\Gamma$-indecomposable face $\sigma$ write
$\sigma=\omega(\sigma)\oplus^\Gamma {\rm r}(\sigma)$. This
implies $\sigma=\omega(\sigma)\in\Omega$. \qed

We conclude the section with some examples. The first two examples
show that for any fan there are a maximal and a minimal embedded
semifan.

\smallskip

\begin{example}\label{ex: orbits} \enspace

If $\Sigma$ is a fan in a lattice $N$, then $\Sigma \subset
\Sigma$ is an embedded semifan. The corresponding stratification of
the toric variety $X$ associated to $\Sigma$ is the decomposition of
$X$ into the orbits of the big torus $T \subset X$. The
orbit stratification is the finest among all toric
stratifications of $X$.

\end{example}

\smallskip

\begin{example} \label{ex: sing} \enspace

For a fan $\Sigma$ in a lattice $N$ and a toric action of
$\Gamma$ on $X_\Sigma$ let 
${\rm Sing}^\Gamma(\Sigma)$ denote the set of all maximal indecomposable 
parts ${\rm sing}^\Gamma(\sigma)$, where $\sigma \in \Sigma$.
Then ${\rm Sing}^\Gamma(\Sigma) \subset \Sigma$ is an embedded
$\Gamma$-semifan, and for any other embedded
$\Gamma$-semifan $\Omega \subset \Sigma$
one has ${\rm Sing}^\Gamma(\Sigma) \subset \Omega$. The toric stratification
corresponding to ${\rm Sing}^\Gamma(\Sigma) \subset \Sigma$ is the
stratification by singularity type on the toric variety
$X_\Sigma$ with the toric action of $\Gamma$. The toric stratification
corresponding to ${\rm Sing}(\Sigma) \subset \Sigma$ is the
coarsest among all toric stratifications of $X$. The toric stratification
corresponding to ${\rm Sing}^\Gamma(\Sigma) \subset \Sigma$
is the
coarsest among all toric stratifications of $X$ which are
compatible with the action of $\Gamma$.   
\end{example}
\bigskip

\begin{example} \label{ex: affine} Let ${\bf A}^2\supset K^*\times K^*$
be a toric variety with the  strata $s_0:=\{(0,0)\}$,
$s_1:=K^* \times \{0\}$, and $s_2:={\bf A}^2\setminus (s_0\cup
s_1)$. Then  ${\bf A}^2$ with the above stratification is a
stratified toric variety. 

We have $\sigma_{0}=\langle(1,0),(0,1)\rangle$,
$\sigma_{1}=\langle(0,1)\rangle,$ $\sigma_{2}=\{(0,0)\},$ 

$\overline{\sigma_{0}}=\{\langle(1,0),(0,1)\rangle\}$, 
$\overline{\sigma_{1}}=\{\langle(0,1)\rangle\}$,
$\overline{\sigma_{2}}=\{(0,0)\}, \langle(1,0)\rangle\}$.
\end{example}

\section{Stratified toroidal varieties} \label{se: stratified toroidal}

\subsection {Definition of a stratified toroidal variety} 

\begin{definition} (see also \cite{Danilov1}). Let $X$ be a variety
or a noetherian scheme over $K$. $X$ is called {\it
toroidal} if for any $x\in X$ there exists an open affine
neighborhood $U_x$ and an \'etale morphism $\phi_x:U_x\rightarrow
X_{\sigma_x}$ into some affine toric variety
$X_{\sigma_x}$. 

\end{definition}
\begin{definition} \label{de: gsmoth} Let $\Gamma$ act on
noetherian schemes $X$ and $Y$.  We say that  a $\Gamma$-equivariant
morphism $f: Y \to X$ is $\Gamma$-{\it smooth} (resp. $\Gamma$-{\it \'etale})
if there is a smooth  (resp. \'etale) morphism 
$f': Y' \to X'$ of varieties with trivial action of
$\Gamma$ and
$\Gamma$-equivariant  fiber square
 \[\begin{array}{rccccccc}
&X \times_{X'}Y' & \simeq &Y&\buildrel f \over \rightarrow
&X&& \\ 

&&& \downarrow & &\downarrow&& \\ 

&&&Y' & \buildrel f' \over \rightarrow  &X' \\
\end{array}\]
\end{definition}

\begin{lemma} \label{le: gsmoth} Let $Y\to X$ be a
$\Gamma$-smooth morphism. If the good quotient $X//\Gamma$ exists then
$Y//\Gamma$ exists and  $Y=Y//\Gamma\times_{X//\Gamma}X$.
\end{lemma}
\noindent{\bf Proof.}
Let $Y_0:=Y'\times_X{X//\Gamma}$. Then
$Y=Y_0\times_{X//\Gamma}X$ and $Y//\Gamma=Y_0$. \qed

\begin{definition} \label{de: G-toroidal variety} Let $X$ be a 
toroidal variety or a toroidal scheme. We say that an action
of a group 
${\Gamma}$ on $X$ is {\it toroidal} if for any
$x\in X$ there exists an open
$\Gamma$-invariant neighborhood $U$ and a
 $\Gamma$-equivariant, $\Gamma_x$-smooth, 
morphism $U\rightarrow X_\sigma$, into a toric variety
$X_\sigma$ with a toric action of 
$\Gamma\supset \Gamma_x$.
\end{definition}

 \begin{definition} \label{de: stratified toroidal}
 Let $X$ be a 
toroidal variety or a toroidal scheme. We say that a  
stratification $S$ of $X$ is {\it toroidal} if
any point $x$ in a  stratum $s \in S$ admits
 an open neighborhood $U_x$  and an \'etale morphism
$\phi_x: U_x\rightarrow
X_{\sigma}$ into a stratified toric variety such that $s\cap U_x$ equals $\phi_x^{-1}(O_{\sigma})$
and the intersections $s'\cap U$, $s'\in S$,  are precisely the inverse
images of strata in $X_\sigma$. 

\end{definition}
\begin{definition}\label{de: G-stratified toroidal}
 Let $X$ be
toroidal variety or a toroidal scheme with a toroidal
action of $\Gamma$. A toroidal stratification $S$ of $X$ is
said to be {\it compatible with the action of $\Gamma$} if 

\begin{enumerate}
\item Strata of $S$ are invariant with respect to the
action of $\Gamma$.
\item All points in one stratum $s$ have the
same isotropy group $\Gamma_s$. 
\item For any  point $x\in s$ there exists a $\Gamma$-invariant
neighborhood $U$ and a $\Gamma$-equivariant,
$\Gamma_s$-smooth morphism
 $U \rightarrow
X_\sigma$ into a $\Gamma$-stratified toric variety, such
that $s\cap U_x$ equals $\phi_x^{-1}(O_{\sigma})$
and the intersections $s'\cap U$, $s'\in S$,  are precisely the inverse
images of strata in $X_\sigma$. 

\end{enumerate}
\end{definition}

If $X$ is a  toroidal variety or a toroidal scheme with a toroidal
stratification $S$ then we shall speak of a {\it stratified
toroidal variety} (resp. a {\it stratified
toroidal scheme}). If $X$ is a toroidal variety or a
toroidal scheme with a
toroidal action of $\Gamma$ and a toroidal
stratification $S$ compatible with the action of $\Gamma$
then we shall speak of a {\it $\Gamma$-stratified
toroidal variety} (resp. {\it $\Gamma$-stratified
toroidal scheme}).

If $(X,S)$ is a $\Gamma$-stratified toroidal variety then for any
stratum $s\in S$ set $\Gamma_s:=\Gamma_x$, where $x\in s$.
Note that if $s\leq s'$ then $\Gamma_s\subseteq \Gamma_{s'}$.

\begin{remark}
 A stratified toroidal variety can be considered as a  
$\Gamma$-stratified toroidal variety with trivial action of
the trivial group $\Gamma=\{e\}$.
\end{remark}

A simple generalization of Example \ref{ex: sing} is the following lemma:

\begin{lemma} \label{le: singularity type}  Let $X$ be a toroidal
variety with a toroidal action of $\Gamma$. Let ${\rm Sing}^\Gamma(X)$ be the
stratification determined by the singularity type ${\rm
sing}^\Gamma(x)$ (see Definition \ref{de: singularity
type}) . Then ${\rm Sing}^\Gamma(X)$ is a toroidal
stratification compatible with  the action of $\Gamma$.\qed
\end{lemma}

\begin{definition}\label{de: toroidal embeddings} (see also \cite{KKMS}). 
A {\it toroidal embedding} is a
stratified toroidal variety such that any point
$p$ in a stratum $s \in S$ admits an open neighborhood
$U_p$  and an \'etale morphism
of stratified varieties $\phi_p: U_p \rightarrow
X_{\sigma}$ into a toric variety with orbit stratification
(see Example \ref{ex: orbits}). 

A {\it toroidal embedding} with $\Gamma$-action is a
$\Gamma$-stratified toroidal variety which is a toroidal 
embedding.

\end{definition}

\begin{remark} The definition of a toroidal embedding differs from the original
definition in \cite{KKMS}. It is equivalent to the definition of a 
toroidal embedding without self-intersections 
(see Section \ref{se: toroidal embeddings}).
\end{remark}

\subsection{Existence of invariant stratifications on smooth
varieties with $\Gamma$-action}

\noindent
\begin{lemma} \label{le: existence2} 
\begin{enumerate}

\item Let $X$ be a  smooth variety with $\Gamma$-action.  
Then  there
exists a toroidal stratification $\Sing^\Gamma(X)$ of $X$
which is compatible with the action of $\Gamma$.

\item Let $X$ be a  smooth variety with $\Gamma$-action and $D$
be a $\Gamma$-invariant divisor with simple normal crossings. Let
$S_D$ be the stratification determined by the components of
the divisor $D$. Set
$$\sing^{\Gamma,D}(x):=(\sing^{\Gamma}(x), D(x)),$$ \noindent
where
$D(x)$ is the set of components of $D$ passing through $x\in X$ 

Then $\sing^{\Gamma, D}$ determines a toroidal stratification
$\Sing^\Gamma(X,D)$ compatible with the action of $\Gamma$. In particular  all the
closures of strata from $\Sing^\Gamma(X,D)$
have normal crossings with components of $D$. 

\end{enumerate}
\end{lemma}

\noindent{\bf Proof.} For any $x$ we can find $\Gamma$-semiinvariant
parameters $u_1,\ldots,u_k$ describing the orbit
$\Gamma\cdot x$  in
some $\Gamma$-invarinat neighborhood $U_x$ of $x$.

Additionally in (2) we require that each component of $D$
through $x\in X$ is described by one of the parameters in
 $U_x$.  We define a smooth morphism $\phi_x:U_x\to {\bf A}^k$ by
$\phi_x(y)=(u_1(y),\ldots,u_k(y))$.  The action of $\Gamma$ on $X$
defines the action of $\Gamma$ on ${\bf A}^k$ with standard
coordinates be $\Gamma$-semiinvariant with corresponding
weights. Then $\Sing^\Gamma({\bf A}^k,D_A)$, where
$D_A=\overline{\phi(D)}$,
defines a  $\Gamma$-invariant stratification on ${\bf A}^k$,
whose strata are determined by standard coordinates.
Therefore $\Sing^\Gamma(X,D)$ is also a stratification on
$X$ which is locally a pull-back of $\Sing^\Gamma({\bf A}^k,D_A)$.

The morphism $\phi_x$ is $\Gamma$-equivariant. Denote
by $u_1,\ldots,u_k,\ldots,u_n$ all local $\Gamma_x$-semiinvariant
parameters at $x$. Let $\widetilde{\phi}_x:U_x\to {\bf A}^n$ be the \'etale
$\Gamma_x$-equivariant morphism  
defined by ${\widetilde{\phi}_x(y)}:=(u_1(y),\ldots,u_n(y))$. By
Luna's fundamental lemma (see Luna \cite{Luna}) we can  find
a $\Gamma_x$-invariant neighborhood $U'_{x}\subset U_x$ for
which  the induced quotient morphism
 ${\widetilde{\phi}_x}/\Gamma_x: U'_x//\Gamma_x \rightarrow
{\bf A}^n//\Gamma_x$ is \'etale
and $U'_x\simeq U'_x//\Gamma_x \times_{{\bf A}^n//\Gamma_x}
{{\bf A}^n}$. By Sumihiro \cite{Sumihiro} we can find a 
$\Gamma$-invariant open affine
$V\subset \Gamma\cdot U'$ such that 
 $V\simeq V//\Gamma_x \times_{V//\Gamma_x}
{{\bf A}^n}$. Let $p: {\bf A}^n\to {\bf A}^k$ denote the standard projection on the first coordinates.
Then $p\widetilde{\phi}_x: V\to {\bf A}^k$ is a
$\Gamma$-equivariant $\Gamma_x$-smooth morphism.
By definition $(X,\Sing^\Gamma(X,D))$ is a $\Gamma$-stratified toroidal variety.
\qed

\subsection{Conical semicomplexes}

Here we generalize the notion of a semifan. For this we first have to
figure out those semifans that describe affine stratified toric
varieties. In analogy to usual cones and fans we shall denote them
by small Greek letters $\sigma, \tau$, etc.:

\smallskip

{\bf Definition.}\enspace
Let $N$ be a lattice. A {\it semicone\/} in $N$ is a semifan $\sigma$
in $N$ such that the support $| \sigma |$ of $\sigma$ occurs as an
element of $\sigma$.

\smallskip

The {\it dimension\/} of a semicone is the dimension of its
support. Moreover, for an injection $\imath \colon N \to N'$ of
lattices, the {\it image\/} $\imath(\sigma)$ of a semicone $\sigma$ in
$N$ is the semicone consisting of the images of all the elements of
$\sigma$.

Note that every cone becomes a semicone by replacing it with the set
of all its faces. Moreover, every semifan is a union of maximal
semicones. Generalizing this observation, we build up semicomplexes
from semicones:

\smallskip

\begin{definition}\label{de: conical semicomplex} \enspace

Let $\Sigma$ be a finite collection of semicones $\sigma$ in
lattices $N_{\sigma}$ with $\dim(\sigma) = \dim(N_{\sigma})$. 
Moreover, suppose that there is a partial ordering ``$\le$'' on
$\Sigma$.

We call $\Sigma$ a {\it semicomplex\/} if for any pair $\tau \le
\sigma$ in $\Sigma$ there is an associated linear injection
$\imath^{\sigma}_{\tau} \colon N_{\tau} \to N_{\sigma}$ such that 
$\imath^{\sigma}_{\tau}(N_{\tau}) \subset N_{\sigma}$ is a saturated
sublattice and
\begin{enumerate}
\item $\imath^{\sigma}_{\tau} \circ
\imath^{\tau}_{\varrho}(\varrho) =
\imath^{\sigma}_{\varrho}(\varrho)$,

\item $\imath^{\sigma}_{\varrho}(|\varrho|) =
\imath^{\sigma}_{\tau}(|\tau|)$ implies $\varrho = \tau$, 

\item $\sigma = \bigcup_{\tau \le \sigma}
\imath^{\sigma}_{\tau}(|\tau|)=\bigcup_{\tau \le \sigma}
\imath^{\sigma}_{\tau}(\tau)$. 
\end{enumerate}
\end{definition}

\smallskip

As a special case of the above notion, we recover back the notion of a
(conical) complex introduced by Kempf, Knudsen, Mumford and Saint-Donat:

\begin{definition}\label{de: conical complex} \enspace
A {\it complex\/} is a semicomplex $\Sigma$ such that every $\sigma
\in \Sigma$ is a fan.
\end{definition}

\begin{definition}\label{de: conical semicomplex2} 

A semicomplex $\Sigma$ is called a {\it $\Gamma$-semicomplex} if
there is a collection of algebraic groups
$\Gamma_\sigma\subset T_\sigma$, where $\sigma\in \Sigma$, such that
\begin{enumerate}

\item Each $\sigma\in \Sigma$ is a  $\Gamma_\sigma$-semicone.

\item For any $\tau\leq \sigma$, $\Gamma_{\tau}\subset\Gamma_\sigma$ and 
the morphism $\imath^{\sigma}_{\tau}$ induces a 
$\Gamma_{\tau}$-equivariant morphism $X_{\tau}\to
X_{\sigma}$. Moreover $\Gamma_{\tau}=(\Gamma_{\sigma})_{\tau}=\{g\in\Gamma_\sigma\mid gx=x, x\in O_\tau\}$.

\end{enumerate}
\end{definition}

\begin{remarks}
\begin{enumerate}
\item Each  semicomplex can be considered as a
$\Gamma$-semicomplex with trivial groups $\Gamma_\sigma$

\item From now on we shall often identify vectors of $\sigma$
with their images under $\imath^{\sigma}_{\tau}$ if $\tau\leq
\sigma$, and simplify
the notation replacing $\imath^{\sigma}_{\tau}$ with set-theoretic inclusions.

\item    
In what follows  we shall write $\sigma\preceq
\sigma'$ means that $\sigma$ is a face of the cone
$\sigma'$, while $\sigma\leq
\sigma'$ to mean that $\sigma$ is a face of the semicone $\sigma'$.

\item
In the case when $\Sigma$ is a conical complex the
relations $\leq$ and $\preceq$ coincide.

\item Note that cones in a semicomplex intersect along a union of
common faces. 
\end{enumerate}
\end{remarks}

Denote by $\Aut(\sigma)$  the group of automorphisms
of the $\Gamma_\sigma$-semicone $\sigma$.

\begin{lemma} \label{le: s2} For any $\Gamma$-semicomplex $\Sigma$ and any
$\varrho\leq \tau\leq \sigma$ in $\Sigma$, there is an automorphism \\
$\alpha_{\varrho}\in\Aut(\varrho)$ such that 
$\imath^{\sigma}_{\tau} \circ
\imath^{\tau}_{\varrho} =\imath^{\sigma}_{\varrho} \alpha_\varrho.$
\end{lemma}
\noindent{\bf Proof.} By definition,  the maps $\imath^{\sigma}_{\tau} \circ
\imath^{\tau}_{\varrho} $ and
$\imath^{\sigma}_{\varrho}$ are both linear isomorphisms of $\varrho$
onto the image
$\imath^{\sigma}_{\tau} \circ
\imath^{\tau}_{\varrho}(\varrho)= 
\imath^{\sigma}_{\varrho}(\varrho)$. Thus $\alpha_{\varrho}:=
{(\imath^{\sigma}_{\varrho})}^{-1}\imath^{\sigma}_{\tau} \circ
\imath^{\tau}_{\varrho}$.
\qed

\begin{remarks} \begin{enumerate}

\item It follows  from Lemma \ref{le: s2} that  vectors in 
$\sigma$ are in general defined up
to automorphisms from $\Aut(\sigma)$. Hence  the notion of support of a 
semicomplex as a
topological space which is the  totality of such vectors is
not well defined.
However if we consider, for instance,  vectors invariant
under  all such  automorphisms then the relevant topological
space can be constructed (see notion of {\it stable
support} \ref{de: stable support}) and plays the key role
for semicomplexes and their subdivisions.

\item In further considerations we shall introduce the
notion of {\it oriented semicomplex} (see Definition
\ref{de: oriented semicomplex}) which allows a
smaller group of automorphisms $\Aut(\sigma)^0
\subset \Aut(\sigma)$. Consequently, the corresponding stable
support of an oriented semicomplex is essentially larger. This
allows one to perform more subdivisions and carry out certain
important birational transformations.   

 \end{enumerate}
\end{remarks}
 
\noindent \begin{lemma}\label{le: semic} 
A $\Gamma$-semifan $\Sigma$ in $N$  determines
a $\Gamma$-semicomplex $\Sigma^{\rm semic}$ where each cone
$\sigma\in\Sigma$ determines the semicone consisting of all faces
$\tau\leq\sigma$ and $\Gamma_\tau=\{g\in \Gamma\mid gx=x, x\in O_\tau\}$.\qed
\end{lemma}

\begin{definition}\label{de: semicomplexes isomorphism}
By an {\it isomorphism} of two $\Gamma$-semicomplexes $\Sigma\to
\Sigma'$ we mean a face bijection
$\Sigma\ni\sigma\mapsto \sigma'\in\Sigma'$ such that
$\Gamma_\sigma=\Gamma_{\sigma'}$,  along with a collection of
face $\Gamma_\sigma$-isomorphisms $j_\sigma:\sigma\to\sigma'$, such that
for any $\tau\leq \sigma$, there is a $\Gamma_\tau$-automorphism 
$\beta_\tau$ of $\tau$ such that $j_{\sigma}
\imath^\sigma_\tau=\imath^{\sigma'}_{\tau'}j_{\tau}\beta_\tau$.
 
\end{definition}

\bigskip
\subsection{Inverse systems  of affine algebraic groups}

\begin{definition}\label{de: affine proalgebraic group}    
By an {\it affine proalgebraic group} we mean an affine
group scheme that is the limit of
 an inverse system $(G_i)_{i\in
 {\bf N}}$ of affine algebraic groups and 
algebraic group homomorphisms $\phi_{ij}:G_i\rightarrow
G_j$, for $j\geq i$.
\end{definition}

\begin{lemma} \label{le: epimorphisms} Consider the natural morphism
$\phi_i:G\to G_i$. Then $H_i:=\phi_i(G)$ is an algebraic subgroup
of $G_i$, all induced morphisms $H_j\to H_i$ for $i\leq j$
are epimorphisms and 
$G=\lim_{\leftarrow}G_i=\lim_{\leftarrow}H_i$. In
particular $K[H_i]\subset K[H_{i+1}]$ and $K[G]=\bigcup K[H_i]$.
\end{lemma}

\noindent {\bf Proof.}  The set $H'_i:=\bigcap_{j>i}
\phi_{ji}(G_j)$ is an intersection of algebraic subgroups
of $G_i$. Hence it is an algebraic subgroup of $G_i$. Note
that the induced homomorphisms $\phi^H_{ij}:H'_j\to H'_i$ are
epimorphisms and
that $G=\lim_{\leftarrow}G_i=\lim_{\leftarrow}H'_i$. Consequently,
$\phi_i(G)=H_i=H'_i$. \qed

\begin{lemma}\label{le: K-points}
 The set $G^K$ of $K$-rational points  of $G$  is an
abstract group which is the inverse limit
$G^K=\lim_{\leftarrow}G^K_i$ in the category of abstract groups.
\end{lemma}

\noindent {\bf Proof.} By the previous lemma we can assume
all morphisms $G_i\to G_j$ to be epimorphisms.
Any $K$-rational point $x$ in $G$ is mapped to $K$-rational points $x_i$ in
$G_i$. This gives an  abstract group homomorphism 
$\phi: G^K\to\lim_{\leftarrow}G^K_i$. 
We also have $K[G]=\bigcup K[G_i]$. We have to show that any 
point $x$ of $\lim_{\leftarrow}G^K_i$ determines a unique
point $y$
in $G^K$. The  point $x$ determines a  sequence of
maximal ideals $m_i\in K[G_i]$ such that $m_i\subset 
m_{i+1}$. Let $m:= \bigcup m_i$. Then $m\cap K[G_i]=m_i$.
Let $f\in K[G]$. Then $f\in K[G_i]$ for some $i$ and $f-k
\in m_i$ for some $k\in K$. Hence $f\in k+m$, which means
that $K[G]/m=K$. The ideal $m$ defines a  $K$-rational point $y$ on
$G^K$ which is mapped to $x$. Note that the point $y$ is
unique since there is a  unique ideal $m$ for which $m\cap K[G_i]=m_i$
. \qed

By abuse of notation we shall identify $G$ with $G^K$. 

\bigskip
\subsection{Group of automophisms of the completion of
a local ring}

Let $X$ be a stratified noetherian scheme over $\Spec(K)$.
Let $S$ denote the stratification  and $x\in X$ be
a closed $K$-rational point. Let $\widehat{\cO}_{X,x}$
denote the completion of the local ring of $x$ on
the  scheme $X$. Let $u_1,\ldots u_k\in
\widehat{\cO}_{X,x}$ generate the maximal ideal 
$m_{x,X}\subset\widehat{\cO}_{X,x}$. Then $\widehat{\cO}_{X,x}=K[[u_1,\ldots,u_k]]/I$, 
where $I$ is the defining ideal. Set $\widehat{X}_x:={\rm
Spec}(\,\widehat{\cO}_{X,x})$ and $X^{(n)}_x:=
{\rm Spec}({\cO}_{X,x}/m_{x,X}^{ n+1})$ for $n\in{\bf N}$.

Let ${\rm Aut}(\widehat{X}_x, S)$ (respectively ${\rm
Aut}(X^{(n)}_x, S)$) denote the automorphism group of 
$\widehat{X}_x$ (resp. $X^{(n)}_x$) preserving all closures
$\overline{s}\ni x$ of strata  $s\in S$. (resp. preserving all subschemes 
\\${\rm Spec}({\cO}_{X,x}/(I_{\overline{s}}+m_{x,X}^{n+1}))$,  for all $s\in
S$)

If  $\Gamma$ acts on $X$ denote by ${\rm
Aut}^\Gamma(\widehat{X}_x,S)$ (respectively ${\rm
Aut}^\Gamma(X^{(n)}_x, S)$)  the group of all
$\Gamma_x$-equivariant atomorphisms of $(\widehat{X}_x,S)$
(resp. $(X^{(n)}_x, S)$).

In further considerations set
$G:={\rm Aut}(\widehat{X}_x,S)$ (resp. $G:={\rm
Aut}^\Gamma(\widehat{X}_x,S))$ in the case of the action of $\Gamma$, and $G_n:={\rm Aut}(X^{(n)}_x,S)$ (resp. $G_n:={\rm
Aut}^\Gamma(X^{(n)}_x,S)$).
\begin{lemma} \label{le: group action} 

\begin{enumerate}

\item $G_n$ is an algebraic group
acting on $X^{(n)}_x$. That is there exists a  co-action
$$\Phi_n^*: {\cO}_{X,x}/m_{x,X}^{n+1}\to 
K[G_n]\otimes
{\cO}_{X,x}/m_{x,X}^{n+1}$$ defining
the action morphism
$\Phi_n: G_n\times X^{(n)}\to X^{(n)}$\\ 
and an action automorphism
$\Psi_n:= 
\pi_n\times\Phi_n: 
G_n\times X^{(n)}\to G_n\times X^{(n)}$\\ 
where $\pi_n$ is the projection on $G_n$.

\item The $G_n$ form an inverse system
of algebraic groups and $G=\lim_{\leftarrow} G_n$ is a proalgebraic group.

\item The co-actions in (1) define a ring homomorphism
$ \Phi^*:\widehat{\cO}_{X,x} \to 
K[G][[\widehat{\cO}_{x,X}]]$\\
and a \underline{formal action morphism}
$ \Phi:
G\widehat{\times}\widehat{X}_x:=\Spec K[G][[\widehat{\cO}_{X,x}]]\to
\widehat{X}_x$\\
such that for any $g\in G^K$ the restriction
$
\Phi_{|\{g\}\widehat{\times}\widehat{X}_x}: 
\{g\}\widehat{\times}\widehat{X}_x\to \widehat{X}_x$ is
given by the action of $g$.
\item The automorphisms $\Psi_n$ of $X^{(n)}$
define a \underline{ formal action automorphism}
$ \Psi$ of \\ 
\centerline{$G\widehat{\times}\widehat{X}_x=\lim_{\leftarrow}G\times X^{(n)}$}
such that for any $g\in G^K$ the restriction
$
\Psi_{|\{g\}\widehat{\times}\widehat{X}_x}: 
\{g\}\widehat{\times}\widehat{X}_x\to
\{g\}\widehat{\times}\widehat{X}_x$ is the automorphism
given by the action of $g$.

\item Let $u_1,\ldots,u_k$ denote the local
parameters on $\widehat{X}_x$. There exist regular functions 
$g_{\alpha\beta}$ generating $K[G]$ such
that the action is given by $u^{\alpha}\mapsto \sum g_{\alpha\beta}u^\beta$, where $\alpha,\beta\in\ZZ_{\geq 0}^n$ are multiindices.
\end{enumerate}

\end{lemma}

\noindent{\bf Proof.} (1) Write 
$K[X^{(n)}_x]=K[u_1,\ldots,u_k]/(I+m_{x}^{n+1})$. 
Set $W^n:=K[u_1,\ldots,u_k]/(m_{x}^{n+1})$,\\ 
$I^n:=(I+m_{x}^{n+1})/m_{x}^{n+1}$, $I_s^n:=(I_s+m_{x}^{n+1})/m_{x}^{n+1}$. \\
Define the product $W^n\times W^n \to W^n$ to be the bilinear
map for which  
$u^{\alpha_1}\cdot u^{\alpha_2}=u^{\alpha_1+\alpha_2}$. Set \\ 
$$H_n:=\{g\in{\rm Gl}(W^n)\mid
g(u^{\alpha_1}\cdot u^{\alpha_2})=g(u^{\alpha_1})\cdot g(u^{\alpha_2})\}.$$

$H_n$ is an algebraic subgroup of ${\rm Gl}(W^n)$. 
Then $$G_n=\{g\in H^n\mid g(I^n)=I^n, g(I_s^n)=I_s^n, 
\quad \mbox{ for all}\quad 
 s\in S, ga=ag \quad \mbox{ for all}\quad 
 a\in \Gamma_x\}$$ \noindent is
algebraic. The embedding $G_n\subset
{\rm Gl}(W^n)$ determines the natural co-action  $K[W^n]\to 
K[G_n]\otimes K[W_n]$,  $$u^\alpha\mapsto
\sum_{\beta \in {\bf Z}_{\geq 0}^n}g_{\alpha\beta}u^\beta,$$ which factors to
$K[X^{(n)}_x]\to K[G_n]\otimes K[X^{(n)}_x]$.

(2) If  the morphism $X^{(n)}_x\hookrightarrow X^{(n+1)}_x$ commutes with 
automorphisms ,i.e. if $\Phi$ is an automorphism of
$X^{(n+1)}_x$,  then $\Phi_{|X^{(n)}}:X^{(n)}\to X^{(n)}\subset
X^{(n+1)}_x$ is an automorphism of $X^{(n)}$. This defines
the morphisms  
${\rm Aut}(X^{n+1}_x,S)\to 
{\rm Aut}(X^{(n)}_x,S)$.

The $K$-rational points of the proalgebraic group $G=\lim_{\leftarrow} 
G_n$ can be identified, by Lemma
\ref{le: K-points}, with ($\Gamma_x$-equivariant)
automorphisms of $\widehat{X}_x$ preserving strata.

(3) The morphisms $\widehat{\cO}_{X,x} \to 
K[G_n][[{\cO}_{x,X}]]/m^n_{x,X}\to
K[G][[\widehat{\cO}_{X,x}]]/m^n_{x,X}$ determine \\
\centerline{$\widehat{\cO}_{X,x}\to 
\lim_{\leftarrow} K[G][[\widehat{\cO}_{X,x}]]/m^n_{x,X}$}

(4) and (5) follow from (1) \qed

\begin{lemma}\label{le: group ideal} Let $G=\lim_{\leftarrow}G_i$ be 
a connected
proalgebraic group acting on $\widehat{X}_x$. Let
$\Phi:G\widehat{\times}\widehat{X}_x\to\widehat{X}_x$ be
the action morphism  and
$\pi:G\widehat{\times}\widehat{X}_x\to\widehat{X}_x$ be the
standard projection. Then for any
ideal $I$ on $\widehat{X}_x$ the following conditions are
equivalent: 
\begin{itemize}
\item $I$ is invariant with respect to the action of the
abstract group  $G^K$ 
\item $\Phi^*(I)=\pi^*(I)$
\end{itemize}
\end{lemma}

 \noindent {\bf Proof.} It suffices to prove the equivalence
of the conditions
for any scheme $X^{(n)}_x$ and ideal $I_n=I\cdot K[X^{(n)}_x]$. In this
situation the assertion  reduces to the well known case of an
algebraic action of $G_n$ on $X^{(n)}_x$. If $\Phi^*(I_n)=\pi^*(I)$ then for any $g\in
G_n^K$ and $f\in I_n $, we have  $g\cdot
f=\Phi^*_g(f)\in\pi_g^*(I_n)=I_n$. Thus $I_n$ is $G_n^K$-invariant.

Now suppose $I_n$ is $G_n^K$-invariant. Let  $\Phi^*(f)=\sum
f_ig_i$ for $f\in I_n$. We can assume that $g_i\in K[G_n]$
are linearly independent. Hence we
can find elements $g^i\in G^K_n$ such that
$v_i:=(g_i(g^1),\ldots, g_i(g^l))$ are linearly independent
vectors for $i=1,\ldots, l$. Then we find coefficients
$\alpha_{ij}$ such that $\sum_j \alpha_{ij}v_j=(0,\ldots,1_j,\ldots,0)$.
Then $f_i=\sum_j \alpha_{ij}{f_i}g_i(g^j)=\sum_j \alpha_{ij}(g^j\cdot f_i) 
\in I_n$. Thus for any $f\in I_n$, $\Phi^*(f)\in \pi^*(I_n)$, which gives
$\Phi^*(I)\subset \pi^*(I)$. 

Denote by ${\rm inv}$ the automorphism of $G\widehat{\times}\widehat{X}_x$ induced
 by taking the inverse $g\mapsto g^{-1}$. Then $\pi^*(I)=\Psi^{-1*}\Phi^*(I)\subset
\Psi^{-1*}\pi^*(I)=\Psi^{*}{\rm inv}^*\pi^*(I)=\Psi^{*}\pi^*(I)=\Phi^{*}(I)$.
\qed 

\bigskip\goodbreak
Let $L$ denote an algebraically closed field containing $K$. 
Set $$ X^{(n)L}_x:=X^{(n)}_x\times_{\Spec{K}}\Spec{L},$$  
$$\widehat{X}^L_x:=\lim_{\leftarrow}X^{(n)L}_x=\Spec{L[[\cO_{X,x}]]}$$
Let ${\rm Aut}^\Gamma_L(X^{(n)L}_x,S)$ denote the group of
$\Gamma_x$-equivariant
$L$-automorphisms of $X^{(n)L}_x$ preserving the relevant
ideals of strata and $G:={\rm Aut}_L(\widehat{X}^L_x,S)$ denote the group of
$L$-automorphisms of $\widehat{X}^L_x$. 

\begin{lemma}\begin{enumerate}\label{le: base}
\item ${\rm Aut}^\Gamma_L(X^{(n)L}_x,S)=
{\rm Aut}^\Gamma_K(X^{(n)}_x,S))\times_{\Spec{K}}\Spec{L}$.   
   
\item 
${\rm Aut}^\Gamma_L(\widehat{X}^L_x,S)={\rm
Aut}^\Gamma_K(\widehat{X}_x,S)\times_{\Spec{K}}\Spec{L}$.

\item Let $G\subset {\rm Aut}_K(\widehat{X}_x,S)$ be a
proalgebraic subgroup. 
An ideal $I\subset \widehat{\cO}_{X,x}$  is 
$G$-invariant iff
$I\cdot L[[\cO_{X,x}]]$ is $G\times_{\Spec{K}}\Spec{L}$-invariant. \qed 
\end{enumerate}
\end{lemma}
\noindent{\bf Proof.}
(1) and (2) follow from the construction of
${\Aut}^\Gamma_K(\widehat{X}_x,S)$. 

(3)($\Rightarrow$) It follows from Lemma 
\ref{le: group ideal} that $I\subset \widehat{\cO}_{X,x}$  is 
$G$-invariant iff
$\Phi^*(I)= \pi^*(I)$  or 
$\Psi^*(\frac{I+m_x^n}{m_x^n})=
\pi^*(\frac{I+m_x^n}{m_x^n})$.
Let $\Psi^L$ and
$\pi$ be pull-backs of the morphisms $\Psi$ and
$\pi^L$ under ${\Spec{L}}\to \Spec{K}$. Then
$\Psi^{L*}(\frac{I+m_x^n}{m_x^n} \cdot L)=  
\pi^{L*}(\frac{I+m_x^n}{m_x^n}\cdot L)\subset L\otimes_{K}\widehat{\cO}_{X,x}/m_x^n$,
which shows , by Lemma \ref{le: group ideal}, that 
$I\cdot L[[\cO_{X,x}]]$ is
$G\times_{\Spec{K}}\Spec{L}$-invariant. 

($\Leftarrow$). Assume the latter
holds then $\Psi^{L*}(\frac{I+m_x^n}{m_x^n}\cdot L)= 
\pi^{L*}(\frac{I+m_x^n}{m_x^n}\cdot L)\subset
L\otimes_{K}\widehat{\cO}_{X,x}/m_x^n$. 
Let $\Ga(L/K)$ be the Galois group of the extension $K\subset
L$. Both homomorphisms $\pi^{L*}$ and
$\Psi^{L*}$ are $\Ga(L/K)$-equivariant.
Considering $\Ga(L/K)$-invariant elements gives $\Psi^*(\frac{I+m_x^n}{m_x^n})=
\pi^*(\frac{I+m_x^n}{m_x^n})$.\qed

\begin{example}\label{ex: differential} Let $\phi_n:{\rm Aut}(\widehat{X}_x,S)\to {\rm Aut}(X^{(n)}_x,
S)$ denote the natural morphisms. For $n=1$ we get 
the differential mapping:
$$d=\phi_1:{\rm Aut}(\widehat{X}_x, S) \longrightarrow
{\rm Aut}(X^{(1)}_x,S)\subset
{\rm Gl}({\rm Tan}_{X,x})).$$  
\end{example}

\bigskip
\subsection{Proof of Demushkin's Theorem}

In what follows we shall use the following
 generalization of
Theorem \ref{th: Dem0}.

\begin{lemma}\label{le: Dem2} 

Let $\sigma$ and $\tau$ be
two  $\Gamma$-semicones in isomorphic lattices $N_\sigma \simeq N_\tau$.

Then  there exists a $\Gamma$-equivariant isomorphism
 $\widehat{X}_{\sigma}
 \simeq \widehat{X}_{\tau}$ preserving strata iff
 there exists an isomorphism of $\Gamma$-semicones 
$\sigma \simeq \tau$. 

\end{lemma}

For a  $\Gamma$-semicone $\sigma$  in a lattice $N_\sigma$ 
denote by
${\rm Aut}(\widehat{X}_{\sigma})$ the group of all $\Gamma$-equivariant automorphisms
 of $\widehat{X}_\sigma$ which preserve strata defined by 
faces   of the $\Gamma$-semicone $\sigma$.

\begin{lemma}(Demushkin \cite{Demushkin}, also Gubeladze
\cite{Gubeladze}) \label{le: Dem4}
Let $\sigma$ and $\tau$ be
two  $\Gamma$-semicones in isomorphic lattices $N_\sigma \simeq N_\tau$.
\begin{enumerate}
\item  
The torus  $T_\sigma$ associated to $
\widehat{X}_{\sigma}$ determines a maximal torus in  the proalgebraic group
 $G:={\rm Aut}(\widehat{X}_{\sigma})$. 

\item 
Let $d:G\to \Gl(\Tan_{X_\sigma,O_\sigma})$ be the differential morphism as
in Example
 \ref{ex: differential}. Then $d(T_\sigma)$ is a
maximal torus in $d(G)$.

\item 
Any $\Gamma$-equivariant isomorphism $\widehat{\phi}:\widehat{X}_{\sigma}
\simeq \widehat{X}_{\tau}$, preserving strata
determines an action of the torus $T_\tau$ on 
$\widehat{X}_{\sigma}$ and a group
embedding $T_\tau\hookrightarrow G$.  
If $d(T_\sigma)=d(T_\tau)$ then
$T_\sigma$ and $T_\tau$  are conjugate in $G$.

\item Any two tori $T_\sigma, T_\tau\subset G$ which are
determined by $\widehat{\phi}$,  are conjugate in $G$. 
\end{enumerate}
\end{lemma}

\noindent {\bf Proof.} 
(1) The torus $T_\sigma$   is maximal in
$G$ since its centralizer consists of
$T_\sigma$-equivariant automorphisms of
$\widehat{X}_\sigma$, preserving strata
and therefore it coincides
with $T_\sigma$ ($T_\sigma$-equivariant
automorphisms multiply characters by constants and are
defined by elements of  $T_\sigma$); (see also Step 4 in \cite{Demushkin}).

(2) By Lemma \ref{le: epimorphisms} we can write $G=\lim_{\leftarrow} G_i$, where
all homomorphisms $G_i\to G_j$ are epimorphisms.
We show that the relevant images of
$T_\sigma$ are maximal tori in
all groups $G_i$. Otherwise we find maximal
tori $T_{i}\subset G_i$ containing the  images of $T_\sigma$,
such that $\lim_{\leftarrow} T_{i}$  is a maximal torus in
$G$ containing $T_\sigma$ as a proper\ subgroup (see
\cite{Borel} Proposition 11.14(1));
 (see also step 9 in \cite{Demushkin}).

(3) Let $T_\sigma\subset G$ be the maximal torus associated to
$\widehat{X}_{\sigma}$. Let $u_1,\ldots,u_k$ be characters of
$T_\sigma$ generating the maximal ideal of $O_\sigma$ and
$x_1,\ldots,x_k$ be semiinvariant functions determined by the  action of 
$T_\sigma$,
induced on the tangent space ${\rm Tan}_{O_{\sigma},X_{\sigma}}$.  
Characters of $T_\sigma$ define the    
 natural $T_\sigma$-equivariant  embedding
$\phi_{\sigma}:\widehat{X}_{\sigma}
\rightarrow {\rm Tan}_{O_{\sigma},X_{\sigma}}$ 
into the tangent space induced by the ring epimorphism
 $\phi^*_{\sigma}:K[[x_1,\ldots,x_k]]\to
K[[u_1,\ldots,u_k]]$ sending $x_i$ to $u_i$.  
 The morphism $\phi_{\sigma}$ determines the
isomorphism $\widehat{\phi}_{\sigma}:
\widehat{X}_{\sigma}\rightarrow
\widehat{\phi}_{\sigma}(\widehat{X}_{\sigma})$ preserving strata corresponding to
 the isomorphism $\widehat{\phi}^*_{\sigma}:K[[x_1,\ldots,x_k]]/I \to
K[[u_1,\ldots,u_k]]$ where  the ideal $I$ is generated by all
elements $x^\alpha-x^\beta$ where $x^\alpha, x^\beta$ are
monomials on which $T_\sigma$ acts with the same weights.

If $d(T_\sigma)=d(T_\tau)$ then  
$\phi_{\sigma}(\widehat{X}_{\sigma})=\phi_{\tau}(\widehat{X}_{\tau})=\Spec K[[x_1,\ldots,x_k]]/I $.
Consequently, $\beta:=\phi_{T_\tau}^{-1}\phi_{T_\sigma}$ 
is a $T_\sigma=d(T_\sigma)$-equivariant 
isomorphism $\widehat{X}_{\sigma}\to
\widehat{X}_{\tau}$ preserving strata. 
It follows that $\beta$ determines a conjugating automorphism with
$T_\sigma=\beta T_\tau\beta^{-1}$.

(4)  By  Borel's theorem (\cite{Borel}) and (2), the tori $d(T_\sigma)$
and $d(T_\tau)$ are conjugate in $d(G)$. Consequently, we
can find a torus $T$ conjugate to  $T_\sigma$
and such that $d(T_\tau)=d(T)$. By (3) $T_\tau$ and $T$ 
are conjugate. \qed 

\bigskip
\noindent {\bf Proof of Lemma \ref{le: Dem2}.}
Let   
$\widehat{\phi}: \widehat{X}_{\sigma}\rightarrow
\widehat{X}_{\tau}$ be a $\Gamma$-equivariant isomorphism
preserving strata. By Lemma \ref{le: Dem4} there is a
$\Gamma$-equivariant automorphism 
$\alpha$ of $\widehat{X}_\sigma$ for which ${
T}_\sigma=\alpha{
T}_\tau\alpha^{-1}$.
Then  $T_\sigma$ acts on $\widehat{X}_{\tau}$ as a big torus and
the $\Gamma$-equivariant isomorphism
$\widehat{\phi}\circ\alpha^{-1}:\widehat{X}_{\sigma}\rightarrow
\widehat{X}_{\tau}$ is $T_\sigma$-equivariant and
preserves strata. Such an isomorphism yields an
isomorphism of the associated $\Gamma$-semicones. \qed

\subsection{Singularity type of nonclosed points}
\bigskip
For a ring $L$ and a cone $\sigma$ in $N_\sigma$ 
 denote by 
$L[[{\sigma}^\vee]]$ the completion of the ring $L[{\sigma}^\vee]$
at the  ideal $m_\sigma$ of $O_\sigma$.

\begin{lemma}\label{le: decomposition} Let $L$ be a ring
containing $K$, and $\sigma$ and
$\tau$ be $\Gamma$-semicones in isomorphic lattices
$N_\sigma\simeq N_\tau$.
 Let $m_\sigma, m_\tau$ denote the
ideals  of $O_{\sigma}$ and $O_{\tau}$ respectively.
 The action of $\Gamma$ on $K[[{\sigma}^\vee]]$ and $K[[{\tau}^\vee]]$
induces an action on
$L[[{\sigma}^\vee]]$ and $L[[{\tau}^\vee]]$ which is
trivial on $L$.

 Any $\Gamma$-equivariant isomorphism 
$\phi: L[[{\sigma}^\vee]]\simeq L[[{\tau}^\vee]]$ 
for which $\phi(m_\sigma)=m_\tau$ can be decomposed as 
$\phi=\phi_0\phi_1$, where $\phi_0$ is a $\Gamma$-equivariant automorphism of 
$L[[{\tau}^\vee]]$ identical on monomials and $\phi_1$ 
is a $\Gamma$-equivariant $L$-isomorphism.

\end{lemma}

\noindent {\bf Proof.} 
(1)  The restriction of $\phi$ to ${L}$
is a ring monomorphism 
$\phi_{|{L}}:{L}\to 
{L}[[{\tau}^\vee]]$. This ring homomorphism
can be extended to a  ring homomorphism 
$\phi_0:{L}[[{\tau}^\vee]]\to
{L}[[{\tau}^\vee]]$ by sending monomials
identically to the same monomials, i.e. 
$\phi_0(\sum_{\alpha\in{\tau}^\vee}
a_\alpha x^\alpha)=\sum_{\alpha\in{\tau}^\vee}
 \phi_0(a_\alpha) x^\alpha$, where $ \phi_0(a_\alpha)=
\sum_{\beta\in{\tau}^\vee} a_{\alpha\beta} x^\beta$. 
Note that in the expression $\sum_{\alpha\in{\tau}^\vee}
\phi_0(a_\alpha) x^\alpha$, the coefficient of $x^\alpha$
is a finite sum since there are finitely many
possibilities to express $\alpha$ as $\alpha'+\beta$ where
$\alpha', \beta\in {\tau}^\vee$. Thus the above expression defines a
formal power series and consequently $\phi_0$ is a ring homomorphism. 
Then $\phi_{0|L}$ determines an  automorphism   ${\phi}_1: {L}\to
{L}[[{\tau}^\vee]]/m_\tau={L}$. Let
${\phi}_2$ denote the
automomorphism of ${L}[[{\tau}^\vee]]$ induced
by the ${\phi}_1$ and identical on monomials.
 The composition 
$\psi=\phi_2^{-1}\phi_0$ is a ring endomorphism which as a
$K$-linear transformation can be written in the form
$\psi=id+p$, where $p(m^i)\subset m^{i+1}$. Such a linear
transformation is invertible (with inverse
$id-p+p^2-\ldots$),  hence $\psi$ is a linear isomorphism and  ring
isomorphism. This implies that $\phi_0$ is an
automorphism and $\phi_1:=\phi_0^{-1}\phi$ is an ${L}$-isomorphism.  

Note that the above morphisms are $\Gamma$-equivariant.
\qed 

\bigskip
The following is a generalization of the Demushkin Lemma.

\begin{lemma}\label{le: Dem3} Let $\sigma$ and
$\tau$ be $\Gamma$-semicones in isomorphic lattices
$N_\sigma\simeq N_\tau$, and $L$ be a field
containing $K$.
Consider the induced  action of $\Gamma$ on
$L[[{\sigma}^\vee]]$ and $L[[{\tau}^\vee]]$ which is
trivial on $L$.
  There exists a $\Gamma$-equivariant isomorphism
$L[[{\sigma}^\vee]]\simeq L[[{\tau}^\vee]]$ over $K$, preserving strata iff 
$\sigma$ and $\tau$ are $\Gamma$-isomorphic.

\end{lemma}

\noindent{\bf Proof.} Denote by $\phi$
an isomorphism 
$L[[{\sigma}^\vee]]
\buildrel \phi \over \simeq  L[[{\tau}^\vee]]$.
 By Lemma \ref{le: decomposition}, $\phi=\phi_0\phi_1$, where
$\phi_1$ is an ${L}$-isomorphism.
Tensoring with the algebraic closure $\overline{L}$ of $L$
and
taking completion determines an $\overline{L}$-isomorphism 
$\overline{\phi}_1: \overline{L}[[{\sigma}^\vee]]
\simeq \overline{L}[[{\tau}^\vee]]$. 
It suffices to apply Theorem \ref{th: Dem} to the above 
$\overline{L}$-isomorphism 

Note that the above isomorphisms are $\Gamma$-equivariant.
\qed

This lemma allows us to extend the notion of singularity type to
any nonclosed points on  toric varieties $X_{\Delta}$ and
even on some other toroidal schemes.

\begin{lemma} \label{le: 1} Let $\Delta$ be a subdivision
of $\sigma$ and
$\widehat{X}_{\Delta}:=X_{\Delta}\times_{X_{\sigma}}\widehat{X}_{\sigma}$
be a toroidal scheme with an  action of a group $\Gamma\subset T_\sigma$.
Let
$p\in Y=\widehat{X}_{\Delta}$ be a point which need not
be closed. Set
$\Gamma_p:=\{g\in\Gamma\mid g(p)=p, g_{|K_p}={\rm
id}_{|K_p}\}$. 

\begin{enumerate}

\item Let $O_{\tau,Y}\in Y$ denote the locally closed
subscheme defined by a toric orbit ${O}_\tau\subset
X_{\Delta}$, where  $\tau\in\Delta$. Then $\widehat{\cO}_{O_{\tau},Y}\simeq
K_{O_{\tau}}[[\underline{\tau}^\vee]]$, where $K_{O_{\tau}}$ is the
residue field of ${O_{\tau}}$.  

\item There is a  $\Gamma_p$-equivariant isomorphism 
$\widehat{\cO}_{p,Y}\simeq K_p[[\sigma_p^\vee]]$,
where $K_{p}$ is the
residue field of ${p}$ and $\sigma_p$ is a
uniquely determined cone. Moreover 
$\sing^\Gamma(p):=(\Gamma_p,
\underline{\sing^{\Gamma_p}(\sigma_p)})=\sing^\Gamma(\tau):=(\Gamma_\tau,
\underline{\sing^{\Gamma_\tau}(\tau)}) $, where $O_{\tau,Y}\in Y$ 
is the minimal orbit scheme 
 containing  $p$.

\end{enumerate}
\end{lemma}

\noindent{\bf Proof.} 
(1)  
Let $X_{\tau,Y}:=X_{\tau}\times_{X_\sigma}\widehat{X_\sigma}\subset
Y$. Then  $K[X_{\tau,Y}]=K[\tau^\vee]\otimes_
{K[\sigma^\vee]}K[[\sigma^\vee]]$. Let $T_\tau\subset T$ be
the torus
corresponding to the sublattice $N_\tau:= N\cap\lin(\sigma)$.
Then $T_\tau$ acts on $K[X_{\tau,Y}]$ with nonnegative
weights. Therefore the subring of ${T_\tau}$-invariant
functions $K[X_{\tau,Y}]^{T_\tau}$ in $K[X_{\tau,Y}]$
equals
$K[\tau^\perp]\otimes_{K[\sigma^\vee\cap
\tau^\perp]}K[[\sigma^\vee\cap
\tau^\perp]]$. The ideal $I=I_{O_{\tau,
Y}}\subset K[X_{\tau,Y}]$ of the orbit
$O_{\tau, Y}$ is generated by all monomials with positive weights
in the set $\sigma^\vee+\tau^\vee\subset
\tau^\vee=\underline{\tau}^\vee\oplus \tau^\perp$ and
consequently it is generated by
$\underline{\tau}^\vee\setminus \{0\}$.
Thus $K[X_{\tau,Y}]^{T_\tau}\simeq K[O_{\tau,Y}]$ and 
$\lim_{\leftarrow}K[X_{\tau,Y}]/I^k\simeq 
K[O_{\tau,Y}][[\underline{\tau}^\vee]]$.
This gives $K[\widehat{Y}_{O_{\tau, Y}}]
\simeq K(O_{\tau,Y})[[\underline{\tau}^\vee]]$.

(2) Let $O_{\tau,Y}\in Y$ 
be the minimal orbit scheme 
 containing  $p$. First we prove that 
$\Gamma_{O_\tau}=\Gamma_p$. We have the obvious inclusion 
$\Gamma_{O_\tau}\subset \Gamma_p$. Now if $g\in \Gamma_p$ then
for any $f\in K[{O_\tau}]$, $(g(f)-f)\in I_p$, where
$I_p\subset K[{O_\tau}]$ describes $p$. Since
$\Gamma_p=\Gamma_{p'}$ for any $p'\in\Gamma\cdot p$, it follows that 
$(g(f)-f)\in \bigcap_{h\in T} h\cdot I_p =\{0\}$.

 The ideal $I_p\subset {\cO}_{p,O_{\tau}}$ of
 $p$ is generated by local parameters
$u_1,\ldots,u_l$. Then there are $\Gamma_p$-equivariant
isomorphisms 
$\widehat{\cO}_{p,O_{\tau}}\simeq K_p[[u_1,\ldots,u_l]]$ and 
$\widehat{\cO}_{p,Y}\simeq K_p[[u_1,\ldots,u_l]][[\tau^\vee]]=
K_p[[(\tau\oplus^{\Gamma_p} \langle e_1,\ldots,e_l \rangle)^\vee]] $
where $\widehat{\cO}_{p,O_{\tau}}=\widehat{\cO}_{p,Y}^{T_\tau}$.
Therefore $\sigma_p=\tau\oplus^{\Gamma_p} \langle e_1,\ldots,e_l
\rangle$ and $\underline{\sing^{\Gamma_p}}(\sigma_p)=
\underline{\sing^{\Gamma_\tau}(\tau)}$.
\qed

By
Lemmas \ref{le: Dem3} and by \ref{le: 1} the singularity type 
of a nonclosed point $p$, 
$ {\rm sing}(p):={\rm sing}(\sigma_p)$ (resp. $ {\rm sing}^{\Gamma_p}(p):=
(\Gamma_p, \underline{\rm sing}^{\Gamma_p}(\sigma_p))$, is  
uniquely determined.
Moreover singularity type determines a stratification ${\rm
Sing}(Y)$ (resp. ${\rm
Sing}^\Gamma(Y)$) on $Y$ such that all points in the same stratum have the
same singularity type. This yields

\begin{lemma}\label{le: singularity type2} Let $\Gamma$
act on 
$\widehat{X}_{\Delta}:=X_{\Delta}\times_{X_{\sigma}}\widehat{X}_{\sigma}$
There is a 
stratification ${\rm Sing}^\Gamma(\widehat{X}_{\Delta})$  of
$\widehat{X}_{\Delta}$ which is determined by singularity
type and therefore preserved by any $\Gamma$-equivariant
automorphism of $\widehat{X}_{\Delta}$.  \qed    
\end{lemma}

\bigskip
\subsection{Semicomplex associated to a stratified toroidal
variety}

For a $\Gamma$-semicone $\sigma$ we denote by $X_\sigma$ the
associated stratified toric variety.

\begin{definition} \label{de: associated semicomplex}  
Let $(X,S)$ be a $\Gamma$-stratified toroidal
variety. We say that a semicomplex $\Sigma$ is {\it associated\/} to
$(X,S)$ if there is a bijection $\Sigma \to S$ with the following
properties: Let $\sigma \in \Sigma$ map to
$s=\strat_X(\sigma) \in S$. Then any $x \in
s$ admits an open $\Gamma$-invariant neighborhood
$U_\sigma \subset X$ and a $\Gamma$-equivariant $\Gamma_s$-smooth morphism 
 $ \phi_\sigma \colon U_\sigma \to X_{\sigma}$ of stratified varieties such that  $s
\cap U$ equals $\phi_\sigma^{-1}(O_{\sigma})$ and the intersections $s' \cap
U$, $s' \in S$, are precisely the inverse images of the strata of
$X_{\sigma}$ and the action of $\Gamma$ 
 on $X_{\sigma}$ extends the action of
$\Gamma_s$.

We call the  smooth morphisms $U_\sigma\to X_{\sigma}$ from the above definition 
 {\it charts}. A collection of charts
satisfying the conditions from the above definition is
called an {\it atlas}.

\end{definition}

\begin{remarks} 
\begin{enumerate}
\item Different
$(X,S)$ may have the same associated semicomplex $\Sigma$. 
Smooth varieties endowed with the trivial stratification have the 
associated semicomplex consisting of one zero cone.
\item
The action of $\Gamma_s=\Gamma_\sigma$ is fixed for
$\sigma$ while the action of $\Gamma$ on $X_\sigma$ depends upon charts.
\end{enumerate}
\end{remarks}

\noindent 

\begin{lemma}\label{le: associated semicomplex} For any 
$\Gamma$-stratified toroidal variety  $(X,S)$ there
exists a unique (up to isomorphism) associated $\Gamma$-semicomplex
$\Sigma$. Moreover $\tau\leq\sigma$ iff
$\overline{\strat_X(\tau)}\supset \strat_X(\sigma)$. 
$(X,S)$ is a toroidal embedding iff
$\Sigma$ is a complex.
\end{lemma}

\noindent {\bf Proof } First we assign to any stratum $s$
a  semicone $\sigma$. 
 
By Definition \ref{de: stratified toroidal} there  is a
$\Gamma_s$-smooth  morphism 
${\phi}_\sigma: U_\sigma \rightarrow X_{\sigma}$ into a
stratifed toric variety $(X_\sigma,S_\sigma)$ preserving strata.
Note that $\Gamma_s$ acts trivially on $O_\sigma$. So we
can assume that $\sigma=\underline{\sigma}$ is of maximal dimension in
$N_\sigma$ by composing $\phi$, if necessary, with a suitable projection. 

We define the $\Gamma_s$-semicone  $\sigma$ to be the $\Gamma_s$-semifan 
associated to the $\Gamma_s$-stratified toric variety $(X_{\sigma},
S_{\sigma})$.

Then there is  an open subset $U$ of $U_\sigma$ 
intersecting $s$ and a { $\Gamma_s$-\'etale}  morphism 
$\phi:=\phi_\sigma\times\psi: U \to X_{\sigma}\times {\bf A}^k=X_{\tau}$
preserving strata,
where  $\tau={\sigma\times \langle
e_1,\ldots, e_{{\rm dim}(s)}}\rangle $,
and $\psi: U\to {\bf A}^k$ is a morphism defined by local
parameters $u_1,\ldots,u_k$ on $s$. The stratification of
$X_\tau$ is defined by the embedded semifan $\sigma\subset
\overline{\tau}$, consisting of the faces of the semicone $\sigma$ in the
lattice $N_\tau$.
Then
$\widehat{\phi}: \widehat{X}_x \to
\widehat{X}_\tau$ is a $\Gamma_s$-equivariant isomorphism
preserving strata.
By Lemma \ref{le: Dem2} 
the embedded $\Gamma_s$-semifan $\sigma\subset \tau$ and
the $\Gamma_s$-semicone $\sigma$
are determined uniquely up to isomorphism. 
We write $s=\strat_X(\sigma)$ and define $\Gamma_\sigma:=\Gamma_s$.

Assume that the closure of a stratum $s=\strat_X(\sigma)$
contains a stratum $t=\strat_X(\tau)$. 
Let $\phi_\tau: U_\tau
\to X_\tau$ denote a chart associated with $\tau$. 
Then the stratum $s\cap U_\tau$  
and the strata of $U_\tau$ having $s\cap U_\tau$ in their closure determine a
semicone $\sigma_s\subset \overline{\tau}$. It follows from the
uniqueness of $\sigma$ that there is a
saturated embedding $i^\tau_{\sigma}$ of the semicone
$\sigma$ into  the semicone $\tau$ with the image 
$i^\tau_{\sigma}(\sigma)=\sigma_s$. Then we shall write $\sigma\leq\tau$.

The $\Gamma$-semicomplex $\Sigma$ is defined as the collection of
the $\Gamma_\sigma$-semicones $\sigma$, and the saturated face
inclusions $i^\tau_{\sigma}$ for $\sigma\leq \tau$.

Now let $\Sigma$ and $\Sigma'$ be two semicomplexes
associated to $(X,S)$. By uniqueness there are
isomorphisms of semicones $j_\sigma: \sigma\to \sigma'$ (see Definition 
\ref{de: semicomplexes isomorphism}). These isomorphisms
induce an isomorphism of semicomplexes.

The second part follows from the   fact that
locally toroidal embeddings correspond to fans consisting
of all faces of some cone. 

 \qed

\begin{lemma} \label{le: asemicomplex} Let $(X_\Sigma,S)$
be a $\Gamma$-stratified
toric variety corresponding to an embedded
$\Gamma$-semifan 
$\Omega\subset\Sigma$. Then $(X_\Sigma,S)$ is a $\Gamma$-stratified
toroidal variety with the associated
$\Gamma$-semicomplex $\Omega^{\rm semic}$. 

There is an atlas  $$\cU^{\rm can}(X_\Sigma,S)=\cU(\Sigma,\Omega)$$ on
$(X_\Sigma,S)$ defined as follows: For any $\sigma$ in $N$ such that 
$\omega(\sigma)=\omega$ in $N_\omega$ 
 there is a chart
$\phi_\sigma: X_{\sigma}\to X_{\omega}$ given
by any projection $\pi^\omega_\sigma: \sigma 
\to \omega$ such that $\pi^\omega_{\sigma|\omega}=\id_{|\omega}$.
\end{lemma}

\noindent {\bf Proof.} Follows from Proposition \ref{le: semifans
correspondence} and from the definition of the associated
semicomplex.\qed

\subsection{Local properties of orientation}

\begin{definition}\label{de: orientation}
 We shall call a proalgebraic group
$G$ {\it connected} if
it is a connected affine scheme. 
For any proalgebraic group $G=\lim_{\leftarrow}G_i$, denote by 
$G^0$  its maximal
connected proalgebraic subgroup $G^0=\lim_{\leftarrow}G^0_i$.
\end{definition}
\begin{lemma} $G=\lim_{\leftarrow}G_i$ is connected if each
$G_i$ is irreducible.
\end{lemma}

\noindent{\bf Proof.} 
By Lemma \ref{le: epimorphisms} we can assume all morphisms $G_i\to
G_j$ to be epimorphisms 
and $K[G]=\bigcup K[G_i]$.
If $g_1,g_2\in K[G]$ and $g_1\cdot g_2=0$ then $g_1,g_2\in
K[G_i]$ for some $i$ and $g_1$ or $g_2$ equals zero. \qed
  
\begin{lemma} \label{le: images} Let
$G=\lim_{\leftarrow}G_i$ be a connected proalgebraic group. Then the image
$\phi_i(G)\subset G_i$ of the
natural homomorphism $\phi_i:G\to G_i$ is connected. \qed 
\end{lemma}

\begin{definition}\label{de: the same orientation2} Let $X$
be a stratified (resp.  $\Gamma$-stratified) toroidal
scheme over $K$. We say
that an automorphism (resp. $\Gamma_x$-equivariant automorphism) 
$\phi$ of $X$ {\it preserves
orientation} at a $K$-rational point $x=\phi(x)$ if
it induces an automorphism $\widehat{\phi}\in
\Aut(\widehat{X}_x,S)^0$ (resp. $\widehat{\phi}\in
\Aut^\Gamma(\widehat{X}_x,S)^0$).
\end{definition}
\begin{definition}\label{de: the same orientation}
We say that two \'etale morphisms (resp. $\Gamma$-\'etale
morphisms) of stratified toroidal
schemes over $K$,
$f_1,f_2:(X,S)\to (Y, T)$, {\it determine the same
orientation} at a closed $K$-rational point
 $x\in X$ if $f_1(x)=f_2(x)$ and
$(\widehat{f}_2)^{-1}\circ \widehat{f_1}\in {\rm
Aut}(\widehat{X}_x,S)^0$ 
(resp. $(\widehat{f}_2)^{-1}\circ \widehat{f_1}\in {\rm
Aut}^\Gamma(\widehat{X}_x,S)^0$).
\end{definition}

In further considerations we shall consider the case of
$\Gamma$-stratified toroidal schemes. The case of 
stratified toroidal schemes corresponds to the situation
when $\Gamma$ is a trivial group.

\begin{definition} Let $(X,S)$ is a $\Gamma$-stratified toroidal
scheme and $x\in X$ be a $K$-rational point. We call
functions $y_1,\ldots,y_k$ {\it locally toric parameters}
if $y_1=\phi^*_x(u_1),\ldots, y_k=\phi^*_x(u_k)$,
where $u_1,\ldots,u_k$ are semiinnvariant generators at the
orbit point $O_\sigma$ on 
a $\Gamma$-stratified toric variety $X_{\sigma}$, and
$\psi_x: U_x\rightarrow
X_{\sigma}$ is a $\Gamma_x$-\'etale morphism from an open
$\Gamma_x$-invariant 
neighborhood $U_x$ of $x$ such that
  $\psi_x(x)=O_\sigma$ 
and  the intersections $s'\cap U$, $s'\in S$,  are precisely the inverse
images of strata of $X_\sigma$. 
\end{definition} 

\begin{lemma} Locally toric parameters exist for any
$K$-rational point of a $\Gamma$-stratified toroidal
scheme $X$.
\end{lemma} 
\noindent{\bf Proof} Let $s$ be the stratum through $x$. 
By Definition \ref{de: G-stratified
toroidal} there is a
$\Gamma_x$-\'smooth  morphism from an open
neighborhood $U_x$ of $x$ such that
  $s=\phi_x^{-1}(O_\sigma)$ 
and  the intersections $s'\cap U$, $s'\in S$,  are precisely the inverse
images of strata of $X_\sigma$. Let $x_1,\ldots,x_k$ be
local paramters of $s$ at $x$. Set $g:U\to {\bf A}^1$,
where $g(x)=(x_1,\ldots,x_k)$. Then the morphism
$\psi_x:=\phi_x\times g$ is $\Gamma_x$-\'etale and defines
locally toric parameters at $x$. \qed

\begin{lemma} \label{le: loc} Let $(X,S)$ be a $\Gamma$-stratified
toroidal scheme and $x\in X$ be a $K$-rational point.
Let $y_1,\ldots,y_n$ be locally toric
paramters at $x\in X$.
The ideals of closures of
the strata $I_s\subset {\cO}_{x,X}$
are generated by subsets of  $\{u_1,\ldots,u_k\}$.  
\end{lemma}
\noindent{\bf Proof.} It suffices to show the lemma for characters
generating $\sigma^\vee$ on a $\Gamma$-stratified toric
variety $X_\sigma$. Let $T$ be the big torus on $X_\sigma$.
Each stratum $s$ is ${T}$-invariant irreducible
hence contains a ${\bf T}$- orbit $O_\tau$, where
$\tau\preceq \sigma$. We conclude that
$\overline{s}=\overline{O_\tau}$. But then
$I_{\overline{O_\tau}}\subset K[X_\sigma ]$ is generated by 
functions corresponding to those generating functionals of
$\sigma^\vee$ which are not zero on $\tau$.\qed

\begin{lemma} \label{le: smooth}\begin{enumerate}
\item The  group ${\rm Aut}^\Gamma(\widehat{X}_x,S)$ 
 is connected
for any smooth  $\Gamma$-stratified toroidal 
scheme $(X,S)$ over $\Spec(K)$ and any  
$K$-rational point $x\in X^\Gamma$. 

\item
 Any two $\Gamma$-\'etale morphisms $g_1: (X,S)\to (Y,R)$, 
$g_2: (X,S)\to (Y,R)$ 
between smooth $\Gamma$-stratified toroidal schemes such that $g_1(x)=g_2(x)$ determine the same orientation at $x$. 

\end{enumerate}
\end{lemma}

\noindent{\bf Proof.} (1)  We can replace   $X$ by ${\bf A}^k$ since 
$(\widehat{X}_x,S)\simeq (\widehat{{\bf A}^k}_0,S_A)$ for a
toric stratification $S_A$ on ${\bf A}^k$. Let $m$ be the
maximal ideal of $0\in {\bf A}^k$. 
The automorphism $g\in G:={\rm
Aut}^\Gamma(\widehat{{\bf A}^k}_0,S_A)$ is defined by
$\Gamma$-semiinvariant functions 
$g^*(x_1),\ldots,g^*(x_k)$, where $x_1,\ldots, x_k$ are the
standard coordinates on ${\bf A}^k$ such that the
$\Gamma$-weights of $g^*(x_i)$ and
$x_i$, where $i=1,\ldots,k$, are equal. There is a
birational map $\alpha: {\bf A}^1\to G$ defined by
$$\alpha(z):=(x_1,\ldots,x_k)\mapsto
((1-z)g^*(x_1)+zx_1,\ldots,(1-z)g^*(x_k)+zx_k).$$ 
Note that $\alpha(z)$
defines automorphisms for the open subset $U$ of ${\bf A}^1$,
where the linear parts of coordinates of $\alpha(z)$ are
linearly independent. By Lemma \ref{le: loc} the closure
of each stratum $s$
is described by a subset of $\{x_1,\ldots,x_k\}$. Since $g$
preserves strata, $\overline{s}$ is described by the
corresponding subset of $\{g^*(x_1),\ldots,g^*(x_k)\}$.
Thus the corresponding coordinates of $\alpha(z)$, $z\in
U$, belong to the ideal $I_{\overline{s}}$. Since they are
linearly independent of order $1$ they generate the ideal
$\frac{I_{\overline{s}}}{{I_{\overline{s}}}\cdot m}$ and
by 
the Lemma  of Nakayama
they generate the ideal $I_{\overline{s}}$. Therefore the
automorphisms  $\alpha(z)$, $z\in
U$, preserve strata.
The morphism $\alpha$ ''connects'' the identity $\alpha(1)={\rm
id}$ to an  arbitrary
element  $\alpha(0)=g\in {\rm Aut}^\Gamma(\widehat{X}_x,S)^0 $.

(2) Follows from (1).  \qed

\begin{lemma} \label{le: extensions} 
Let $f: (X,S)\to (Y,R)$ be  a $\Gamma$-smooth morphism of relative
dimension $k$ between $\Gamma$-stratified toroidal schemes. Let $x\in X^\Gamma$ be
a closed $K$-rational point. Let $\Gamma$ act trivially on
${\bf A}^k$ and $g_1,g_2:
X\to {\bf A}^k$ be any two $\Gamma$-equivariant morphisms such that $g_1(x)=g_2(x)=0$ and $f\times
g_i: (X,S)\to (Y,R) \times {\bf A}^k$ are $\Gamma$-\'etale for $i=1,2$. 
Then $f\times
g_1$ and $f\times g_2$ determine the same orientation at
 $x$.

\end{lemma}

\begin{definition}\label{de: extensions} 
We shall call such  a morphism $f\times g$ an {\it
\'etale extension} of $f$ and denote it by $\widetilde{f}$. 
\end{definition}

\noindent {\bf Proof of \ref{le: extensions}}.  Let $y=f(x)$ and
$y_1,\ldots,y_l$ be locally toric parameters at $y\in Y$. Let
$x_1,\ldots,x_k$ be standard coordinates at  $0\in {\bf A}^k$. Set
$v_i=f^*(y_i)$ for $i=1,\ldots,l$ and $w^1_j=g_1^*(x_j)$, 
$w^2_j=g_2^*(x_j)$ for $j=1,\ldots,k$. Then
$v_1,\ldots,v_l,w^1_1,\ldots,w^1_k$ and
$v_1,\ldots,v_l,w^2_1,\ldots,w^2_k$ are locally toric parameters at $x$. The
automorphism $\widehat{(f\times
g_1)}_x\circ \widehat{(f\times g_2)}_x^{-1}: \widehat{X}_x\to
\widehat{X}_x$  maps $w^1_i$ to $w^2_i$.
We can find a linear 
automorphism $\alpha_1\in \{{\rm id_v}\} \times \Gl(k)$, which
preserves $v$-coordinates and
acts nontrivially on $w^2_i$-coordinates so that
$\alpha_1^*(v_i)=v_i$ and $\alpha_1^*(w^2_i)=w^1_i+z_i$,
 where $z_i$ are some functions from the ideal 
$(w^2_1,\ldots,w^2_k)^2+ (v_1,\ldots,v_l)$. Clearly
$\alpha_1\in {\rm Aut}^\Gamma(\widehat{X}_x,S)^0$. 
 Now consider the morphism ${\bf A}^1\to
{\rm Aut}^\Gamma(\widehat{X}_x,S)^0$ such that $ {\bf A}^1\ni t \mapsto 
\phi_t$, where $\phi^*_t(v_i)=v_i$
and $\phi^*_t(w^1_j)=w^1_j+t\cdot z_i$. This shows that
$\phi_1\in{\rm Aut}^\Gamma(\widehat{X}_x,S)^0$. Finally $(\widehat{f\times
g_1}_x)\circ (\widehat{f\times g_2}_x)^{-1} =
\phi^{-1}_1\circ\alpha_1\in {\rm Aut}^\Gamma(\widehat{X}_x,S)^0$. \qed

\bigskip
\noindent \begin{definition}
We say that two $\Gamma$-smooth morphisms $f_1,f_2:(X,S)\to (Y,R)$
of dimension $k$ 
{\it determine the same orientation} at a closed $K$-rational point $x\in X^\Gamma$ 
if $f_1(x)=f_2(x)$ and
there exist \'etale extensions	
$e{f_1},\widetilde{f_2}: (X,S)\to (Y,R)\times {\bf A}^k$
which determine the same orientation at
 $x$.
\end{definition}

\begin{lemma}\label{le: group of automorphisms}

Let $f:(X,S)\to (Y,R)$ be a $\Gamma$-smooth morphism of $\Gamma$-stratified toroidal
schemes over $K$. Let $x\in X^\Gamma$ and $y=f(x)\in Y^\Gamma$ be
$K$-rational points and $s=f^{-1}(y)\in
S$ be the stratum through $x$. Let $y_1,\ldots,y_l$ be 
locally toric $\Gamma_x$-semiinvariant parameters at $y$ and 
$v_1=f^*(y_1),\ldots,v_l=f^*(y_l), w_1,\ldots,w_k$ be locally toric 
$\Gamma_x$-semiinvariant parameters at $x$.
Set $R:=K[[w_1,\ldots,w_k]]$.
Then
$\widehat{\cO}_{X,x}=R[[v_1,\ldots,v_l]]/I$, where $I$ is the ideal span by binomial relations in $v_i$.
${\rm Aut}^\Gamma_{R}(\widehat{X}_x,S)$ is a proalgebraic group of 
$\Gamma_x$-equivariant automorphisms
 preserving strata and the  subring $R$,
and  the monomorphism   
$\beta:{\rm Aut}^\Gamma_{R}(\widehat{X}_x,S)\to {\rm
Aut}^\Gamma(\widehat{X}_x,S)$ induces an isomorphism $$\overline{\beta}:{\rm
Aut}^\Gamma_{R}(\widehat{X}_x,S)/{\rm
Aut}^\Gamma_{R}(\widehat{X}_x,S)^0 \to 
{\rm Aut}^\Gamma(\widehat{X}_x,S)/{\rm Aut}^\Gamma(\widehat{X}_x,S)^0.$$ 

\end{lemma}

\bigskip
\noindent
{\bf Proof.} It suffices to construct a surjective morphism $\alpha:{\rm Aut}^\Gamma(\widehat{X}_x,S)\to
{\rm Aut}^\Gamma_{R}(\widehat{X}_x,S)$ with
connected fibers such that $\alpha\circ\beta=\id_{{\rm
Aut}^\Gamma_{R}(\widehat{X}_x,S)}$.  
 Any automorphism $\phi$ from ${\rm
Aut}^\Gamma(\widehat{X}_x,S)$ maps  $w_1,\ldots,w_k,v_1,\ldots,v_l$ to 
$w'_1,\ldots,w'_k,v'_1,\ldots,v'_l$. We put
$\alpha(\phi)^*(w_i)=w_i$,
$\alpha(\phi)^*(v_j)=v'_j$. Then $\alpha(\phi)$ is an endomorphism
since $w_i$ are algebraically independent of $v_j$. It is an
automorphism since we can easily define an inverse homomorphism. 
 The fiber  $\alpha^{-1}(\alpha(\phi))$ consists of all elements of the
type $\phi=\alpha(\phi)\circ\phi_1$,  where $\phi_1^*(w_i)=w'_i$ and 
$\phi_1^*(v'_j)=v'_j$.
 It is connected: for
any two elements $\phi=\alpha(\phi)\circ\phi_1$ and $
\phi'=\alpha(\phi)\circ\phi'_1$ we can define a rational map
   $\Phi:{\bf A}^1\to {\rm
Aut}^\Gamma(\widehat{X}_x,S)$ by
$\Phi(t)^*(v'_j)=v'_j$, $\Phi(t)^*(w_i)=tw_i+(1-t)w'_i$.
The latter
is defined on an open subset $U\subset {\bf A}^1$ containing $0$ and $1$. 
\qed

\begin{lemma}\label{le: sections} 
Let $\phi_i: (X,S)\to (Y,R)$  for $i=1,2$ be two 
$\Gamma$-smooth morphisms of $\Gamma$-stratified toroidal schemes such that
$\phi_1(x)=\phi_2(x)=y\in Y$ for a $K$-rational point $x\in
X^\Gamma$ and   
 strata in $S$ are preimages of
strata in $T$ and $\phi^{-1}_1(y)=\phi^{-1}_2(y)$. Assume
there exists a smooth scheme $V$ with a trival action of
$\Gamma$ and a $\Gamma$-equivariant morphism $g:
X\to V$ such that  $\phi_i\times
g: (X,S)\to (Y,R) \times V$ are $\Gamma$-smooth. Define
$X':=g^{-1}(g(x)), S':=\{s\cap X'\mid s\in S\}$. 

Then $\phi_1$ and $\phi_2$ determine the
same orientation iff their restrictions $\phi'_i: (X',S')\to
(Y,R)$ do. 
\end{lemma}

\noindent{\bf Proof.}
We can assume that $g$ is $\Gamma$-\'etale by replacing, if necessary, $(Y,R)$ with
$(Y,R)\times {\bf A}^m$, where $m=\dim(X)-\dim(Y)-\dim(V)$,  
and $\phi_i: (X,S)\to (Y,R)$ and  $\phi_i\times
g: (X,S)\to (Y,R) \times V$ with $\Gamma$-smooth morphisms $\widetilde{\phi_i}:
(X,S)\to (Y,R)\times {\bf A}^m$ and $\Gamma$-\'etale morphisms $\widetilde{\phi_i}\times
g: (X,S)\to (Y,R)\times {\bf A}^m \times V$. 

Let $y=f(x)$ and
$y_1,\ldots,y_l$ be locally toric $\Gamma$-semiinvariant parameters at $y\in Y$. Let
$x_1,\ldots,x_k$ be locally toric parameters at  $p:=g(x)\in V$. Set
$v^1_i=\phi_1^*(y_i)$, $v^2_i=\phi_2^*(y_i)$ for
$i=1,\ldots,l$ and  
$w_j=g^*(x_j)$ for $j=1,\ldots,k$. Then
$v^1_1,\ldots,v^1_l,w_1,\ldots,w_k$ and
$v^2_1,\ldots,v^2_l,w_1,\ldots,w_k$ are locally toric parameters at $x$. The
automorphism $\alpha:=\widehat{(\phi_1\times
g)}_x^{-1}\circ \widehat{(\phi_2\times g)}_x: \widehat{X}_x\to
\widehat{X}_x$  maps the first set of
parameters onto the second one. Thus $\alpha$ belongs to the
group of automorphisms ${\rm
Aut}^\Gamma_{R}(\widehat{X}_x,S)$ preserving $R=K[[w_1,\ldots,w_k]]$.
The restriction of each automorphism from ${\rm
Aut}^\Gamma_{R}(\widehat{X}_x,S)$ to 
$X'={\rm Spec}(\widehat{\cO}_{X,x}/(w_1,\ldots,w_k))$ is an
automorphism. On the other hand we can write
$\widehat{\cO}_{X,x}=\widehat{\cO}_{X',x}[[w_1,\ldots,w_k]]$.
Hence each automorphism in ${\rm
Aut}^\Gamma_{K}(\widehat{X'_x},S')$ determines an automorphism in 
${\rm Aut}^\Gamma_{R}(\widehat{X}_x,S)$. We come to
a natural epimorphism of proalgebraic
 groups:
$${\rm res}: {\rm Aut}^\Gamma_{R}(\widehat{X}_x,S)\to 
{\rm Aut}^\Gamma_{K}(\widehat{X'_x},S')
.$$

The kernel of ${\rm res}$ is a proalgebraic
group $H$ 
 consisting of all automorphisms
$\beta\in 
{\rm Aut}^\Gamma_{R}(\widehat{X}_x,S)$ which can
be written in the form  $\beta(v^1_i)=v^1_i+r_i$, $\beta(w_j)=w_j$, where
$r_i\in (w_1,\ldots,w_k)\cdot\widehat{\cO}_{X,x}$. For any fixed $\beta$ and $t\in
{\bf A}^1$ yield
$r_i^t:=r_i(v^1_1,\ldots,v^1_l,t\cdot w_1,\ldots,t\cdot w_k)$ and
$\beta^t(v^1_i)=v^1_i+r^t_i$, $\beta^t(w_j)=w_j$. This
gives a morphism $t:{\bf A}^1\to H$ such that 

$t(1)=\beta^1=\beta$ and $t(0)={\rm id}$. Consequently, $H$ is 
connected and $$res^{-1}({\rm Aut}^\Gamma_{K}(\widehat{X'}_x,S')^0=
H\cdot({\rm Aut}^\Gamma_{K}(\widehat{X'}_x,S')^0=
{\rm Aut}^\Gamma_{R}(\widehat{X}_x,S)^0.$$ 
By the
above and Lemma \ref{le: group of automorphisms} the homomorphisms
$$ {\rm Aut}^\Gamma_{K}(\widehat{X}_x,S)
\buildrel\beta\over\longleftarrow 
{\rm Aut}^\Gamma_{R}(\widehat{X}_x,S)
\buildrel {\rm res}\over\longrightarrow 
 {\rm Aut}^\Gamma_{K}(\widehat{X'}_x,S')$$
define isomorphisms
$$ {\rm Aut}^\Gamma_{K}(\widehat{X}_x,S)/
{\rm Aut}^\Gamma_{K}(\widehat{X}_x,S)^0
\buildrel \overline{\beta}^{-1} \over \longrightarrow 
{\rm Aut}^\Gamma_{R}(\widehat{X}_x,S)/{\rm Aut}^\Gamma_{R}(\widehat{X}_x,S)^0
\buildrel\overline{\rm res}_1\over\longrightarrow 
{\rm Aut}^\Gamma_{K}(\widehat{X'}_x,S')/{\rm
Aut}^\Gamma_{K} (\widehat{X'}_x,S')^0.$$
Finally we see that $\widehat{(\phi_1\times
g)}_x\circ \widehat{(\phi_2\times g)}_x^{-1}\in {\rm Aut}^\Gamma_{K}(\widehat{X}_x,S)^0$ iff
$\widehat{\phi'_1}_x\circ \widehat{\phi'_2}_x^{-1}\in {\rm Aut}^\Gamma_{K}(\widehat{X'}_x,S')^0$
 \qed

\subsection{Orientation on stratified toroidal varieties}

\begin{definition} \label{de: orientation3} 
 Let $(X,S)$ be a 
$\Gamma$-stratified toroidal variety (or $\Gamma$-stratified toroidal scheme)
with an associated $\Gamma$-semicomplex
$\Sigma$. Let $\tau\leq \sigma$ and $x\in
s$. Let $\imath^\sigma_\tau: \tau\to\sigma$ denote the standard
embedding   
and $\phi_{\sigma}: U_{\sigma}\to X_{\sigma}$ be
a chart.  For any $\sigma'\preceq \sigma$ in $N_\sigma$ 
such that $\omega(\sigma')=\tau$ in $N_{\tau}$ we denote
by $\phi^{\sigma'}_{\sigma}:
U_\sigma^{\sigma'}:=\phi^{-1}_{\sigma} (X_{\sigma'})\to
X_{\sigma'}$  the restriction of $\phi_{\sigma}$ to $U_\sigma^{\sigma'}$. 

Let $\pi^\tau_{\sigma'}:X_{\sigma'}\to X_{\tau}$
denote the toric morphism induced by
any projection 
$\overline{\pi}^\tau_{\sigma'}: 
{\sigma'}\to \tau$ such that
$\overline{\pi}^\tau_{\sigma'}\circ \imath^{\sigma}_{\tau}= \id_{|\tau}$. 
\begin{enumerate}
\item
We say that the
$\Gamma$-stratified 
toroidal variety $(X,S)$ with atlas $\cU$ is {\it oriented}
if for any two charts $\phi_i: U_i\to X_{\sigma_i}$, where $i=1,2$, and 
 any $\sigma'_i\preceq\sigma_i$
such that $\omega(\sigma'_i)=\sigma \leq\sigma_i$
the $\Gamma_\sigma$-smooth morphisms 
$\pi^\sigma_{\sigma'_i}\phi^{\sigma'_i}_{\sigma_i}:
U_{\sigma_i}^{\sigma'_i}\to X_{\sigma}$ 
determine the same orientation at
any $x\in U_{\sigma_1}^{\sigma'_1}\cap U_{\sigma_2}^{\sigma'_2}
\cap \strat_X(\sigma)$. 

\item
Let $(X,S)$ be a $\Gamma$-stratified
toroidal scheme with  atlas $\cU$ which
contains a reduced subscheme $W$ of finite type over $K$.
We say that
$(X,S)$ is {\it oriented} along $W$  if for any two
 charts $\phi_i: U_i\to X_{\sigma_i}$, where $i=1,2$, and 
 any $\sigma'_i\preceq\sigma_i$
such that $\omega(\sigma'_i)=\sigma \leq\sigma_i$
the $\Gamma_\sigma$-smooth morphisms 
$\pi_{\sigma'_i}\phi_{\sigma_i}^{\sigma'_i}:
U_{\sigma_i}^{\sigma'_i}\to X_{\sigma}$ 
determine the same orientation at
any $x\in U_{\sigma_1}^{\sigma'_1}\cap U_{\sigma_2}^{\sigma'_2}
\cap \strat_X(\sigma) \cap W$.

\item We shall call such $W$ a {\it $K$-subscheme}.
\end{enumerate}
\end{definition}

The following lemma is a generalization of Lemma \ref{le: asemicomplex}
\begin{lemma} \label{le: asemicomplex2} Let $(X_\Sigma,S)$
be a $\Gamma$-stratified toric variety corresponding to an embedded
semifan $\Omega\subset \Sigma$. Then $(X_\Sigma,S)$ is an oriented
$\Gamma$-stratified toroidal variety with the associated
oriented $\Gamma$-semicomplex $\Omega^{\rm semic}$ and atlas 
$\cU^{\rm can}(\Omega,\Sigma)$. 
\qed
\end{lemma}

\begin{remark} By Lemma \ref{le: extensions} the above
definition does not depend upon the choice of the projection $\pi_\sigma$. 
\end{remark}

\begin{definition} Let $(X,S)$  be a $\Gamma$-stratified
toroidal variety  (respectively   a $\Gamma$-stratified
toroidal scheme with a $K$-subscheme $W$) with two atlases $\cU_1$ and $\cU_2$
such that $(X,S,\cU_1)$ and $(X,S,\cU_2)$ are oriented
(resp. oriented along W).
Then  $\cU_1$ and $\cU_2$ on
$(X,S)$ are
{\it compatible} (resp. {\it compatible along $W$}) if
$(X,S,\cU)$, where $\cU:=\cU_1\cup\cU_2$, is oriented
(resp. oriented along $W$). 
\end{definition}

\begin{lemma}\label{le: induce}
Let $f:(X,S)\to (Y,R) $ be a $\Gamma$-smooth morphism of $\Gamma$-stratified
toroidal schemes such that the strata in $S$ are preimages of
strata in $R$ and all strata in $R$ are dominated by strata
in $S$. Assume that $(Y,R,\cU)$ with associated
$\Gamma$-semicomplex $\Sigma$ is oriented along a
$K$-subscheme $W$. Define $f^*(\cU):= \{\phi f\mid \phi \in
\cU\}$. Then $(X,S,f^*(\cU))$ with the associated
$\Gamma$-semicomplex $\Sigma$ is oriented along any
$K$-subscheme $W'\subset f^{-1}(W)$. 
\end{lemma}
\noindent{\bf Proof.} Let $\phi:U\to X_\sigma$ be a chart
on $Y$. Then $\phi$ is $\Gamma$-equivariant,
$\Gamma_\sigma$-smooth morphism such that the intersections
$r\cap U$ are inverse images of strata of $X_\sigma$. Then
$\phi f:\phi^{-1}(U)\to X_\sigma$ has exactly the same
properties. Since $X_\sigma//\Gamma_\sigma$ exists it
follows from Lemma \ref{le: gsmoth} that $\phi f$ is $\Gamma_\sigma$-smooth.
\qed

\begin{lemma}\label{le: sections2} 
Let $\phi_i: (X,S)\to (Y,R)$  for $i=1,2$ be two 
$\Gamma$-smooth morphisms of $\Gamma$-stratified toroidal schemes such that 
the  strata in $S$ are preimages of
strata in $R$. Assume  $(Y,R,\cU)$ is oriented
and 
there exists  a $\Gamma$-equivariant morphism $g:
X\to V$ into a smooth scheme with trivial action of
$\Gamma$,  such that $\phi_i\times
g: (X,S)\to (Y,R) \times V$ are $\Gamma$-smooth. Set
$X':=g^{-1}(g(x)), S':=S\cap X'$. Let $W\subset X'$ be a $K$-subscheme.
Let $\psi_i: (X',S')\to
(Y,R)$  denote the restrictions of $\psi_i$.

Then $\phi^*_1(\cU) $ and $\phi^*_2(\cU)$  are compatible on
$(X,S)$ along $W$ iff 
$\psi^*_1(\cU) $ and $\psi^*_2(\cU)$  are compatible on
$(X',S')$ along $W$. 

\end{lemma}
\noindent{\bf Proof.} Follows from Lemma \ref{le: sections}.\qed

\begin{example} \label{ex: toroidal embeddings} 
{\it The group of automorphisms of the formal completion of a
toroidal embedding at a closed point is connected}.

The formal completion of a
toroidal embedding at a closed point is isomorphic to
$\widehat{X}_{\sigma\times \reg(\sigma)}$, where
$\reg(\sigma)$ is a regular cone and strata are defined by all
faces of $\sigma$. 
 Let $x_1,\ldots, x_k$ be semiinvariant coordinates on
$\widehat{X}_\sigma$ corresponding to characters
$m_1,\ldots,m_k$. Denote by $x_{k+1},\ldots,z_{k+l}$ the standard
coordinates on $X_{\reg(\sigma)}$ corresponfing to
characters $m_{k+1},\ldots,m_{k+l}$. 
Let $\phi$ be an automorphism  of 
$\widehat{X}_{\sigma\times\reg(\sigma)}$ which preserves the torus
orbit stratification. 
 Since all the divisors in the orbit
stratification  are preserved, the divisors
 of $x_i$, where $i=1,\ldots,k$, are also preserved, so we have 
 $\phi^*(x_i)=x_is_i$,
where
$s_i=a_{i0}+a_{i\alpha}x^\alpha+\ldots$ is invertible
and $\alpha=(\alpha_1,\ldots,\alpha_{k+l})$ 
denote  multiindices. 
By Lemma \ref{le: extensions},
 we may assume that $\phi^*(x_i)=x_i$ for $i=k+1,\ldots,l$.
An integral vector 
$v:=(v_1,\ldots,v_{l+k})\in\inte(\sigma\times \reg(\sigma))$ defines  
the $1$-parameter
subgroup $t\mapsto t^v$ of the ''big'' torus $T$. 
For any $t\in T$ let $t^v$ denote the automorphism defined as
$$x_i\mapsto x_i\cdot t^{\langle v,m_i\rangle}.$$ 
consider the morphism
$\phi_t:=t^{-v}\phi
t^v$.  In
particular $\phi_1=\phi$ and $$x_i\circ\phi_t=
x_i(a_{i0}+t^{\alpha_1m_1v_1+\ldots+\alpha_{l+k}m_{l+k}v_{l+k}}
a_{i\alpha}x^\alpha+\ldots),
$$ \noindent where $\alpha_1m_1v_1+\ldots+\alpha_{l+k}m_{l+k}v_{l+k}= 
\langle v, \alpha_1m_1+\ldots+\alpha_{l+k}m_{l+k} \rangle \geq
0$. Then $\phi_0$
is the well defined automorphism $x_i\mapsto a_{i0}x_i$. The latter
automorphism belongs to the  torus $T$. Since we can ''connect'' any
automorphism $\phi$ with an element of the torus   the
connectedness of the group of automorphisms  follows. 

 {\it Any toroidal embedding is an oriented
stratified variety}. 

By Lemma \ref{le: sections} and
by the connectedness of the automorphism group any two charts
determine the same orientation.

\end{example}

\begin{example} \label{ex: isolated} Let $X:=\{x\in {\bf A}^4\mid x_1x_2=x_3x_4\}$ be a
toric variety with isolated singularity. Then $X$ is a
stratified toroidal variety with the stratification
consisting of the singular point $p$ and of its complement.
The blow-up of the ideal of the point is a resolution of singularities
with the exceptional divisor determining the valuation. Consequently
the 
valuation determined by the point is preserved by all
automorphisms in $\Aut(\widehat{X}_x,S)$
. We
have the natural homomorphism
$d:\Aut(\widehat{X}_x,S)\to {\rm Gl}(\Tan_{X,x})$. The kernel of
$d$ is a connected proalgebraic group (see Lemma \ref{le:
group structure}(3)). The image
$d(\Aut(\widehat{X}_x,S))$ consists of linear automorphisms
preserving the ideal of $x_1x_2-x_3x_4$, that is, multiplying the
polynomial by a nonzero constant. Hence 
$d(\Aut(\widehat{X}_x,S))=K^*\cdot O$, where $K^*$ acts by
multiplying the coordinates by $t\in K^*$ and $O$ is the group of
linear automorphisms preserving $x_1x_2-x_3x_4$. By the linear
change of  coordinates $x_1=y_1-iy_2$, $x_2=y_1+iy_2$, 
$x_3=y_3-iy_4$, $x_4=y_3+iy_4$ we transform the polynomial into 
$y_1^2+y_2^2+y_3^2+y_4^2$ (${\rm char}(K)\neq 2$). This shows that $O$ is conjugate
to the group of orthogonal matrices with $K$-rational coefficients $O(4,K)$.
$O(4,K)$ consists of two components with $O(4,K)^0=SO(4,K)$.
Therefore $\Aut(\widehat{X}_x,S)=ker(d)\cdot K^*\cdot O$ consists
of two components. Fix any isomorphism $\phi:X\to
X_\sigma$, which can be considered as a chart. 
An oriented atlas consists of charts
compatible with $\phi$. Consequently, there are two
orientations corresponding to any two incompatible
isomorphisms $\phi:X\to
X_\sigma$.  
\end{example}

\subsection{Subdivisions of oriented semicomplexes}

\begin{definition} \label{de: oriented semicomplex} For any
$\Gamma$-semicone $\sigma$ denote by ${\rm
Aut}(\sigma)^0$  the group of all its automorphisms
which define the $\Gamma$-equivariant automorphisms of
${X}_{\sigma}$ preserving orientation.

By an {\it oriented semicomplex} (resp. {\it oriented
$\Gamma$-semicomplex}) we mean a semicomplex (resp. 
$\Gamma$-semicomplex) $\Sigma$ for
which for any $\sigma\leq\tau\leq\gamma$ there is 
$\alpha_\sigma\in {\rm Aut}(\sigma)^0$ for which 
$\imath^\gamma_{\tau}\imath^\tau_{\sigma}=\imath^{\gamma}_\sigma\alpha_\sigma$.
\end{definition}

\begin{definition}\label{de: semicomplexes isomorphism2}
By an {\it isomorphism} of two oriented $\Gamma$-semicomplexes $\Sigma\to
\Sigma'$ we mean a bijection of faces
$\Sigma\ni\sigma\mapsto\sigma'\in\Sigma'$ along with a collection of
face isomorphisms $j_\sigma:\sigma\to\sigma'$
such that
for any $\tau\leq \sigma$, there is an automorphism 
$\beta_\tau\in \Aut(\tau)^0$    such that 
$j_{\sigma}\imath^\sigma_{\tau}=\imath^{\sigma'}_{\tau'}j_{\tau}\beta_\tau$.
 
\end{definition}
\begin{remarks} 
\begin{enumerate}
 
\item An oriented semicomplex can be viewed as an oriented 
$\Gamma$-semicomplex
with trivial groups $\Gamma_\sigma$. 
\item 
We will show that the $\Gamma$-semicomplex associated
to an oriented $\Gamma$-stratified toroidal variety is oriented (Lemma
\ref{le: orientation4}).

\item The description of the group
$\Aut(\sigma)^0$ is given in Lemma \ref{le: orientation group}.
\end{enumerate}
\end{remarks}

\begin{lemma}
 Let $\Delta$ be a subdivision of a fan $\Sigma$ and
$\sigma\in \Sigma$. Then \\$\Delta{|\sigma}:=\{\delta\in
\Delta \mid \delta\subset\sigma\}$ is a subdivision of
${\sigma}$.\qed   
\end{lemma}

\bigskip
\noindent \begin{definition} \label{de: subdivision of
an oriented semicomplex} A {\it subdivision }  of an oriented
$\Gamma$-semicomplex $\Sigma$ is a collection $\Delta$ of fans 
$\Delta^\sigma$ in $N_\sigma$ where $\sigma\in \Sigma$
such that

\begin{enumerate}
\item For any $\sigma\in \Sigma$, $\Delta^\sigma$ is a
subdivision of $\overline{\sigma}$ which is
$\Aut(\sigma)^0$-invariant.  

\item For any $\tau\leq \sigma$,
$\Delta^{\sigma}|\overline{\tau}= \imath^{\sigma}_\tau(\Delta^\tau)$.

\end{enumerate}
\end{definition}

\bigskip

\begin{remarks}\begin{enumerate} 
\item By abuse of terminology we shall understand by
a subdivision of a semicone $\sigma$ a subdivision of the
fan $\overline{\sigma}$.
\item A subdivision  of an
oriented semicomplex can be viewed as a subdivision of 
an oriented $\Gamma$-semicomplex
with trivial groups $\Gamma_\sigma$. 

\item By definition vectors in faces
$\sigma$ of an oriented $\Gamma$-semicomplex $\Sigma$ are defined up
to automorphisms from $\Aut(\sigma)^0$. Consequently, the
faces of subdivisions $\Delta^\sigma$ are defined up
to automorphisms from $\Aut(\sigma)^0$ and in
general don't give a structure of semicomplex. Only
invariant faces of $\Delta^\sigma$ may define  a
semicomplex. 

\item The condition on $\Delta^\sigma$ to be
$\Aut(\sigma)^0$-invariant is for {\it canonical}
subdivisions replaced with a  somewhat stronger condition
of similar nature which says
that the induced morphism
$\widehat{X}_{\Delta^\sigma}:={X}_{\Delta^\sigma}\times_{{X}_{\sigma}}
\widehat{X}_{\sigma}\to \widehat{X}_{\sigma}$
is $\Aut^\Gamma(\widehat{X}_{\sigma})^0$-equivariant.
The latter condition can be translated into the 
condition that all "new" rays of the subdivision are in the {\it stable support} of $\Sigma$.  

\end{enumerate}
\end{remarks}

 \begin{definition}\label{le: starsub}
Let $\Delta$ be a subdivision of an oriented $\Gamma$-semicomplex
$\Sigma$. Let $\varrho_\sigma$ be an  ray passing through
the interior of the cone $\sigma\in\Sigma$ such that defining a collection of 
rays $\varrho:=\{\varrho_\tau\in\tau\mid\tau\geq\sigma\}$
such that for any $\tau\geq\sigma$, the ray
$\varrho_\tau=\imath^\tau_\sigma(\varrho_\sigma)$ is 
$\Aut(\tau)^0$-invariant . 
By the {\it star
subdivision} of $\Delta$ at $\rho$ we mean the subdivision
$$\varrho\cdot\Delta:=\{\varrho_\tau\cdot\Delta^\tau \mid \tau\geq\sigma\}\cup\{\Delta^\tau \mid 
\tau\not\geq\sigma\}.$$
\end{definition}

\bigskip
\subsection{Toroidal morphisms of stratified toroidal varieties}

\begin{definition}\label{de: toroidal modification}  Let $(X, S)$ be a
a stratified toroidal variety 
(resp. a $\Gamma$-stratified toroidal variety) 
with an atlas ${\cU}$. We say that $Y$ is a {\it toroidal
modification} of $(X, S)$ if

\begin{enumerate} 

\item There is given a proper morphism $f:Y\to X$ (resp.
proper $\Gamma$-equivariant morphism) such that 
for any $x\in s=\strat_X(\sigma)$ there exists a chart $x\in U_\sigma\to
X_{\sigma}$, a subdivision $\Delta^\sigma$ of $\sigma$,
and a fiber square  
\[\begin{array}{rcccccccc}
&& & U_\sigma & \buildrel \phi_\sigma \over\longrightarrow & X_{\sigma}&&&\\

&&&\uparrow f & & \uparrow  &&&\\

U_\sigma \times_{X_{\sigma}}
X_{\Delta^\sigma} &&\simeq& f^{-1}(U_\sigma)
 &\buildrel \phi^f_\sigma \over \longrightarrow &
X_{\Delta^\sigma} &&& \\

\end{array}\] 

\item
For any  point $ x$ in a stratum $s$   every
automorphism (resp. $\Gamma_s$-equivariant automorphism) 
$\alpha$ of $\widehat{X}_x$ preserving strata  and orientation
 can be lifted
to an automorphism (resp. $\Gamma_s$-equivaiant automorphism)
$\alpha'$ of $Y\times_X\widehat{X}_x$.
\end{enumerate}
\end{definition}

\begin{definition}
\label{de: canonical stratification}

Let $Y$ be a toroidal modification of a stratified toroidal
variety $(X,S)$.
By  a {\it canonical  stratification $R$}   on $Y$ we mean   
the finest stratification on $Y$ satisfying the conditions:
\begin{enumerate}

\item Strata of $R$ are mapped onto strata of $S$.

\item
For any chart $\phi_\sigma: U_\sigma\to X_\sigma$ on $(X,S)$ there is an
embedded subdivision
$\Omega^\sigma\subset\Delta^\sigma$ such that the strata of 
$ R\cap f^{-1}(U_\sigma)$ are defined as inverse images 
$(\phi^f_\sigma)^{-1}(\strat(\omega))$ of strata on
$(X_{\Delta^\sigma}, S_{\Omega^\sigma})$. 
\item
For any  point $x$ in a stratum $s$   every
$\Gamma_s$-equivariant automorphism $\alpha$ of
$\widehat{X}_x$ preserving strata and
orientation   
 can be lifted
to an $\Gamma_s$-equivariant automorphism $\alpha'$ of 
$Y\times_X\widehat{X}_x$ preserving strata. 

\end{enumerate}
\end{definition} 

If $Y$ is a toroidal modification of $(X,S)$ and $R$ is a
canonical stratification on $Y$ then 
we shall also speak of a {\it toroidal
morphism} $(Y,R)\to (X,S)$ of stratified toroidal varieties.

\begin{remarks} \begin{enumerate}
\item The definition of a toroidal morphism of 
stratified toroidal varieties is a generalization of
the definition of a toroidal morphism of toroidal embeddings.

\item The condition (2) of Definition 
\ref{de: toroidal modification} is similar to
Hironaka's ''allowability'' condition, used by  Mumford in the theory of toroidal
embeddings (see \cite{KKMS} and section 
\ref{se: toroidal embeddings}). It means that a toroidal morphism, which
by condition (1) is induced locally by charts, in fact
does not depend on the charts i.e.,
 if a morhism which satisfies
the Hironaka condition is induced locally by the diagram (1) for
some chart, then it is also induced locally by any other
chart and the given subdivision  
(see Lemma \ref{le: isomorphisms}).

\item In the above definition $\Omega^\sigma$ consists of those
faces of $\Delta^\sigma$ whose relative interior intersect the {\it
stable support} (see Lemma \ref{le: stab-subdivisions}).

\item Only some special subdivisions $\Delta^\sigma$ of $\sigma$
define (locally)
toroidal modifications (see example \ref{ex: example}). 
These subdivisions,  called {\it
canonical}, induce modifications which are independent
of the charts. They are characterized by the property that
"new" rays of the subdivision are contained in the stable
support (Proposition \ref{pr: simple}).

\end{enumerate}
\end{remarks}

\begin{example} \label{ex: example}
Let $X= {\bf A}^2$ and $S$ be the stratification on $X$ that
consists of the  point $s_0=(0,0)$ and its complement $s_1=
{\bf A}^2\setminus \{(0,0)\}$. Then $(X,S)$ is an oriented stratified
toroidal variety corresponding to the semicomplex $\Sigma=
\{\sigma_{0},\sigma_{1}\}$
where $\sigma_{0}=\langle(1,0),(0,1)\rangle$, $\sigma_{1}=\{(0,0)\}$.  

\begin{itemize}
\item Let $Y\to X$ be a toric morphism corresponding to
the normalized blow-up of $I=(x^2,y)$. This morphism of
toric varieties corresponds to the subdivision 
$\Delta:=\langle (1,2)\rangle\cdot\Sigma$ of $\Sigma$.
Therefore condition (1) of Definition \ref{de: toroidal
modification} is satisfied
for the identity chart.

The condition (2) of the definition is not satisfied. The
automorphism $\alpha$ defined by $\alpha(x)=y$, $\alpha(y)=x$
doesn't lift.

\item Let $Y\to X$ be a blow-up of $I=(x,y)$. Then $Y$ is a
toric variety and let $S_Y$ be
the stratification consisting of $s_1$ and the toric orbits
in the exceptional curve ${\bf P}^1$: $s_2=\{0\}$,
$s_3=\{\infty\}$ and $s_4={\bf
P}^1\setminus\{0\}\setminus\{\infty\}$. Then $Y$ is a
toroidal \underline{modification} of $(X,S)$. The
automorphism $\alpha$ defined by $\alpha(x)=y$, $\alpha(y)=x$
lifts to the automorphism permuting $s_2$ and $s_3$. Therefore
$(Y,S_Y)\to (X,S)$ is NOT a toroidal \underline{morphism}.

\item Let $Y\to X$ be a blow-up of $I=(x,y)$ and $S_Y$ be
the stratification consisting of the exceptional curve $s_2=D$
and its complement $s_1$. Then $(Y,S_Y)\to (X,S)$ is a
toroidal morphism corresponding to the 
subdivision $\Delta$ of $\Sigma$, where
$\Delta=\{\langle(1,0),(1,1)\rangle, \langle(1,1)\rangle,
\langle(1,1),(0,1)\rangle, \{(0,0)\}$. The stratification
$S_Y$ is described by the semicomplex
$\Omega=\{\langle (1,1) \rangle,\{(0,0)\}\}\subset
\Delta$. After blow-up some faces "disappear" from
$\Omega\neq \Delta$. The remaining faces of
$\Omega$ are $\Aut(\sigma)^0$-invariant.
\end{itemize}
\end{example}

\begin{lemma} \label{le: fibers} Let $ Y\to X$ be a
toroidal modification of an oriented
$\Gamma$-stratified toroidal variety $(X,S)$. Then all
fibers $f^{-1}(x)$, where $x\in s$, are isomorphic.
\end{lemma}
\noindent{\bf Proof}. It follows from  condition (1) of
Definition \ref{de: toroidal modification} that all fibers
$f^{-1}(x)$, where $x\in s=\strat_X(\sigma)\cap U_\sigma$, are
isomorphic for any chart $U_\sigma \to X_\sigma$. \qed

\begin{lemma} \label{le: dominate} Let $ Y\to X$ be a
toroidal modification of an oriented
$\Gamma$-stratified toroidal variety $(X,S)$. Then the exceptional divisors
of the  morphism $f: Y\to X$ dominate strata.
\end{lemma}
\noindent{\bf Proof.} Suppose an exceptional divisor $D$
does not dominate a stratum. The generic point of $f(D)$ is in a
stratum $s$ and $\dim(f(D)) < \dim(s)$. By Lemma \ref{le: fibers},
 $s$ is not a generic open
stratum in $X$. The dimension of a generic fiber $f^{-1}(x)$,
where $x\in f(D)$, is greater than or equal to 
$\dim(D)-\dim(f(D))=n-1-\dim(f(D)$. The dimension of a 
generic fiber $f^{-1}(x)$, where $x\in s$, is at most
$\dim(f^{-1}(s))-\dim(s)< n-1-\dim(f(D)$. This contradicts
 Lemma \ref{le: fibers}. \qed

\begin{definition}\label{de: vertices}
If $\Sigma$ is a fan then we denote by $\Ver(\Sigma)$  the set of all 1-dimensional
faces (rays) in  $\Sigma$.
\end{definition} 

\begin{lemma}
 Let $\sigma=\tau_1\oplus\tau_2$ and $\Delta$ 
be a subdivision of $\tau_1$. Then \\
$\Delta\oplus \tau_2:=\{\delta\oplus
\sigma 
\mid \delta\in \Delta,\,\, \sigma\preceq \tau_2 \}$ 
is a subdivision of $\tau$ called the \underline{joint}
of $\Delta$ and $\tau_2$. \qed
\end{lemma}

\begin{lemma} \label{le: sum} Let $\Omega\subset\Sigma$ be an
embedded $\Gamma$-semifan. Let $\Delta$ be  a subdivision of
$\Sigma$ such that any ray $\varrho\in\Ver(\Delta\setminus\Sigma)$ 
passes through the interior of 
$\omega\in\Omega$. Then for any cone
$\sigma=\omega(\sigma)\oplus {\rm r}(\sigma)$ in $\Sigma$, 
 $\Delta|{\sigma}=(\Delta|\omega(\sigma))\oplus {\rm r}(\sigma)$. 
\end{lemma}

\noindent{\bf Proof.} 
Since $\omega:=\omega(\sigma)$ is a maximal face of $\sigma$ 
which belongs to $\Omega$,  all ''new'' rays which are
contained in $\sigma$ are in
$\omega$ and therefore the maximal cones in
$\Delta|\sigma$ are of the form 
$\tau=\omega\oplus {\rm r}(\sigma)$ and $\tau\subset \omega$. 
This shows that $\Delta_1:=\Delta|{\sigma}$ is a subfan of 
$\Delta_2:=(\Delta|\omega)\oplus {\rm r}(\sigma)$. By properness we get
$|\Delta_1|=|\Delta_2|$ and hence $\Delta_1=\Delta_2$. \qed

\begin{lemma}. \label{le: iso} Let $\psi_1: Y_1\to X$ and
$\psi_2: Y_2\to X$ be proper
birational  morphisms of normal noetherian schemes.
\begin{enumerate}

\item Suppose $\widehat{X}_x\times_{X}Y_1$ and
$\widehat{X}_x\times_{X}Y_2$ are isomorphic over
$\widehat{X}_x$ for a
$K$-rational point $x\in X$. Then there exists an open
neighborhood $U$ of $x$ such that the proper birational
map $\psi_1^{-1}(U)\to \psi_2^{-1}(U)$ is an isomorphism. 

\item Suppose $\widehat{X}_x\times_{X}Y_1\to
\widehat{X}_x\times_{X}Y_2$ is a  proper  morphism over
$\widehat{X}_x$ for a
$K$-rational point $x\in X$. Then there exists an open
neighborhood $U$ of $x$ such that the proper birational
map $\psi_1^{-1}(U)\to \psi_2^{-1}(U)$ is a morphism.
\end{enumerate}
\end{lemma}
\noindent{\bf Proof} (1) Suppose that the sets 
$\psi_1^{-1}(U)$ and $\psi_2^{-1}(U)$
are not isomorphic for any $U$. Then the closed subset of $Y_2$ where
 $\psi_1\psi_2^{-1}$ is not an isomorphism interesect
the fiber $\psi_2^{-1}(x)$.

Let $Y_3$ be a component in $Y_1\times_{X}Y_2$ which
dominates $X$. Since all varieties are normal,   
at least one of the proper morphisms $Y_3\to Y_i$ is not an
isomorphism for $i=1,2$ over $\psi_i^{-1}(U)$ for any $U\ni x$. Then it
collapses a curve to a point over $x$.  

Let $\widehat{X}_3$ be a  component in 
$\widehat{X}_1\times_{\widehat{X}_x}\widehat{X}_2$ which
dominates $\widehat{X}_x$. The morphisms $\widehat{X}_3\to
\widehat{X}_i$ for $i=1,2$ are pull-backs of the morphisms
$Y_3\to Y_i$ via \'etale morphisms. This shows that they are
not both isomorhisms, and neither is the induced birational map 
$\widehat{X}_1\to\widehat{X}_2$, which contradicts the
assumption. 

(2) The same reasoning. \qed

\begin{lemma}\label{le: fiber square}
Let $\phi_\sigma:U_{\sigma}\to
X_{\sigma}$  denote a chart on an oriented
$\Gamma$-stratified toroidal variety $(X,S)$. 
Let $\sigma'\preceq\sigma$ denote a face for which  
$\omega(\sigma')=\varrho\leq\sigma$.

Consider the smooth morphism
$\psi_\sigma:=\pi^\varrho_{\sigma'}\phi^{\sigma'}_{\sigma} :
U_{\sigma'} \to X_\varrho$ induced by the chart $\phi_\sigma$.

Let $\Delta^\sigma$ be a subdivision of $\sigma$ and 
 $\Delta^\varrho:=\Delta^\sigma|\varrho$ be a subdivision of $\varrho$. 

Suppose 
${f}_\sigma:
\widetilde{U}_{\sigma}:=U_{\sigma}\times_{X_{\sigma}} X_{\Delta^\sigma}\to
U_{\sigma}$  satisfies  Condition (2) of Definition \ref{de: toroidal modification}.
Then the following diagram consists of
two fiber squares of smooth morphisms:

\[\begin{array}{rcccccccc}
&& & U_\sigma^{\sigma'} & \longrightarrow & X_{\sigma'}=X_{\varrho\oplus
{\rm r}(\sigma')}& \longrightarrow &X_{\varrho} &\\

&&&\uparrow f & & \uparrow  &&\uparrow&\\

U_\sigma:=U\times_{X_{\sigma'}} X_{\Delta^\sigma}
&&\simeq&  f^{-1}(U_\sigma^{\sigma'})
 & \rightarrow &
X_{\Delta^\sigma|\sigma'}=X_{\imath^\sigma_\varrho(\Delta^{\varrho})\oplus
{\rm r}(\sigma')} &
\rightarrow& X_{\Delta^\varrho}& \\

\end{array}\] 
\end{lemma}

\noindent{\bf Proof}
It follows from Lemma \ref{le: dominate} that
the assumption of Lemma \ref{le: sum} is satisfied for $\Delta^\sigma$.
By Lemma \ref{le: sum}, 
${\Delta^\sigma|\sigma'}=
{\Delta^{\sigma}|\varrho\oplus {\rm r}(\sigma')}=\imath^\sigma_\varrho(\Delta^{\varrho})\oplus {\rm r}(\sigma')$.
\qed

\begin{lemma}\label{le: isomorphisms}   
Let $\phi_\sigma:U_{\sigma}\to
X_{\sigma}$ and $\phi_\tau:U_{\tau}\to
X_{\tau}$ denote two charts on an oriented
$\Gamma$-stratified toroidal variety $(X,S)$. 
Let $\Delta^{\sigma}$ be a subdivision of $\sigma$ and $\Delta^{\tau}$ denote
a subdivision of $\tau$.
Assume  that for any $\varrho\leq \sigma, \tau$, 
$\Delta^\sigma|\varrho=\Delta^\tau|\varrho$. 
Suppose 
${f}_\sigma:
\widetilde{U}_{\sigma}:=U_{\sigma}\times_{X_{\sigma}} X_{\Delta^\sigma}\to
U_{\sigma}$ and ${f}_\tau: 
\widetilde{U}_{\tau}:=U_{\tau}\times_{X_{\tau}}X_{\Delta^\tau}\to
U_{\tau}$ satisfy Condition (2) of Definition \ref{de: toroidal modification}.

Then ${\phi}_\sigma^{-1}(U)$
and ${\phi}_\tau^{-1}(U)$
are isomorphic over $U:=U_{\sigma}\cap U_{\tau}$.

\end{lemma}
\noindent{\bf Proof.}
Suppose that the relevant
sets are not isomorphic over  $\Spec({\cO}_x)$ for $x\in
\strat_X(\varrho)\cap U$. Then there are faces
$\sigma'\preceq\sigma$ and $\tau'\preceq\tau$ such that 
${\omega(\sigma')}={\omega(\tau')=\varrho}$.

Set by $\reg(\varrho)=\langle e_1,\ldots,e_{\dim(\strat_X(\varrho))}\rangle$. 
The  diagram from Lemma \ref{le: fiber square} gives the following fiber square diagram of
\'etale extensions:
\[\begin{array}{rcccccccccc}
 & U & \longrightarrow &X_{\varrho\times reg(\varrho)}&\,&& & U & \longrightarrow &X_{\varrho\times reg(\varrho)}&\\

&\uparrow f_\sigma &&\uparrow&\,&&&\uparrow f_\tau &&\uparrow&\\

&f_\sigma^{-1}(U)&\longrightarrow&
X_{\Delta\times reg(\varrho)}&\,&& &f_\tau^{-1}(U)

 & \rightarrow &X_{\Delta\times \reg(\varrho)}& \\

\end{array}\]
\noindent where the horizontal morphisms are \'etale.

Consider the following  fiber square diagram: 

\[\begin{array}{rcccccc}
&\widehat{X}\times_U U_\sigma&\buildrel \alpha' \over 
\longrightarrow&
 \widehat{X}\times_U U_\tau &\buildrel
\widehat{\psi}'_{\sigma} \over \longrightarrow&
{X}_{\Delta^\varrho\times \reg(\varrho)} 
\times_{{X}_{\varrho\times \reg(\varrho)}} 
\widehat{X}_{\varrho\times \reg(\varrho)}\\
&\downarrow&&\downarrow&&\downarrow &\\

&\widehat{X}_x&\buildrel\alpha \over
\longrightarrow&\widehat{X}_x&\buildrel
\widehat{\psi}_{\sigma} \over \longrightarrow&  
\widehat{X}_{\varrho\times \reg(\varrho)},
\end{array}\]

\noindent where $\widehat{\psi}_{\sigma}$ and 
$\widehat{\psi}_{\tau}$  are isomorphisms of the completions of the local schemes
induced by \'etale extensions of the smooth morphisms
${\psi}_{\sigma}$  and ${\psi}_{\tau}$,
 $\alpha$ is an isomorphism such that 
$\widehat{\psi}_{\tau}:=\widehat{\psi}_{\sigma}\circ
\alpha$, and $\widehat{\psi}'_{\sigma}$,
$\widehat{\psi}'_\tau$ and $\alpha'$  are liftings 
of $\widehat{\psi}_{\sigma}$,
$\widehat{\psi}_\tau$ and $\alpha$.

 It follows from
the diagram that $\alpha$ is an automorphism of
$\widehat{X}_x$ which can be lifted to an automorphism
$\alpha'$ of $\widehat{X}_1$. Note that $\alpha$ 
 preserves strata and orientation. It follows that 
$\widehat{X}\times_U U_\sigma$ and $\widehat{X}\times_U U_\tau$.
are isomorphic.
By Lemma \ref{le: iso},  $\widetilde{U}_\sigma$ and
$\widetilde{U}_\tau$ are isomorphic over a neighborhoood of
$x\in X$, which
contradicts  the choice of $x\in X$. \qed

\bigskip
\subsection{Canonical subdivisions of oriented semicomplexes}\label{se: -semicomplexes}

As a corollary from the proof of Lemma \ref{le:
isomorphisms} we obtain the following lemma.
 \begin{lemma}\label{le: tilde} Let $(X,S)$ be an oriented 
$\Gamma$-stratified toroidal
variety of dimension $n$ with associated
oriented $\Gamma$-semicomplex $\Sigma$ and let $f:Y\to (X,S) $ be a toroidal modification.
Let $x$ be a
closed point in stratum $\strat_X(\sigma)\in S$, $\phi_\sigma: U\to
X_{\sigma}$ be a chart of a neighborhood $U$ of $x$
and $\Delta^\sigma$ be a subdivision of $\sigma$ for which
there is a fiber square 

\[\begin{array}{rcccccc}
&& & U & \buildrel \phi_\sigma \over \longrightarrow &X_{\sigma} &\\

&&&\uparrow f &&\uparrow&\\

&&&f^{-1}(U)

 & \buildrel \phi^f_\sigma \over \longrightarrow &X_{\Delta^\sigma}&, \\

\end{array}\]
\noindent where the horizontal morphisms are $\Gamma_\sigma$-smooth.
Set $$\reg(\sigma):=\langle
e_1,\ldots,e_{n-\dim(N_\sigma)}\rangle=\langle
e_1,\ldots,e_{\dim(\strat_X(\sigma))}\rangle.$$

Then 
\begin{enumerate}
\item
there is a fiber square of \'etale extensions
\[\begin{array}{rcccccc}
&& & U & \buildrel \widetilde{\phi}_\sigma \over \longrightarrow &X_{\widetilde{\sigma}} &\\

&&&\uparrow f &&\uparrow&\\

&&&f^{-1}(U)

&\buildrel \widetilde{\phi}^f_\sigma \over \longrightarrow 
&X_{\widetilde{\Delta}^\sigma},& \\

\end{array}\]
\noindent where the horizontal morphisms are $\Gamma_\sigma$-\'etale and where
$$\widetilde{\sigma} := \overline{\sigma}\times
\reg(\sigma),\quad\mbox{and} \quad \widetilde{\Delta}^\sigma:=\Delta^\sigma\times \reg(\sigma).
$$

\item $X_{\widetilde{\sigma}}$ is a
$\Gamma_\sigma$-stratified toric variety with the strata
described by embedded $\Gamma_\sigma$-semifan
$\sigma\subset \widetilde{\sigma}$. Moreover the strata on
$U$ are exactly inverse images of strata of $X_{\widetilde{\sigma}}$.

\item There is a fiber square of isomorphisms
\[\begin{array}{rcccccc}
&& & \widehat{X}_x & \buildrel \widehat{\phi}_\sigma \over \simeq &{\widetilde{X}_{\sigma}} &\\

&&&\uparrow \widehat{f}_x &&\uparrow&\\

&&&Y\times_X \widehat{X}_x

 & \buildrel \widehat{\phi}^f_\sigma \over \simeq &{\widetilde{X}_{\Delta^\sigma}} & \\

\end{array}\]
where

 $$\widetilde{X}_{\sigma}:=
\widehat{X}_{\widetilde{\sigma}},\quad \mbox{and} \quad
 \widetilde{X}_{\Delta^\sigma}:=
X_{\widetilde{\Delta^\sigma}}\times_{{X}_{
\widetilde{\sigma}}}
\widetilde{X}_{\sigma}.$$
\item $\widetilde{X}_{\sigma} $ is a
$\Gamma_\sigma$-stratified toroidal scheme with the strata 
described by embedded $\Gamma_\sigma$-semifan
$\sigma\subset \widetilde{\sigma}$. The isomorphism
$\widehat{\phi}_\sigma$ preserves strata.
\item
The morphism $\widehat{f}_x:Y\times_X \widehat{X}_x\to
\widehat{X}_x$ is
 ${\rm Aut}^\Gamma(\widehat{X}_x,S)^0$-equivariant 
.

\item
The morphism $\widetilde{X}_{\Delta^\sigma}\to
\widetilde{X}_{\sigma} $  is 
${\rm Aut}(\widetilde{X}_{\sigma})^0$-equivariant.
 \qed

\end{enumerate}
\end{lemma}
\begin{remark}
In what follows we assume that $n$ in the definition of $\widetilde{\sigma}$,
$\widetilde{\Delta}$, $\widetilde{X}_{\sigma}$, and
$\widetilde{X}_{\Delta}$ denotes the dimension of relevant 
stratified toroidal varieties. If the 
semicomplexes are not associated to stratified toroidal
varieties then $n$ will  denote some, a priori given,
natural number, satisfying  the  condition that
all simplices and semicomplexes have dimension
$\leq n$. Thus the definition of $\widetilde{\sigma}$,
$\widetilde{\Delta}$, $\widetilde{X}_{\sigma}$, and
$\widetilde{X}_{\Delta}$ will make sense for any 
cones
$\sigma$ and any fans $\Sigma$ considered.
\end{remark}

\begin{lemma} \label{le: isomorphisms2}
 Let $\Sigma$ be an oriented $\Gamma$-semicomplex and
$\sigma\in \Sigma$.
Let $\phi_1, \phi_2: (U, S)\to (X_{\sigma},S_{\sigma})$ be  smooth
morphisms of $\Gamma$-stratified toroidal schemes which determine the same
orientation at a closed point $x\in \strat_X(\sigma)$. Let $\Delta^\sigma$
be a subdivision of $\sigma$ such that $\widetilde{X}_{\Delta^\sigma}\to 
\widetilde{X}_{\sigma}$ is ${\rm
Aut}(\widetilde{X}_{\sigma},S_{\sigma})^0$-equivariant. 
Let
$\widetilde{U}_i=U\times_{X_{\sigma}}X_{\Delta^\sigma}$
 denote the fiber products for $\phi_i$, where $i=1,2$. Then
$\widetilde{U}_i$ are isomorphic over some neighborhood
$V\subset U$ of $x\in X$.
\end{lemma}
\noindent {\bf Proof.}  Identical with the proof of Lemma \ref{le:
isomorphisms}.\qed

We shall assign to the faces of an oriented
$\Gamma$-semicomplex $\Sigma$ the collection of
connected proalgebraic groups $$G^{\sigma}=
{\rm
Aut}(\widetilde{X}_{\sigma})^0.$$

\begin{definition} \label{de: canonical} A proper
subdivision $\Delta=\{\Delta^\sigma\mid\sigma\in\Sigma\}$ of an oriented $\Gamma$-semicomplex $\Sigma$ is
called {\it canonical}  if
 for any $\sigma \in \Sigma$, 
$G^{\sigma}$ acts
on $\widetilde{X}_{\Delta^\sigma}$ (as an abstract group) and the morphism
$\widetilde{X}_{\Delta^\sigma}\to
\widetilde{X}_{\sigma} $ is $G^{\sigma}$-equivariant.
\end{definition}

\bigskip
\begin{remark} This technical definition is replaced by
a simple combinatorial definition by using the notion of
stable support (see section  \ref{se: simple}). Also
$\widetilde{X}_{\sigma}$ and
$\widetilde{X}_{\Delta^\sigma}$ are replaced with 
$\widehat{X}_{\sigma}$ and
$\widehat{X}_{\Delta^\sigma}$ (Proposition \ref{pr: ssimple}). 
\end{remark}

\bigskip
\subsection{Correspondence between canonical subdivisions
of semicomplexes and toroidal modifications}

\begin{theorem}\label{th: modifications} Let  $(X,S)$ be
 an oriented $\Gamma$-stratified toroidal variety 
with the associated oriented $\Gamma$-semicomplex    
$\Sigma$. There exists a 1-1 correspondence
between the toroidal modifications $Y$ of  $(X,S)$ and 
the canonical subdivisions $\Delta$ of $\Sigma$. 
\begin{enumerate}
\item

If $\Delta$ is a canonical subdivision of $\Sigma$ then
the
toroidal modification associated to it is defined locally by 
\[\begin{array}{rcccccccc}
&& & U_\sigma &  \longrightarrow & X_{\sigma}&&&\\

&&&\uparrow f & & \uparrow  &&&\\

U_\sigma \times_{X_{\sigma}}
X_{\Delta^\sigma} &&\simeq& f^{-1}(U_\sigma)
 & \rightarrow &
X_{\Delta^\sigma} &&& \\

\end{array}\] 

\item
If $Y^1\to X$,
$Y^2\to X$ are toroidal modifications associated to
canonical subdivisions $\Delta_{1}$ and $\Delta_{2}$ of $\Sigma$ then 
the natural birational map $Y^1\to Y^2$ is a morphism iff
$\Delta_{1}$ is a subdivision of $\Delta_{2}$.
\end{enumerate}
\end{theorem}

\noindent{\bf Proof}. 
We shall  construct the canonical
subdivision $\Delta$ of $\Sigma$ associated to a given toroidal
modification $Y$. 

By 
Definition \ref{de:
toroidal modification} for any point $x$ there is a chart
and a subdivision $\Delta^\sigma$ of
$\sigma$ and the fiber square stated in  property (1). 
For two charts $\phi_{\sigma,1}: U_{\sigma,1}\to
X_\sigma$ and $\phi_{\sigma,2}: U_{\sigma,2}\to
X_\sigma$ we
may get two subdivisions $\Delta^\sigma_1$ and
$\Delta^\sigma_2$ of $\sigma$. 
These two charts induce by Lemma \ref{le: tilde}
isomorphisms $\widetilde{\phi}_{\sigma,i}: {\widehat{X}_x}\to
\widetilde{X}_\sigma$ and their liftings
$\widetilde{\phi}^f_i: Y\times_{X}{\widehat{X}_x}\simeq 
\widetilde{X}_{\Delta^\sigma_1}\simeq
\widetilde{X}_{\Delta^\sigma_2}$. The isomorphisms
$\widetilde{\phi}_{\sigma,i}$ define an automorphism
$\alpha:=\widetilde{\phi}_{\sigma,2}\widetilde{\phi}^{-1}_{\sigma,1}$
of $\widetilde{X}_\sigma$
preserving orientation and strata, and its  lifting  isomorphism
$\alpha_Y:=\widetilde{\phi}^f_{\sigma,2}({\widetilde{\phi}_{\sigma,1}}^{f})^{-1}:
\widetilde{X}_{\Delta^\sigma_1}\to
\widetilde{X}_{\Delta^\sigma_2}$. On the other hand
$\alpha$ lifts to an automorphism $\alpha^f$ of
$\widetilde{X}_{\Delta^\sigma_1}\simeq
Y\times_{X}{\widehat{X}_x}$. Finally the isomorphism
$\alpha'_Y\alpha_Y(\alpha^f)^{-1} :\widetilde{X}_{\Delta^\sigma_1}\to
\widetilde{X}_{\Delta^\sigma_2}$ is a lifting of the identity
morphism on $\widetilde{X}_\sigma$. Thus it is an identity
morphism on $\widetilde{X}_{\Delta^\sigma_1}=
\widetilde{X}_{\Delta^\sigma_2}$ and
$\Delta^\sigma_1 = \Delta^\sigma_2$. This implies 
the uniqueness of $\Delta^\sigma$. 
By  Lemma \ref{le: tilde}(6) , $i_x:\widetilde{X}_{\Delta^\sigma}\to
\widetilde{X}_{\sigma} $  is 
${\rm Aut}(\widetilde{X}_{\sigma})^0$-equivariant.
Thus the fan $\Delta^\sigma$ is $\Aut(\sigma)^0$-invariant.
It follows from the uniqueness that if $\tau\leq \sigma$
then
$\Delta^\sigma{|\tau}=\imath^\sigma_\tau(\Delta^\tau)\simeq
\Delta^\tau$. Therefore the
subdivisions $\Delta^\sigma$ define the canonical
subdivision $\Delta$ of $\Sigma$.  

Now let $\Delta$ be a canonical subdivision of the oriented
$\Gamma$-semicomplex $\Sigma$.
We shall construct 
the toroidal modification of $(X,S)$ associated to $\Delta$. 
  For any  chart $U_{\sigma,a}\to X_{\sigma}$ set 
$$\widetilde{U}_{\sigma,a}:=U_{\sigma,a}
\times_{X_{\sigma}}X_{\Delta^\sigma}$$

For any two charts $\phi_{\sigma,a}: U_{\sigma,a}\to X_{\sigma}$ and
$\phi_{\tau,b}:U_{\tau,b}\to X_{\tau}$ set 
$$\widetilde{U}_{\sigma,\tau,a,b} :=
(U_{\sigma,a}\cap U_{\tau,b})\times_{X_{\sigma}}X_{\Delta^\sigma},$$
\noindent where the fiber product is defined for the
restriction of $\phi_{\sigma,a}: U_{s,a}\to X_{\sigma}$ to 
$U_{\sigma,a}\cap U_{\tau,b}$. 

It follows from Lemma \ref{le: isomorphisms} that
 $\widetilde{U}_{\sigma,\tau,a,b}$ is
isomorphic over $U_{\sigma,a}\cap U_{\tau,b}$ to
$\widetilde{U}_{\tau,\sigma, b,a}$.
Thus we can  glue all the sets $\widetilde{U}_{\sigma,a}$
along $\widetilde{U}_{\sigma,\tau,a,b}$ and get a separated
variety $Y$. Note that the action of $\Gamma$ lifts to any 
subset $\widetilde{U}_{\sigma,a}$ and hence to an action
on the whole $Y$. 

A chart $\phi_\sigma:U_\sigma\to
X_{\sigma}$ defines  a $\Gamma_\sigma$-equivariant isomorphism  
$\widehat{X}_x\buildrel \alpha_1\over
\to\widetilde{X}_{\sigma}$.
For any $x\in\strat_X(\sigma)$ any  automorphism 
$\alpha$ of $\widehat{X}_x\simeq \widetilde{X}_{\sigma}$
preserving strata and orientation can be lifted
to a $\Gamma_\sigma$-equivariant automorphism $\alpha'$ 
of $\widehat{X}_y\times_X Y\simeq \widetilde{X}_{\Delta^\sigma}$.

(1)  Follows from the construction of a toroidal modification
from a canonical subdivision. 

(2) Let  $Y^1\to X$,
$Y^2\to X$ be toroidal modifications associated to
canonical subdivisions $\Delta_1$ and $\Delta_2$. 
If $\Delta_1$ is a subdivision of $\Delta_2$  then 
the natural birational map $Y^1\to Y^2$ is a morphism since
it is a morphism for each chart. On the other hand if 
$Y^1\to Y^2$ is a morphism then $Y^1\times_{X}{\widehat{X}_x}\simeq 
\widetilde{X}_{\Delta^\sigma_{1}}\to Y^2\times_{X}{\widehat{X}_x}\simeq 
\widetilde{X}_{\Delta^\sigma_{2}}$ is a toric morphism. By
Lemma \ref{le: iso} the latter morphism is induced by a 
toric morphism ${X}_{\Delta^\sigma_{1}}\to
X_{\Delta^\sigma_{2}}$. Hence  ${\Delta^\sigma_{1}}$  is a
subdivision of ${\Delta^\sigma_{2}}$ for each $\sigma \in \Sigma$
and consequently $\Delta_1$ is a subdivision of
$\Delta_2$. This completes the proof of Theorem \ref{th: modifications}. \qed

\section{Stable valuations}
For simplicity we consider
only valuations with integral values.
\subsection{Monomial valuations}
\begin{definition} 
Let $R$ be a local noetherian
ring, $u_1,\ldots,u_k$ be  generators of its maximal ideal,
and let $a_1,\ldots,a_k$ be nonnegative integers. For any $a\in
{\bf Z}$ set
$$J_a:=(u^{i_1,\ldots,i_k}\mid i_1a_1+\ldots+i_ka_k\geq a)\subset R,$$ 
\noindent  where ${i_1,\ldots,i_k}$, is a
multiindex in ${\bf Z}_{\geq 0}^k$.
We call a valuation $\nu$ {\it monomial} with  basis $u_1,\ldots,u_k$ 
and weights $a_1,\ldots,a_k$ if for any $f\in R\setminus\{0\}$, 
$$\nu(f)={\rm max}\{a\in {\bf Z}\mid f\in J_a\}.$$
\end{definition}

\begin{lemma} \label{le: monomial0} Let $\nu$ be a
nonnegative monomial valuation of  a local noetherian
ring $R$. Then for any $a\in {\bf Z}$, $I_{\nu,a}=J_a$. \qed
\end{lemma}

\begin{lemma} \label{le: monomial} Let $\nu$ be a
nonnegative valuation of  a local noetherian
ring $R$, 
and $u_1,\ldots,u_k$ be generators of the maximal
ideal of $R$. Let $\nu_0$ be a monomial valuation with  basis
$u_1,\ldots,u_k$ such that $\nu_0(u_i)=\nu(u_i)\geq 0$. Then 
$\nu_0(f)\leq \nu(f)$ for any
$f\in R$.  
\end{lemma}

\noindent{\bf Proof.} $I_{\nu_0,a}:=\{f\in R\mid 
\nu_0(f)\geq a\}=J_{a}\subset I_{\nu,a}.$
\qed   

\begin{lemma} \label{le: monomial2} Let $\nu$ be a
nonnegative monomial valuation  of  a local noetherian
ring $R$ with basis $u_1,\ldots,u_k$.
Let $g$ be an automorphism of $R$ such that
$\nu(g(u_i))\geq \nu(u_i)$ for any $i=1,\ldots,k$. 
Then $\nu(g(f))=\nu(u_i)$  for any $f\in R$.
\end{lemma}
\noindent{\bf Proof.} Note that
$\nu(g(u^{i_1,\ldots,i_k}))\geq
\nu(g(u^{i_1,\ldots,i_k}))$. By Lemma \ref{le: monomial0},
$g(I_{\nu,a})= g(J_a)\subset I_{\nu,a}$.
If $g(I_{\nu,a})\subset\neq I_{\nu,a}$ then we obtain an
infinite ascending chain of ideals $I_{\nu,a}\subset
g^{-1}(I_{\nu,a})\subset g^{-2}(I_{\nu,a})\subset\ldots$. \qed

\subsection{Toric and locally toric valuations}
 Let $X$ be an
algebraic variety and 
 $\nu$ be a valuation of the field $K(X)$ of rational functions. By the valuative criterion of
separatedness and properness the valuation ring of $\nu$
dominates the local ring of a uniquely determined  point (in general nonclosed)
$c_{\nu}$ on a complete variety $X$. (If $X$ is not
complete such a point may not exist). We call the
closure of $c_\nu$ the {\it center of the valuation } $\nu$ and
denote it by $\ce(\nu)$ or $\ce(\nu,X)$. For any $x\in \ce(\nu)$
and $a\in {\bf Z}_{\geq 0}$, 
$$I_{\nu,a,x}:=\{f\in {\cO}_{X,x}\mid 
\nu(f)\geq a\}$$ \noindent is an ideal in ${\cO}_{X,x}$. 
For  fixed $a$  these ideals define a coherent
sheaf of ideals ${\cI}_{\nu,a}$ supported at $\ce(\nu)$.  

\begin{lemma} \label{le: center}
Let $v$ be an integral vector in the support of the fan
$\Sigma$. Then the toric valuation $\val(v)$ is centered on
$\overline{O}_\sigma$, where $\sigma$ is the cone whose
relative interior contains $v$.
\end{lemma}
\noindent{\bf Proof.} Let $I_{O_\sigma}\subset K[X_\sigma]$
be the ideal of the orbit $O_\sigma$. Let $f=\sum a_ix^{m_i}$,
$m_i\in\sigma^\vee$, be a regular function on $X_\sigma$.
Then $\val(v)(f)>0$ iff $(v,m_i)>0$ for all $a_i\neq 0$.
But $(v,m_i)\geq 0$ by definition and $(v,m_i)=0$ iff
$m_i\in\sigma^\perp$. Thus we have $(v,m_i)> 0$ iff
$m_i\in\sigma^\vee\setminus \sigma^\perp$. Finally
$I_{O_\sigma}=\{f\in K[X_\sigma]\mid \val(v)(f)>0\}$ and the
valuation ring dominates the local ring of $O_\sigma$.\qed

Let $L$ be a ring containing $K$, $\sigma$ be a cone of the
maximal dimension in $N_\sigma$ and $v$ be an integral
vector of  $\sigma$.
Toric valuation $\val(v)$ can be defined on any ring
$L[[\sigma^\vee]]$ and its arbitrary subring $R$ containing
$L[\sigma^\vee]$. For any $f=\sum a_ix^{m_i}$, set
$$\val(v)(f)=\min\{(m_i,v)\mid a_i\neq 0\}.$$
 
The valuation $\val(v)$ of $R$ will be sometimes denoted by
$\val(v,R)$ or $\val(v,\Spec(R))$.

\begin{lemma} \label{le: invar} Let $\sigma$ be a cone of a
maximal dimension and $T_\sigma$ denote the
''big torus'' acting on $\widehat{X}_\sigma$. Let $\mu$ be a
$T_\sigma$-invariant valuation on $ \widehat{X}_\sigma$. Then
there exists an integral vector $v\in\sigma$ such that
$\val(v)=\mu$. 
\end{lemma}
\noindent{\bf Proof.} The valuation $\mu$ defines a linear
function on $\sigma^\vee$ corresponding to an integral vector $v\in\sigma$.
Since $\mu$ is $T_\sigma$ invariant so are the ideals
$I_{\mu,a}$. Therefore the ideals $I_{\mu,a}$ are generated
by characters $x^m$, where $(m,v)\geq a$. Consequently,
$I_{\mu,a}= I_{\val(v),a}$ and $\mu=\val(v)$.\qed 

\begin{lemma} \label{le: ind} Let $\sigma$ be a cone of the
maximal dimension in a lattice $N_\sigma$ and $v\in\sigma$ be its
integral vector. Let $R$ be a
subring of $K[\widehat{X}_\sigma]$ that contains
$K[{X}_\sigma]$. Then for any $a\in {\bf Z}$ we have 
$I_{\val(v,{X}_\sigma),a}\cdot R=I_{\val(v,R),a}$
\end{lemma}
\noindent{\bf Proof.} Denote by $R'$ the localization of $R$
at the maximal ideal of $O_\sigma$.
Then
$I_{\val(v,R'),a}=I_{\val(v,\widehat{X}_\sigma),a}\cap R'=
I_{\val(v,{X}_\sigma),a}\cdot K[\widehat{X}_\sigma]\cap R'$. 
Since $K[\widehat{X}_\sigma]$ is faithfully flat over $R'$
we obtain $I_{\val(v,{X}_\sigma),a}\cdot
K[\widehat{X}_\sigma]\cap R'= I_{\val(v,{X}_\sigma),a}\cdot
R'$. Thus we have $I_{\val(v,R'),a}=I_{\val(v,{X}_\sigma),a}\cdot
R'$. Then $I_{\val(v,R),a}=I_{\val(v,R'),a}\cap R= 
I_{\val(v,{X}_\sigma),a}\cdot
R' \cap R=I_{\val(v,{X}_\sigma),a}\cdot
R$.
\qed

\begin{lemma} \label{le: ind2} Let $X_\Sigma$ be a toric
variety with a toric action of $\Gamma$. Let $f: U\to X_\Sigma$ be a
smooth $\Gamma$-equivariant morphism. Let $v\in \inte(\sigma)$, where $\sigma\in
\Sigma$, be an integral vector. Assume that the inverse
image of $\overline{O}_\sigma$ is irreducible. Then there
exists a $\Gamma$-invariant valuation $\mu$ on $U$, such that $
{\cI}_{\mu,a}=f^{-1}({\cI}_{\val(v),a})\cdot{\cO}_{U}$.
\end{lemma}

\noindent{\bf Proof.} Let $x$ belong to the support $\supp
(f^{-1}({\cI}_{\val(v),a})\cdot{\cO}_{U})=f^{-1}(\overline{O}_\sigma)$. 
Then $f(x)\in O_\tau$, where $\tau\geq\sigma$. Let $\pi:
X_\tau\to X_{\underline{\tau}}$ denote the natural
projection. Then ${\cI}_{\val(v,X_\tau),a}=
\pi^{-1}({\cI}_{\val(v,X_{\underline{\tau}},a})\cdot{\cO}_{X_\tau}$.
Thus it suffices to prove the lemma replacing $X_\Sigma$
with $X_{\underline{\tau}}$ and $U$ with the inverse image
$U'$ of $X_\tau$. Let
$U_x\subset U'$ be an open neighborhood of $x$ 
for which there exists
an \'etale extension $\widetilde{f}: U_x\to
X_{\underline{\tau}\times \reg(\underline{\tau})}$.
It defines an isomorphism $\widehat{f}: \widehat{X}_x\simeq
\widehat X_{\underline{\tau}\times
\reg(\underline{\tau})}$. The toric valuation
$\val(v,X_{\underline{\tau}})$ defines the toric valuation $\val(v,X_{\underline{\tau}\times
\reg(\underline{\tau})})$ and the corresponding valuation
$\overline{\mu}_x$ on $\widehat{X}_x$. By Lemma \ref{le: ind} , $\overline{\mu}_x$ determines the valuation
${\mu}_x$ of ${\cO}_{X,x}$ such that 
$({\cI}_{\mu_x,a})_x=f^{-1}({\cI}_{\val(v),a})\cdot{\cO}_{X,x}$.
Thus $\mu_x$ is a valuation on $U_x$ supported on the
irreducible set $\overline{f^{-1}(O_\sigma)}$. Note that
both sheaves of ideals 
${\cI}_{\mu_x,a}$ and
$f^{-1}({\cI}_{\val(v),a})\cdot{\cO}_{X}$ are equal in  some
open neighborhood $V_x$ of $x$. Therefore the valuation $\mu_y$
is the same for all closed point $y\in
\overline{f^{-1}(O_\sigma)}\cap V_x$. Since
$\overline{f^{-1}(O_\sigma)}$ is irredducibe the valuations
$\mu_x$ are the
same for all $x\in \overline{f^{-1}(O_\sigma)}$ and define
a unique valuation $\mu$. Since the sheaves of ideals 
$f^{-1}({\cI}_{\val(v),a})\cdot{\cO}_{X}={\cI}_{\mu,a}$ are
$\Gamma$-invariant $\mu$ is also $\Gamma$-invariant.\qed

\bigskip
Let $f:X\rightarrow Y$ be a dominant morphism and $\nu$ be  a
valuation of $K(X)$. Then $f_*(\nu)$ denotes the valuation
which is the restriction
of $\nu$ to $K(Y)\simeq f^*(K(Y))\subset K(X)$. 

Let $\nu$
be a valuation on $X$ and ${\cI}_{\nu,a}$ be the associated sheaves of ideals.  
Let $f:Y\to X$ be a
morphism for which $f^{-1}({\cI}_{\nu,a})\cdot{\cO}_{Y}$
determine a  unique
valuation $\mu$ on $Y$ such that for sufficiently divisible
integer $a$,  
${\cI}_{\mu,a}=f^{-1}({\cI}_{\nu,a})\cdot{\cO}_{Y}$.  We
denote this valuation by $\mu=f^{*}(\nu)$.

\begin{definition} By a {\it locally toric valuation} on an $\Gamma$-toroidal
variety $X$ we mean a valuation $\nu$ of $K(X)$ such that for any
point $x\in \ce(\nu,X) $ there exists a
$\Gamma$-equivarinat $\Gamma_x$-smooth morphism $\phi: U\to
X_\sigma$ from an open $\Gamma$-invariant neighborhood $U$
of $x$ to a toric variety $X_\sigma$ with a toric action of
$\Gamma$ such that $\mu=\phi^*(\val(v))$, where $v$ is an
integral vector of $\sigma$. 

If $\mu$ is a locally toric valuation on $X$ then by 
$\mu_{|\widehat{X}_x}$ we denote the induced valuation on 
$\widehat{X}_x\simeq\widehat{X}_\sigma$.
 \end{definition}

\begin{definition} \label{de: blow}
By the
{\it blow-up} $\bl_{\nu}(X)$ of $X$ at a locally toric
valuation $\nu$ we mean the normalization of 
$$\Proj({\cO\oplus \cI}_{\nu,1}\oplus {\cI}_{\nu,2}\oplus\ldots),$$
\end{definition}

\begin{lemma} For any natural $l$, $\Proj({\cO\oplus \cI}_{\nu,1}\oplus
{\cI}_{\nu,2}\oplus\ldots)=\Proj({\cO\oplus
\cI}_{\nu,l}\oplus {\cI}_{\nu,2l}\oplus\ldots)$ 
.\qed
\end{lemma}

  Denote by $\bl_{\cJ}(X)\to X$ the
blow-up of  any coherent sheaf of
ideals ${\cJ}$.

\noindent \begin{lemma}\label{le: blow-up valuation} Let $X_\Sigma$ be a toric variety
associated to a fan $\Sigma$ in $N$ 
and $v\in |\Sigma|\cap N$. Then $\bl_{\val(v)}(X_\Sigma)$ is the toric
variety associated to the subdivision $\langle
v\rangle\cdot\Sigma$. Moreover for any sufficiently
divisible integer $d$, $\bl_{\val(v)}(X_\Sigma)$ is the normalization of
the blow-up of  ${\cI}_{\nu,d}$.
\end{lemma}

\noindent {\bf Proof.} It follows from  Definition \ref{de: blow}  that 
$\bl_{\val(v)}(X_\Sigma)$
is a normal toric variety. By Lemma \ref{le: center}, the sheaf of ideals
${\cI}_{\val(v),d}$ is supported on $\overline{O}_\sigma$,
where $v\in \inte(\sigma)$, $\sigma\in\Sigma$.  
By \cite{KKMS}, Ch.I, Th. 9, it determines a convex piecewise
linear function $f:=\ord_{{\cI}_{\val(v),d}}:|\Sigma|\to {\bf R}$ such that
that $f=0$ on any $\tau\in\Sigma$ that does not contain $v$, and 
 $f(u)={\rm min}\{ (m,u)\mid m \in \tau^{\vee}\cap M, (m,v) \geq
d\}$ for any $\tau$ containing $v$. Assume that $\tau$
contains $v$ and
let $\tau'$ be a face
of $\tau$ that does not contain $v$. Let $m_{\tau',\tau} \in \tau^\vee$ be a
nonnegative integral functional on $\tau$ which equals $0$
exactly on $\tau'$. Assume that $d$ divides all $(m_{\tau',\tau},v)$. 
Then $f_{|\tau'+ \langle v \rangle} =
\frac{d}{(m_{\tau',\tau},v)}\cdot m_{\tau',\tau}$.
Consequently, ${\cI}^k_{\val(v),d}={\cI}_{\val(v),kd}$ and
$\bl_{\val(v)}(X_\Sigma)=\bl_{{\cI}_{\val(v),d}}(X_\Sigma)$. 
By \cite{KKMS} Ch II Th. 10, the blow-up of ${\cI}_{\val(v),d}$ corresponds to
the minimal subdivision $\Sigma'$ of $\Sigma$ such that $f$ is 
linear on each cone in $\Sigma'$. By the above, 
$$\Sigma'=\Sigma \setminus {\rm Star}(\sigma,\Sigma)
\cup \bigcup_{\tau'\in \overline{{\rm
Star}(\sigma,\Sigma)}\setminus {{\rm
Star}(\sigma,\Sigma)}}
\tau'+ \langle v\rangle =
 \langle v \rangle\cdot\Sigma.$$ \qed

\begin{proposition}\label{pr: blow} For any locally toric
valuation $\nu$ on $X$ there exists an
integer $d$ such that 

$\bullet$ $\bl_\nu(X)=\bl_{{\cI}_{\nu,d}}(X)$. 

$\bullet$ The   valuation $\nu$ is induced by an irreducible $Q$-Cartier
divisor on $\bl_\nu(X)$.
\end{proposition}

\noindent{\bf Proof.} By quasicompactness of $X$ one can 
find a finite open affine covering $U_i$ of $X$ of charts
$\phi_i: U_i\to X_{\sigma_i}$ such that the valuation $\nu$
on each $U_i\subset X$ is equal to $\phi^*_i(\val(v_i)$. It
follows from Lemma \ref{le: blow-up valuation} that for any
$i$ there
exists $d_i$ that
$\bl_\nu(U_i)=\bl_{{\cI}_{\nu,d_i}}(U_i)$. It suffices to
take $d$ equal the product of all $d_i$. 
\qed

\bigskip
\subsection{Stable vectors and stable valuations}

\begin{definition} \label{de: stable valuations} Let
$(X,S)$ be an oriented 
$\Gamma$-stratified toroidal variety. 
 A  locally toric $\Gamma$-invariant valuation $\nu$
is {\it stable} on $X$
 if
\begin{enumerate}

\item $\ce(\nu,X)=\overline{s}$ for some stratum $s \in S$. 
\item For any $x\in \overline{s}$, $\nu_{|\widehat{X}_x}$ is invariant with
respect to any $\Gamma_s$-equivariant automorphism of
$(\widehat{X}_x,S)$, 
  preserving orientation.
\end{enumerate}
Let $f:Y\to (X,S)$ be a toroidal modification of an
oriented $\Gamma$-stratified toroidal variety. 
Then a valuation $\nu$ on $Y$ is called $X$-{\it stable}
if $f_*(\nu)$ is stable on $(X,S)$. 
\end{definition}

\begin{definition}  \label{de: stable vectors} Let  $\Sigma$  be an 
 oriented $\Gamma$-semicomplex.  
 A  vector $v_0\in {\rm
int}({\sigma}_0)$, where $\sigma_0\in \Sigma$,  
is $\Sigma$-{\it stable} on $\Sigma$ if for any $\tau\geq
\sigma_0$ the corresponding
valuation ${\rm val}(\imath^\tau_{\sigma_0}(v))$ on $\widetilde{X}_{\tau}$
is $G^\tau$-invariant. A vector $v\in \sigma$, where $\sigma\in\Sigma$ is $\Sigma$-{\it
stable} if there is a stable vector $v_0\in {\rm
int}({\sigma_{0}})$, where $\sigma_0\leq \sigma$, for which
$v=\imath^\sigma_{\sigma_0}(v_0)$.

\end{definition}

Let  $\Sigma$ be an oriented $\Gamma$-semicomplex. For any $\sigma \in \Sigma$ set
$${\rm stab}(\sigma):=
\{a\cdot v\mid  v\in\sigma, \,\,\,\, 
v \, \mbox{
\,is a
stable vector}\,\, , a\in {\bf Q}_{\geq 0}\}.$$

\begin{lemma}\label{le: glueing}

\begin{enumerate} 
\item Stable vectors $v\in\sigma$ are $\Aut(\sigma)^0$-invariant. 

\item 
For any $\varrho\geq \tau \geq \sigma$ and any vector $v\in {\rm
stab}(\sigma)$,
$\imath^\varrho_{\tau}\imath^\tau_{\sigma}(v)=\imath_\sigma^\varrho(v)$.

\item
$\imath_\sigma^{\tau}({\rm stab}(\sigma))
=\imath_\sigma^{\tau}(\sigma)\cap {\rm stab}(\tau)$ for any
$\tau\geq\sigma$.

\end{enumerate} 

\end{lemma}
\noindent{\bf Proof.}(1) Follows from definition.

(2) By definition there is $v_0\in
\inte(\sigma_{0})$ for which  $v=\imath^\sigma_{\sigma_0}(v_0)$. First we prove
that for any $\tau\geq \sigma$, we have
$\imath^\tau_{\sigma}\imath^\sigma_{\sigma_0}(v_0)=\imath^\tau_{\sigma_0}(v_0)$. 

By Definition \ref{de: oriented semicomplex} 
there is an automorphism $\alpha_{0}\in\Aut(\sigma_{0})^0$
which preserves the orientation
such that $\imath^\tau_{\sigma}\imath^\sigma_{\sigma_0}
=\imath^\tau_{\sigma_0}\alpha_{0}$.
By Definition \ref{de: stable vectors} the 
valuation $\val(v_0)$ is $G^{\sigma_0}$-invariant on
$\widetilde{X}_{\sigma_0}$ and hence $v_0$ is
preserved by $\alpha_{0}$. This gives 
$\imath^\tau_{\sigma}\imath^\sigma_{\sigma_0}(v_0)=\imath^\tau_{\sigma_0}\alpha_{0}(v_0)=
 \imath^\tau_{\sigma_0}(v_0)$
.
Now let $\varrho\geq\tau\geq\sigma$.
By the above we have
$\imath^\varrho_{\sigma}\imath^\sigma_{\sigma_0}(v_0)=\imath^\varrho_{\sigma_0}(v_0)$
Then $\imath^\varrho_{\tau}\imath^\tau_{\sigma}(v)=\imath^\varrho_{\tau}\imath^
\tau_{\sigma_0}(v_0)=\imath^\varrho_{\sigma_0}(v_0)$.

Finally
$\imath^\varrho_{\sigma_0}(v_0)=\imath^\varrho_{\sigma}\imath^{\sigma}_{\sigma_0}(v_0)
=\imath^\varrho_{\sigma}(v)$.

(3) 
If $v$ is a stable vector on
$\sigma$ then by Definition \ref{de: stable vectors}  
there is $\sigma_0\leq \sigma$ and a stable vector $v_0\in \inte (\sigma_{0})$ such
that $v=\imath^{\sigma}_{\sigma_0}(v_0)$. Then for $\tau\geq\sigma$,
$\imath^\tau_{\sigma_0}(v_0)=\imath^\tau_{\sigma}\imath^\sigma_{\sigma_0}(v_0)=
\imath^\tau_{\sigma}(v)$ is stable on
$\tau$. 
If $w\in \imath^\tau_{\sigma}(\sigma)\cap {\rm stab}(\tau)$ then 
there is $\sigma_0\leq \sigma$ and a stable vector $v_0\in \inte (\sigma_{0})$ such
that
$w=\imath^\tau_{\sigma_0}(v_0)=\imath^\tau_{\sigma}\imath^\sigma_{\sigma_0}(v_0)=
\imath^\tau_{\sigma}(v)$,
where $v=\imath^\sigma_{\sigma_0}(v_0)$.
 \qed

\begin{example}
Let $\Sigma$ be a regular semicomplex and 
$ \sigma=\langle e _1,\ldots, e_k\rangle \in \Sigma$. Then
$e:=e_1+\ldots+e_k$ is stable since it corresponds to the
valuation of the stratum $\strat_X(\sigma)$ on $\widetilde{X}_{\sigma}$. The subset $s$ defines a smooth
subvariety in a smooth neighborhood and therefore it
determines a valuation. If we blow-up the subvariety $\strat_X(\sigma)$  then this valuation
coincides with the one defined by the exceptional divisor.
\end{example}
\begin{example} \label{ex: affine2} Consider the  stratified
toroidal variety ${\bf A}^2$ from Example \ref{ex: affine}. The valuation of the point $(0,0)$ corresponds to the
vector $(1,1)\in N$, the valuation of $\overline{s} =
{\bf A}^1\times \{0\}$ corresponds to the vector $(0,1)\in N$. 
Then $|{\rm Stab(\Sigma)}|=\langle(1,1),(0,1)\rangle$. The valuations from
$\langle(1,0),(0,1)\rangle \setminus \langle(1,1),(0,1)\rangle$ are not stable: Let
$\widetilde{{\bf A}^2}\to {\bf A}^2$ be
the blow-up of the point $(0,0)\in {\bf A}^2$. The valuations from
$\langle(1,0),(0,1)\rangle \setminus \langle(1,1),(0,1)\rangle$
are centered on $\widetilde{{\bf A}^2}$ at $l_2\cap D$ or $l_2$ (using notation from
Example \ref{ex: affine}) and the automorphism
 $\phi:{\bf A}^2\to {\bf A}^2$, 
$\phi(x_1,x_2)=(x_1+x_2,x_2)$, preserves  strata
and moves the subvarieties $l_2$ and $l_2\cap D$. In particular it changes the relevant valuations. 
\end{example}

\begin{lemma} \label{le: stable valuations} 
 Let $(X,S)$ be a $\Gamma$-stratified toroidal
variety with the associated oriented $\Gamma$-semicomplex $\Sigma$. 
The following conditions are equivalent:
\begin{enumerate}
\item $\nu$ is stable on $X$.

\item There exists a stable vector $v\in {\rm
int}({\sigma})$, where $\sigma\in \Sigma$, 
such that for any $\tau\geq \sigma$  and any chart $\phi:U\to X_{\tau}$ 
we have 
$\nu = {\phi}^*({\rm {val}}(\imath^\tau_{\sigma}(v)))$.  
\end{enumerate}
\end{lemma}

\noindent{\bf Proof.} Let $x$ denote a closed point of 
$\overline{\strat_X(\sigma)}$. Then $x\in
\strat_X(\tau)$, where $\tau\geq \sigma$. 
The valuation $\nu$ on $X$ determines a $G^\tau$-invariant valuation on
$\widehat{X}_x\simeq\widetilde{X}_{\tau}$ which is 
in particular a toric valuation so, by Lemmas \ref{le:
invar} and \ref{le: center}
it corresponds to a unique $v\in \inte(\sigma)$. Thus
 $v$ is stable and ,by Lemma \ref{le: ind2}, $\nu = {\phi}^*({\rm {val}}(v))$.

Let $v$ be a stable vector. The valuation ${\rm val}(v)$  
determines a monomial  valuation 
$\nu:=\widehat{\phi}^*({\rm val}(v))$ on 
$\widehat{X}_x$, where $x\in \overline{\strat_X(\sigma)}\cap
U$. Since $\val(v)$ is
$G^\tau$-invariant on
$\widetilde{X}_{\tau}$ for any $\tau\geq \sigma$ 
we see that $\nu_{|\widehat{X}_x}$ is invariant with
respect to any $\Gamma_x$-equivariant automorphism
preserving strata and orientation. Moreover
$\nu_{|\widehat{X}_x}$ and $\nu$ do not depend on $\phi$.\qed

 \begin{lemma}\label{le: locally monomial} 
Let $f: Y\to X$ be a toroidal modification  of a $\Gamma$-stratified toroidal
variety $(X,S)$ associated to a canonical subdivision
$\Delta$ of $\Sigma$. 
The following conditions are equivalent:
\begin{enumerate}
\item $\nu$ is $X$-stable on $Y$.

\item There exists a $\Sigma$-stable vector $v\in {\rm
int}({\sigma})$, where  $\sigma\in\Sigma$, 
such that for any $\tau\geq \sigma$ and morphism $\phi^f_\tau:
f^{-1}(U_\tau) \to X_{\Delta^\tau}$, induced
by a chart $\phi_\tau:U_\tau\to X_{\tau}$ 
we have 
$\nu = \phi^{f*}_\tau({\rm {val}}(v))$.  
\end{enumerate}
\end{lemma} 

\noindent
{\bf Proof.} Let $ y\in \ce(\nu,Y)$, $x=f(y)\in X$ and
$\phi:U\to X_{\sigma}$ be a chart at $x$. We can extend
this $\Gamma^\sigma$-smooth morphism to an \'etale morphism 
$\widetilde{\phi}:U\to X_{\widetilde{\sigma}}$. 
By Lemma \ref{le: stable valuations} we have 
$\nu = \widetilde{\phi}^*({\rm {val}}(v))$ for some stable
vector $v\in\sigma$.  

Consider the fiber squares

\[\begin{array}{rccccccccc}

&f^{-1}(U) & \buildrel \widetilde{\phi}_f \over \longrightarrow  & 
X_{\widetilde{\Delta}}&&&&Y\times\widehat{X}_x
& \buildrel \widehat{\phi}_f \over \longrightarrow  & 
\widetilde{X}_{{\Delta}}\\

&\downarrow f&  &\downarrow i  
&&&&\downarrow \widehat{f}&  &\downarrow \widehat{i}\\

&{U} 
 &\buildrel \widetilde{\phi}\over \longrightarrow &
{X_{\widetilde{\sigma}}}&&&&\widehat{X}_x 
 &\buildrel \widehat{\phi}\over \longrightarrow &
\widetilde{X}_{\sigma}, 

\end{array}\] 

\noindent where $\widetilde{\phi}^f_\sigma$ is the induced \'etale morphism and $i$ is a
toric morphism. 
Then $\widehat{\phi}$ as well as $\widehat{\phi^f_\sigma}$ are
isomorphisms 
and  we have
$\nu=
\widehat{f}_*^{-1}\widehat{\phi}^*i_*({\rm val}(v))=
\widehat{\phi}_f^*({\rm val}(v))$ on $Y\times_X\widehat{X}_x$. 
Consequently
$\nu=
f_*^{-1}\widetilde{\phi}^*i_*({\rm val}(v))=
\widetilde{\phi}_f^*({\rm val}(v))$ 
since the above valuations coincide on $f^{-1}(U)$. \qed

\begin{lemma} \label{le: minimal vectors2} Let $\sigma$ and
$\tau$ be $\Gamma$-semicones in isomorphic lattices
$N_\sigma\simeq N_\tau$, and $L$ be a field
containing $K$.
Consider the induced  action of $\Gamma$ on
$L[[{\sigma}^\vee]]$ and $L[[{\tau}^\vee]]$ which is
trivial on $L$.
 Let $\psi: L[[\sigma^\vee]]\simeq L[[{\tau}^\vee]]$ be a 
$\Gamma$-equivariant isomorphism over $K$. 
  Let $v_\sigma\in \inte(\sigma) $ be a minimal internal vector of
$\sigma$. Then
there exists a minimal internal vector  $v_\tau\in {\rm
int}(\tau)$ such
that $\psi_*({\rm val}(v_\sigma))={\rm val}(v_\tau)$. 
\end{lemma}

\noindent{\bf Proof.} Note that ${\rm val}(v_\sigma)(f_\sigma)> 0$ and
$\psi_*({\rm val}(v_\sigma))(f_\tau)> 0$,  for 
any functions $f_\sigma\in m_\sigma$, $f_\tau\in m_\tau$ 
from the maximal ideals
$m_\sigma\subset L[[\sigma^\vee]]$ and $m_\tau\subset
L[[\tau^\vee]]$ respectively. 
The dual cone $\tau^{\vee}$ of regular characters  
 defines a subgroup in the
multiplicative group of rational
functions. The valuation $\psi_*({\rm val}(v_\sigma))$ determines
a group homomorphism of that subgroup into ${\bf Z}$ and
hence defines 
a linear function on $\tau^{\vee}$ corresponding to  a vector
$v_\tau \in \inte(\tau) $. We have to show that
$\val(v_\tau)=\psi_*({\rm val}(v_\sigma))$ (it is not clear whether 
$\psi_*({\rm val}(v_\sigma))$ is monomial with respect to characters). By
Lemma \ref{le: monomial}, $\psi_*({\rm val}(v_\sigma))\geq {\rm
val}(v_\tau)$. Let $v_\sigma'\in \inte(\sigma)$ correspond to the
linear function determined by  $\psi^{-1*}({\rm val}(v_\tau))$
 on the cone 
$\sigma^{\vee}$. Then 
${\rm val}(v_\sigma)=\psi^{-1*}\psi_*({\rm val}(v_\sigma))\geq
\psi^{-1*}({\rm val}(v_\tau))\geq
{\rm val}(v'_\sigma)$. Finally, ${\rm val}(v_\sigma)\geq {\rm
val}(v'_\sigma)$. By the
minimality of $v_\sigma$ it follows that $v_\sigma=v'_\sigma$ and 
${\rm val}(v_\sigma)=\psi^{-1*}({\rm val}(v_\tau))={\rm val}(v'_\sigma)$, or equivalently
$\psi^{*}({\rm val}(v_\sigma))={\rm val}(v_\tau)$. 

We need to show that $v_\tau$ is a minimal internal vector.
Suppose it isn't. Write $v_\tau=v_{\tau 0}+v_{\tau 1}$, where
$v_{\tau0}\in {\rm int(\tau)}$ and $v_{\tau 1}\in {\rm \tau}$.   
 By the
above find $v_{\sigma 0}\in {\rm int(\sigma)}$ 
such that 
${\rm val}(v_{\sigma 0})=\psi^{-1*}({\rm val}(v_{\tau 0}))$.
Then
 ${\rm val}(v_\sigma)=\psi^{-1*}({\rm val}(v_\tau)) \geq \psi^{-1*}({\rm
val}(v_{\tau 0}))= {\rm val}(v_{\sigma 0})$. By the minimality of $v_\sigma$, we
get $v_\sigma=v_{\sigma 0}$ and $v_\tau=v_{\tau 0}$. 
\qed

\begin{lemma} \label{le: Hilbert} Let $\sigma$ be a
$\Gamma$-semicone and
 $v\in \sigma$ be a vector such
that for any
$\phi\in G^{\sigma}$ there exists $v_\phi\in \sigma$ such that 
$\phi^*(\val(v))=\val(v_\phi)$
. Then $\val(v)$ is $G^{\sigma}$-invariant.
\end{lemma}
\noindent {\bf Proof.}
 Set $X:=\widetilde{X}_\sigma,\,\,  S:=S_{\widetilde{\sigma}} ,\,\,
x:= O_{\widetilde{\sigma}}$. Let $W$ denote the set of all
vectors ${v_\phi} \in \sigma$ for all $\phi\in {{\rm
Aut}}(X,S)$. For any natural $a$, the ideals 
$I_{\val(v_\phi),a}:=\{f\in\widehat{{\cO}}_{X,x}\mid
 \val(v_\phi)(f)\geq a\}$ are
generated by monomials. They define the same Hilbert-Samuel
function
$k\mapsto \dim_K(K[X]/(I+m^k))$, where $m
$ denotes the maximal ideal. It follows that the set $W$ is finite.

On the other hand
since the ideals $I_{\val(v_\phi),a}$
 are
generated by monomials  they can be distinguished by
 the ideals ${\rm gr}(I_{\val(v_\phi),a})$ in the graded ring
$${\rm gr}(\widehat{{\cO}}_{X,x_{\sigma}})=
\widehat{{\cO}}_{X,x}/m_{x,X}\oplus
m_{x,X_{\sigma}}/m_{x, X_{\sigma}}^2\oplus\ldots,$$ \noindent where
$m_{x,X}\subset \widehat{{\cO}}_{X,x}$ is the
maximal ideal  of the point $x$. 

Let  
$d:{{\rm Aut}}(\widetilde{X}_{\sigma},
S_{\sigma})^0\to\Gl(\Tan_{X,x})$ be the differential morphism (see Example 
\ref{ex: differential}). Then the connected algebraic group $d({{\rm Aut}}(\widetilde{X}_{\sigma},
S_{\sigma})^0)$ acts  algebraically on
the connected component of the Hilbert scheme of graded ideals with fixed Hilbert polynomial 
(see Example \ref{ex: differential}) . In
particular it acts trivially on the finite subset ${\rm gr}(I_{\val(v_\phi),a})$,
and consequently ${{\rm Aut}}(\widetilde{X}_{\sigma},
S_{\sigma})^0$
preserves all
$I_{\val(v_\phi),a}$.\qed

\begin{lemma} \label{le: divisors} Let $\sigma$ be a
$\Gamma$-semicone and
 $\widetilde{X}_{\Delta^\sigma}\to
\widetilde{X}_{\sigma}$ be a $G^{\sigma}$-equivariant morphism.
Then  all its
exceptional divisors are $G^{\sigma}$-invariant.
\end{lemma}
\noindent {\bf Proof.} Any automorphism $g\in G^{\sigma}$ maps an
exceptional divisor $D$ to another exceptional divisor
$D'$. Hence the valuation $\val_D$ defined by $D$ satisfies
 the condition of the previous lemma, and therefore it is $G^{\sigma}$-invariant.  
\qed

\begin{definition} Let $\sigma$ be a $\Gamma$-semicone.
A valuation $\val(v)$, where $v\in 
\sigma$, is {\it $G^{\sigma}$-semiinvariant} if for any $\phi\in 
G^{\sigma}$ there exist $v'\in \sigma$ such that $\phi^*(\val(v))=\val(v')$.

\end{definition}
By abuse of terminology a vector $v\in \sigma$ will be
called 
{\it $G^{\sigma}$-semiinvariant} (resp. {\it
$G^{\sigma}$-invariant}) if the corresponding valuation  
$\val(v)$ is $G^{\sigma}$-semiinvariant (resp. $G^{\sigma}$-invariant).

\begin{definition} Let $\Sigma$ be an oriented $\Gamma$-semicomplex.
 $v\in {\rm
int}(\sigma)$ is {\it semistable} if for any $\tau\geq
\sigma$,  $\val(v)$ is $G^\tau$-semiinvariant on $\widetilde{X}_\tau$.
\end{definition}

Let $\Ver(\Sigma)$ denote the set of all $1$-dimensional
rays in the fan $\Sigma$.

\begin{definition} Let $\sigma$ be a
$\Gamma$-semicone and
 $\widetilde{X}_{\Delta^\sigma}\to
\widetilde{X}_{\sigma}$ be a $G^{\sigma}$-equivariant
morphism. A cone $\delta\in \Delta^\sigma$ is {\it
$G^\sigma$-invariant}
if 
\begin{enumerate}
\item there is $\tau\leq\sigma$ such that 
$\inte(\delta)\subset \inte(\tau)$,

\item $\overline{O}_\delta\subset 
\widetilde{X}_{\Delta^\sigma}$ is $G^\sigma$-invariant. 
\end{enumerate}
\end{definition}

\begin{definition} Let 
$\Delta$ be a 
canonical subdivision of an oriented $\Gamma$-semicomplex
$\Sigma$. A cone  $\delta\in \Delta^\sigma$, $\sigma\in\Sigma$ is called a {\it $\Sigma$-stable
face} of $\Delta$ if 
\begin{enumerate}
\item 
$\inte(\delta)\subset \inte(\sigma)$,
\item $\imath^\tau_{\sigma}(\delta)\in \Delta^\tau$ is
$G^\tau$-invariant for any
$\tau\geq \sigma$.
\end{enumerate}
\end{definition}

\begin{lemma} \label{le: minimal vectorss} Let
 Let $\sigma$ be a
$\Gamma$-semicone and
 $\widetilde{X}_{\Delta^\sigma}\to
\widetilde{X}_{\sigma}$ be a $G^{\sigma}$-equivariant
morphism.
\begin{enumerate} 
\item All $G^{\sigma}$-semiinvariant vectors $v\in 
\sigma$ are $G^{\sigma}$-invariant.  
\item $\Ver(\Delta^\sigma)\setminus
\Ver(\overline{\sigma})$ are $G^{\sigma}$-invariant. 
\item Let $\delta$ be a $G^\sigma$-invariant face of
$\Delta^\sigma$. Then any minimal
internal vector $v\in \inte(\delta)$ are $G^{\sigma}$-invariant.
 \item Let $\delta$  be a 
\underline{$\Gamma$-indecomposable} face of $\Delta^\sigma$.
Then  $\delta$ is $G^\sigma$-invariant and all its 
 minimal internal vectors $v$ are $G^{\sigma}$-invariant.
In particular all minimal generators of
any face $\delta\in\Delta^\sigma$ are $G^{\sigma}$-invariant.

\item Let $\delta\in \Delta^\sigma$ be a minimal  face such that $
\inte(\delta)\subset \inte(\tau)$, where $\tau\neq
\{0\}$. Then $\delta$ is $G^{\sigma}$-invariant and
$\inte(\delta)$ contains a $G^{\sigma}$-invariant
vector.
\end{enumerate}
\end{lemma}
\noindent{\bf Proof.}
(1)  Follows from Lemma \ref{le: Hilbert}. 
(2) Follows from Lemma \ref{le: divisors}. 
(3) By Lemma \ref{le: minimal vectors2}, $v$ is
$G^{\sigma}$-semiinvariant. 
By (1), $v$ is $G^{\sigma}$-invariant. 

(4) Let $\tau\preceq\sigma$ be a face of $\sigma$ such
that $\inte(\delta)\subset\inte(\tau)$. By Lemmas \ref{le: sum}
and \ref{le: dominate}, $\Delta^\sigma|\tau=
\Delta^\sigma|\omega(\tau)\oplus r(\tau)$ for some
regular cone $r(\tau)$. Since $\delta$ is
$\Gamma$-indecomposable we conclude that $\delta\in 
\Delta^\sigma|\omega(\tau)$ and
$\tau=\omega(\tau)\leq \sigma$. 
By definition $G^\sigma$ acts on 
$\widetilde{X}_{\Delta^\sigma}$ and consequently on
strata from
 the stratification 
${\rm Sing}^\Gamma(\widetilde{X}_{\Delta^\sigma})$ (see
Lemma \ref{le: singularity type2}).

Let ${\rm sing_0}\in {\rm
Sing}(\widetilde{X}_{\Delta^\sigma})$ denote the stratum
corresponding to indecomposable face $\delta_0\in\Delta^\sigma$.

The images of the $G^\sigma$ action on ${\rm sing_0}$ form a finite invariant
subset of strata
${\rm Sing_0}\subset {\rm Sing}(\widetilde{X}_{\Delta^\sigma})$ 
. All strata from
${\rm Sing_0}$  correspond to some isomorphic
cones $\delta_i\in \Delta^\sigma$ for $i=0,\ldots,l$.
Let $v\in \inte(\delta)$ be 
a minimal internal vector. 
Then by Lemma \ref{le: minimal vectors2} for any 
$\widetilde{\phi}\in G^\sigma$ the image 
$\widetilde{\phi}_*({\rm val}(v))$
is  a valuation on $\widetilde{X}_{\Delta^\tau}$
corresponding to a minimal internal vector of $\delta_i$.
Hence $v$ is $G^\sigma$-semiinvariant and by (1) it is
$G^\sigma$-invariant. 
 By Lemma \ref{le: minimal0}  each minimal
generator  is a minimal internal vector of some
indecomposable cone $\delta$, hence by the above it is $G^\sigma$-invariant.

(5) For any $\tau\geq\sigma$ the cone $\delta$ corresponds
to  the maximal component 
$O_{{\delta}}$ in $\widetilde{X}_{{\Delta^\sigma}}$
which dominates $O_{{\tau}}\subset 
\widetilde{X}_{\sigma}$. Then
$G^{\sigma}$ permutes the set of maximal components of the
inverse image of $O_{{\tau}}$ 
dominating $O_{{\tau}}$ and we repeat the
reasoning from (4).
\qed 
\begin{remark}
The above lemma holds if we replace $G^{\sigma}$ with its
connected proalgebraic subgroup.
\end{remark}

As a direct corollary of Lemma \ref{le: minimal vectorss}
is the following
\begin{lemma} \label{le: minimal vectors} Let 
$\Delta$ be a 
canonical subdivision of an oriented $\Gamma$-semicomplex $\Sigma$. 
\begin{enumerate} 
\item   
 All the semistable vectors in $\Sigma$ are stable. 

\item For any $\sigma \in \Sigma$ all elements of $\Ver(\Delta^\sigma)\setminus
\Ver(\overline{\sigma})$ are $\Sigma$-stable.
\item Let $\delta$ be a $\Sigma$-stable face of $\Delta$. Then any minimal
internal vector $v\in \inte(\delta)$ is $\Sigma$-stable.  

\item Let $\delta$  be a 
\underline{$\Gamma$-indecomposable} face of $\Delta^\sigma$.
Then  $\delta$ is $\Sigma$-stable and all its 
 minimal internal vectors are stable. 
In particular all minimal generators of
any face
$\delta\in\Delta^\sigma$ are $\Sigma$-stable.

\item Let $\delta\in \Delta^\sigma$ be a minimal  face such that $
\inte(\delta)\subset \inte(\sigma)$, where $\sigma\neq
\{0\}$. Then $\delta$ is $\Sigma$-stable and
$\inte(\delta)$ contains a $\Sigma$-stable
vector. \qed
\end{enumerate}
\end{lemma}

\bigskip
\subsection{Star subdivisions and blow-ups of stable valuations}

\noindent
\begin{proposition} \label{pr: blow-ups}

\begin{enumerate}
\item Let $\Sigma$ be an oriented  $\Gamma$-semicomplex and 
$v_1,\ldots,v_k$ be $\Sigma$-stable vectors. Then 
$\Delta:=\langle v_k\rangle \cdot\ldots\cdot \langle v_1
\rangle\cdot\Sigma$ is a canonical subdivision of $\Sigma$.

\item Let $(X,S)$ be a  $\Gamma$-stratified toroidal
variety with an associated  oriented $\Gamma$-semicomplex
$\Sigma$. Let $v_1,\ldots,v_k$ be stable vectors
and $\nu_1,\ldots,\nu_k$ be the associated 
stable valuations. Then the composite of blow-ups 
$\bl_{\nu_k}\circ\ldots\circ \bl_{\nu_1}:
X_k\rightarrow X$ is the  toroidal modification of $X$
associated to the canonical subdivision
 $\Delta:=\langle v_k\rangle \cdot\ldots\cdot \langle v_1
\rangle\cdot\Sigma$ of $\Sigma$.
\end{enumerate}
\end{proposition}

\noindent{\bf Proof.} Induction on $k$. Let 
$\Delta_k:=\langle v_k\rangle \cdot\ldots\cdot \langle v_1
\rangle\cdot\Sigma$ be a canonical subdivision
corresponding to $f:=\bl_{\nu_k}\circ\ldots\circ \bl_{\nu_1}:
X_k\rightarrow X$. Then $v_{k+1}\in\inte(\sigma)$ defines an
$X$-stable valuation $\nu_{k+1}$ on $X_k$. It follows from \ref{le:
locally monomial} that for any chart $\phi_\sigma: U_\sigma\to X_{\sigma}$ the 
blow-up of $f^{-1}(U_\sigma)\subset X_j$ at $\nu_{k+1}=\phi^{f*}(\val(v_{k+1}))$ corresponds 
to the star subdivision of $\Delta^\sigma_{k}$ at $v_{k+1}$. \qed

\section{Correspondence between toroidal morphisms and
canonical subdivisions}

\subsection{Definition of stable support}

By Lemma \ref{le: glueing} we are in a position to glue pieces ${\rm
stab}(\sigma)$ into one topological space. 

\begin{definition}\label{de: stable support}
The {\it stable support}  of a $\Gamma$-semicomplex $\Sigma$  is the
topological space ${\rm Stab}(\Sigma) :=\bigcup_{\sigma\in
\Sigma} {\rm stab}(\sigma)$. 
\end{definition}

 Let 
$G\subset 
{\rm Aut}(\widetilde{X}_\sigma)$ be any abstract
algebraic group. 
Let $I_G$ denote the set of ${G}$-invariant toric
valuations on $\widetilde{X}_\sigma$. Set
$${\rm Inv}(G,\sigma):=\{a\cdot v\mid \val(v)\in I_G,\,\,
a\in {\bf Q}_{\geq 0}\},\quad {\rm Inv}(\sigma):={\rm Inv}(G^\sigma,\sigma).$$  

\begin{lemma}\label{le: stable support} Let $\Sigma$ be a $\Gamma$-semicomplex.
\begin{enumerate}
\item ${\rm stab(\sigma)}\cap \inte(\sigma)=\bigcap_{\tau\geq \sigma}{\rm
Inv}(\tau)\cap \inte(\sigma)$.

\item ${\rm stab(\sigma)}=
\bigcup_{\tau\leq \sigma} ({\rm stab}(\tau)\cap \inte(\tau))$.

\end{enumerate}
\end{lemma}
\noindent{\bf Proof.} Follows from the definitions of stable
valuation and stable support.\qed

\bigskip
\subsection{Convexity of the stable support}

\begin{lemma} \label{le: convex} 
\begin{enumerate} 

\item For any abstract subgroup 
$G\subset {\rm Aut}(\widehat{X}_\sigma)$ the cone
  ${\rm Inv}(G,\sigma)$ is  convex.

\item Let $\Sigma$ be a $\Gamma$-semicomplex. Then for any   
 $\sigma \in \Sigma$, ${\rm stab}(\sigma)$ is  convex.

\end{enumerate}
\end{lemma}

\noindent {\bf Proof.}
(1). Let $v_1,v_2\in I_G$ and 
 $\bl_{{\rm val}(v_1)}\circ \bl_{{\rm
val}(v_2)}: X_{\langle v_1\rangle\cdot \langle v_2
\rangle\cdot{{\sigma}}} \rightarrow
 {X}_{\sigma} $ be the toric morphism. The induced morphism 
$(\bl_{{\rm val}(v_1)}\circ \bl_{{\rm
val}(v_2)})\widehat{}: X_{\langle v_1\rangle\cdot \langle v_2
\rangle\cdot{{\sigma}}}\times_{X_{\sigma}}
\widehat{X}_{\sigma}\rightarrow
 \widehat{X}_{\sigma}$
is $G$-equivariant.
The exceptional divisors $D_1, D_2$ of
$\bl_{{\rm val}(v_1)}\circ \bl_{{\rm
val}(v_2)}$ correspond to $
v_1,v_2 \in \sigma$. The
$D_1$ and $D_2$ intersect along a stratum $O_{\tau}$
corresponding to the cone 
$\tau:=\langle v_1, v_2 \rangle$. 
Then $\widehat{X}_{\tau}$ 
is a ${G}$-invariant local scheme of toric variety at the
generic point of the orbit $O_{\tau}$
 and  $\{D_1,D_2,O_{\tau}\}$ is the  
orbit stratification on $(\widehat{X}_{\tau})$
which is also $G$-invariant.

Let $u_1,\ldots,u_k$ denote semiinvariant generators of 
$\widehat{{\cO}}_{X_{\tau},O_{\tau}} \simeq
K(O_{\tau})[[u_1,\ldots,u_k]]$
Each automorphism from $G$  preserves the orbit stratification and
multiplies the generating monomials $u_i$ by invertible
functions. Therefore it does not change the valuations
${\rm val}(v)$, where $ v\in\tau$,
on $\widehat{{\cO}}_{X_{\tau},O_{\tau}} \simeq
K(O_{\tau})[[u_1,\ldots,u_k]]$.

(2) Follows from (1) and from  Lemma \ref{le: stable support} 
\qed

\bigskip
\subsection{Toroidal embeddings and stable support}

\begin{lemma} \label{le: tembeddings} Let  $(X,S)$ be a toroidal
embedding. Then 
\begin{enumerate}

\item All integral
vectors in the faces $\sigma$ of $\Sigma$ are stable. 

\item
All subdivisions of $\Sigma$ are canonical. 
\end{enumerate}

\end{lemma}

\noindent{\bf Proof.} (1) Each automorphism $g$ from $G^{\sigma}$  
preserves divisors, hence
multiplies the generating monomials by invertible
functions. Consequently, it  does not change 
the valuations ${\rm val}(v)$, where $ v\in\sigma$,
on $\widetilde{X}_{\sigma}$. 

(2) Let $\Delta^\sigma$ be a
subdivision of $\sigma$. 
For any $\delta\in \Delta^\sigma$
each automorphism $g$ of $\widetilde{X}_{\sigma}$
lifts to an automorphism $g'$ of
$\widetilde{X}_{\delta}=\widetilde{X}_{\sigma}
\times_{X_{\widetilde{\sigma}}}
{X}_{\widetilde{\delta}}$ which also multiplies monomials
by suitable invertible functions. Therefore $g$ lifts
to the scheme $\widetilde{X}_{\Delta^\sigma}$.  \qed

\bigskip
\subsection{Minimal vectors and stable support}

\begin{lemma}\label{le: parr} 
 Let $\sigma$ be a
$\Gamma$-semicone and
 $\widetilde{X}_{\Delta^\sigma}\to
\widetilde{X}_{\sigma}$ be a $G^{\sigma}$-equivariant
morphism.

\begin{enumerate}

\item Let $\delta\in\Delta^\sigma$ be a simplicial  cone 
where $\sigma\in\Sigma$. Then all vectors in
${\rm par}(\delta)$ are $\Sigma$-stable.

\item
Let $\delta\in \Delta^\sigma$ be a $G^{\sigma}$-invariant simplicial face.
Then all vectors from
$\overline{{\rm par}(\delta)}\cap {\rm int(\delta)}$ are $G^{\sigma}$-invariant.

\end{enumerate}
\end{lemma}
\noindent{\bf Proof.} 
(1) Each $v\in {\rm par}(\delta)$ is a nonnegative
integral combination of
minimal generators of $\delta$. Minimal generators of
$\delta$ are $G^{\sigma}$-invariant by Lemma 
\ref{le: minimal vectorss}(4). Their linear combination is
$G^{\sigma}$-invariant by Lemma
\ref{le: convex}.

(2) Let $v\in\overline{{\rm par}(\delta)}\cap {\rm
int(\delta)}$. Write $\delta=\langle
w_1,\ldots,w_k\rangle$ and $v=\sum \alpha_iw_i$, where
$0\leq \alpha_i\leq 1$. If $v$ is a minimal internal vector
then it is $G^{\sigma}$-invariant by Lemma \ref{le: minimal vectorss}(3). If not
then $v=v'+v''$, where $v'=\sum \beta_iw_i$, $0<
\beta_i\leq \alpha_i\leq 1$, is a minimal internal vector and
$v''=\sum (\alpha_i-\beta_i)w_i \in {\rm par}(\delta)$. By
(1) and Lemma \ref{le: convex}, $v'$ is $G^{\sigma}$-invariant.  
\qed

 A direct corollary from the abo Lemma is the following Lemma.
\begin{lemma}\label{le: par} Let 
$\Delta$ be a 
canonical subdivision of an  oriented $\Gamma$-semicomplex
$\Sigma$.  
\begin{enumerate}

\item Let $\delta\in\Delta^\sigma$ be a simplicial  cone 
where $\sigma\in\Sigma$. Then all vectors in
${\rm par}(\delta)$ are $\Sigma$-stable.

\item
Let $\delta\in \Delta^\sigma$ be a $\Sigma$-stable simplicial cone.
Then all vectors from
$\overline{{\rm par}(\delta)}\cap {\rm int(\delta)}$ are $\Sigma$-stable.
\qed

\end{enumerate}
\end{lemma}

\subsection{Canonical stratification on a toroidal
modification} \label{se: canonical stratification}

\begin{lemma} \label{le: stab-subdivisions} 
Let $\Delta$ be a canonical subdivision of  an oriented 
$\Gamma$-semicomplex $\Sigma$. For any $\sigma\in\Sigma$ set $$
\Delta^\sigma_\stab:=\{\delta\in\Delta^\sigma\mid 
\inte(\delta)\cap {\rm stab}(\sigma)\neq\emptyset\}.$$ 
Then 
\begin{enumerate}
\item $\Delta^\sigma_\stab\subset\Delta^\sigma$ is an
embedded $\Gamma_\sigma$-semifan in $N_\sigma$.

\item  The subset $\Delta^\sigma_\stab\subset\Delta^\sigma$
consists of all $\Sigma$-stable faces in $\Delta^\sigma$.

\item For any $\tau\leq\sigma$,
$\Delta^\sigma_\stab|\imath^\sigma_\tau(\tau) = 
\imath^\sigma_\tau(\Delta^\tau_\stab)$.

\item  For any $\sigma\in\Sigma$, the
stratification $S^{\stab}_\sigma$ on
$\widetilde{X}_{\Delta^\sigma}$ determined by the embedded $\Gamma$-semifan 
$\Delta^\sigma_\stab\subset\Delta^\sigma$ is $G^{\sigma}$-invariant.

\end{enumerate}
\end{lemma}

\noindent{\bf Proof.}  
(1)  By Lemma \ref{le: convex},
 $\stab(\sigma)$ is convex and there is a
unique maximal face $\omega\preceq\delta$ whose relative
interior $\inte(\omega)$ intersects the stable support. 
Thus $\omega \in \Delta^\sigma_\stab$.
By Lemma 
\ref{le: minimal vectors}(4), either $\sing^\Gamma(\delta)=\{0\}$
or $\inte(\sing^\Gamma(\delta))\cap\stab(\sigma)\neq \emptyset $.
In both cases $\sing^\Gamma(\delta)\leq \omega$ and consequently
$\delta=\omega\oplus^\Gamma r(\delta)$, for some
regular cone $r(\delta)$.

(2) If $\delta\in\Delta^\sigma$ is $\Sigma$-stable then by Lemma \ref{le:
minimal vectors}(3) its minimal internal
vector is $\Sigma$-stable. Therefore $\delta\in\Delta^\sigma_\stab$.
If the relative interior of $\delta\in\Delta^\sigma_\stab$ 
contais a stable vector
$v\in \inte(\tau)$, where $\tau\leq\sigma$, then
$\inte(\delta)\subset\inte(\tau)$, 
$\delta\in\Delta^\tau_\stab$, and for any $\varrho\geq\tau$,
the closure of $O_{\imath^\varrho_\tau(\delta)}\subset
\widetilde{X}_{\Delta^\varrho}$ is exactly the center of $\val(v)$,
and therefore is $G^\varrho$-invariant.

(3) (4) follow from (2).

\qed

As a simple corollary of the above we get
\begin{lemma} There is a $\Gamma$-semicomplex $\Delta_\stab$
obtained by glueing the semicomplexes $\Delta^{\sigma\,\,\semic}_\stab$
along $\Delta^{\tau\,\,semic}_\stab$, where $\tau\leq\sigma$.\qed
\end{lemma}
\noindent{\bf Proof} For any
$\omega\in\Delta^{\sigma\,\,\semic}_\stab$ denote by
$\sigma(\omega\leq\sigma$ the semicone in $\Sigma$ for which
$\inte(\omega)\subset\inte(\sigma(\omega))$. Then 
$\omega\in\Delta^{\sigma(\omega)\,\,semic}_\stab$.
For $\omega\in \Delta^{\sigma(\omega)\,\,semic}_\stab$ and 
$\gamma\in \Delta^{\sigma(\gamma)\,\,semic}_\stab$,
write $\omega\leq \gamma$ if
$\sigma(\omega)\leq\sigma(\gamma)$ and
$\imath^{\sigma(\gamma)}_{\sigma(\omega)}(\omega)\leq
\gamma$. Then for $\omega\leq \gamma$ we set
$\imath^{\gamma}_{\omega}:=
\imath^{\sigma(\gamma)}_{\sigma(\omega)}$ and
$\Gamma_\omega:= (\Gamma_{\sigma(\omega)})_\omega$.
\qed
\begin{proposition} \label{pr: canonical stratification} 
Let $(X,S)$ be an oriented $\Gamma$-stratified
toroidal variety with an associated oriented $\Gamma$-semicomplex
$\Sigma$. Let $Y$ be a toroidal modification of $(X,S)$
corresponding to a canonical subdivision $\Delta$ of $\Sigma$.
Then there is a canonical stratification $R$ on $Y$ with  the
following properties:
\begin{enumerate}

 \item Let $\phi^f_\sigma: f^{-1}(U_\sigma)\to
{X}_{\Delta^\sigma}$ denote the morphism 
induced by a chart $\phi_\sigma:U_\sigma\to
{X}_{\sigma}$. The intersections $r\cap f^{-1}(U_\sigma)$,
$r\in R$  are precisely the inverse images of
the strata associated to the  embedded 
$\Gamma$-semifan $\Delta^\sigma_\stab \subset\Delta^\sigma$.

\item   The closures $\overline{r}$ of strata  $r\in R$ are  centers
of $X$-stable valuations  on $Y$.  

\item The  $\Gamma$-semicomplex  
associated to $(Y,R)$ is equal to 
$\Sigma_R=\Delta_\stab$ and the atlas is given by
$\bigcup_{\sigma\in\Sigma}\cU(\Delta^\sigma,\Delta^\sigma_\stab)$.  

\item The stratification $R$ is the finest stratification
on $Y$ satisfying the conditions:

\begin{enumerate}
\item The morphism $f$ maps strata in $R$ onto strata in $S$.
\item
There is an embedded semifan
$\Omega^\sigma\subset\Delta^\sigma$ such that the intersections 
$ r\cap f^{-1}(U_\sigma)$, $r\in R$, are inverse images of
strata of $S_{\Omega^\sigma})$ on
$(X_{\Delta^\sigma}$. 
\item
For any  point $ x$    every $\Gamma_x$-equivariant
automorphism $ \alpha$ of $\widehat{X}_x$
  preserving strata and
orientation   
 can be lifted
to an automorphism $\alpha'$ of $Y\times_X\widehat{X}_x
\to Y$ preserving strata. 
\end{enumerate}
\end{enumerate}
\end{proposition}

\noindent{\bf Proof.}
For any face $\omega\in \Delta^\sigma_\stab$ we find a
$\Sigma$-stable vector $v_\omega$ in its relative interior. By Lemma
\ref{le: stable valuations} the vector $v_\omega$
corresponds to an  $X$-stable valuation
$\nu_\omega$  on $Y$.
We  define the closure of a stratum $r\in R$
asociated to $\omega$ to be 
$$\overline{r}:=\overline{\strat_Y(\omega)}:=\ce(\nu_\omega).$$
\noindent Then define the stratum $r$ as
$$ r:=\strat_Y(\omega):=\overline{\strat_Y(\omega)}\setminus \bigcup_{\omega'< \omega}
\overline{\strat_Y(\omega')}.$$

By Lemma \ref{le: strata} the strata $\strat(\omega)$ 
of the stratification associated to the embedded semifan
$\Delta^\sigma_\stab \subset \Delta^\sigma$ satisfy  the
condition
$\overline{\strat(\omega)}=\ce(\val(v_\omega),X_{\Delta^\sigma})$, 
$\strat(\omega):=\overline{\strat(\omega)}\setminus \bigcup_{\omega'< \omega}
\overline{\strat(\omega')}$. It follows by the above that
the sets $r$ define stratification $R$ satisfying 
conditions (1), (2) and (3). 

(4) Since the strata of $R$ on $Y$ are defined by centers of
$X$-stable valuations,  conditions (a), (b) and (c) are satisfied.

Let $Q$ be a stratification on $Y$ satisfying
conditions (4a), (4b) and (4c). 

For the morphism $\phi^f_\sigma$ induced by a chart
$\phi_\sigma$ and for any $x\in
U_\sigma$, let 
$\widetilde{\phi}^f_\sigma: Y\times_X\widehat{X}_x\to 
\widetilde{X}_{\Delta^\sigma}$
be a $G^\sigma\simeq \Aut(\widehat{X}_x,S)^0$-equivariant
isomorphism of formal completions mapping the strata of $Q$
isomorphically onto the strata of
$\widetilde{S}_{\Omega^\sigma}$ on
$\widetilde{X}_{\Delta^\sigma}$ defined by the embedded
semifan $\Omega^\sigma\subset \widetilde{\Delta}^\sigma$

By (4a) the strata of $
\widetilde{S}_{\Omega^\sigma}$ 
on $\widetilde{X}_{\Delta^\sigma}$ are
$G^{\sigma}$-invariant. 
 Hence for any two isomorphisms $
\widetilde{\phi}^f_{\sigma,i}: Y\times_X\widehat{X}_x\to 
\widetilde{X}_{\Delta^\sigma}$, where $i=1,2$,
induced by 
charts $\phi_{\sigma,i}$  the induced  automorphism 
 $\widetilde{\phi}^f_{\sigma,1}(\widetilde{\phi}_{\sigma,2}^{f})^{-1}\in G^\sigma$
of
$\widetilde{X}_{\Delta^\sigma}$ maps
$\widetilde{S}_{\Omega_1^\sigma}$ to  $\widetilde{S}_{\Omega_2^\sigma}$.
Since both stratifications are $G^\sigma$-invariant we get
$\widetilde{S}_{\Omega_1^\sigma}= \widetilde{S}_{\Omega_2^\sigma}$,
$\Omega_1^\sigma=\Omega_2^\sigma$.  
Hence $\Omega^\sigma$ does not depend upon a chart.
Consequently, if $\tau\leq\sigma$ then
$\Omega^\sigma|\tau=\imath^\sigma_\tau(\Omega^\tau)$.
By condition (4a), the relative interior of a face
$\omega\in \Omega^\sigma$ is contained in the interior of
a face $\tau \leq \sigma$. Moreover the closure $O_{\imath^\varrho_\tau(\omega)}\subset
\widetilde{X}_{\Omega^\varrho}$ is $G^\varrho$-invariant
for any $\varrho\geq\tau$. This
shows that all faces $\omega\in \Omega^\sigma$  are
$\Sigma$-stable. Consequently, 
 $\Omega^\sigma\subset\Delta_\stab^\sigma$ and finally by Lemma \ref{le: strata2}.
the corresponding stratification $S_{\Delta^\sigma}$ is finer
than $S_{\Omega^\sigma}$. \qed 

\subsection{Correspondence between toroidal
 morphisms  and canonical subdivisions} 

Proposition \ref{pr: canonical stratification} can be
rephrased as follows:
\begin{theorem} \label{th: correspondence} Let  $(X,S)$ be an
oriented stratified toroidal variety
(resp. $\Gamma$-stratified toroidal variety) 
with the associated oriented semicomplex (resp. $\Gamma$-semicomplex) 
$\Sigma$.  There is a $1-1$ correspondence
between the toroidal morphisms of stratified 
(resp. $\Gamma$-stratified) toroidal
varieties $f:(Y,R)\to (X,S)$ and 
canonical  subdivisions $\Delta$ 
of the oriented semicomplex  
(resp. orineted $\Gamma$-semicomplex) $\Sigma$. 
Moreover the semicomplex (resp. $\Gamma$-semicomplex) 
associated to $(Y,R)$ is given by
$\Sigma_R= \Delta_{\stab}$. \qed
\end{theorem}

\begin{remark} 
 In particular, if $(X,S)$ is a
toroidal embedding (see Lemma \ref{le: tembeddings}) then
$\Sigma$ is a complex and all its subdivisions are
canonical (Lemma \ref{le: tembeddings}). We get a
1-1-correspondence between the subdivisions of the
complex $\Sigma$ and the toroidal morphisms $(Y,R)\to (X,S)$ (see \cite{KKMS}).
\end{remark}

\bigskip
\begin{example}\label{ex: isolated 3} Let $X$ be a variety with isolated
singularity of type $x_1x_2=x_3x_4$ with the
stratification consisting of the singular point and its
complement (as in Example \ref{ex: isolated}). Then the 
associated semicomplex $\Sigma$ consists of the  cone over a square and
its vertex. As  follows from Example \ref{ex: isolated4}
the stable support consists of the ray over the center of this
square. Consequently, by Proposition \ref{pr: simple} there are three canonical nontrivial
subdivisions of $\Sigma$. One is the star subdivision at
the stable ray corresponding to the blow-up of the point. Then
$\Delta_\stab$ consists of the ray and its vertex and
corresponds to the toroidal embedding defined by the
exceptional divisor and its complement. The other two are subdivisions defined
by diagonals corresponding to two small resolutions of
singularities. $\Delta_\stab$ consists of the cone over the
diagonal and its vertex and corresponds to the smooth
stratified toroidal variety with the stratification defined by the exceptional
curve-the preimage of the singular point and its complement.  
\end{example}

The language of stratified toroidal varieties allows a
combinatorial description of the \underline{Hironaka twist}. 

\bigskip
\begin{example}(\cite{Hironaka1})\label{ex: Hironaka}  Let $X={\bf P}^3$ be a projective
$3$-space containing  two curves $l_1$ and $l_2$ intersecting transversally in two
points $p_1$ and $p_2$. These data define a stratified
toroidal variety $(X,S)$. The associated semicomplex
$\Sigma$ consists of two
$3$-dimensional cones $\sigma_1$ and $\sigma_2$
corresponding respectively to the points $p_1$ and $p_2$,
and sharing two $2$-dimensional faces $\tau_1$ and $\tau_2$
corresponding to $l_1$ and $l_2$. The Hironaka` twist $Y$
is obtained by glueing  the consecutive blow-ups of
$X\setminus\{p_i\}$ at  $l_1$ and $l_2$ taken in two
different orders, along the isomorphic open subsets over
$X\setminus\{p_1\}\setminus\{p_2\}$. Then $Y$ is a
stratified toroidal variety. The preimage of $p_i$ consists
of two irreducible curves $l_{i1}$ and $l_{i2}$
intersecting at $p_{i0}$. The stratification $T$ on $Y$ is  determined by
the above four curves, two points and two exceptional
divisors-preimages of curves. Let $v_i$ denote
the sum of the two generators of $\tau_i$. Let $v^i$
denote the sum of the generators of $\sigma_i$.

$\Delta$ is the subdivision
of $\Sigma$ obtained by glueing  the consecutive star subdivisions of
$\sigma_i$ at  $\langle v_1\rangle $ and $\langle v_2\rangle $ taken in two
different orders, along the star subdivisions of $\tau_i$
at $\langle v_i\rangle$. 
It follows from Example \ref{ex: Hironaka2} that 
$\Stab(\Sigma)$ is the union of the cones 
$\langle v_1,v_2,v^1\rangle\subset \sigma_1$ and 
$\langle v_1,v_2,v^2\rangle\subset \sigma_2$.
Then $\Delta_{\stab}$ consists of the relevant
cones in $\Delta$ whose relative interior intersects
$\Stab(\Sigma)$. These cones correspond to the above
mentioned strata on $Y$.

\end{example}

\section{Canonical subdivisions and stable support}

\subsection{Isomorphisms of local rings and linear
transformations of stable supports}
\begin{lemma} \label{le: linear} 
\begin{enumerate} 

\item
Let $\sigma$ and $\tau$ be two cones of the same dimension.
 Let $\psi:\widehat{X}_{\sigma}\simeq 
\widehat{X}_{\tau}$
be an isomorphism  preserving the closures of the toric orbits. Then
there is  a linear isomorphism
$L_\psi:{\sigma}
\to \tau$ such that 
$\psi_*(\val(v))=\val(L_\psi(v))$  for any $v\in\sigma$.

\item Let $\sigma$ and $\tau$ be two $\Gamma$-semicones of the same
dimension. Let $\psi:\widetilde{X}_{\sigma}\simeq 
\widetilde{X}_{\tau}$ be a $\Gamma$-equivariant isomorphism 
preserving strata. Then there is 
a linear isomorphism
$L_\psi:{\Inv(\sigma)}
\to {\Inv(\sigma')}$ such that 
$\psi_*(\val(v))=\val(L_\psi(v))$  for any $v\in\Inv(\sigma)$. 
\end{enumerate}
\end{lemma}

\noindent{\bf Proof.}
  
(1)  The isomorphism $\psi$ maps the divisors of 
the orbit stratification on $\widetilde{X}_{\sigma}$
to the divisors of the orbit stratification on
$\widetilde{X}_{\tau}$.  Consequently, the local
semiinvariant parameters at $O_{\sigma}$ are mapped to local
parameters at $O_{\tau}$ which differ
from semiinvariant ones by invertible functions. This defines a
linear isomorphism of cones $L:{\sigma}
\to {\tau}$ mapping faces of
${\sigma}$ corresponding to
strata to suitable faces of ${\tau}$
 Let ${\psi}_L: \widehat{X}_\sigma\simeq 
\widehat{X}_{\tau}$
denote the induced morphism. Then $\phi:=\psi^{-1}_L\psi$ is an
automorphism of
$\widehat{X}_{\sigma}$ preserving strata.
 By Lemma \ref{le: tembeddings}
the automorphism $\phi$ preserves all valuations, hence
$\psi_*(\val(v))= \psi_{L*}(\val(v))=\val(L(v))$.

(2)   Let  $v_1, v_2\in \Inv(\sigma)$. Then
$\psi_*(\val(v_1))$ and $\psi_*(\val(v_2))$ are
$G^\tau$-invariant, and correspond to  
$v_1',v_2'\in \Inv(\tau)$. 
Set
$\Delta^\sigma=\langle v_2 \rangle\langle v_1 \rangle\cdot\sigma$
and $\Delta^\tau=\langle v'_2 \rangle\langle v'_1 \rangle\cdot\tau$.

The isomorphism $\psi:\widetilde{X}_{\sigma}\simeq 
\widetilde{X}_{\tau}$ lifts to an isomorphism 
$\Psi:
\widetilde{X}_{\Delta^\sigma}\simeq 
\widetilde{X}_{\Delta^\tau}$.  Let $D_1$ and $D_2$ (resp. $D'_1$ and $D'_2$) 
denote the exceptional divisors on
$\widetilde{X}_{\Delta^\sigma}$ (resp. $\widetilde{X}_{\Delta^\tau}$)
corresponding to $v_1$ and $v_2$ (resp. $v'_1$ and $v'_2$)
The generic point of intersection $p:=D_1\cap D_2$ is mapped 
to the point of intersection $p':=D'_1\cap D'_2$. Set 
$\delta_1:=\langle v_1, v_2 \rangle\in \Delta^\sigma$ and 
$\delta_2:=\langle v'_1, v'_2 \rangle\in \Delta^\tau$
. Then 
$(\widetilde{X}_{\delta_1},\{D_1,D_2,p\})$ is a formal completion of a toroidal
embedding which is mapped  to $(\widetilde{X}_{\delta_2},\{D'_1,D'_2,p'\})$.

 By (1) $\Psi_*$ determines a linear transformation on 
$\langle v_1, v_2 \rangle$ and consequently a linear
transformation  of $\Inv(\sigma)$ which is a linear isomorphism. 
 \qed

\begin{definition} Let $\sigma$ be a
$\Gamma$-semicone. A vector  $v_\sigma \in {\rm
int}( \widetilde{\sigma})$ is {\it absolutely invariant}
if $\val(v_\sigma)$ is
invariant with respect to ${\rm Aut}
(\widetilde{X}_{\sigma})$. 
\end{definition}
\begin{lemma} \label{le: invariant} Let $\sigma$ be a
$\Gamma$-semicone. There
exists an absolutely invariant  vector $v\in {\rm
int}( \widetilde{\sigma} )$. 
\end{lemma}

\noindent{\bf Proof.} We can assume that
$\sigma=\widetilde{\sigma}$ replacing $\sigma$, if
necessary, by the $\Gamma$-semicone
$\widetilde{\sigma}:=\sigma\cup |\widetilde{\sigma}|$.
Let $v\in {\rm
int}( \widetilde{\sigma} )$
be a minimal internal vector which  by Lemma \ref{le:
minimal vectors}(3) is a stable vector. 
By Lemma \ref{le: minimal vectors2} for any $g\in {\rm Aut}
(\widetilde{X}_{\sigma})$, $g_*(\val(v))=\val(v')$, where $v'$,
is also minimal. The set $W$ of minimal vectors  is contained in
$\overline{\rm par}(\sigma)$ and therefore it is finite. By
the above, ${\rm Aut}(\widetilde{X}_{\sigma})$ acts on $W$.
Let $v_\sigma$ denote the sum of all the minimal vectors. By
Lemma \ref{le: linear}, ${\rm
Aut}(\widetilde{X}_{\sigma})$ acts on $\Inv(\sigma)\supset W$ and therefore it
acts trivially on $v_\sigma$.\qed

\subsection{Initial terms defined by a monomial valuation}
\begin{lemma} \label{le: valuation} Let
$F(x_1,\ldots,x_k)\in K[x_1,\ldots,x_k]$
be a quasihomogeneous polynomial of weight $w(F)$ with
respect to weights
$w(x_1),\ldots,w(x_k)$. Let $\nu$ be a valuation of a ring
$R$ and $u_1,\ldots,u_k\in R$ be such that
$\nu(u_i)\geq w(x_i)$. Then $\nu(F(u_1,\ldots,u_k))\geq
w(F)$. If $F$ belongs to the ideal $ (x_l,\ldots,x_k)$ and $\nu(u_i)> w(x_i)$
for $i=l, \ldots, k$ then $\nu(F(u_1,\ldots,u_k)) >
w(F)$.
\end{lemma}
\noindent{\bf Proof.} The proof can be reduced to the
situation, in which case $F$ is a monomial when it immediately
follows from the definition of valuation. \qed

\begin{definition} Let $\nu$ be a monomial valuation of $K[[x_1,\ldots,x_k]]$
with a basis $x_1,\ldots,x_k$ and weights $a_1,\ldots,
a_k>0$. For  $f\in K[[x_1,\ldots,x_k]]$ let
$f=f_0+f_1+\ldots$ be a decomposition of $f$
into an infinite sum of quasihomogeneous polynomials such
that $\nu(f_0)<\nu(f_1)<\ldots$.
Then by the  {\it initial term} of  $f$ we mean
$\ini_\nu(f)=f_0$. Let $G=(G_1,\ldots,G_k)$ be an  endomorphism
of $K[[x_1,\ldots,x_k]$ such that $x_i\mapsto G_i$. Then
$$\ini_\nu(G):=(\ini_\nu(G_1),\ldots,\ini_\nu(G_k))$$
\noindent is the 
endomorphism determined by the initial terms of $G_i$.
\end{definition}  

\begin{lemma} \label{le: initial}
Let $G=(G_1,\ldots,G_k)$ be an automorphism of
$K[[x_1,\ldots,x_k]]$ preserving a monomial valuation
$\nu$: $G_*(\nu)=\nu$. 
Then $\ini_\nu(G)$
is an automorphism of $K[[x_1,\ldots,x_k]]$ preserving $\nu$.
Moreover for any $f\in K[[x_1,\ldots,x_k]]$ we have 
$\ini_\nu(G(f))=\ini_\nu(G)(\ini_\nu(f))$.
\end{lemma}
\noindent{\bf Proof.} Write
$G_i=\ini_\nu(G_{i})+G_{ih}$, where $\nu(G_{ih})>\nu(G_i)$. 

For any $f\in K[[x_1,\ldots,x_k]]$ write
$f=\ini_\nu(f)+f_{h}$,  where $\nu(f_{h})>\nu(f)$. Then we have
$G(f)=f(G_1,\ldots,G_n)=\ini_\nu f(G_1,\ldots,G_n)+f_{h}(G_1,\ldots,G_n)$,
 where by Lemma \ref{le: valuation}, 
$\nu(f_{h}(G_1,\ldots,G_n))>\nu(f)$.

$\ini_\nu(f)(x_1+y_1,\ldots,x_k+y_k)$ is a quasihomogeneous
polynomial in $x_1,\ldots,x_k,y_1,\ldots,y_k$ with weights
$\nu(x_1)=\nu(y_1) 
,\ldots,\nu(x_k)=\nu(y_k)$. Write
$\ini_\nu f(x_1+y_1,\ldots,x_k+y_k)= \ini_\nu(f)(x_1,\ldots,x_n)+f_h(x_1,y_1,\ldots,x_n,y_n)$ where
$f_h$ is a quasihomogeneous polynomial in the ideal
$(y_1,\ldots,y_n)\cdot K[[x_1,\ldots,x_k,y_1,\ldots,y_k]]$. 
Applying Lemma \ref{le: valuation} to $x_i=\ini_\nu(G_i)$ and ,
$y_i=G_{ih}$ we obtain  
\[\begin{array}{rc}
&(\ini_\nu(f)(G_1,\ldots,G_n))=
\ini_\nu(f)(\ini_\nu(G_1)+G_{1h},\ldots,\ini_\nu(G_n)+G_{nh}))=\\
&\ini_\nu(f)(\ini_\nu(G_1),\ldots,\ini_\nu(G_n)))+ 
f_h(\ini_\nu(G_1),\ldots,\ini_\nu(G_n),G_{1h},\ldots,G_{nh})),
\end{array}\]
\noindent where by Lemma \ref{le: valuation}, 
$\nu(f_h(\ini_\nu(G_1),\ldots,\ini_\nu(G_n),G_{1h},\ldots,G_{nh}))>\nu(f)$,
which gives $$\ini_\nu(G(f))=\ini_\nu(G)(\ini_\nu(f)).$$

Now let $G^{-1}=(G'_1,\ldots,G'_n)$. Then by the above,
$$ G(x_i)\circ
G^{-1}(x_i)=G'_{i}(G_{1},\ldots,G_{n})=x_i\quad \mbox{and}$$ 
$$\ini_\nu(G'_{i})(\ini_\nu(G_{1}),\ldots,\ini_\nu(G_{n}))=x_i.$$
This shows that $\ini_\nu(G)$ and $\ini_\nu(G^{-1})$ are automorphisms.  
They are
defined by quasihomogeneous polynomials of $\nu$-degrees 
$\nu(x_1),\ldots,\nu(x_k)$.
It follows from Lemma \ref{le: valuation} that 
$\nu(\ini_\nu(G(f)))\geq \nu(f)$
and $\nu(\ini_\nu(G^{-1}(f)))\geq \nu(f)$ for any $f\in
K[[x_1,\ldots,x_k]]$. Consequently, $\nu(\ini_\nu(G(f)))= \nu(f)$.\qed

\begin{lemma}\label{le: valuation2}
 Let $u_1,\ldots,u_n$ be local
semiinvariant parameters on  $\widetilde{X}_{\sigma}$. Write
$K[\widetilde{X}_{\sigma}]= K[[x_1,\ldots,x_k]]/I$. Let
$\nu$ be a toric valuation on 
$\widetilde{X}_{\sigma}$ invariant with respect to the group of  
all automorphisms of $\widetilde{X}_{\sigma}$ centered at 
$O_{\widetilde{\sigma}}$. 
Then $\nu$ defines a monomial valuation on
$K[[x_1,\ldots,x_k]]$ with basis $x_1,\ldots,x_n$ and
weights $\nu(u_1),\ldots,\nu(u_n)$.
Moreover:  
\begin{enumerate}
\item Any automorphism $g$ of
$K[\widetilde{X}_{\sigma}]$ extends to an automorphism $G$
of $K[[x_1,\ldots,x_k]]$ preserving $\nu$  and $I$.

\item The automorphism $\ini_\nu(G)$ of
$K[[x_1,\ldots,x_k]]$ defines an automorphism 
$\ini_\nu(g)$ of $K[\widetilde{X}_{\sigma}]$ preserving $\nu$.
\end{enumerate} 
\end{lemma}
\noindent{\bf Proof.} Any automorphism $g$ of
$K[\widetilde{X}_{\sigma}]$ sends parameters
$u_1,\ldots,u_n$ to $g_1=g^*(u_1),\ldots,g_n=g^*(u_n)$. Moreover
$\nu(g_i)= \nu(u_i)$ so we may assume that
$g_i=G_i(u_1,\ldots, u_n)$ are given
by some formal power series  $G_i\in K[[x_1,\ldots,x_k]]$
of $\nu$-weights  $\nu(u_i)$. Let $h=g^{-1}$, and let 
$h_1=h^*(u_1),\ldots,h_n=h^*(u_n)$  be given
by some formal power series  $H_i\in K[[x_1,\ldots,x_k]]$
of $\nu$-weights $\nu(u_i)$.
Then $\nu(H_i(G_1,\ldots,G_k))\geq \nu(u_i)$. On the other
hand
$H_i(G_1(u_1,\ldots,u_n),\ldots,G_k(u_1,\ldots,u_n))=u_i$.
Hence $(H_i(G_1,\ldots,G_k)-x_i)\in I$. In other words
$H_i(G_1,\ldots,G_k)=x_i+F_i(x_1,\ldots,x_k)$, where
$\nu(F_i)\geq \nu(u_i)$ and $F_i\in I$. Consequently, for
any $F\in K[[x_1,\ldots,x_k]]$, $\nu(F)\geq \nu(F(G_1,\ldots,G_k))$. 

Let $m\subset K[[x_1,\ldots,x_k]]$ be the
maximal ideal. By definition $m/(m^2+I)=m/m^2$ and
thus $I\subset m^2$. This gives $F_i\in m^2$.
Then the ring automorphism $\Phi:=H\circ G$ of 
$K[[x_1,\ldots,x_k]]$ can be written as a $K$-linear
transformation $\Phi=\id+F$, where
$F(K[[x_1,\ldots,x_k]])\subset I$ and $F(m^i)\subset m^{i+1}$.
It follows from Lemma \ref{le: valuation} that for any
formal power series $x\in
K[[x_1,\ldots,x_k]]$, $\nu(\Phi(x))\geq\nu(x)$ and hence
$\nu(F(x))\geq \min\{\nu(\phi(x)),\nu(x)\}=\nu(x)$.
The inverse of $\Phi$ can be written as
$\Phi^{-1}=\id+F'$, where $F'=-F+F^2-F^3+\ldots$. By the above 
$F'(K[[x_1,\ldots,x_k]])\subset I$ and 
$\nu(F'(x))\geq \nu(x)$, which shows that
$\nu(\Phi^{-1}(x))\geq \nu(x)$ and finally
$\nu(\Phi^{-1}(x))= \nu(x)$. 
Hence $G$ and $G^{-1}=\Phi^{-1}\circ H$ are automorphisms
satisfying $\nu(G^{-1}(x)\geq \nu(x)$.
By Lemma \ref{le: valuation} and the above,  $\nu(G(x))
\geq \nu(x)$ for any $x\in K[[x_1,\ldots,x_k]]$. 
Hence $\nu(G(x))=\nu(x)$.
We have shown that any automorphism $g$ of
$K[\widetilde{X}_{\sigma}]$ is defined by an automorphism
$G$ of
$K[[x_1,\ldots,x_k]]$ preserving the valuation $\nu$ and the ideal $I$.
Then by Lemma \ref{le: initial}, $\ini_\nu(G)$ is an
automorphism of $K[[x_1,\ldots,x_k]]$ preserving the valuation 
$\nu$ and the ideal $\ini_\nu(I)=I$. Therefore it defines an
automorphism $\ini_\nu(G)$ of $K[\widetilde{X}_{\sigma}]$. \qed

\subsection{Structure of the group of automorphisms of the
completion of a local ring of a toric variety}

\begin{lemma} 
 
\label{le: group structure} 
Let $\sigma$ be a $\Gamma$-semicone.
Let $\nu$ be the valuation on 
$\widetilde{X}_{\sigma}$ associated to the absolutely
invariant vector $v_\sigma$. Let $G^\sigma_{\nu}:={\rm Aut}_\nu
(\widetilde{X}_{\sigma})$ denote the group of all
quasihomogeneous automorphisms with respect to
semiinvariant coordinates on $\widetilde{X}_{\sigma}$
with  weights detrmined by $\nu$.  Then 

\begin{enumerate} 
\item $G^\sigma_{\nu}$ is an \underline{algebraic group}.

\item There is a
surjective morphism 
$\Phi:{\rm Aut}(\widetilde{X}_{\sigma})\to {\rm Aut}_\nu
(\widetilde{X}_{\sigma})$ defined by taking the lowest degree
quasihomogeneous part.

\item The  kernel ${\rm ker}_\sigma$ of $\Phi$ is a
connected proalgebraic group.

\item ${\rm Aut}(\widetilde{X}_{\sigma}) =G^\sigma_{\nu}\cdot {\rm ker}_\sigma$. 

\item For any $g\in {\rm ker}_\sigma$ there exists a morphism $i_{g}:{\bf A}^1\to {\rm ker}_\sigma$
such that $g= i_{g}(1)$, $\id=i_{g}(0)$. 

\item For any $g\in {\rm Aut}
(\widetilde{X}_{\sigma})^0$ there exists a morphism
$i_g: W :={\bf A}^1\times (G^{\sigma}_{\nu})^0\to 
{\rm Aut}(\widetilde{X}_{\sigma})^0$
such that $\{g,\id\}\subset i_g(W)$.

\end{enumerate}
\end{lemma}

\noindent {\bf Proof.} 
Let $u_1,\ldots,u_k$ be semiinvariant parameters on
$\widetilde{X}_{\sigma}$ generating the maximal ideal of
$O_{\widetilde{\sigma}}$. Denote by $m_1,\ldots,m_k\in 
{\widetilde{\sigma}}$ the
corresponding characters. 

(1)  Let  ${\rm Aut}_\nu
({{\bf A}^k})$ denote the group of all
quasihomogeneous automorphisms of ${\bf A}^k$. By Lemma \ref{le: valuation2}, ${\rm Aut}_\nu
(\widetilde{X}_{\sigma})$  is isomorphic to the subgroup of ${\rm Aut}_\nu
({{\bf A}^k})$ consisting of the automorphisms preserving
$I$.
The automorphisms of   ${\rm Aut}_\nu
({{\bf A}^k})$ 
are described by sets of $k$ quasihomogeneous polynomials with
weights $\nu(u_1),\ldots,\nu(u_k)$ and having linearly
independent linear terms. Thus  ${\rm Aut}_\nu
({{\bf A}^k})$ and ${\rm Aut}_\nu
(\widetilde{X}_{\sigma})$ are algebraic.  
   
(2) This follows from Lemma $\ref{le: valuation}$ and $\ref{le: valuation2}$

(3)(4)(5) The $1$-parameter
subgroup $t\mapsto t^{\nu_\sigma}$ defines the action  on 
$\widetilde{X}_{\sigma}$ such that
$\psi_t(u_i)=t^{(m_i,v_\sigma)}u_i=t^{\nu(u_i)}u_i=t^{\alpha_i}$,
where $\alpha_i=\nu(u_i)$ is a $\nu$-weight of $u_i$. 
The action defines an embedding $K^* < {\rm
Aut}(\widetilde{X}_{\sigma})$. 

Write $g\in {\rm ker}_\sigma$ in terms of coordinate functions,
$g=(g_1,\ldots,g_k)$, where $g_i=g^*(u_i)$. Let
$\alpha_{i,0}=\alpha_i$ denote the $\nu$-weight of $g_i$ and  let 
$g_i=u_i+g_{\alpha_{i,1}}+\ldots+g_{\alpha_{i,j}}+\ldots$
be the decomposition according to $\nu$-weights
$\alpha_{i,0}<\alpha_{i,1}<\ldots<\alpha_{i,j}<\ldots$.  Define the morphism
$i_g: K^* \to {\rm ker}_\sigma$ by $i_g(t):={\psi_t}^{-1}\circ g \circ
\psi_t$. Then $(i_g(t))^*(u_i)=({\psi_t}^{-1}\circ g \circ
\psi_t)^*(u_i)=(g \circ
\psi_t)^*t^{-\alpha_i}(u_i)=\psi_t^*(t^{-\alpha_i}g_i)=
t^{-\alpha_i}(t^{\alpha_{i,0}}u_i+t^{\alpha_{i,1}}g_{\alpha_{i,1}}+\ldots+
t^{\alpha_{i,j}}
g_{\alpha_{i,j}}+\ldots)=
=u_i +t^{\alpha_{i,1}-\alpha_{i,0}}
g_{\alpha_{i,1}}+\ldots+
t^{\alpha_{i,j}-\alpha_{i,0}}g_{\alpha_{i,j}}+\ldots$. 
The above morphism extends  to a  morphism
$i_g: {\bf A}^1\to {\rm ker}_\sigma$
''connecting'' $\id$ to $g$. This shows that ${\rm ker}_\sigma$ is
connected. Note that $\Psi$ maps the connected component ${\rm Aut}_\nu
(\widetilde{X}_{\sigma})^0$  to the
connected component  ${\rm
Aut}(\widetilde{X}_{\sigma})^0$. 

(6) Write $g\in {\rm Aut}
(\widetilde{X}_{\sigma},S)^0$ as $g=g_1g_2$, where $g_1\in
{\rm ker}_\sigma$, $g_2\in G^\sigma_{\nu}$. Then for $t\in {\bf A}^1$ and
$h\in G^\sigma_{\nu}$, set
$i_g(t,h):=i_{g_1}(t)\cdot h$. \qed

\begin{lemma}(\cite{Demushkin})\label{le: surjection})  
Let $\sigma$ be a $\Gamma$-semicone. 
Then ${\rm Aut}(\widetilde{X}_{\sigma})^0\subset 
{\rm Aut}(\widetilde{X}_{\sigma})$ is a normal subgroup and
there is a natural surjection 
${\rm Aut}(\sigma) \longrightarrow 
{\rm Aut}(\widetilde{X}_{\sigma})/{\rm Aut}(\widetilde{X}_{\sigma})^0.$
\end{lemma}

\noindent{\bf Proof.} 
Set $G:={\rm Aut}_{\nu}(\widetilde{X}_{\sigma})=
{\rm Aut}_{\nu}({X}_{\widetilde{\sigma}})$ 
(notation of Lemma \ref{le: group structure}). 
It
follows from  Lemma \ref{le: group structure} that the
natural inclusion $G\subset {\rm
Aut}(\widetilde{X}_\sigma)$ determines a group isomorphism  $G/G^0\simeq
{\rm Aut}(\widetilde{X}_\sigma)/ ({\rm Aut}(\widetilde{X}_\sigma))^0$. Then we
get a  surjection $N_G(T)/T\to G/G^0$. Let $H:={\rm
Aut}_T(X_{\sigma})$ be the group of the
$\Gamma$-equivariant automorphisms of $X_\sigma$,
preserving the big torus $T$. It
suffices to show that $N_G(T)=H$ and $H/T\simeq
\Aut(\sigma)$. If $g\in N_G(T)$ then $g{
T}g^{-1}=T$. Let $x\in T\subset X_{\widetilde{\sigma}}$. Then $g{
T}g^{-1}x=Tx$ or equivalently ${
T}g^{-1}x=g^{-1}Tx$. Since the latter subset is open,
$g^{-1}x$ is in the big open orbit ${
T}g^{-1}x={
T}x$, in other words, $g^{-1}x\in T$ which shows that
$g\in H$ and $N_G(T)\subset H$. 
Let $H^0$ be the connected component of $H$
containing $\id$. Then $H_0$ acts trivially on irreducible
components of the complement of $T$,that is, on ${T}$-invariant Weil and Cartier divisors. Hence it
multiplies characters by invertible functions on $X_{\widetilde{\sigma}}$.
Invertible functions on $T$ are monomials. Invertible
monomials on $X_{\widetilde{\sigma}}$ are constants. Thus $H_0$ acts on
characters multiplying them by nonzero constants. This
shows that $H_0=T$. Hence for  any $g\in H$, $g{T}g^{-1}=T$, which shows 
that $N_G(T)\supset H$.
$H/T$ can be identified with the subgroup of $H$
preserving $1\in T\subset X_{\widetilde{\sigma}}$, that is, the automorphisms
mapping characters to characters. The latter subgroup is
equal to the group ${\rm Aut}(\widetilde{\sigma})$ of
automorphisms of embedded semifan $\sigma\subset
\widetilde{\sigma}$.  

The restriction of an automorphism 
in ${\rm Aut}(\widetilde{\sigma})$ is an automorphism
of $\sigma$. Any automorphism of $\sigma$ defines an
automorphism of
$\widetilde{\sigma}=\sigma\times\reg(\sigma)$. Thus there
are natural group homomorphisms $i: {\rm
Aut}({\sigma})\to {\rm Aut}(\widetilde{\sigma})$  and 
$p:{\rm Aut}(\widetilde{\sigma})\to {\rm Aut}({\sigma})$
such that $ p\circ i=\id_{\sigma}$. 
Moreover $p$ corresponds to the restriction of 
automorphisms of $X_{\widetilde{\sigma}}$ to the subvariety 
 $X_{\sigma}\times O_{\reg(\sigma)}\subset X_{\widetilde{\sigma}}$. 
By Lemma \ref{le: sections}, $p({\rm
Aut}(\widetilde{\sigma})^0={\rm Aut}({\sigma})^0$.   
\qed   

\subsection{Invariant and semiinvariant valuations}
\begin{lemma} \label{le: semi}
Let $v\in\sigma$ be an integral vector. Then $\val(v)$
is  an ${\rm
Aut}(\widetilde{X}_{\sigma})^0$-invariant valuation on
$\widetilde{X}_{\sigma}$ iff  it is ${\rm
Aut}(\widetilde{X}_{\sigma})$-semiinvariant.
\end{lemma}
\noindent{\bf Proof.} $(\Leftarrow)$. If $\val(v)$ is ${\rm
Aut}(\widetilde{X}_{\sigma})$-semiinvariant then it is ${\rm
Aut}(\widetilde{X}_{\sigma})^0$-semiinvariant. Hence by
Lemma \ref{le: minimal vectorss}(1) it is ${\rm
Aut}(\widetilde{X}_{\sigma})^0$-invariant. 
$(\Rightarrow)$. If $\val(v)$ is  ${\rm
Aut}(\widetilde{X}_{\sigma})$-invariant then by Lemma
\ref{le: surjection}, ${\rm
Aut}(\widetilde{X}_{\sigma})$ acts on $\val(v)$ as the finite
group ${\rm Aut}(\sigma)\longrightarrow 
{\rm Aut}(\widetilde{X}_{\sigma})/ {\rm Aut}(\widetilde{X}_{\sigma})^0 $. Therefore $\val(v)$ is ${\rm
Aut}(\widetilde{X}_{\sigma})$-semiinvariant. \qed

\begin{lemma} \label{le: semi2} Let $L$ be any
algebraically closed field containing $K$.
Then $\val(v)$ defines ${\rm
Aut}(\widetilde{X}_{\sigma})^0$-invariant valuation on 
$\widetilde{X}_{\sigma}$ iff it defines an ${\rm
Aut}_L(\widetilde{X}_{\sigma})^0 $-invariant valuation on 
$\widetilde{X}^L_{\sigma}$.
\end{lemma} 
\noindent {\bf Proof.} Follows from Lemma \ref{le: base}(3). 
\qed

For a $\Gamma$-semicone $\sigma$ denote by
$\widehat{\Inv}(\sigma)$ the cone generated by
vectors $v\in |\sigma|$, for which valuations
$\val(v,\widehat{X}_{\sigma})$ are invariant with
respect to ${\rm
Aut}(\widehat{X}_{\sigma})^0.$

\begin{lemma}\label{le: hat}
For a $\Gamma$-semicone $\sigma$  the two  cones are equal
$\widehat{\Inv}(\sigma)={\Inv}(\sigma)$.
\end{lemma}
\noindent {\bf Proof.} Note that
$K[\widetilde{X}_{\sigma}]=
L[\widehat{X}_{\sigma}]$ for suitable
$L:=K[[x_1,\ldots,x_k]]$. Let $(L)$ denote the algebraic
closere of the quotient
field of $L$.
Each automorphism of $K[[\sigma^\vee]]$ defines the
automorphism of $L[[\sigma^\vee]]$ which is constant on $L$.
Each automorphism in ${\rm^\vee
Aut}(\widehat{X}_\sigma)^0$ determines an  automorphism in ${\rm
Aut}(\widetilde{X}_\sigma)^0$, so that there is a
monomorphism of proalgebraic groups $\imath: {\rm
Aut}(\widehat{X}_\sigma)^0 \to {\rm
Aut}(\widetilde{X}_\sigma)^0$. Therefore if
$\val(v,\widetilde{X}_\sigma)$ is invariant with respect to ${\rm
Aut}(\widetilde{X}_\sigma)^0$ then it is $\imath({\rm
Aut}(\widehat{X}_\sigma)^0)$-invariant. 
Since the subring $K[\widehat{X}_{\sigma}]\subset
K[\widetilde{X}_{\sigma}]$ is $\imath({\rm
Aut}(\widehat{X}_\sigma)^0)$-invariant then the 
restriction $\val(v,\widehat{X}_\sigma)$ of
$\val(v,\widetilde{X}_\sigma)$ to $\widehat{X}_\sigma$ is ${\rm
Aut}(\widehat{X}_\sigma)^0$-invariant.

Now, if $\val(v,\widehat{X}_\sigma)$ is ${\rm
Aut}(\widehat{X}_\sigma)^0$-invariant then  by
Lemma \ref{le: semi2}, it defines an ${\rm
Aut}_{(L)}(\widehat{X}^{(L)}_\sigma)^0$-invariant valuation.
By Lemma \ref{le: semi},  $\val(v,\widehat{X}^{(L)}_\sigma)$ is ${\rm
Aut}_{(L)}(\widehat{X}^{(L)}_\sigma)$-semiinvariant.
The proalgebraic group ${\rm
Aut}_{L}(\widetilde{X}_\sigma)={\rm
Aut}_{L}(\widehat{X}^L_\sigma)$  is a subgroup (as
an abstract group) of  ${\rm
Aut}_{(L)}(\widehat{X}^{(L)}_\sigma)$.  Thus the
restriction $\val(v,\widetilde{X}_\sigma)$ of 
$\val(v,\widehat{X}^{(L)}_\sigma)$ is  ${\rm
Aut}_{L}(\widetilde{X}_\sigma)$-semiinvariant.
By Lemma \ref{le: decomposition} any
automorphism $\phi$ in  ${\rm
Aut}(\widetilde{X}_\sigma)$ can be
decomposed as $\phi=\phi_0\phi_1$, where $\phi_0$  preserves
monomials and $\phi_1\in {\rm
Aut}_{L}(\widetilde{X}_\sigma)$.
 Therefore $\val(v,\widetilde{X}_\sigma)$ is ${\rm
Aut}(\widetilde{X}_\sigma)$-semiinvariant and consequently
by Lemma \ref{le: semi} it is ${\rm
Aut}(\widetilde{X}_\sigma)^0$-invariant. \qed

\begin{lemma}\label{le: stableincl} If $\sigma\leq \tau$ then
\begin{enumerate}
\item $\Inv(\sigma)\subset \Inv(\tau)$.
\item $\stab(\sigma)\subset\stab(\tau)\cap \sigma$.
\end{enumerate}
\end{lemma}
\noindent {\bf Proof.}
(1) Let $v\in \Inv(\sigma)$. Any automorphism
$\phi\in \Aut(\widetilde{X}_{\tau})$ preserves the
stratum $\strat(\sigma)$ and thus defines an automorphism
$\widehat{\phi}$ of the completion
of $Y:=\widetilde{X}_{\tau}$ at $\strat(\sigma)$ which  by Lemma
\ref{le: 1} is isomorphic to the spectrum of the ring
$K[\widehat{Y}_{\strat(\sigma)}]=K_\sigma[\widehat{X}_{\sigma}]$,
where $K_\sigma$ is the residue field of the generic orbit
point $O_\sigma$. By Lemma \ref{le: decomposition} any
automorphism $\widehat{\phi}$ of $K_\sigma[\widehat{X}_{\sigma}]$
 can be
written as $\widehat{\phi}=\phi_0\phi_1$, where $\phi_0$ preserves
monomials and $\phi_1\in {\rm
Aut}_{K_\sigma}(\widehat{X}_{\sigma})$. By Lemma \ref{le:
hat}, $\val(v,\widehat{X}_{\sigma})$ is  ${\rm
Aut}(\widehat{X}_{\sigma})^0$-invariant. Let
$\overline{K_\sigma}$ be the algebraic closure of $K_\sigma$.
By Lemma \ref{le: semi2},
$\val(v,\widehat{X}^{\overline{K_\sigma}}_{\sigma})$ is  ${\rm
Aut}_{K_\sigma}(\widehat{X}^{\overline{K_\sigma}}_{\sigma})^0$-invariant.
By Lemma \ref{le: semi}, 
$\val(v,\widehat{X}^{\overline{K_\sigma}}_{\sigma})$ is  ${\rm
Aut}_{K_\sigma}(\widehat{X}^{\overline{K_\sigma}}_{\sigma})$-semiinvariant. 
Then its restriction
$\val(v,\widehat{X}^{K_\sigma}_{\sigma})$ is  ${\rm
Aut}_{K_\sigma}(\widehat{X}^{K_\sigma}_{\sigma})$-semiinvariant.
Thus $\phi_{1*}(\val(v,\widehat{X}^{K_\sigma}_{\sigma}))$  can
be one of the finitely many toric valuations on 
$\widehat{X}^{K_\sigma}_{\sigma}$ for any $\phi_1\in {\rm
Aut}_{K_\sigma}(\widehat{X}_{\sigma})$. But then 
$\widehat{\phi}_{*}(\val(v,\widehat{Y}_{\strat(\sigma)}))=\phi_{1*}(\val(v,
\widehat{Y}_{\strat(\sigma)}))$
can be one of the finitely many toric valuations on 
$K[\widehat{Y}_{\strat(\sigma)}]$  for all
$\phi\in \Aut(\widetilde{X}_{\tau})$. The restrictions of these
valuations to the local ring of $\widetilde{X}_{\tau}$ at
$\strat(\sigma)$ define finitely many valuations on
$\widetilde{X}_{\tau}$. Hence $\val(v,\widetilde{X}_{\tau})$ is semiinvariant and
consequently invariant on $\widetilde{X}_{\tau}$. 
(2) Follows from (1) and from Lemma \ref{le: stable support}. 
\qed

\subsection{Group of divisor classes of the completion of a local
ring of a toric variety}

\begin{lemma} \label{le: divisor classes} Let $\sigma$ be a
cone in a lattice $N$ and $\Delta$ be a  suddivision $\sigma$.
Let $\widehat{X}_\Delta:={X}_\Delta\times_{X_\sigma}\widehat{X}_\sigma$.
The following groups of divisor
classes are isomorphic. (The isomorphisms are determined by the
natural morphisms). 
\begin{enumerate}
\item $\Cl(\widehat{X}_\sigma)\simeq \Cl(X_\sigma)
$, $\Pic(\widehat{X}_\sigma)\simeq\Pic(X_\sigma)
$.

\item $\Cl(\widehat{X}_\Delta)\simeq\Cl({X}_\Delta)$, 
$\Pic(\widehat{X}_\Delta)\simeq\Pic({X}_\Delta)$ 
.

\item For the affine 
toric variety $X_\sigma$ set
${\bf A}^1\widehat{\times}\widehat{X}_\sigma:=\lim_{\to} {\bf A}^1\times
X^{(n)}_\sigma=\Spec(K[t][\widehat{X}_\sigma])$. Then \\ 
$\Cl({\bf A}^1\widehat{\times}\widehat{X}_\sigma)\simeq
\Cl({\bf A}^1{\times} {X}_\sigma)\simeq \Cl({X}_\sigma)$
, $\Pic({\bf A}^1\widehat{\times}\widehat{X}_\sigma)\simeq
\Pic({\bf A}^1{\times} {X}_\sigma)\simeq \Pic({X}_\sigma)$.

\item For any 
 subdivision $\Delta$ of $\sigma$ set 
${\bf A}^1\widehat{\times}\widehat{X}_\Delta:=
({\bf A}^1\widehat{\times}\widehat{X}_\sigma)\times_{X_\sigma}{X}_\Delta$
. Then \\ $\Cl({\bf A}^1\widehat{\times}\widehat{X}_\Delta)\simeq
\Cl({\bf A}^1{\times}\widehat{X}_\Delta)\simeq \Cl({X}_\Delta)$, 
$\Pic({\bf A}^1\widehat{\times}\widehat{X}_\Delta)\simeq
\Pic({\bf A}^1{\times}\widehat{X}_\Delta)\simeq \Pic({X}_\Delta)$.

\end{enumerate}
\end{lemma}

\noindent {\bf Proof.}  Let $v_1,\ldots,v_k$ denote the
generators of $\sigma\cap N_\sigma$ and let $\tau=\langle e_1,\ldots,e_k
\rangle$ denote the regular
$k$-dimensional cone. Let $\pi:\tau\to \sigma$ be the
projection defined by
$\pi(e_i)=v_i$. Then $\pi$ defines a surjective morphism of
lattices $\pi:N_\tau\to N_\sigma$ whose kernel is a
saturated sublattice $N\subset N_\tau$. Thus the projection
$\pi$
defines the quotient map
$X_\tau \to X_\sigma=X_\tau//T$ for the subtorus $T\subset
T_\tau$ corresponding to the sublattice $N\subset N_\tau$.
This also defines the quotient morphism $\alpha: \widehat{X}_\tau\to \widehat{X}_\sigma$, 
where $ K[\widehat{X}_\sigma] = K[\widehat{X}_\tau]^T$.

Denote by $T_\tau\subset X_\tau $, $T_\sigma\subset
X_\sigma $ the relevant tori and set
$\widehat{T}_\tau:=\widehat{X}_\tau\times_{X_\tau}
T_\tau\subset \widehat{X}_\tau$,
$\widehat{T}_\sigma:=\widehat{X}_ \sigma\times_{X_\sigma}
T_\sigma\subset \widehat{X}_\sigma$. Then both schemes are nonsingular and let
$\alpha': \widehat{T}_\tau \to \widehat{T}_\sigma $ be the restriction of $\alpha$. We also have 
$ K[\widehat{T}_\sigma] = K[\widehat{T}_\tau]^T$.

The morphism $\alpha$ defines a group homomorphism
$\Cl(\widehat{T}_\sigma)=\Pic(\widehat{T}_\sigma)
\buildrel \alpha^*\over \longrightarrow 
\Pic(\widehat{T}_\tau)=\Cl(\widehat{T}_\tau)=0$. The last
group is $0$ since the ring $K[\widehat{X}_\tau]$ as well as
its localization $K[\widehat{T}_\tau]$ are UFD. Let $D$ be
an effective Cartier divisor on $\widehat{T}_\sigma$. Then $\alpha^*(D)$ is $T$-
invariant and principal on
$\widehat{T}_\tau$. This  means that the ideal $I$ of $\alpha^*(D)$ is
generated by $f\in K[\widehat{T}_\tau]$. By multiplying
by a suitable monomial we can asume that $f\in
K[\widehat{X}_ \tau]^T$. Note that all $t\cdot f$, $t\in T $
also belong to $I$. Let $M_T$ denote the lattice of
characters of $T$. Let $f_\beta$ denote the component of $f$
with weight $\beta\in M_T$. By considering
$K[\widehat{X}_\tau]/m^k$, where $m\subset K[\widehat{X}_\tau]$ is the maximal ideal,
we see that all $f_\beta+m^k\in I\cdot
K[\widehat{X}_\tau]/m^k$. Moreover $I\cdot
K[\widehat{X}_\tau]/m^k$ is generated by all
$f_\beta+m^k$. This yields $I=(f_\beta)_{\beta\in M_T}$ or
$(1)=(f_\beta/f)$ showing that there is a $\beta_0$ such
that $f_{\beta_0}/f$ is invertible, which means
$I=(f_{\beta_0})$. By multiplying by
monomials we can assume that $f_{\beta_0}$ is $T$-invariant.
Thus $f_{\beta_0}\in
K[\widehat{T}_\sigma]^T=K[\widehat{X}_\tau]$  and generates the
ideal of $D$. But then $D$ is principal on
$\widehat{X}_\tau$ and finally $\Cl(\widehat{T}_\sigma)=0$.

Now any divisor $D$ in $\Cl(\widehat{X}_\sigma)$
(respectively in $\Cl(\widehat{X}_\Delta)$) is linearly
equivalent to the  $T$-invariant one 
$D':=D-(f_D)$, where $D_{|\widehat{T}_\sigma}=
(f_D)_{|\widehat{T}_\sigma}$.

(3)(4) We repeat the reasoning from (1) and (2) and use the
fact that $K[{\bf A}^1\widehat{\times}\widehat{X}_\tau]=K[t][[x_1,\ldots,x_k]]$ is UFD
(see Bourbaki\cite{Bourbaki}).
\qed

\begin{lemma}\label{le: trivial action}
\begin{enumerate}
\item $G^{\sigma}$ acts trivially on
$\Cl(\widetilde{X}_{\Delta^\sigma})$. 

\item Let $\Delta^\sigma$ be a subdivision
of $\sigma$ such that $\widetilde{X}_{\Delta^\sigma}\to 
\widetilde{X}_{\sigma}$ is $G^{\sigma}$-equivariant. Then
$G^{\sigma}$ acts trivially on $\Cl(\widetilde{X}_{\Delta^\sigma})$.
\end{enumerate}  
\end{lemma}

\noindent {\bf Proof.} We shall use the notation and
results from Lemma \ref{le: group structure}. 

(1) The natural morphism 
$\widetilde{X}_{\Delta^\sigma}\to {X}_{\Delta^\sigma}$ is
$G^\sigma_{\nu}$-equivariant.  By Lemma \ref{le: divisor classes} this
morphism induces a $G^\sigma_{\nu}$-equivariant isomorphism  
$\Cl(\widetilde{X}_{\sigma})\to \Cl({X}_{\sigma})$.
By Sumihiro \cite{Sumihiro}  the algebraic linear group acts trivially on
$\Cl({X}_{\sigma})$, which yields a trivial action on 
$\Cl(\widetilde{X}_{\sigma})$. By Lemma \ref{le: group
structure} it suffices to show that 
${\rm ker}_\sigma$ acts trivially on $\Cl(\widetilde{X}_{\sigma})$. 

Let $\widetilde{X}_\sigma^{\rm ns}\subset 
\widetilde{X}_\sigma$ denote the open subset of
nonsingular points. Its complement is of codimension 2.
We get a  $G^{\sigma}$-equivariant isomorphism
$\Pic(\widetilde{X}_\sigma^{\rm ns})= 
\Cl(\widetilde{X}_\sigma^{\rm ns})\simeq 
\Cl(\widetilde{X}_\sigma)$. It suffices to prove
that ${\rm ker}_\sigma$ acts trivially on
$\Pic(\widetilde{X}_\sigma^{\rm ns})$. Fix $g\in G$.
 Let $i_{g,{\bf A}^1}:{\bf A}^1\to {\rm ker}_\sigma$ be the morphism from Lemma 
\ref{le: group structure}. 
Let $\Phi:
G^{\sigma}\widehat{\times}\widetilde{X}_\sigma^{\rm ns}\to 
\widetilde{X}_\sigma^{\rm ns} $ be the
action morphism from Lemma \ref{le: group action} and set
$\Phi_{{\bf A}^1}:=i_{g,{\bf A}^1}\circ 
\Phi: 
{\bf A}^1\widehat{\times}\widetilde{X}_\sigma^{\rm ns}\to 
\widetilde{X}_\sigma^{\rm ns} $. Denote by 
$p: {\bf A}^1\widehat{\times}\widetilde{X}_\sigma^{\rm ns}\to 
\widetilde{X}_\sigma^{\rm ns} $  the standard
projection. By Lemma \ref{le: divisor classes}(4), 
$p^*:\Pic(\widetilde{X}_\sigma^{\rm ns})\to
\Pic({\bf A}^1\widehat{\times}\widetilde{X}_\sigma^{\rm ns})$ is
an isomorphism . Let $D\in \Pic(\widetilde{X}_\sigma^{\rm ns})$.
Let $j_g:\widetilde{X}_\sigma^{\rm ns}\to  
\{g\}\widehat{\times}\widetilde{X}_\sigma^{\rm ns}
\subset {\bf A}^1\widehat{\times}\widetilde{X}_\sigma^{\rm ns}$
denote the standard embedding. Then $p\circ
j_g=\id_{|\widetilde{X}_\sigma^{\rm ns}}$. Since
$p^*$ is an isomorphism there is $D'\in
\Pic(\widetilde{X}_\sigma^{\rm ns})$ such that   
$\Phi_{{\bf A}^1}^*(D)\simeq p^*(D')$.   
Therefore
$j_g^*\Phi_{{\bf A}^1}^*(D)=g\cdot D\simeq
j_g^*p^*(D')=D'$. This implies $h\cdot D\simeq D$ for any
$h\in G={\rm ker}_\sigma \cdot\Aut_\nu(\widetilde{X}_{\sigma})$

(2) Let $L\subset \widetilde{X}_{\Delta^\sigma}$ denote the
complement of the set where the birational
morphism $\psi:\widetilde{X}_{\Delta^\sigma}\to 
\widetilde{X}_{\sigma}$ is an isomorphism. 
It is a $G^{\sigma}$-invariant subset of
codimension 2. Therefore we get $G^{\sigma}$-equivariant
homomorphisms $\Cl(\widetilde{X}_{\sigma})\simeq
\Cl(\widetilde{X}_{\sigma}\setminus L) \hookrightarrow 
\Cl(\widetilde{X}_{\Delta^\sigma})$. Note also that 
$\Cl(\widetilde{X}_{\Delta^\sigma})$  is generated by
$\Cl(\widetilde{X}_{\sigma})$ and by the exceptional divisors of $\psi$ 
which are $G^{\sigma}$-invariant by Lemma \ref{le: divisors}. 
Finally, all elements of
$\Cl(\widetilde{X}_{\Delta^\sigma})$ are $G^{\sigma}$-invariant.\qed

\bigskip
\subsection{Simple definition of a canonical subdivision of a
semicomplex}\label{se: simple}
 
Recall that for any fan  $\Sigma$  we denote 
by $\Ver(\Sigma)$  the set of all $1$-dimensional
rays in $\Sigma$.   

\begin{proposition}\label{pr: simple} Let $\Sigma$ be an oriented
$\Gamma$-semicomplex. A  subdivision $\Delta$ of $\Sigma$
is canonical if for any $\sigma \in \Sigma$ $$\Ver(\Delta^\sigma)\setminus
\Ver(\overline{\sigma}) \subset{\rm
Stab}(\Sigma).$$
\end{proposition}
\noindent{\bf Proof.} This is an immediate
consequence of Lemma \ref{le: stable support} and the following.

\begin{proposition}\label{pr: ssimple} Let $\sigma$ be a
$\Gamma$-semicone and $\Delta^\sigma$ be a
 subdivision of $\sigma$. Then
the following conditions are equivalent:
\begin{enumerate}
\item $\widetilde{X}_{\Delta^\sigma}\to \widetilde{X}_{\sigma}$ is
$G^{\sigma}$-equivariant (where $G^{\sigma}=\Aut(\widetilde{X}_{\sigma})^0$). 
\item $\Ver(\Delta^\sigma)\subset
\Ver(\overline{\sigma})\cup {\rm Inv}({\sigma})$.   
\item $\widehat{X}_{\Delta^\sigma}\to \widehat{X}_{\sigma}$ is
$\Aut(\widehat{X}_{\sigma})^0$-equivariant.
\end{enumerate}
\end{proposition}

Before the proof of Propositions \ref{pr: ssimple} and
\ref{pr: simple} we need to show a few lemmas below.

\begin{lemma} \label{le: multiplication} Let $\sigma$ be a
$\Gamma$-semicone. If $\Delta_1$ and $\Delta_2$ are
 subdivisions of $\sigma$  for which 
$\widetilde{X}_{\Delta_i}\to \widetilde{X}_{\sigma}$ is
$G^{\sigma}$-equivariant then for the subdivision
$$\Delta_1\cdot\Delta_2:= \{\sigma_1\cap \sigma_2\mid
\sigma_1\in \Delta_1, \sigma_2\in \Delta_2\}$$\noindent 
the morphism 
$\widetilde{X}_{\Delta_1\cdot\Delta_2}\to \widetilde{X}_{\sigma}$ is
$G^{\sigma}$-equivariant.

Let $I_1$ and $I_2$ be
$G^{\sigma}$-invariant on $\widetilde{X}_{\sigma}$. Let 
$\Delta_1$ and $\Delta_2$ be subdivisions of $\sigma$
corresponding to the normalized blow-ups at $I_1$ and
$I_2$. Then $I_1\cdot I_2$ corresponds to 
$\Delta_1\cdot\Delta_2$.
\end{lemma}

\noindent {\bf Proof.} It follows from the universal
property of the fiber product that $\Delta_1\cdot\Delta_2$ corresponds
to the normalization of an irreducible component in 
$\widetilde{X}_{\Delta_1}\times_{\widetilde{X}_{\sigma}}
{\widetilde{X}_{\Delta_2}}$. 

The second part of the lemma is an
immediate consequence of  the relations between ideals
and $\ord$-functions: $\ort(I_1\cdot
I_2)=\ort(I_1)+\ort(I_2)$ corresponds to $\Delta_1\cdot\Delta_2$(\cite{KKMS}). \qed

\begin{lemma}\label{le: vertices}  
Let $\sigma$ be a $\Gamma$-semicone.
Let $\tau\subset |\sigma|$ be a cone
with all rays $\Ver(\tau)$ in $\Ver(\overline{\sigma})\cup{\rm Inv}(\sigma) 
$. 
\begin{enumerate} 

\item There exists a fan subdivision 
$\Delta_\tau$ of $\sigma$ which
contains $\tau$ as its face and for which 
$\widetilde{X}_{\Delta^\tau}\to \widetilde{X}_{\tau}$ is
$G^{\sigma}$-equivariant. 

\item If all rays of $\tau$ which are not in ${\rm
Inv}(\sigma)$ determine a face $\varrho$ of $|\sigma|$
(and $\tau$) then there exists
a $G^{\sigma}$-invariant ideal $I$ such that  the normalized  
blow-up of $I$ corresponds to the subdivision $\Delta_\tau$ of
$\sigma$ containing  $\tau$.

\end{enumerate}
\end{lemma}

\noindent {\bf Proof.}  Write $\tau=\langle v_1,\ldots,v_l\rangle$, 
where $v_1,\ldots,v_k\in \Inv(\sigma)$ and
$v_{k+1},\ldots,v_{l}$  are in $\Ver(\overline{\sigma})$.

(1) Let $\{\tau_i\mid i\in J_0\}$ denote the set of all
one-codimensional faces of $\tau$.
For any one-codimensional
face $\tau_i$ of $\tau$, where $i\in J_0$,  
 find an integral functional $H_i$ such that
$H_{i|\tau_i}=0$ and $H_{i|\tau\setminus\tau_i}  >0$.
Find functionals  $H_j$, $ j\in J_1$  with
common zeros
exactly on ${\rm lin}(\tau)$. Then $\tau=\{x\in
\sigma\mid H_i(x)\geq 0, H_j(x)=0, i\in J_0, j\in J_1\}$.  

Let $F_0$ be an integral functional such that for any $i\in
J_0$ and $j\in J_1$, $F_0+H_i$ and $F_0-|H_j|$ are strictly greater than
zero on $\sigma\setminus \{0\}$.  
Set $F_0(v_1)=n_1,\ldots,F_0(v_l)=n_l$.
 Let $\nu_1,\ldots,\nu_l$ be
valuations corresponding to $v_1,\ldots,v_l$. 
Set $$I:=
I_{\nu_1,n_1}\cap\ldots\cap I_{\nu_l,n_l}.$$ Let 
$p:\widetilde{X}_{\Delta_0}\to \widetilde{X}_{\sigma}$ be the
composition of the $G^{\sigma}$-equivariant blow-ups $\bl_{\nu_k}\circ\ldots\circ \bl_{\nu_1}$
corresponding to  $G^{\sigma}$-invariant valuations $\nu_1$,\ldots,$\nu_k$.
Then all valuations $\nu_1,\ldots,\nu_l$ correspond to Weil
divisors $D_1,\ldots,D_l$ on $\widetilde{X}_{\Delta_0}.$ 
Set $D=n_1D_1+\ldots+n_lD_l$. Then $$p_*(I_D)=\{f\in p_*(
{\cO}_{\widetilde{X}_{\Delta_\tau}})={\cO}_{\widetilde{X}_{\sigma}}\mid
\nu_{i}(f)\geq n_i,
i=1,\ldots,l\}= I.$$ By Lemma \ref{le: trivial action}, 
$G^{\sigma}$ acts
trivially on $\Cl(\widetilde{X}_{\Delta_0})$; then for any
$g\in G$ we have $g_*(D)=D+(f_g)-(h_g)$, where $f_g, h_g\in
{\cO}_{\widetilde{X}_{\sigma}}$. In other words for any
$g\in G$, $$f_g\cdot I_D=h_g\cdot I_{g_*(D)},$$ which implies
$$f_g\cdot I=h_g\cdot {g_*(I)}.$$
This means that the action of $G^{\sigma}$ lifts to the blow-up of
$I$, and to its  normalization. The normalized blow-up of $I$
corresponds to a $G^{\sigma}$-equivariant subdivision $\Delta_\tau$ of $\sigma$ into 
maximal cones, where the piecewise linear function
$\ord(I)={\rm min}\{F\in \sigma^\vee\mid F(v_i)\geq n_i\}$ is
linear. We have to show that $\tau\in \Delta_\tau$. 
Note that $\ord(I)$ is linear on $\tau$ since by
definition $\ord(I)_{|\tau}=F_{0|\tau}$.

 Let $x\not\in \tau$. Consider first the case $x\in{\rm lin}(\tau)$.  There exists  
a one-codimensional face $\tau_{i_0}$, $i_0\in J_0$, of
$\tau$ such that $H_{i_0}(x)<0$. Then $F_1:=F_0+H_{i_0} \in \sigma^\vee$ 
and  by definition $F_1\geq
\ord(I)$, since $H_{i_0}(v_i)\geq 0$ and $F_1(v_i)=F_0(v_i)+H_{i_0}(v_i)\geq
n_i$. This implies $\ord(I)(x)\leq F_1(x)<F_0(x)$. Now let 
$x\not\in {\rm lin}(\tau)$. Then there is $H_{j_0}$, $j_0\in
J_1$ such that $H_{j_0}(x)\neq 0$. Set
$F_1:=F_0+H_{j_0}$ if $H_{j_0}(x)<0$ and $F_1:=F_0-H_{j_0}$ otherwise. 
The definition of $F_0$ shows that  $F_1$ is strictly greater than
zero on $\sigma\setminus \{0\}$. We have  $F_1(v_i)=F_0(v_i)+H_{j_0}(v_i)=
n_i$ and consequently $\ord(I)\leq F_1$. Again
$\ord(I)(x)\leq F_1(x)<
F_0(x)$.  This implies $\ord(I)< F_0$ 
off ${\rm lin}(\tau)$.   

(2)  The reasoning is very similar. Let $J_0$ denote the set
of all  one-codimensional faces $\tau_i$  of
$\tau$ which contain  $\varrho$.  For any $i\in J_0$
 find an integral functional $H_i$ such that
$H_{i|\tau_i}=0$ and $H_{i|\tau\setminus\tau_i}  >0$.
Let $H_j$, $ j\in J_1$,  be functionals with
common zeros
exactly on ${\rm lin}(\tau)$. Then $\tau=\{x\in
\sigma\mid H_i(x)\geq 0, H_j(x)=0, i\in J_0, j\in J_1\}$.  

Let $F_0$ be an integral functional such that
$F_{0|\tau}=0$, and for any $i\in
I$ and $j\in J$, $F_0+H_i$, $F_0-|H_j|$ are strictly greater than
zero on $\sigma\setminus \varrho$ for $j\in J_1$. 
Set $F_0(v_1)=n_1,\ldots,F_0(v_k)=n_k$. Then $F_0(v_{k+1})=F_0(v_l)=0$.
 Let $\nu_1,\ldots,\nu_k$ be
valuations corresponding to $v_1,\ldots,v_k$. 
Set $$I:=
I_{\nu_1,n_1}\cap\ldots\cap I_{\nu_k,n_k}.$$ By definition $I$
is $G^{\sigma}$-invariant. The normalization of the blow-up of $I$
corresponds to a $G^{\sigma}$-equivariant subdivision
$\Delta_\tau$ of $\sigma$ into 
maximal cones, where the piecewise linear function
$\ord(I)={\rm min}\{F\in \sigma^\vee\mid F(v_i)\geq n_i\}$ is
linear. We have to show that $\tau\in \Delta_\tau$. 
Note that $\ord(I)$ is linear on $\tau$ since by
definition $\ord(I)_{|\tau}=F_{0|\tau}$.

The rest of the proof is the same as in (1).\qed

\begin{lemma} \label{le: Ssigma} Let $\sigma$ be a
$\Gamma$-semicone and $  \tau \subset
|\sigma|$ be a cone all of whose  rays are in $\Ver({\overline{\sigma}})\cup
{\rm Inv}(\sigma)$. Let $\Delta_\tau$ be a 
subdivision of $\sigma$ containing $\tau$ as its face
and for which 
$\widetilde{X}_{\Delta_\tau}\to \widetilde{X}_{\sigma}$ is
$G^{\sigma}$-equivariant. 
Then 
$$S(\tau, \Delta_\tau):=\{\varrho\in \Delta_\tau \mid (\varrho \setminus
\tau) \cap {\rm Inv}(\sigma) = \emptyset \}$$     
\noindent is a subfan of $\Delta_\tau$ for which the open subset
$\widetilde{X}_{S(\tau, \Delta_\tau)}={X}_{S(\tau,
\Delta_\tau)}\times_{X_\sigma}\widetilde{X}_\sigma$ of
$\widetilde{X}_{\Delta_\tau}$  is
$G^{\sigma}$-invariant. 
\end{lemma}

\noindent {\bf Proof } It suffices to show that 
each $G^{\sigma}$-orbit intersecting $\widetilde{X}_{S(\tau, \Delta_\tau)}:=
{X}_{S(\tau, \Delta_\tau)}\times_{{X}_\sigma}\widetilde{X}_{\sigma}$
is contained in this subscheme. 
Since $\widetilde{X}_{S(\tau, \Delta_\tau)}$ is an open subset
of $\widetilde{X}_{\Delta_\tau}$
any $G^{\sigma}$-orbit contains the generic point of
$O_\varrho$ for $\varrho\in S(\tau, \Delta_\tau)$. By Lemma 
\ref{le: minimal vectorss}(3), $\inte(\varrho)$  intersects 
${\rm Inv}(\sigma)$.  
The other toric orbits in the $G^{\sigma}$-orbit
 correspond to the cones $\varrho'\in {\rm Star}(\varrho,
\Sigma)$ for which $\varrho'\setminus\varrho =\emptyset$. Thus
$\varrho$ is a face of $\tau$ and all $\varrho'$ belong to $S(\tau, \Delta_\tau)$. \qed

\begin{lemma} \label{le: Ssigma2} Let $\sigma$ be a
$\Gamma$-semicone. Then
$ S(\tau):=S(\tau, \Delta)$ does not depend
upon a subdivision $\Delta$ of $\sigma$ containing $\tau$ and all of whose 
 rays are in $\Ver(\overline{\sigma})\cup
{\rm Inv}(\sigma)$.
\end{lemma}

\noindent {\bf Proof.} Let $\Delta_1$ and $\Delta_2$ be two
 subdivisions containing $\tau$ and for
which  all rays are in $\Ver(\overline{\sigma})\cup
{\rm Inv}(\sigma)$. By Lemmas
\ref{le: multiplication} and \ref{le: vertices} we may find
subdivisions 
$\Delta_\varrho$, where
${\varrho\in\Delta_1 \cup \Delta_2}$ and  
the common subdivision 
$\Delta:=\prod_{\varrho\in\Delta_1 \cup \Delta_2}\Delta_\varrho$
such that $\widetilde{X}_{\Delta}\to
\widetilde{X}_{\sigma}$ is $G^{\sigma}$-equivariant.
Then $\Ver(\Delta_1\cdot\Delta_2)\subset
\Ver(\Delta)\subset\Ver({\sigma})\cup
{\rm Inv}(\sigma)$. Hence we can assume that
$\Delta_1$ is a subdivision of $\Delta_2$ replacing $
\Delta_1$ with $\Delta_1\cdot\Delta_2$ if necessary. Let
$\varrho\in S(\tau, \Delta_2)$. By Lemma \ref{le:
stab-subdivisions}(1) applied to $\Delta_\varrho$ we can
represent $\varrho$ as 
$\varrho:=\varrho' \oplus \langle e_1,\ldots,e_k \rangle$, where
 $ \inte( \varrho') \cap {\rm
Inv}(\sigma)\neq \emptyset$, $(\varrho\setminus\varrho')\cap {\rm
Inv}(\sigma)= \emptyset$ and thus $\varrho'\preceq
\tau$. Since all new rays of $\Delta_{1}|\varrho$
are in ${\rm
Inv}(\sigma)$ we find that
$\Delta_{1}|\varrho=\Delta_{1}|\varrho'\oplus \langle e_1,\ldots,e_k \rangle$.
But $\Delta_{1}|\varrho'=\varrho'$  since
$\varrho'\preceq\tau$.
Therefore $\Delta_{1}|\varrho=\varrho$ and
$\varrho\in\Delta_1$,  and consequently 
$\varrho\in S(\tau, \Delta_1)$.\qed 

\begin{lemma} \label{le: Ssimple3} Let $\sigma$ be a $\Gamma$-semicone.
Let $\Delta$ be any subdivision of $\sigma$ with all
''new '' rays belonging to ${\rm Inv}(\sigma)$. Then for
any  $\tau\in\Delta$,
$\widetilde{X}_{S(\tau)}$ is $G^{\sigma}$-invariant.
\end{lemma}
\noindent{\bf Proof} By Lemma \ref{le: vertices}, we find a
subidvision $\Delta_\tau$ of $\sigma$. By Lemma \ref{le: Ssigma2},
$\widetilde{X}_{S(\tau)}=\widetilde{X}_{S(\tau,\Delta_\tau)}$.
By Lemma \ref{le: Ssigma}, the latter  scheme is $G^{\sigma}$-invariant.
\qed

\bigskip
\noindent {\bf Proof of Proposition \ref{pr: ssimple}.}
$(1)\Rightarrow (2)$ All exceptional divisors determine
$G^{\sigma}$-semiinvariant and therefore $G^{\sigma}$-invariant
valuations. The latter  correspond to vectors in 
$\Ver(\Sigma)\subset\Ver({\sigma})$.

$(1)\Leftarrow (2)$
Let $\Delta$ be any subdivision of $\sigma$ containing
$\tau$ with all
''new '' rays belonging to ${\rm Inv}(\sigma)$. 
Then by Lemmas \ref{le: Ssimple3}, 
$\widetilde{X}_{\Delta}$ is the
union of open $G^{\sigma}$-invariant neighborhoods
$\widetilde{X}_{S(\tau)}$, where $\tau\in\Delta$, and
consequently,  $\widetilde{X}_{\Delta}$ is $G^{\sigma}$-invariant. 

$(2)\equiv (3)$ By Lemma \ref{le: hat} vectors defining $G^{\sigma}$-invariant
valuations are the same as those defining \\ 
$\Aut(\widehat{X}_{\sigma})^0$-invariant valuations.
Then the proof of the equivalence is the same as for
$(1)\equiv (2)$. We have to
replace  $\widetilde{X}_{\sigma}$ with $\widehat{X}_{\sigma}$.

\qed

\begin{lemma}\label{le: ideal} Let $\sigma$ be a
$\Gamma$-semicone and  $G^{\sigma}$ be a connected
proalgebraic group acting on $X:=\widetilde{X}_{\sigma}$. 
Let $\Delta^\sigma$ be a
subdivision of $\sigma$ such that $G^K$ acts on
$Y:=\widetilde{X}_{\Delta^\sigma}$ and the toric morphism $\psi: Y\to X$
is a proper birational $G^K$-equivariant
morphism. Then there exists a toric ideal $I$ on $X$ such
that the normalization of the blow-up of $I$,
$Z\to X$, factors as
$Z\to Y\to X$.
\end{lemma}

\noindent {\bf Proof.} First subdivide $\Delta^\sigma$ so that
all its faces satisfy the condition of  Lemma \ref{le: vertices}(2).
We can make $\Delta^\sigma$ simplicial  by applying 
additional star subdivisions at all
its $1$-dimensional rays (see \cite{Wlodarczyk1}). 
Let  $\tau$ be a  face of $\Delta^\sigma$ which does not 
satisfy the condition of  Lemma \ref{le: vertices}(2). 
Denote by $\tau'$ a minimal face of $\tau$ with rays
in $\Ver(\overline{\sigma})$ and which is not a face of the cone $\sigma$. Then
there is $\sigma'$ which is a
face of $\sigma$ and for which
$\inte(\tau')\subset \inte(\sigma')$.
Write $\sigma'=\sing^\Gamma(\sigma')\oplus^\Gamma\langle
e_1,\ldots,e_k\rangle$. Then $\tau'\cap \sing^\Gamma(\sigma')$ is a
face of $\tau'$ but it is not a face of $\sigma$.
(Otherwise  $\tau'$ was a face of $\sigma$.)
By minimality  $\tau'=\tau'\cap \sing^\Gamma(\sigma')\subset
\sing^\Gamma(\sigma')$ and $\sigma'=\sing^\Gamma(\sigma')$ is
indecomposable and since $\inte(\tau')\subset
\inte(\sigma')$ we can find by Lemma 
\ref{le: minimal vectorss}(5) a vector
 $v\in\inte(\tau)\cap \Inv(\sigma)$.
By applying  star subdivision at $\langle
v\rangle$ to $\Sigma$ we eliminate the cone $\tau$ (and
some others cones not satisfying the condition of  Lemma 
\ref{le: vertices}(2)). After
a finite number of steps we arrive at a subdivision $\Delta^\sigma$
with all faces satisfying the condition of  Lemma \ref{le: vertices}(2).
 
Then by Lemma \ref{le: vertices}(2), for any $\tau\in\Delta$, we  
construct the ideal $I_\tau$ such that the normalization of
the blow-up of $I_\tau$
corresponds to the subdivision $\Delta^\sigma_\tau$ of $\sigma$
containing $\tau$.  
Then by Lemma \ref{le: multiplication}  the $G^{\sigma}$-invariant
ideal $I:=\prod_{I_\tau\in
\Sigma}I_\tau$  correponds to the subdivision $\Delta_I:=\prod_{\tau\in
\Sigma}\Delta^\sigma_\tau$ of $\sigma$. \qed

\begin{lemma}\label{le: extendsub}
Let $\sigma$ be a face of a fan $\Sigma$. Let $\Delta_1$ be
a subdivision of  $\sigma$ such that $\Ver(\Delta_1)=\Ver(\sigma)$. Then
there is a subdiision $\Delta_2$ of $\Sigma$ such that
$\Delta_2|\sigma=\Delta_1$ and $\Ver(\Delta_2)=\Ver(\Sigma)$.
\end{lemma}
\noindent{\bf Proof.}
Set $\{v_1,\ldots,v_k\}=\Ver(\Sigma)\setminus
\Ver(\sigma)$, $\Delta_0:=\langle v_1 \rangle\cdots
\langle v_k \rangle\cdot\Sigma$. Then $\sigma\in
\Delta_0$. We shall construct
a subdivision $\Delta_2$ for any cone in $\Delta_0$. If
$\tau\in\Delta_0$ does not intersect $\sigma$ then we put 
$\Delta_2|\tau=\tau$.
If $\tau\in\Delta_0$ intersects $\sigma$ along a common face
$\sigma_0$ then all rays $\Ver(\tau\setminus \sigma_0)$ are
centers of star subdivisions and therefore are
linearly independent of the other rays of $\sigma$
and in particular of $\lin(\sigma_0)$ (After a star
subdivision at $\langle v \rangle$, $v$ becomes linearly
independent of all other rays of cones containing $v$. Hence the
subdivision $\Delta_2|\sigma_0$ extends to a unique
subdivision of $\tau$ with no new rays,
$\Delta_2|\tau:=\{\delta+\sum_{\varrho\in(\Ver(\tau)\setminus
\Ver(\sigma_0))}\varrho \mid \delta\in \Delta|\sigma_0\}$.\qed

\begin{lemma} \label{le: subdivision} Let $\tau\subset\sigma$ be a cone such that
$\Ver(\tau)\subset\Ver(\sigma)$. Then there is a proper
subdivision $\Delta$ of $\sigma$ containing $\tau$ and such 
that $\Ver(\Delta)=\Ver(\sigma)$. 
\end{lemma}
\noindent{\bf Proof.} By Lemma \ref{le: extendsub} 
it suffices to prove the lemma for the situation when
$\Ver(\sigma)\setminus \Ver(\tau)$ consists of one ray.
Then for any $\tau$ we find $\sigma'$ containing $\tau$ and
for which $\card(\Ver(\sigma)\setminus \Ver(\tau))=1$. By induction we can
find a subdivision of $\sigma'$ and then by Lemma \ref{le:
extendsub} extend it to a subdivision of $\sigma$.
 
Let $\tau=\langle v_1,\ldots,v_k\rangle$ and 
$\sigma=\langle v_1,\ldots,v_k,v_{k+1}\rangle$. 
If $\tau$ is simplicial and $\tau$ is not a face of
$\sigma$ then there is a unique 
linear relation between $v_1,\ldots,v_k,v_{k+1}$.
Without loss of generality we can assume that it is
$a_rv_r+a_{r+1}v_{r+1}+\ldots+a_lv_l=a_{l+1}v_{l+1}+\ldots+a_{k+1}v_{k+1}$,
where $r\geq 1$ and all  coefficients are positive. Then
there are exactly two
simplicial subdivisions of 
$\sigma$ with no
''new'' rays and with maximal cones $\langle
v_1,\ldots,\check{v}_i,\ldots,v_k,v_{k+1}\rangle$, where
$i=r,\ldots l$ and $i=l+1,\ldots,k+1$ (see
\cite{Wlodarczyk1}). The second subdivision
contains $\tau$.
If $\tau=\langle v_1,\ldots,v_k\rangle$ is not simplicial
then for the cone $\tau^\vee$ we
find a simplicial cone $(\tau_0^\vee,M_{\tau_0})$ and a linear epimorphism 
$\phi^\vee: (\tau_0^\vee,M_{\tau_0}\to(\tau^\vee,M_\tau)$ 
mapping rays of $\tau_0^\vee$ to 
rays of $\tau_0^\vee$. Moreover let $H^\vee\subset
M_{\tau_0}$ denote the vector
subspace which is the kernel of
$\phi^\vee$. Then $\tau^\vee\simeq(\tau_0^\vee+H^\vee)/H^\vee$.
The dual morphism $\phi: (\tau,N_\tau) \to
(\tau_0,N_{\tau_0})$ is a monomorphism which maps
$\tau$ isomorphically to $\tau_0\cap H$, where $H=\{v\in N_{\tau_0}\mid v_{|H^\vee}=0\}$.
Then also $v_{k+1}\in H\setminus \tau_0$. Note that $v_{k+1}\not\in -\tau_0$
since otherwise $v_{k+1}\in (-\tau_0)\cap H=-\tau$. Thus
we can apply the previous case to 
$\sigma_0:=\langle v_{k+1}\rangle + \tau_0\supset \tau_0$ and obtain a
subdivision $\Delta_0$ of $\sigma_0$ containing
$\tau_0$ with no ''new '' rays. Intersecting $\Delta_0$
with $H$ defines a subdivision $\Delta$ of $\sigma$ containing
$\tau$. All rays in $\Sigma$ except  $v_{k+1}\in H$ are
obtained by intersecting some faces of $\tau_0$ with $H$,
hence belong to $\Ver(\tau)$.\qed

\section{Orientability of stratified toroidal varieties 
and resolution of singularities}
\bigskip
\subsection{Orientation group of an affine toric variety}

\begin{definition}\label{de: orientation group}
Let $(X,S)$ be a stratified noetherian scheme over $K$ and
$x\in X$ be a $K$-rational point fixed under the
$\Gamma$-action. By the {\it orientation group} of $(X,S)$ at $x\in
X$ we mean
$\Theta^\Gamma(X,S,x):=\Aut(\widehat{X}_{X,x})/\Aut(\widehat{X}_{X,x})^0$.
\end{definition}

\begin{lemma}\label{le: orientation group}
Let $\sigma$ be a $\Gamma$-semicone. Let $\widehat{\Inv}(\sigma)=\Inv(\sigma)$ 
be the
set of the vectors $v\in\sigma$ corresponding to  
$\Aut(\widehat{X}_{\sigma})^0$-invariant valuations $\val(v)$ on
$\widehat{X}_{\sigma}$ (see Lemmma \ref{le: hat}).
 Let 
$\Aut(\sigma)_{\Inv}$ be the group of all elements $g$ of 
$\Aut(\sigma)$ satisfying the following conditions:
\begin{enumerate}

\item $g_{|\Inv(\sigma)}=\id_{|\Inv(\sigma)}$.

\item For any subdivision $\Delta$ of $\sigma$ such that 
$\Ver(\Delta)=\Ver(\overline{\sigma})$, $g$ lifts to an automorphism
of $\Delta$.

\end{enumerate}

Then
$\Aut(\sigma)_{\Inv}=Aut(\sigma)^0$ and
$\Theta^\Gamma(X_{\sigma},S,O_{\sigma})=\Aut(\sigma)/\Aut(\sigma)_{\Inv}$.
  
\end{lemma}
\noindent{\bf Proof.} By Lemma \ref{le: surjection} there is
an epimomorphism $\Aut(\sigma)\longrightarrow 
\Theta^\Gamma(X_{\sigma},S,O_{\sigma})$ with  kernel $\Aut(\sigma)^0$. 
 By definition all $g\in
\Aut(\sigma)^0$ preserve $\Inv(\sigma)$. By 
Proposition \ref{pr: simple} the subdivisions $\Delta$ satisfying 
condition (2) of the lemma determine
$\Aut(\widehat{X}_{\sigma})^0$-equivariant
morphisms $\widehat{X}_{\Delta}\to\widehat{X}_{\sigma}$ which
commute with morphisms defined by $g\in \Aut(\sigma)^0$. Thus
automorphisms $g\in \Aut(\sigma)^0$ lift to automorphisms of $\Delta$. 
Therefore $\Aut(\sigma)^0\subseteq
\Aut(\sigma)_{\Inv}$. We need to prove that $\Aut(\sigma)^0\supseteq
\Aut(\sigma)_{\Inv}$.

Let $v_1,\ldots,v_k$ denote the generators of $\sigma$. Let
$\sigma_0$ denote the regular cone spanned by the standard
basis $e_1,\ldots,e_k$ of some $k$-dimensional lattice.
Let $\pi:\sigma_0\to\sigma$ be the projection defined by
$\pi(e_i)=v_i$. Each automorphism of $\sigma$ lifts to a unique
automorphism of $\sigma_0$. This defines a group 
homomorphism $\pi^*:\Aut(\sigma)_{\Inv}\to \Aut(\sigma_{0})$.
Let $\Aut(\sigma_{0})_\Inv$ denote the group of automorphisms of
the cone $\sigma_{0}$ preserving $\pi^{-1}({\Inv(\sigma)})$.

\begin{lemma}\label{le: sigma}
$\Aut(\sigma_{0})_\Inv\simeq \Aut(\sigma)_\Inv$.
 
\end{lemma}
\noindent{\bf Proof.}
Any automorphism $g\in \Aut(\sigma)_{\Inv}$ preserves all vectors from
$\Inv(\sigma)$ and consequently it preserves the vectors
from $\overline{\Inv(\sigma)}:={\rm
lin}(\Inv(\sigma))\cap\sigma=\lin(\overline{\Inv(\sigma)})\cap
\sigma$.
Set
$\overline{\Inv(\sigma_0)}:=\pi^{-1}(\overline{\Inv(\sigma)})=
\langle w_1,\ldots, w_m \rangle$. Then, by definition
$\overline{\Inv(\sigma_0)}=\lin(\overline{\Inv(\sigma_0)})\cap
\sigma_0$.
 
Let $g':=\pi^*(g)$. We shall  show that
$g'$ preserves $\overline{\Inv(\sigma_0)}$.

Let $w_i\in\inte(\tau_i)$, where
$\tau_i=\langle e_{ij_1}, \ldots, e_{ij_l}\rangle\preceq
\sigma_0$. Then , by definition, $w_i$ defines a unique ray
of the cone $\tau_i\cap
\overline{\Inv(\sigma_0)}=\tau_i\cap
\lin(\overline{\Inv(\sigma_0)})$. Then 
$\ker(\pi)\subset \lin(\overline{\Inv(\sigma_0)})$ and
$\ker(\pi)\cap
\lin(\tau_i)=\ker(\pi)\cap\lin(\overline{\Inv(\sigma_0)})\cap
\lin(\tau_i)=\ker(\pi)\cap\lin\{w_i\}=0$. Consequently
$\pi$ maps  $\tau_i$ isomorphically onto a cone
$\pi(\tau_i)=\langle v_{ij_1}, \ldots, v_{ij_l}\rangle$
containing a unique ray $\pi(w_i)\in\overline{\Inv(\sigma)}$.
By Lemma \ref{le: subdivision} such a cone is a face of some 
subdivision $\Delta$ of $\sigma$ satisfying  condition (2)
of the lemma.  Thus
$\pi(\tau_i)$ is a unique face of $\Delta$ containing
$\pi(w_i)$ in its interior. Therefore it is preserved by
$g$. Consequently, $g'$ preserves
$\tau_i$ and $w_i$, and hence all points from
$\overline{\Inv(\sigma_0)}$. 

Now, let $g'$ be an automorphism of
$\sigma_0$ preserving $\overline{\Inv(\sigma_0)}$. Then 
$g'$ defines a permutation $g$ of $v_1,\ldots,v_k$. 

We will show that such a
permutation preserves linear relations between
$v_1,\ldots,v_k$. 

Let $p:=a_{i_1}v_{i_1}+\ldots +a_{i_r}v_{i_r}=
a_{i_{r+1}}v_{i_{r+1}}+\ldots +a_{is}v_{i_s}$ be a minimal
linear relation, i.e. $s-1=\dim\{v_{i_1}, \ldots, v_{i_s}\}$. 
Then both cones $\sigma_1:=\langle v_{i_1}, \ldots,
v_{i_r}\rangle $ and $\sigma_2:=\langle v_{i_{r+1}}, \ldots,
v_{i_\sigma}\rangle $ are simplicial and by Lemma \ref{le: subdivision} 
they are faces of two
subdivisions $\Delta_{1}$ and $\Delta_{2}$ of $\sigma$
satisfying the condition 2 of the lemma. Then
$\langle p \rangle \in \Ver(\Delta_1\cdot\Delta_2)$ and
by Lemma \ref{le: multiplication} and Proposition 
\ref{pr: ssimple},  $p 
\in {\Inv(\sigma)}$. Consequently
the point $p':=a_{i_1}e_{i_1}+\ldots +a_{i_r}e_{i_r}\in 
\overline{\Inv(\sigma_0)}$. This implies
$p'=g'(p')=a_{i_1}g'(e_{i_1})+\ldots +a_{i_r}g'(e_{i_r})$ and
$p=a_{i_1}v_{i_1}+\ldots +a_{i_r}v_{i_r}=a_{i_1}g(v_{i_1})+\ldots
+a_{i_r}g(v_{i_r})$. Analogously $p=a_{i_{r+1}}g(v_{i_{r+1}})+ 
\ldots +a_{is}g(v_{i_s})$. Finally, $a_{i_1}g(v_{i_1})+\ldots
+a_{i_r}g(v_{i_r})=a_{i_{r+1}}g(v_{i_{r+1}})+ 
\ldots +a_{is}g(v_{i_s})$ and $g$ defines a linear
automorphism preserving vectors from
$\overline{\Inv(\sigma)}$. 

We need to show that $g$ satisfies condition (2). 
Let $\Delta$ be a 
subdivision satisfying  (2). If $\tau$ is any face
of $\Delta$ whose relative interior intersects 
${\Inv(\sigma)}\subset \overline{\Inv(\sigma)}$
then it lifts to a face $\tau'$ of $\sigma_0$ whose relative interior
intersect $\overline{\Inv(\sigma_0)}$. Thus $\sigma'$ is
preserved by $g'$ and $\sigma$ is
preserved by $g$. Therefore all
 faces of $\Delta$ whose relative interior intersect
${\Inv(\sigma)}$ are preserved. If the relative interior of $\tau\in
\Delta$ does not intersect ${\Inv(\sigma)}$ then $\tau\in
S(\sigma')\subset \Delta$, where $\sigma'\in \Delta$,  is the
maximal face of $\tau$ whose relative interior intersects
${\Inv(\sigma)}$ (see Lemma \ref{le: Ssigma}). By the above,
$\sigma'$ is preserved by $g$. Thus by
\ref{le: Ssigma2} $g(\tau)\in
S(\sigma',g(\Delta))=S(\sigma',\Delta)=S(\sigma') $. Consequently,  all 
faces of $\Delta$ are mapped to  faces of $\Delta$ and
finally $g$ defines 
an automorphism of $\Delta$. Lemma \ref{le: sigma} is
proven. \qed

Write $X_{\sigma}=X_{\sigma_0}//T_0$ for the torus
$T_0$ corresponding to the kernel of $\pi$. Also
$\widehat{X}_{\sigma}=\widehat{X}_{\sigma_0}//T_0$ and
$K[\widehat{X}_{\sigma}]=K[\widehat{X}_{\sigma_0}]^{T_0}$. 

Let ${\Aut(\widehat{X}_{\sigma_{0}})}_{\Inv}$ denote the
group of all $T_0$-{equivariant} automorphisms $g$
in $\Aut(\widehat{X}_{\sigma_{0}})$ which preserve invariant
valuations, i.e $g_*(\val(v))=\val(v)$
for
any $ v\in\overline{\Inv(\sigma_0)}$. By definition
each automorphism $g\in
{\Aut(\widehat{X}_{\sigma_{0}})}_{\Inv}$ defines an
automorphism of $K[\widehat{X}_{\sigma_0}]^{T_0}$. This
gives  the homomorphism
$p: {\Aut(\widehat{X}_{\sigma_{0}})}_{\Inv}\longrightarrow
\Aut(\widehat{X}_{\sigma})$. 
It now suffices to prove the following Lemma:
\begin{lemma}\label{le: connect}
$\Aut(\widehat{X}_{\sigma_{0}})_{\Inv}$ is
connected.
\end{lemma}
Consequently, $\Aut(\sigma)_{\Inv}\subset 
p(\Aut(\widehat{X}_{\sigma_{0}})_{\Inv})\subset
\Aut(\widehat{X}_{\sigma})^0$ and 
$\Aut(\sigma)_{\Inv}\subset\Aut(\sigma)^0$. \qed

\noindent {\bf Proof of Lemma \ref{le: connect}}.
The same reasoning as in Lemma \ref{le: smooth}. Let
$x_1,\ldots, x_k$ be the standard coordinates on $X_{\sigma_0}={\bf A}^k$.
Any automorphism $g\in
\Aut(\widehat{X}_{\sigma_{0}})_{\Inv}$ is given by
$\Gamma$-semiinvariant functions 
$g^*(x_1),\ldots,g^*(x_k)$ with
 the corresponding $\Gamma$-weights and such that
$\val(v)(g^*(x_i))=\val(v)((x_i))$, where 
$v\in \overline{\Inv}({\sigma_{0}})\cap N_{\sigma_0}$. There is a
birational map $\alpha: {\bf A}^1\to G$ defined by
$$\alpha(z):=(x_1,\ldots,x_k)\mapsto
((1-z)g^*(x_1)+zx_1,\ldots,(1-z)g^*(x_k)+zx_k).$$ 
 $\alpha(z)$
defines a $\Gamma$-equivarinat automorphism of
$\widehat{X}_{\sigma_{0}}$ for all $z$ in  the open subset $U$ of ${\bf A}^1$,
where the linear parts of coordinates of $\alpha(z)$ are
linearly independent. Note that for any integral vector 
$v\in \overline{\Inv}({\sigma_{0}})$ and $z\in U$,
$\val(v)((1-z)g^*(x_i)+zx_i)\geq \val(v)(x_i)$. By Lemma
\ref{le: monomial2}, $\alpha(z)$ preserves $\val(v)$. 
\qed
\begin{example}\label{ex: isolated2} Let $X_\sigma\subset
{\bf A}^4$ be a toric variety described by $x_1x_2=x_3x_4$. Then
$\sigma=\langle v_1,v_2,v_3,v_4 \rangle\subset N^{\bf Q}_\sigma$ is a cone 
over a
square, where $N^{\bf Q}_\sigma:=\{(x_1,x_2,x_3,x_4)\mid
x_1+x_2-x_3-x_4\}\subset {\bf Q}^4$, $v_1=(1,0,0,1)$, $v_2=(0,1,1,0)$, 
$v_3=(1,0,1,0)$ $v_4=(0,1,0,1)$. 
Then $\Aut(\sigma)={\bf D}_8$ consists of all
isometries 
of $\sigma$. $\Aut(\sigma)_{\Inv}$ consists of all isometries preserving
diagonals: the reflections with respect to the diagonals and
rotation through $\pi$. Consequently
$$\Theta(X_\sigma,S,O_\sigma)=\Aut(\sigma)/\Aut(\sigma)_{\Inv}
\simeq {\bf Z}_2.$$
\end{example}

\bigskip
\subsection{Stratified toroidal varieties are orientible
}\label{se: existence}

 \begin{lemma}\label{le: orientation0}
Let 
$\phi_i: X \to X_\sigma$ for $i=1,2$ be
$\Gamma$-smooth morphisms of a $\Gamma$-stratified toroidal
scheme $(X,S)$ to a 
$\Gamma$-stratified toric variety $(X_\sigma,S_\sigma)$
such that strata of $S$ are precisely the inverse images of
strata of $S_\sigma$.
Denote by $\Delta_1,\ldots,\Delta_l$ all
subdivisions of $\sigma$ for which
$\Ver(\Delta_i)=\Ver(\overline{\sigma})$. Let $v_1,\ldots,v_r$ be stable
vectors such that $\lin(v_1,\ldots,v_r)\supseteq
\Inv(\sigma)$. Denote by $\Delta_{l+1},\ldots,\Delta_{l+r}$
the star subdivisions of
$\sigma$ at
$\langle v_1 \rangle,\ldots, \langle v_r \rangle$. 
The following conditions are equivalent:
\begin{enumerate}
\item $\phi_1$ and $\phi_2$ determine the same orientation at 
 a $K$-rational point $x\in s=\phi_i^{-1}(\strat)\sigma))$. 
\item
 There is an open
neighborhood $U$ of $x$ such that
$U_{j1}:=U\times_{X_\sigma}X_{\Delta_j}$, defined via
$\phi_1$, and $U_{j2}:=U\times_{X_\sigma}X_{\Delta_j}$, defined via
$\phi_2$, are isomorphic over $U$ for any $j=1,\ldots,{l+r}$.
\end{enumerate}
\end{lemma} 

\noindent {\bf Proof.}   
$(\Rightarrow)$ Suppose $\phi_1$ and $\phi_2$ determine the same
orientation at $x$. By Proposition \ref{pr: simple} the  
morphism $\widetilde{X}_{\Delta_j}\to \widetilde{X}_\sigma$
is $G^\sigma$-equivariant. Therefore by Lemma \ref{le: isomorphisms2} for
any  subdivision $\Delta_j$ of $\sigma$ there is an open
neighborhood $U_j$ of $x$ such that
$U'_{j1}:=U_j\times_{X_\sigma}X_{\Delta_j}$, defined via
$\phi_1$, and $U'_{j2}:=U_j\times_{X_\sigma}X_{\Delta_j}$, defined via
$\phi_2$, are isomorphic over
$U_j$. It suffices to put $U=\bigcap U_j$.

$(\Leftarrow)$
Suppose $\phi_1$ and $\phi_2$ do not determine the same
orientation at $x$.
 Then by Lemma \ref{le: orientation group} 
we find an automorphism $\phi$ of $\sigma$  such
that for the corresponding automorphism
$\overline{\phi}\in\Aut(X_\sigma)$, $\phi_1$ and
$\overline{\phi}\phi_2$ determine the same
orientation at $x$. By the above there is an open $U\ni x$ for which
$U_{j,\phi_1}$ and $U_{j,\phi\phi_2}$ are isomorphic.
But $\overline{\phi}$ does not preserve orientation so by
 \ref{le: orientation group}, $\phi$ does not satisfy one of
 the conditions (1)-(2)
of the lemma and does not
 lift to an automorphism of some
$\Delta_j\neq \phi(\Delta_j)$. Consequently, 
$U_{j,\phi\phi_2}\simeq U_{j,\phi_1}\simeq U\times_{X_\sigma}X_{\Delta_j}$ is
not isomorphic to $U_{j,\phi_2}\simeq U\times_{X_\sigma}X_{\phi(\Delta_j)}$ for any open neighborhood
$U$ of $x$. \qed

 \begin{lemma}\label{le: orientation1}
Let $(X_\Delta,S_\Omega)$ be a $\Gamma$-stratified toric
variety corresponding
to an embedded semifan $\Omega\subset\Delta$.
Let 
$\phi_i: X \to X_\Delta$ for $i=1,2$ be
$\Gamma$-smooth morphisms of a $\Gamma$-stratified toroidal
scheme $(X,S)$ to a $\Gamma$-stratified toric
variety $(X_\Delta,S_\Omega)$ such that the strata of $S$
are precisely the inverse images of strata of $S_\Omega$.
 
For a cone $\omega\in\Omega$ denote by
${\omega}(\Delta)$ the subset $\{\tau\in \Delta\mid
\omega(\tau)=\omega\}$ of $\Delta$.
Let $\Delta_1,\ldots,\Delta_l$ be all 
subdivisions of $\omega$
satisfying  condition (2) of Lemma
\ref{le: orientation group}.  Let $\Delta_{l+1},\ldots,\Delta_{l+r}$
denote the star subdivisions of
$\omega$ at stable vectors
$\langle v_1 \rangle,\ldots, \langle v_r \rangle$ for which
$\lin(v_1,\ldots,v_r)\supseteq
\Inv(\sigma)$.
Denote by ${\Delta_j}(\Delta)$ the subdivisions of
${\omega}(\Delta)$ induced by ${\Delta_j}$, i.e. for
any $\delta=\omega(\delta)\oplus r(\delta)\in
{\omega}(\Delta)$, 
${\Delta_j}(\Delta)|\delta=\Delta_j\oplus r(\delta)$.
Let $x\in \phi_i^{-1}(\strat(\omega))$
be a $K$-rational point. 
The following conditions are equivalent:
\begin{enumerate}
\item $\phi^*_1(\cU(\Delta,\Omega))$ and
$\phi^*_2(\cU(\Delta,\Omega))$  are compatible at $x$.
 
\item
 There is an open
neighborhood $U$ of $x$ such that
$U_{j1}:=U\times_{X_{{\omega}(\Delta)}}X_{{\Delta_j}(\Delta)}$, defined via
$\phi_1$, and $U_{j2}:U\times_{X_{{\omega}(\Delta)}}X_{{\Delta_j}(\Delta)}
$, defined via
$\phi_2$, 
are isomorphic over $U$ for any $j=1,\ldots,{l+r}$
\end{enumerate}
\end{lemma} 
\noindent {\bf Proof.} Let $\phi_i(x)\in\delta_i$, where
$\delta_i\in\Delta$ and $i=1,2$. Denote by $U_i$ the
inverse image $\phi^{-1}_i(X_{\delta_i})$. Since
$x\in \phi_i^{-1}(\strat(\sigma))$,  we have $\omega{(\delta_i)}=\omega$.
The atlases $\phi^*_i(\cU^\can(X_\Delta,S_\Delta))$
are compatible at $x$ if the morphisms
$\pi_\sigma^{\sigma_i}\phi_i: U_i \to X_\omega$ determine the same
orientation at $x$. This is equivalent by Lemma \ref{le:
orientation0} to the fact
that $U\times_{X_{\sigma}}X_{\Delta_j}$, defined via
$\pi^\sigma_{\sigma_i}\phi_{i}$, are isomorphic for some
$U$. The latter schemes are isomorphic to 
$U\times_{X_{{\sigma}_i}}X_{{\Delta_j}(\Delta)|{\delta_i}}=
U\times_{X_{{\omega}(\Delta)}}X_{{\Delta_j}(\Delta)}$.\qed

 \begin{lemma}\label{le: orientation2} 
Let 
$\phi_i: (X,S) \to (X_\Delta,S_{\Omega})$ for $i=1,2$ be
$\Gamma$-smooth morphisms  of a  $\Gamma$-stratified toroidal scheme to a 
$\Gamma$-stratified toric variety $(X_\Delta,S_{\Omega})$
such that the strata of $S$ are precisely the inverse images
of strata of $S_\Omega$.
Assume that a stratum $s\in S$ is a variety over $K$ and
 $\phi^*_1(\cU(\Delta,\Omega))$ and 
$\phi^*_2(\cU(\Delta,\Omega))$ are compatible at some
point $x\in s$. Then they are compatible along $s$.

\end{lemma}
 
\noindent {\bf Proof.} By Lemma \ref{le:
orientation1},  $\phi^*_1(\cU(\Delta,\Omega))$ and
$\phi^*_2(\cU(\Delta,\Omega))$ are compatible at any
$x\in s$ iff they are compatible at all points of $s\cap
U_x$, where $U_x$  
is an open neighborhood of $x$. Since $s$ is a variety the
sets  $s\cap
U_x$  and  $s\cap U_y$ intersect for any $x,y\in s$.\qed

\begin{lemma} \label{le: openorient} 
Let
$\phi_1, \phi_2: X \to X_{\sigma}$ be two
$\Gamma$-smooth morphisms of a $\Gamma$-stratified toroidal 
variety $(X,S)$ into the $\Gamma$-stratified toric variety $X_{\sigma}$, such that the strata of $S$
are precisely the inverse images of strata of $S_\sigma$.
Assume that $\phi_1, \phi_2$
determine the same orientation at a point $x\in \phi_i^{-1}(\strat(\sigma))$. Then $\phi^*_i(\cU^\can(X_{\sigma},S_{\sigma}))$ and
$\phi^*_2(\cU^\can(X_{\sigma},S_{\sigma}))$ are compatible on $X$.

\end{lemma}
\noindent{\bf Proof.} Let
$\Delta_{1},\ldots,\Delta_{l},\ldots,\Delta_{l+r}$
be subdivisions of $\tau\leq \sigma$ as in
Lemma \ref{le: orientation0}. By Lemma \ref{le: extendsub},
the subdivisions $\Delta_{1},\ldots,\Delta_{l}$ can be
extended to subdivisions ${\Delta}_{1}(\sigma),\ldots,{\Delta}_{l}(\sigma)$
 of $\sigma$ for which $\Ver(\Delta_i)=\Ver(\overline{\sigma})$.
 By Lemma
\ref{le: stableincl} $\Inv(\tau)\subset \Inv(\sigma)$ and
we can extend the star subdivisions
$\Delta_{l+1},\ldots,\Delta_{l+r}$  to the star
subdivisions $\Delta_{l+1}(\sigma),\ldots,\Delta_{l+r}(\sigma)$ of
$\sigma$ at $v_1,\ldots,v_r\in
\Inv(\tau)\subset \Inv(\sigma)$.
By Lemma \ref{le: isomorphisms2} there is an open
neighborhood $U$ of $x$ 
for which $U_{j,\phi_i}=U\times_{X_\sigma}X_{\Delta_j}$
definded for $\phi_1$ and $\phi_2$ are pairwise isomorphic for any $j$.
Let ${\tau}(\sigma)=\{\varrho\preceq \sigma\mid \omega(\varrho)=\tau\}$.
Then there is an open subset
$U'=\phi_1^{-1}(X_{{\tau}(\sigma)})\cap
\phi_2^{-1}(X_{{\tau}(\sigma)})\cap U $ intersecting
$\phi_i^{-1}(\strat(\tau))$ such that 
$U'\times_{X_\sigma}X_{\Delta_j}=
U'\times_{X_{{\tau(\sigma)}}} X_{{\Delta_j(\sigma)}|{{\tau(\sigma)}}}$
defined by the restrictions of $\phi_{i}$ are isomorphic.
This implies by Lemma \ref{le: orientation1} that 
$\phi^*_1(\cU^\can(X_{\sigma},S_{\sigma})$ and
$\phi^*_2(\cU^\can(X_{\sigma},S_{\sigma}))$ are compatible
at $x\in \strat(\tau) \cap U'$ and
consequently by  Lemma \ref{le: orientation2} 
they are compatible at all the points of $\phi_i^{-1}(\strat(\tau))$.\qed

\begin{lemma} \label{le: orientation4} The $\Gamma$-semicomplex
$\Sigma$ 
associated   to an oriented
$\Gamma$-stratified toroidal variety $(X,S)$ 
 is oriented. 
\end{lemma}
\noindent{\bf Proof.}
For any $\sigma\leq \tau$, the inclusion morphism
$\imath^\tau_{\sigma}: \sigma\to \tau$
identifies $\sigma$ with a face of $\tau$. Let
$\overline{\pi}^\sigma_{\tau}: N_\tau\to N_\sigma$ denote a
projection such that $\overline{\pi}^\sigma_{\tau}\circ \imath^\tau_{\sigma}=
\id_{|\sigma}$. Let ${\pi}^\sigma_{\tau}: X_\sigma\to
X_\tau $ denote the
induced toric morphism. Consider faces $\sigma\leq \tau
\leq \varrho$ of the semicomplex $\Sigma$. Let 
$\phi_{\sigma}: U_\sigma \to X_\sigma, \phi_\tau: U_\tau
\to X_\tau, \phi_\varrho: U_\varrho \to X_\varrho$ denote 
any charts on $(X,S)$ corresponding to the cones $\sigma,
\tau$ and $\varrho$.  Since $\tau\leq \varrho$, the open subset
$U_\varrho$ intersects $\tau$. For a generic point $x_\tau$
of $\strat(\tau)\cap U_\varrho $,  the morphism
 $\phi_\varrho$ sends $x_\tau$ to the toric orbit $O_\tau\subset X_\varrho$.
Since $(X,S)$ is oriented we find that
$\pi^\tau_{\varrho}\phi_{\varrho}: U^\tau_\varrho \to X_\tau$
and $\phi_{\tau}: U_\tau \to X_\tau$
determine the same orientation at $x_\tau$. Analogously 
$\pi^\sigma_{\varrho}\phi_{\varrho}$, $\phi_{\sigma}$ and
$\pi^\sigma_{\tau}\phi_{\tau}$ determine the same orientation 
for a generic point $x_\sigma$ of $\strat(\sigma)$. 
Applying Lemma \ref{le: openorient} to $\pi^\tau_{\varrho}\phi_{\varrho}$ and $\phi_{\tau}$ we see that
$\pi^\sigma_{\tau}\pi^\tau_{\varrho}\phi_{\varrho}$ and
$\pi^\sigma_{\tau}\phi_{\tau}$ determine the same orientation for
a generic point of $\strat(\sigma)$. Finally
$\pi^\sigma_\tau\pi^\tau_\varrho\phi_{\varrho}$ and
$\pi^\sigma_\varrho\phi_{\varrho}$ determine the same orientation and
thus so do $\pi^\sigma_\tau\pi^\tau_\varrho$ and 
$\pi^\sigma_\varrho$.
 Write $\imath^\varrho_\tau \imath^\tau_\sigma=\imath^\varrho_\sigma
\alpha_{\sigma}$,
where $\alpha_{\sigma}\in {\rm Aut}(\sigma)$.  
 Then $\alpha_{\sigma}{\pi}^\sigma_\tau{\pi}^\tau_\varrho 
\imath^\varrho_\sigma= 
\alpha_{\sigma}{\pi}^\sigma_\tau{\pi}^\tau_\varrho
\imath^\varrho_\tau \imath^\tau_\sigma\alpha_{\sigma}^{-1}=
\id_{|\sigma}$.
By Lemma \ref{le: sections} the toric morphisms induced by
the projections
$\alpha_{\sigma}{\pi}^\sigma_\tau{\pi}^\tau_\varrho$ and
${\pi}^\sigma_\varrho$ determine the same orientation. Thus the 
morphisms induced by
$\alpha_{\sigma}{\pi}^\sigma_\tau{\pi}^\tau_\varrho$ and 
${\pi}^\sigma_\tau{\pi}^\tau_\varrho$ determine the same
orientation, which gives $\alpha_{\sigma}\in {\rm Aut}(\sigma)^0$. 
\qed 
\bigskip

\begin{lemma} \label{le: existence} 
On any $\Gamma$-stratified toroidal
 variety $(X,S)$ there exists an atlas  ${\cU}$ such that
$(X,S,{\cU})$ is oriented.
\end{lemma}

\noindent{\bf Proof.} We shall improve the charts from
${\cU}$ and inclusion maps in the associated $\Gamma$-semicomplex
$\Sigma$.
For any $\sigma \in \Sigma$ fix a chart $\phi_{\sigma,a_\sigma}:
U_{\sigma,a_\sigma}\to X_{\sigma}$.
 By Lemma \ref{le: surjection} for any chart $\phi_{\sigma,a}:U_{\sigma,a}\to X_{\sigma}$
we can find $\psi_{\sigma,a}$ of $X_{\sigma})$ induced by
an automorphism of $\sigma$ 
such
that $\phi_{\sigma,a_\sigma}$ and
$$\overline{\phi}_{\sigma,a}:=\psi_{\sigma,a}\circ\phi_{\sigma,a}$$
determine the same orientation for some points of $\strat(\sigma)\cap U_{\sigma,a_\sigma}\cap
U_{\sigma,a}$. 
By Lemma \ref{le: orientation2} the above morphisms determine 
the same orientation for all points from $\strat(\sigma)\cap U_{\sigma,a_\sigma}\cap
U_{\sigma,a}$. 

Let $\sigma\leq \tau$. Denote by $\sigma'\preceq\tau$ the
face cone in $N_\tau$ such that $\underline{\sigma'}=\sigma$. 
By Lemma \ref{le: surjection}  find 
$\alpha_{\sigma}\in {\rm Aut}(\sigma)$ such that
for the inclusion map 
$\overline{\imath^\tau_{\sigma'}}:={\imath^\tau_{\sigma}}
\circ\alpha_\sigma$
 any two charts $\phi_{\sigma,a_\sigma}: U_{\sigma,a_\sigma}\to 
X_{\sigma}$, 
$\pi_{\sigma'}^{\sigma}\phi^{\sigma'}_{\tau,a_{\tau}}: 
U^{\sigma'}_{\tau,a_{\tau}}\to
X_{\sigma}$ define the same orientation at some point of  
$x\in \strat(\sigma)\cap
U_{\sigma,a_\sigma}\cap U^{\sigma'}_{\tau,a_{\tau}} 
$. 
By Lemma \ref{le: orientation2} they determine the same
orientation at all points of $\strat(\sigma)\cap
U_{\sigma,a_\sigma}\cap U^{\sigma'}_{\tau,a_{\tau}} 
$. We need to verify that for any two charts $\phi_{\sigma,a}$ and
$\phi_{\tau, b}$, where $\sigma\leq\tau$, and any
$\tau'\preceq \tau$ such that $\omega(\tau')=\sigma$ the morphisms
$\phi_{\sigma,a}$ and
$\phi^{\tau'}_{\tau,b}$ determine the same orientation
for all $x\in \strat(\sigma)\cap U_{\sigma,a}\cap U^{\tau'}_{\tau,b}$. 
 By Lemma \ref{le: orientation2} it suffices to show that the
above morphisms  determine the same orientation for some $x\in 
\strat(\sigma)\cap U_{\sigma,a}\cap
U_{\sigma,a_\sigma} \cap U^{\tau'}_{\tau,b}\cap U^{\sigma}_{\tau,a_\tau}$. 
By the above $\phi_{\sigma,a}$ and $\phi_{\sigma,a_\sigma}$
determine the same orientation at all $x\in
\strat(\sigma)\cap U_{\sigma,a}\cap
U_{\sigma,a_\sigma}$. Also $\phi_{\tau,b}$ and $\phi_{\tau,a_\tau}$
determine the same orientation at all  $x\in \strat(\tau)\cap U_{\tau,b}\cap
U_{\tau,a_\tau}$. By Lemma \ref{le: openorient},  
$\phi_{\tau,b}^*(\cU^\can(X_\tau,S_\tau))$ and 
$\phi_{\tau,a_\tau}^*(\cU^\can(X_\tau,S_\tau))$ are
compatible on $U_{\tau,b}\cap
U_{\tau,a_\tau}$. Hence $\phi^{\tau'}_{\tau,b}$ and 
$\phi^{\sigma'}_{\tau,a_\sigma}$ determine the same
orientation at points $x\in 
\strat(\sigma) \cap U^{\tau'}_{\tau,b}\cap U^{\sigma'}_{\tau,a_\tau}$.
Finally $\phi_{\sigma,a}$ and
$\phi^{\tau'}_{\tau,b}$ determine the same orientation for
all $x\in \strat(\sigma)\cap U_{\sigma,a}\cap U^{\tau'}_{\tau,b}$.
\qed

\bigskip
\subsection{Resolution of singularities
of  toroidal varieties in arbitrary characteristic}\label{se: resolution} 
We call a 
$\Gamma$-semicomplex $\Sigma$ {\it regular} if all cones
$|\sigma|$, where $\sigma\in\Sigma$ are regular. 
\begin{proposition} \label{pr: Danilov} For any oriented 
$\Gamma$-semicomplex $\Sigma$ there exists a sequence of
stable vectors 
$v_1,\ldots,v_m\in {\rm Stab}(\Sigma) $  such that the subdivision
 $\Delta:=\langle v_k\rangle \cdot\ldots\cdot \langle v_1
\rangle\cdot\Sigma$ of $\Sigma$ is a regular
$\Gamma$-semicomplex. 
\end{proposition}

\noindent{\bf Proof.} 
Let $\sigma_1,\ldots,\sigma_k$ denote all the faces of $\Sigma$.
We shall construct by induction on $i$ the canonical
subdivision $\Delta_i$ of $\Sigma$ such that the
fan ${\Delta_i^{\sigma_j}}$ is regular for all
$j\leq i$.  Suppose $\Delta_i$ is
already constructed. 
By \cite{KKMS} or \cite{Danilov1} 
there exist centers $v^i_1,\ldots,v^i_{k_i}$  for the fan
$\Delta^{\sigma_{i+1}}_i$ in $N_{\sigma_{i+1}}$ such that
$\langle v^i_{k_i} \rangle\cdot\ldots\cdot\langle v^i_1
\rangle\cdot\Delta^{\sigma_{i+1}}_i$ is regular.
 The subsequent centers $v^i_j$ of star subdivisions coul be chosen to be  minimal internal vectors in indecomposable faces of
 $\langle v^i_{j-1}\rangle \cdot\ldots\cdot \langle v^i_1
\rangle\cdot\Delta^{\sigma_{i+1}}_i$. 
Therefore by Lemma \ref{le: minimal vectors}   all centers
in the desingularization process are $\Sigma$-stable. Hence by
Proposition \ref{pr: blow-ups}, $\Delta_{i+1}:
=\langle v^i_k \rangle\cdot\ldots\cdot\langle v^i_1
\rangle\cdot\Delta_{i}$ is canonical. By construction 
$\Delta^{\sigma_{i+1}}_{i+1}$ is regular. Note that 
$\Delta^{\sigma_{j}}_{i+1}=\Delta^{\sigma_j}_{i}$ for
$j<i+1$ since during the desingularization regular cones
remain unaffected. \qed

\begin{theorem} \label{th: resolution}
 For any $\Gamma$-stratified toroidal variety $(X,S)$
(respectively a $\Gamma$-toroidal variety $X$) there
exists a sequence of $\Gamma$-equivariant blow-ups at locally monomial valuations  
$\bl_{\nu_k}\circ\ldots\circ \bl_{\nu_1}(X)$ which is a
 resolution of singularities.
\end{theorem}

\noindent{\bf Proof.} In the case of a nonstratified
toroidal variety  consider $X$ with stratification ${\rm
Sing}^\Gamma(X)$. By Lemma \ref{le: existence} we can assume that
the $\Gamma$-stratified toroidal variety is oriented. 
By Lemma \ref{le: orientation2} the associated semicomplex
is oriented.
By  
Proposition \ref{pr: Danilov} and Proposition \ref{pr:
blow-ups}, $\bl_{\nu_k}\circ\ldots\circ \bl_{\nu_1}(X)$ is smooth.  \qed 

\section{Orientation of toroidal modifications}
\subsection{Lifting group actions}
Let $U$ be an affine variety. Let $\widehat{X}_x$ be a
completion of a variety $X$ at its closed point $x$. Let
$\Theta: Y\to X$ be a proper morphism and 
$\widehat{Y}_x:=\widehat{X}_x\times_XY$. 

Set $$U\widehat{\times}\widehat{X}_x:=\Spec(\lim_{\leftarrow}
K[U]\otimes \cO_x^n/m_{X,x}),$$  
$$U\widehat{\times}\widehat{Y}_x:={(U\widehat{\times}\widehat{X}_x)}\times_XY.
$$
\begin{lemma}\label{le: normal0}
If $U$, $X$ are normal varieties then the scheme
$U\widehat{\times}\widehat{X}_x$ is normal. If $U$, $Y$ are
normal varieties then the scheme
$U\widehat{\times}\widehat{Y}_x$ is normal.

\end{lemma}
\noindent{\bf Proof.} 
(1)  $\pi_X: U\widehat{\times}\widehat{X}_x \to
U\times \widehat{X}_x$ is \'etale at $(u,x)$ for $u\in U$ 
since it determines an
isomorphism of the completions of the local rings at $(u,x)$. $U\times
\widehat{X}_x$ is normal therefore $U\widehat{\times}\widehat{X}_x$ is
normal at $(u,x)\in U\times\{x\}$. But all closed points of 
$U\widehat{\times}\widehat{X}_x$ are in $
U\times\{x\}\subset U\widehat{\times}\widehat{X}_x$. The morphism 
$\pi$ is \'etale in a neighborhood of
$U\widehat{\times}\widehat{X}_x$. Such a neighborhood is
equal to $U\widehat{\times}\widehat{X}_x$ since its
complement if nonempty would contain some closed points.

(2) $\pi_Y: U\widehat{\times}\widehat{Y}_x \to
U\times \widehat{Y}_x$ is \'etale at $(u,y)$ for $u\in U$
and $y\in \theta^{-1}(x)$ since $\pi_Y$ is a pull-back of
$\pi_X$. The rest of the reasoning is the same. \qed

\begin{definition}
\label{de: lifting}
Let a proalgebraic group $G$ act on $\widehat{X}_x$. Let 
$\Phi:G\widehat{\times}\widehat{X}_x\to\widehat{X}_x$ be
the action morphism and 
$\Psi:G\widehat{\times}\widehat{X}_x\to 
G\widehat{\times}\widehat{X}_x$
be the action automorphism (see Lemma \ref{le: group action}). 
Let $\theta: U\to G$
be a morphism from a normal algebraic variety $U$  to $G$. 

Then by the {\it action morphism with
respect to $U$} or simply {\it  action morphism} we mean the
induced morphism
$\Phi_{U,\widehat{X}_x}:U\widehat{\times}\widehat{X}_x\to\widehat{X}_x$.

By the {\it action automorphism with
respect to $U$} or simply {\it action automorphism} we mean the
induced morphism
$\Psi_{U,\widehat{X}_x}:U\widehat{\times}\widehat{X}_x\to 
U\widehat{\times}\widehat{X}_x$.

Let $\alpha: Y\to\widehat{X}_x$ be a $G^K$-equivariant
birational morphism. 
By the {\it lifting} of the action morphism (resp. automorphism)
we mean the morphism
$\Phi_{U,Y}:U\widehat{\times}Y\to Y$ (resp. the automomophism
$\Psi_{U,Y}:U\widehat{\times}Y\to U\widehat{\times}Y$)
 commuting with the induced morphism $\id\widehat{\times}\alpha:U\widehat{\times}Y\to
U\widehat{\times} \widehat{X}_x$ and with $\alpha$.

\end{definition}

\begin{lemma}\label{le: blow-up} Let $G$ be a 
proalgebraic group acting on $\widehat{X}_x$. Let $U\to G$
be any morphism from a normal algebraic variety $U$.
Let $I$ be
a $G$-invariant ideal on $\widehat{X}_x$. Let $\phi: Y\to
\widehat{X}_x$ be the normalization of the blow-up at
$I$.  Then there is a lifting of the action automorphism
(and action morphism)
$\Psi_{U,\widehat{X}_x}: U\widehat{\times} \widehat{X}_x \to
U\widehat{\times} \widehat{X}_x$ to $\Psi_{U,Y}: U\widehat{\times}
Y \to U \widehat{\times} Y$.
\end{lemma}

\noindent {\bf Proof.} Let $\Phi:=\Phi_{U,\widehat{X}_x}: U\widehat{\times}
\widehat{X}_x \to \widehat{X}_x$ denote the action morphism. 
Set $\Psi:=\Psi_{U,\widehat{X}_x}$. By Lemma \ref{le: group
ideal} we know  
that $J=:p^*(I)=\Phi^*(I)=\Psi^*p^*(I)=\Psi^*(J)$ is
preserved by the action automorphism.
  Therefore $J$ lifts to the blow-up of 
$U\widehat{\times} \widehat{X}_x$ at $J$ and to its
normalization. 
Since $U\widehat{\times} \widehat{X}_x\to
\widehat{X}_x$ is \'etale, the normalization of the blow-up of 
$U\widehat{\times} \widehat{X}_x$ at $J$ is isomorphic to 
$U\widehat{\times} Y$. 
\qed

\begin{proposition} \label{pr: liftings} Let $\sigma$ be a cone of 
maximal dimension in $N_\sigma$. Let $G$ be a connected
proalgebraic group acting on $X:=\widehat{X}_\sigma$. 
Let $g:U\to G$ be a morphism from a normal algebraic
variety $U$. Let $G^K$ act on
$Y:=\widehat{X}_\Sigma$ such that $\psi: Y\to X$ is a
$G^K$-equivariant  proper birational
morphism. Then there is a lifting of the action automorphism
$\Psi_{U,X}: U\widehat{\times} X \to U\widehat{\times} X$ to
$\Psi_{U,Y}: U\widehat{\times} Y \to U\widehat{\times} Y$.
\end{proposition}

\noindent {\bf Proof.} Let $\Psi_{U,Y}: U\times Y {{-}\to} U\times Y$ be
a birational map which is a lifting of $\Phi_{U,Y}$. By
Lemmas \ref{le: ideal} and \ref{le: blow-up} we can find a
factorization $Z\to Y\to X$ giving the diagram of proper morphisms

\[\begin{array}{rcccccc}
&U\widehat{\times} Z&&\buildrel \Psi_{U,Z}\over\longrightarrow && 
U\widehat{\times} Z&\\  
&&\searrow&&\swarrow&\\
&\downarrow&&W&&\downarrow&\\
&&\swarrow f_1&&\searrow f_2&&\\
&U\widehat{\times} Y&&\buildrel \Psi_{U,Y} \over {{-}\to} && 
U\widehat{\times} Y&\\  
&\downarrow\phi_{U\widehat{\times} Y}&&&&
\downarrow\phi_{U\widehat{\times} Y}&\\
&U\widehat{\times} X&&\buildrel \Psi_{U,X} \over \longrightarrow&&
U\widehat{\times} X& 
\end{array}\]

\noindent  where $W$ is the irreducible component of the
fiber product of $U\widehat{\times}
Y\times_{U\widehat{\times} X} U\widehat{\times} Y$  which
dominates $U\widehat{\times} Y$. (The fiber product is 
given by the morphisms $\phi_{U\widehat{\times} Y}$ and 
$\Psi_{U,X}\phi_{U\widehat{\times} Y}.)$
Let $u\in U$ be a closed point. Take a pull-back of the above
diagram via the closed embedding: 
$\{u\}\widehat{\times} X\to U\widehat{\times} X$
\[\begin{array}{rcccccc}
&\{u\}\widehat{\times} Z&&\buildrel \Psi_{u,Z}\over\longrightarrow && 
\{u\}\widehat{\times} Z&\\  
&&\searrow&&\swarrow&\\
&\downarrow&&W_u&&\downarrow&\\
&&\swarrow &&\searrow &&\\
&\{u\}\widehat{\times} Y&&\buildrel \Psi_{u,Y} \over {\longrightarrow} && 
\{u\}\widehat{\times} Y&\\  
&\downarrow&&&&\downarrow&\\
&\{u\}\widehat{\times} X&&\buildrel \Psi_{u,X} \over \longrightarrow&&
\{u\}\widehat{\times} X& 
\end{array}\]
The birational map $\Psi_{u,Y}$ is the isomorphism defined
by the
action of $g(u)\in G$ on $Y$. $W_u$ consists of some
components of  the
fiber product  $\{u\}\widehat{\times}
Y\times_{\{u\}\widehat{\times} X} \{u\}\widehat{\times} Y$ 
. In particular it contains the irreducible component dominating
$\{u\}\widehat{\times} Y$ which is isomorphic to $Y$. On
the other hand $\{u\}\widehat{\times} Z\to W_u$ is a proper
surjective morphism. Consequently, $W_u$ is irreducible and
isomorphic to $Y$. By Lemma \ref{le: normal0},
$U\widehat{\times}Y$ is normal. Thus $f_i:W\to
U\widehat{\times}Y$ is a proper birational
morphism onto a  normal scheme which is bijective 
on the set of closed ($K$-rational) points. By the Zariski theorem
$f_i$ is an isomorphism and
$\Psi_{U,Y}=f_1^{-1}f_2$ is an action automorphism. \qed

\bigskip
\subsection{Orientation of toroidal modifications}

\begin{proposition} \label{pr: canonical orientation} Let  $(X,S)$ be an
oriented $\Gamma$-stratified toroidal variety 
with the associated $\Gamma$-semicomplex $\Sigma$. Let $f:(Y,R)\to (X,S)$
be the toroidal morphism associated to a 
canonical  subdivision $\Delta$  of $\Sigma$.
Then
\begin{enumerate}
\item The stratified toroidal variety $(Y,R)$ is \underline{oriented}
with charts from
$\bigcup_{\sigma\in\Sigma}\phi_\sigma^{f*}(\cU(\Delta^\sigma,\Delta^\sigma_\stab))$ 
induced by the charts $\phi_\sigma: U\to X_\sigma$ on $X$
and the associated oriented semicomplex $\Sigma_R=\Delta_\stab$.
\item 
For any  point $x$ in a stratum $s\in S$   every
$\Gamma_s$-equivariant automorphism 
$\alpha$ of $\widehat{X}_x$  preserving strata  and orientation
 lifts
to a  $\Gamma_s$-equivariant automorphism
$\alpha'$ of $Y\times_X\widehat{X}_x$ such that the atlases
\\
${\alpha'}^*\hat{\phi}_\sigma^{f*}
(\cU(\Delta^\sigma,\Delta^\sigma_\stab))$ and $\hat{\phi}_\sigma^{f*}
(\cU(\Delta^\sigma,\Delta^\sigma_\stab))$ are
compatible along $f^{-1}(x)\subset Y\times_X\widehat{X}_x$.
Here $\hat{\phi}_\sigma^{f}: Y\times_X\widehat{X}_x\to
X_{\Delta^\sigma}$ denotes the morphism induced by ${\phi}_\sigma^{f}$.

\item The stable support of the oriented $\Gamma$-semicomplex
$\Sigma_R=\Delta_\stab$ is contained in the stable support of $\Sigma$. 
 
\end{enumerate}
\end{proposition}

\noindent{\bf Proof.}$(\Leftarrow)$ 
Let $\phi_{\sigma,1}:U_{\sigma,1}\to
X_{\sigma}$ and $\phi_{\tau,2}: U_{\tau,2}\to
X_{\tau}$ be two charts on $X$, where $\sigma\leq
\tau$. For (1) we have to show that
$\phi_{\sigma,1}^{f*}(\cU(\Delta^\sigma, \Delta^\sigma_\stab))$
and $\phi_{\tau,2}^{f*}(\cU(\Delta^{\tau}, \Delta^\tau_\stab))$
are compatible on $f^{-1}(U_{\sigma,1}\cap U_{\tau,2})$ along
$f^{-1}(x)$ for any point $x\in
\strat(\sigma)\cap(U_{\sigma,1} \cap U_{\tau,2})$. 

Let $\tau'\preceq \tau$ denote a face of $\tau$ such that
${\omega}(\tau')=\sigma$ and  $\phi_{\tau,2}(x)\in O_{\tau'}$.
Then $\tau'=\sigma\times r(\tau')$, where
$r(\tau')$ denotes a regular cone and
$\tau'\setminus\sigma$ is disjoint from $\Stab(\Sigma)$.
By Lemma \ref{le: sum},
$\Delta^\tau|\tau'=\Delta^\tau|\sigma
\oplus r(\tau')=\Delta^\sigma\times r(\tau')$. Also,
$\Delta^\tau_\stab|\tau'=\Delta^\sigma_\stab\times\{0\}$.   

Denote by $\phi^{\tau'}_{\tau,2}: U^{\tau'}_{\tau,2}\to
X_{\tau'}$,   the
restriction of $\phi_{\tau,2}$ and by
$\phi^{\tau'f}_{\tau,2}$ its lifting.
Then
$\phi^{\tau'f*}_{\tau,2}(\cU(\Delta^\tau|\tau', \Delta_\stab^\tau|\tau'))
$ is the restriction of $\phi^{f*}_{\tau,2}
(\cU(\Delta^\tau, \Delta_\stab^\tau))
$ to the open set $f^{-1}(U^{\tau'}_{\tau,2})$ containing
$f^{-1}(x)$. 
Set $U_{\sigma,2}:=U^{\tau'}_{\tau,2}$, 
$\phi_{\sigma,2}:={\pi}^\sigma_{\tau'}\phi^{\tau'}_{\tau,2}:
U_{\sigma,2}\to
X_{\sigma}$, where ${\pi}^{\sigma}_{\tau'}: X_\tau'
\to
X_\sigma$ is defined by any projection which is identical on
$\sigma\preceq \tau'$.  Let
$\phi^f_{\sigma,2}:={\pi}^{\sigma f}_{\tau'}\phi^{\tau'
f}_{\tau,2}: f^{-1}(U_{\sigma,2})\to
X_{\sigma}$ be the lifting of $\phi_{\sigma,2}$. Then
$$\phi^{\tau'f*}_{\tau,2}(\cU(\Delta^\tau|\tau',\Delta_\stab^\tau|\tau'))=
\phi^{\tau'f*}_{\tau,2}\cU(\Delta^\sigma\times
r(\tau'), \Delta_\stab^\sigma\times
\{0\})$$ \noindent is compatible along $f^{-1}(x)$ with 
$$\phi^{\tau'f*}_{\tau,2}
{\pi}^{\sigma f*}_{\tau'}(\cU(\Delta^\sigma, \Delta_\stab^\sigma))=
\phi^{f*}_{\sigma,2}(\cU(\Delta^\sigma, \Delta_\stab^\sigma)).$$ 
It suffices to show that
$\phi^{f*}_{\sigma,1}(\cU(\Delta^\sigma,
\Delta_\stab^\sigma))$ and 
$\phi^{f*}_{\sigma,2}(\cU(\Delta^\sigma,
\Delta_\stab^\sigma))$ are compatible along $f^{-1}(x)$.

Let $U\subset U_{\sigma,1}\cap U_{\sigma,2}$,  be an open
subset for which there exist \'etale
extensions $\widetilde{\phi}_{i}: U \to
X_{\widetilde{\sigma}}$ of $\phi_{\sigma,i}$, $i=1,2$.
Assume that
$\widetilde{\phi}_{i}(x)=O_{\widetilde{\sigma}}=
O_{{\sigma}\times
\reg(\sigma)}$, where $\reg(\sigma)$ is a regular cone of
dimension $k=\dim(\strat(\sigma))$.
Let $\widetilde{\phi}_{fi}$ be an \'etale extension 
of ${\phi}_{fi}$ which is a lifting of $\widetilde{\phi}_{i}$.

Consider the fiber squares

\[\begin{array}{rccccccccccccc}

&\widetilde{\phi}_i:&U&\to&
X_{\widetilde{\sigma}}&=X_{\sigma\times
\reg(\sigma)}&&&&
\widehat{\phi}_i:&\widehat{X}_x&\to& \widetilde{X}_{\sigma}&\\
&&\uparrow&&\uparrow f_{\widetilde{\Delta}}
&&&&&&\uparrow&&\uparrow \widetilde{f}_{{\Delta}}&\\
&\widetilde{\phi}_{fi}:&f^{-1}(U)&\to&
X_{\widetilde{\Delta}^\sigma}&=X_{\Delta^\sigma\times
\reg(\sigma)}&&&&\widehat{\phi}_{fi}:
&Y\times_X \widehat{X}_x&\to & \widetilde{X}_{\Delta^\sigma}&

\end{array}\]

The morphism $\widetilde{f}_{{\Delta}}$   
is a pull-back of the morphism $f_{\widetilde{\Delta}}$. Therefore
$\widetilde{f}_{{\Delta}}^{-1}(O_{\widetilde{\sigma}})$ is
a  $K$-subscheme of $ \widetilde{X}_{\Delta^\sigma}$. 

The automorphism
$\widehat{\phi}:=\widehat{\phi}_{2}\widehat{\phi}^{-1}_{1}$
of $\widetilde{X}_{\sigma}\simeq \widehat{X}_x$
preserves strata and orientation.
For the proof of conditions (1) and (2) we need to show
that any such automorphism $\widehat{\phi}$ preserving
strata and orientation 
lifts to the automorphism
$\widehat{\phi}^f=\widehat{\phi}^f_{2}(\widehat{\phi}_{1}^{f})^{-1}$
of $\widetilde{X}_{\Delta^\sigma}\simeq Y\times_X \widehat{X}_x$
such that
$
\widehat{\phi}^{f*}(\cU(\Delta^\sigma,
\Delta^\sigma_\stab)$ and $\cU(\Delta^\sigma, \Delta^\sigma_\stab)$ are compatible along 
$\widetilde{f}_{{\Delta}}^{-1}(O_{\widetilde{\sigma}_\sigma})$.

By Lemma \ref{le: group structure} there is  
an action morphism
$\Phi: W \widehat{\times}\widetilde{X}_{\sigma} \to
\widetilde{X}_{\sigma}$ such that

\begin{enumerate}
\item $\Phi_{e}:=\Phi_{|e\widehat{\times}\widetilde{X}_{\sigma}}=
\id_{X_{\sigma}}$. 

\item There exists $g\in W$ such that
$\Phi_{g}:=\Phi_{|g\widehat{\times}\widetilde{X}_{\sigma}}
=\widehat{\phi}$.

By Proposition \ref{pr: liftings}, $\Phi$ can be lifted to  ${\Phi}^f
: W\widehat\times 
\widetilde{X}_{\Delta^\sigma}\to
\widetilde{X}_{\Delta^\sigma}$ such that

\item
${\Phi}^f_{e}=\id_{\widetilde{X}_{\Delta^\sigma}}$.

\item 
${\Phi}^f_{g}=\widehat{\phi}^f$.

\end{enumerate}

Let ${\Pi}: W\widehat\times 
\widetilde{X}_{\sigma}\to
\widetilde{X}_{\sigma}$ and ${\Pi}^f: W\widehat\times 
\widetilde{X}_{\Delta^\sigma}\to
\widetilde{X}_{\Delta^\sigma}$ denote the natural projections. 

Note that
${\Pi}^{-1}(O_{\widetilde{\sigma}})={\Phi}^{-1}(O_{\widetilde{\sigma}})=
W\widehat{\times}{O}_{\widetilde{\sigma}}\simeq
W{\times}{O}_{\widetilde{\sigma}}\simeq W$.
The morphisms $\Pi^f$ and $\Phi^f$ are pull-backs of 
$\Pi$ and $\Phi$ induced by $\widetilde{f}_\Delta$. Thus
$({\Pi}^{f})^{
-1}(f_{\widetilde{\Delta}}^{-1}(O_{\widetilde{\sigma}}))=
({\Phi}^{f})^{ -1}(f_{\widetilde{\Delta}}^{-1}(O_{\widetilde{\sigma}})=
W\widehat{\times}f_{\widetilde{\Delta}}^{-1}(O_{\widetilde{\sigma}})
 \simeq
W{\times}f_{\widetilde{\Delta}}^{-1}(O_{\widetilde{\sigma}})$ 
is a $K$-subscheme of $W\widehat{\times}\widetilde{X}_{\Delta^\sigma}$.

For the sake of our considerations we shall enrich
the stratification $R$ by adding to strata
$r\in R$ their intersections with $f^{-1}(x)$.
Set $\overline{R}=\{r\setminus f^{-1}(x)\mid r\in R\}
\}\cup\{r \cap f^{-1}(x)\mid r\in R \}$.

Strata of $R$ on $f^{-1}(U)$ correspond to
the embedded semifan $\Delta^\sigma_\stab\subset
\Delta^\sigma\times \reg(\sigma)$. Strata of $\overline{R}$
correspond to the embedded semifan $\Omega\subset \Delta^\sigma\times
\reg(\sigma)$, where $\Omega:= \Delta^\sigma_\stab \cup 
\{\omega\times\reg(\sigma)\mid \omega\in
\Delta^\sigma_\stab, 
\inte(\omega)\subset\inte(\sigma)\}$.

For any stratum $\strat_Y(\omega)\in R$, $\omega\in\Delta^\sigma_\stab$,
 dominating the stratum 
$\strat_X(\sigma)$ let $\strat(\omega)$ denote the
corresponding stratum on $\widetilde{X}_{\Delta^\sigma}$.
Let $\widetilde{\omega}:=
\omega\times\reg(\sigma)\in\Omega$. Then
$\strat(\omega)\cap
f_{\widetilde{\Delta}}^{-1}(O_{\widetilde{\sigma}})= 
\strat(\widetilde{\omega})\subset 
f_{\widetilde{\Delta}}^{-1}(O_{\widetilde{\sigma}})
\subset \widetilde{X}_{\Delta^\sigma}$ is a
$G^{\sigma}$-invariant stratum corresponding to
$\widetilde{\omega}\in\Omega$
. Since $\strat(\widetilde{\omega})$ is an irreducible locally closed subset of 
$f_{\widetilde{\Delta}}^{-1}(O_{\widetilde{\sigma}})$ it is an
algebraic variety. 

The natural isomorphism $(\Phi^{f})^{-1}(\strat(\widetilde{\omega}))=
(\Pi^{f})^{-1}(\strat(\widetilde{\omega}))= W\widehat{\times}
\strat(\widetilde{\omega})\longrightarrow
W{\times}\strat(\widetilde{\omega}) $ is a pull-back of the
isomorphism 
$W\widehat{\times}{O}_{\widetilde{\sigma}}\longrightarrow
W{\times}{O}_{\widetilde{\sigma}}$ induced by $
f_{\widetilde{\Delta}}^{-1}(O_{\widetilde{\sigma}})\to 
O_{\widetilde{\sigma}}$.

Thus $W\widehat{\times}\strat(\widetilde{\omega})=
W{\times}\strat(\widetilde{\omega})$ is an algebraic variety.

By definition we have $\Psi^f_{e}= \Pi^f_{e}=
\id_{|\strat(\widetilde{\omega})}$. Hence by Lemma \ref{le: sections2}
applied to the natural projection of 
$W\widehat{\times}\strat(\widetilde{\omega})$
 on $W$
the collections of charts $\Psi^{f*}
(\cU^\can(\widetilde{X}_{\Delta^\sigma},S_\Omega))$
and $\Pi^{f*}(\cU^\can(\widetilde{X}_{\Delta^\sigma}, S_\Omega))$ are compatible along
$\{e\}\widehat{\times}\strat(\widetilde{\omega})$.

Hence by Lemma \ref{le: orientation2},  
$\Psi^{f*}(\cU^\can(\widetilde{X}_{\Delta^\sigma}, S_\Omega))$
and $\Pi^{f*}(\cU^\can(\widetilde{X}_{\Delta^\sigma}, S_\Omega))$ are compatible along
$\{g\}\widehat{\times}\strat(\widetilde{\omega})$. Again by
Lemma \ref{le: sections2},
$\Psi_{g}^{f*}(\cU^\can(\widetilde{X}_{\Delta^\sigma}, S_\Omega))=
\widehat{\psi}^{f*}(\cU^\can(\widetilde{X}_{\Delta^\sigma},
S_\Omega))$
and
$\Pi_{g}^{f*}(\cU^\can(\widetilde{X}_{\Delta^\sigma}, 
S_\Omega))=
\cU^\can(\widetilde{X}_{\Delta^\sigma},
S_\Omega)$ are compatible
along $\strat(\widetilde{\omega})$. The charts
$\phi_{\omega'}$, $\omega'\geq\omega$,  of 
$\cU^\can(\widetilde{X}_{\Delta^\sigma},
S_{\Delta^\sigma_\stab})$ 
are obtained by composing the charts
$\phi_{\widetilde{\omega'}}$ of 
$\cU^\can(\widetilde{X}_{\Delta^\sigma},
S_\Omega)$  with the natural projection
$X_{\widetilde{\omega'}}\to X_{\omega'}$. Thus  the atlases
$\cU^\can(\widetilde{X}_{\Delta^\sigma},
S_{\Delta^\sigma_\stab})$ and 
$\widehat{\psi}^{f*}(\cU^\can(\widetilde{X}_{\Delta^\sigma},
S_{\Delta^\sigma_\stab})$
are compatible
along $\strat({\omega})\cap 
f_{\widetilde{\Delta}}^{-1}(O_{\widetilde{\sigma}})=\strat(\widetilde{\omega})$.

(3) By Lemma \ref{le: orientation4}, $\Sigma_R=\Delta_\stab$ is an
oriented $\Gamma$-semicomplex and the notion of stability
makes sense for $\Sigma_R$. Let $v\in \inte(\omega)$, 
where $\omega\in\Delta^\sigma_\stab$, be a
$\Sigma_R$-stable vector. 
Assume that $\inte(\omega)\subset \inte(\sigma)$. We need
to show that for any $\tau\geq\sigma$, 
any automorphism $\alpha$ of $\Aut(\widetilde{X}_\tau)^0$ preserves
$\val(v,\widetilde{X}_\tau)$.  Find the cone $\omega'\in
{\Delta^\tau}$  such that $\omega'\geq\omega$ and
$\inte(\omega')\subset\inte(\tau)$. By the convexity of
$\stab(\tau)$, $\inte(\omega')$ intersects $\stab(\tau)$
and $\omega'\in {\Delta^\tau}_\stab$.

Denote by
$\pi:\widetilde{X}_{\Delta^\tau}\to {X}_{\Delta^\tau}$ the
standard projection.
By (2) $\alpha$ lifts to an
automorphism $\alpha^f$ of $\widetilde{X}_{\Delta^\tau}$
such that $\alpha^{f*}\pi^*(\cU(\Delta^\tau,
\Delta^\tau_\stab)$ and $\pi^*(\cU(\Delta^\tau,
\Delta^\tau_\stab)$ are compatible along
$f_\Delta^{-1}(O_{\widetilde{\tau}})$. 

Let $p:=O_{\omega'\times\reg(\tau)}=O_{\omega'}\cap f_\Delta^{-1}
(O_{\widetilde{\tau}}) \subset 
\widetilde{X}_{\Delta^\tau}$. Set
$Y:=\widetilde{X}_{\Delta^\tau}$. By Lemma \ref{le: 1},
$\widehat{Y}_p=\widehat{X}_{\omega'\times reg(\tau)}=
\widetilde{X}_{\omega'}$. By assumption
$\val(v,\widehat{Y}_p)$ is invariant with respect to any 
$\Gamma_{\omega'}=\Gamma_\tau$-equivariant automorphisms of
$\widehat{Y}_p=\widetilde{X}_{\omega'}$ preserving
orientation and strata.
The automorphism $\alpha^{f}$ preserves $p$ and induces an
automorphism $\widehat{\alpha}_p^{f}$ of $\widehat{Y}_p$.
Denote by $\widehat{\pi}: \widehat{Y}_p\to X_{\omega'}$ the
morphism induced by $\pi$. 
Since the atlases $\alpha^{f*}_p\pi^*(\cU(\Delta^\tau,
\Delta^\tau_\stab)$ and $\pi^*(\cU(\Delta^\tau,
\Delta^\tau_\stab))$ are compatible the morphisms
$\widehat{\pi}\widehat{\alpha}_p^{f}:\widehat{Y}_p\to
X_{\omega'}$ and $\widehat{\pi}$ determine the same
orientation at $p$. Thus $\widehat{\alpha}_p^{f}$ preserves
orientation at $p$ and $\widehat{\alpha}_{p*}^{f}(\val(v,Y_p))=\val(v,Y_p)$.
Then also
${\alpha}_{*}^{f}(\val(v,\Spec(\cO_{Y,p}))=\val(v, \Spec(\cO_{Y,p})$.
Consequently,  
$\alpha_{*}^{f}(\val(v,\widetilde{X}_{\Delta^\tau})=
\val(v,\widetilde{X}_{\Delta^\tau})$
and  $\alpha_{*}(\val(v,\widetilde{X}_{\tau})=
\val(v,\widetilde{X}_{\Delta^\tau})$. 
\qed

\section{The $\pi$-desingularization lemma of Morelli}

\subsection{Local projections of $\Gamma$-semicomplexes}
\begin{definition} Let $\Sigma$ be a simplicial
$\Gamma$-semicomplex and  $\Delta$ be its canonical subdivision.
 Let $\pi_{\sigma}: \sigma\rightarrow 
\sigma^\Gamma$ be the projection defined by the quotient map $X_{\sigma}\to
X_{\sigma}/\Gamma_\sigma=X_{\sigma^\Gamma}$, for any $\sigma\in
\Sigma$. For any 
$\delta\in\Delta^\sigma$ set
$\Gamma_{\delta}=(\Gamma_{\sigma})_{\delta}=
\{g\in\Gamma_\sigma\mid gx=x, x\in O_\delta\}$. 
By  $\pi_{\delta}: \delta\rightarrow 
\delta^\Gamma$ denote the projection defined by the quotient map $X_{\delta}\to
X_{\delta}/\Gamma_\delta=X_{\delta^\Gamma}$.

We say that $\Sigma$ is {\it strictly
$\pi$-convex} if for any $\sigma \in \Sigma$, 
$\pi_{\sigma}(\sigma)= \sigma^\Gamma$ is strictly
convex (contains no line). We call a semicomplex $\Sigma$
{\it simplicial} if all cones $|\sigma|$, $\sigma\in\Sigma$, are simplicial.

\end{definition}
\begin{lemma} \label{le: projections2} 
\begin{enumerate}
\item Let $\Sigma$ be a   
strictly convex $\Gamma$-semicomplex.   
Then for any $\tau\leq \sigma$, there exists an inclusion 
$\imath^{\sigma\Gamma}_\tau : \tau^\Gamma \hookrightarrow
\sigma^\Gamma $ and the commutative diagram of inclusions:
\[\begin{array}{rccc}
\imath^{\sigma}_\tau:&\tau &\hookrightarrow&\sigma\\  
&\pi_{\tau}\downarrow&&\pi_{\sigma}\downarrow\\
\imath^{\sigma\Gamma}_\tau:&\tau^\Gamma&\hookrightarrow& 
\sigma^\Gamma 
\end{array}\]

\item Let $\Delta$ be a canonical subdivision of $\Sigma$.
Then for any $\delta\in\Delta^\sigma_\stab$ 
we have the commutative diagram of inclusions
\[\begin{array}{rccc}
&\delta &\subset&\sigma\\  
&\pi_{\delta}\downarrow&&\pi_{\sigma}\downarrow\\
&\delta^\Gamma&\hookrightarrow& 
\sigma^\Gamma 
\end{array}\] 
\end{enumerate}
\end{lemma}

\noindent{\bf Proof.}   
(1) Consider the commutative diagram
\[\begin{array}{rccccc}
&
X_{(\tau,N_\sigma)}&\simeq& X_{\tau}\times T&&\\
 &\downarrow &&\downarrow&&\\
&X_{(\tau,N_\sigma)}//\Gamma_\tau&
\buildrel i\over\simeq&X_\tau//\Gamma_\tau
\times T&&\\
 &\downarrow &&\downarrow&&\\
&X_{\pi_{\sigma}(\tau,N_{\sigma})}=
X_{(\tau,N_{\sigma})}/\Gamma_\sigma&\buildrel j\over\simeq
&X_{(\tau,N_\sigma)}//\Gamma_\sigma
\times T/\Gamma_{\sigma}&\simeq& X_{(\tau^\Gamma,N_{\sigma}^\Gamma)}\\

\end{array}\]
\noindent  where $T$ is the relevant torus.
Note that $\Gamma_\tau$ acts trivially on
$\strat(\tau)\subset X_\sigma$, hence in particular on the torus $T$ in
$X_{\tau}\times T$. This gives the isomorphism
$i$. Since  $\Gamma_\sigma/\Gamma_\tau$ acts freely on 
$X_{(\tau,N_\sigma)}//\Gamma_\tau=X_{\tau}//\Gamma_\sigma
\times T$ we get the
isomorphism $j$. Consequently,  
$(\tau^\Gamma,N_{\sigma}^\Gamma)\simeq\pi_{\sigma}(\tau,N_{\sigma})
\subset \sigma^\Gamma$. 

(2) Repeat the reasoning from (1) with appropriate  inclusions.
\qed

\bigskip
\subsection{Dependent and independent cones}
In further considerations we shall assume $\Gamma= K^*$.

\begin{lemma} Let $\Sigma$ be a simplicial strictly
convex $K^*$-semicomplex.
Let $\Dep(\Sigma):=\{\sigma \in \Sigma\mid
\Gamma_\sigma=K^*\}$ and $\Ind(\Sigma):=\{\sigma \in \Sigma\mid
\Gamma_\sigma\neq K^*\}$. Then

\begin{enumerate}
\item For $\sigma \in \Ind(\Sigma)$, $\pi_{\sigma}$ is an immersion into the lattice
$N^\Gamma_\sigma$ of the same dimension.

\item For $\sigma \in \Dep(\Sigma)$, $\pi_{\sigma}$ is  a
submersion onto the lattice $N^\Gamma_\sigma$ such that 
$\dim(N_\sigma)= \dim(N^\Gamma_\sigma)+1$.
Moreover there is a vector $v_\sigma\in N_\sigma$ determined by a
$K^*$-action which spans the kernel of $\pi_{\sigma}$.
\end{enumerate}
\end{lemma}
\noindent {\bf Proof.} If $\sigma \in \Ind(\Sigma)$ then
$\Gamma_\sigma$ is finite and there is an immersion
$M_\sigma^\Gamma\to M_\sigma$, of lattices of the same dimension, where
$M^\Gamma$ is a lattice of $\Gamma$-invariant characters.
The dual morphism $N_\sigma\to N^\Gamma_\sigma$ is also an immersion of
lattices of the same dimension. 
If $\sigma \in \Dep(\Sigma)$ then $\Gamma_\sigma=K^*$ is a
subtorus corresponding to the sublattice $N:=({\bf Q}\cdot
v_\sigma)\cap N_\sigma$ of $N_\sigma$. The projection
$\pi_\sigma$ is defined by the quotient map $N_\sigma\to
N_\sigma/N\simeq N_\sigma^\Gamma$.
\qed

\begin{lemma} Let $\Sigma$ be a simplicial strictly
convex $K^*$-semicomplex. Then for any faces
$\tau\leq\sigma$ in $\Dep(\Sigma)$, 
$\imath^\sigma_\tau(v_\tau)=v_\sigma$.
\end{lemma}
\noindent{\bf Proof.} The induced toric morphism $X_\tau\to X_\sigma$  
is $K^*$-equivariant. \qed

\begin{definition} (see Morelli \cite{Morelli1})
Let $\Delta$ be a canonical subdivision of a
$K^*$-semicomplex $\Sigma$. Let $\sigma\in \Sigma$.
\begin{enumerate}

\item A cone $\delta \in \Delta^\sigma$ is called  {\it independent}
if the restriction of $\pi_\sigma$ to $\delta$ is a lattice
immersion.
Otherwise $\delta \in \Delta^\sigma$  is called  {\it
dependent}. 

\item
A minimal dependent face of $\Delta^\sigma$ is called a
{\it circuit}. 

\item We call an independent face $\tau$ 
{\it up-definite} (respectively
{\it down-definite}) with
respect to  a dependent face $\delta$ if  
 there exists a nonzero functional $F$ on
$\delta\subset N_\sigma^{\bf Q}$ such that 
$F(v_\sigma)> 0$ (respectively $F(v_\sigma)< 0$), and $\tau=\{v \in
\sigma\mid  F(v)=0\}$.  
\end{enumerate}
\end{definition}
\begin{lemma}
  Each dependent cone 
$\delta=\langle v_1,\ldots,v_k\rangle \in \Delta^\sigma$ defines a
unique linear dependence relation
$$r_1\pi_{\sigma}(v_1)+\ldots+r_k\pi_{\sigma}(v_k)=0
\eqno{*}.$$ \noindent This
relation is determined up to proportionality. 
\end{lemma}
\noindent {\bf Proof.} $\pi_\sigma(\delta)$ is $k-1$-
dimensional cone spanned by $k$ rays.\qed

\begin{definition}

\begin{enumerate}
\item The  relation (*) is {\it positively normalized} if
the rays
$\langle v_i \rangle$ 
for which $r_i>0$ form an up-definite face

\item The rays for which $r_i>0$ (in a positively normalized
relation) are called {\it positive}, the rays for which $r_i<0$ are called
{\it negative}, the rays for which $r_i=0$ are called {\it
null} rays. The face $\tau$ of a cone $\delta$ 
is called {\it codefinite} if it contains only positive or only
 negative rays.

\item For any cone $\delta$ we denote by $\delta_{-}$ (respectively
$\delta_{+}$) the fan
consisting of all faces of $\delta$ which are
down-definite (respectively up-definite). By $\delta^-$ (respectively
$\delta^+$) we denote the face of $\delta$ 
spanned by all negative rays and null rays 
(respectively by all positive rays and null rays).
\end{enumerate}
\end{definition}

\begin{lemma} \label{le: weights} The relation (*) is
positively normalized iff
there exists $\alpha>0$, such that
$$r_1v_1+\ldots+r_kv_k=\alpha v_\delta$$
\end{lemma}
\noindent{\bf Proof.} Let $F$ be a nonnegative
functional on $\delta$ which is $0$ exactly on $\langle
v_i\mid r_i>0 \rangle$. Then $\langle
v_i\mid r_i>0 \rangle$ is up-definite iff $F(v_\delta)>0$.
The latter is equivalent to 
$\alpha>0$. \qed

\begin{lemma} 
Let $\tau=\langle v_1,\ldots,\check{v}_i,\ldots,v_k \rangle$ be a  face
of $\delta$ of codimension $1$. Then $\tau$ is up-definite
iff $r_i>0$.
\end{lemma}

\noindent{\bf Proof.} Let $F$ be a nonnegative
functional on $\delta$ which is $0$ exactly on $\tau$.
Then by Lemma \ref{le: weights}, $r_i>0$ iff
$F(v_\delta)>0$. The latter means that
 $\tau$ is up-definite.\qed

\begin{lemma} If $\tau\preceq\delta$ is up-definite
(respectively down-definite) then any face
$\tau'\preceq\tau$ is up-definite (respectively
down-definite) with respect to $\delta$. \qed
\end{lemma}
\noindent{\bf Proof.} Let $F$ be a nonnegative
functional on $\delta$ which is $0$ exactly on $\tau$ and
such that $F(v_\delta)>0$. Let $F'$ be a nonnegative
functional on $\tau$ which is $0$ exactly on $\tau'$. Then 
$nF+F'$, where $n>>0$, is a nonnegative 
functional on $\delta$ which is $0$ exactly on $\tau'$ and 
$(nF+F')(v_\delta)>0$.\qed

\begin{lemma} \label{le: projection} The set $\delta_+$ (respectively $\delta_-$)
is a subfan of $\overline{\delta}$ with maximal faces of
the form $\langle v_1,\ldots,\check{v}_i,\ldots,v_k
\rangle$, where $r_i>0$ (respectively $r_i<0$). In
particular, each boundary face is up-definite or down-definite.
The projection $\pi_\delta$ maps $\delta_+$ (respectively
$\delta_-$) onto the
subdivision $\pi_\delta(\delta_+)$ (respectively
$\pi_\delta(\delta_-)$ of $\pi_\delta(\delta)$. Moreover
the restriction of $\pi$ to any boundary face is a linear isomorphism.
\end{lemma}
\noindent{\bf Proof.} Let $v$ be a vector in
$|\delta^\Gamma|$. Then either the line $\pi^{-1}(v)=\{v\}+\lin(v_\delta)$
intersects the relative interiors of exactly two boundary
faces at one point each and $\pi^{-1}(v)\cap \delta$ is an interval or
$\pi^{-1}(v)$ intersects the relative interior of one
boundary face at a point
and $\pi^{-1}(v)\cap \delta$ is the point. In the first case
one of the faces is up-definite and one is down-definite.
In the second case let $\tau$ be the boundary face
containing $p=\pi^{-1}(v)\cap \delta$. Let $F$ be a
nonnegative 
functional on $\delta^\Gamma$ which is zero exactly on
$\tau$. Let $F'$ be any functional on $N^{\bf Q}_\delta$
which equals $0$  on $\tau$ and
such that $F'(v_\delta)\geq 0$ . Then $(nF\circ\pi)\pm F'$, where
$n>>0$, is nonnegative on $\delta$ and equals $0$ exactly
on $\tau$. Moreover $(nF\circ\pi+F')(v_\delta)>0$ and 
$(nF\circ\pi-F')(v_\delta)<0$. Hence $\tau$ is up-definite
and down-definite.\qed

\begin{lemma} \label{le: weights2} Let
$\delta\in\Delta^\sigma$ be a dependent cone. Set 
$\delta^+:=\langle v_i\mid r_i> 0 \rangle$.
The
ideal $I_{\overline{O}_{\delta^+}}$ of the closure of
${O}_{\delta^+}$ is generated by the set of all characters
with positive $K^*$-weights.
\end{lemma}
\noindent{\bf Proof.} Let $x^m$, $m\in \delta^\vee\cap
M_\delta$ be the character with the positive weight. Then
$(m,v_\delta)> O$. By Lemma \ref{le: weights}, there is
$r_i> 0$ in a positively normalized relation such that $(m,v_i)>0$.
Then $x^m$ is zero on the divisor corresponding to $v_i$
and on ${O}_{\delta^+}$. On the other hand for any $r_i>0$ we can find a
functional $m_i$ such that $(m_i,v_i)>0$,
$(m_i,v_\delta)>O$, $(m_i,v_j)=0$ , where $j\neq i$. Then
$\overline{O}_{\delta^+}$ is the set of
the common zeros of $x^{m_i}$, $r_i> 0$.  \qed

\bigskip
 We associate with a cone $\delta \in
\Delta^\sigma$ and an integral vector $ v \in
\delta^\Gamma$ a vector ${\rm Mid}(v,\delta)\in \delta$ (\cite{Morelli1}.
If $\delta$
is independent face, then ${\rm Mid}(v,\delta)$ is defined to
be the primitive
vector of the ray ${\bf Q}_{\geq
0}\cdot\pi_{\sigma}^{-1}(v) \subset \delta$. 
 If $\delta$ is dependent let $\pi_{\delta_{-}}=\pi_{\sigma|\delta_-}$ and
$\pi_{\delta_{+}}= \pi_{\sigma|\delta_+}$ be the
restrictions of $\pi_\sigma$ to $\delta_{-}$ and
$\delta_{+}$. Then
${\rm
Mid}(v,\delta)$ is the primitive vector of the ray spanned by
$\pi_{\delta_{-}}^{-1}(v)+\pi_{\delta_{+}}^{-1}(v)
\in \delta.$

\begin{remark} By Lemma \ref{le: projections2} the local
projections $\pi_{\sigma}:N_\sigma\to N^\Gamma_\sigma$   
commute with subdivisions and face restrictions  therefore
the above notions are
coherent.
\end{remark}

\bigskip
\subsection{Stable vectors on simplicial
$K^*$-semicomplexes}\label{se: pi-stable} 

In the section $\Sigma$ denotes an oriented 
strictly convex $K^*$-semicomplex.

\begin{lemma} \label{le: tild1} Let $\sigma$ be a semicone
in  $\Sigma$. Then
 $G_\sigma$ acts on 
$$\widetilde{X}_{\sigma^\Gamma}:=\widetilde{X}_\sigma/\Gamma_\sigma$$
\noindent  and the morphism 
$\widetilde{X}_{\sigma}\to
\widetilde{X}_{\sigma^\Gamma}$
is $G_\sigma$-equivariant. Moreover 
$\widetilde{X}_{\sigma^\Gamma}=\widehat{X}_{\sigma^\Gamma\times\reg(\sigma)}$.       
\end{lemma}
\noindent{\bf Proof.}
Follows from the fact that the action of $G^\sigma$
commutes with $\Gamma_\sigma$.\qed

\begin{lemma} \label{le: tild2} Let $\sigma$  be a
dependent cone in $\Sigma$.
Then \begin{enumerate} 
\item   
$\widetilde{X}_{\sigma_{-}}:= {X}_{\sigma_{-}}
\times_{{X}_{\sigma}}\widetilde{X}_{\sigma} \subset 
\widetilde{X}_{\sigma}$ is $G_\sigma$-invariant.

\item
$\widetilde{X}_{{\pi_{\sigma_-}}({\sigma_{-}})}:=
{X}_{{\pi_{\sigma_-}}({\sigma_{-}})}
\times_{{X}_{{\pi(\sigma}})}
\widetilde{X}_{\pi({\sigma})}\to 
\widetilde{X}_{{\pi}({\sigma})}$ is proper $G_\sigma$-equivariant.
\end{enumerate}
\end{lemma}

\noindent{\bf Proof}
(1) By Lemma \ref{le: projection}, $\sigma_-=
\overline{\sigma}\setminus\Star(\sigma^+,\overline{\sigma})$.
By
Lemma \ref{le: weights}, the ideal of
$\overline{O}_{{\sigma}^+}=\widetilde{X}_{\sigma}\setminus
\widetilde{X}_{\sigma_{-}}$
 is 
generated by all
semi-invariant functions with positive weights (see also
\cW, Example 2). 

(2) Follows from Proposition \ref{pr: simple} since 
$\Ver(\pi_{\sigma_-}({\sigma_{-}}))=\Ver(\pi({\sigma}))$ \qed

\begin{lemma} Let $\sigma$ be a cone of the maximal
dimension. Let $F$ denote the functional on 
$\tau:=\sigma^{\vee}$ defined by some integral vector
$v\in N_\sigma$ which is not in $\sigma$. For
any integral $k$
set $\tau_k:=\{v\in 
\tau\mid F(v)=k\}$. Then there are a finite number of vectors 
$w_{k1},\ldots,w_{kl_k}\in \tau_k$ such that $$\tau_k=
\bigcup_{i=1,\ldots,l_k} (w_{ki}+\tau_0).$$
\end{lemma}
\noindent {\bf Proof.} We can replace $\tau$ by a regular
cone $\tau_0$ by considering the epimorphism 
$p:\tau_0\to \tau$ mapping generators to  generators. Then
$F$ defines a functional on $\tau_0$. Let 
$\tau_{0k} :=\{v\in \tau_0\mid F(v)=k\}$. Let
$x_1,\ldots,x_m$ define the standard coordinates on $\tau\simeq 
{\bf Z}^m_{\geq 0}$.
Without loss of generality we can write $F=n_1x_1+\ldots+
n_lx_l- n_{l+1}x_{l+1}-\ldots-n_{r}x_{r}$, where $l<r\leq m$ and
all $n_i>0$. Then $\tau_{0k}=(\tau_{0k}\cap ({\bf Z}^r_{\geq 0}\times \{0\}))+
(\{0\}\times {\bf Z}^{m-r})$, where $
\{0\}\times {\bf Z}^{m-r}\subset \tau_{00}$.

By a {k-\it minimal} vector we shall mean a vector  
$v\in \tau_{0k}\cap ({\bf Z}^r_{\geq 0}\times \{0\})$ such that
there is no $w\in \tau_{00}\cap ({\bf Z}^r_{\geq 0}\times
\{0\})$ for which $v-w\in\tau_{0k}$. 
It suffices to show that the number of $k$-minimal vectors
is finite.
Suppose 
$x_{i_0}(v)\geq (n_1+\ldots +n_l)(n_{l+1}+\ldots+n_r)+|k|$ for some
$i_0 \leq l$ and let $x_{j_0}:=\max\{x_i(v)\mid l+1\leq i\leq
r\}$. Then $x_{j_0}\geq
\frac{n_1x_1+\ldots+n_lx_l}{n_1+\ldots+n_l}=
\frac{n_{l+1}x_{l+1}+\ldots+n_{r}x_{r}+k}{n_{l+1}+\ldots+n_r}\geq
(n_1+\ldots +n_l)(n_{l+1}+\ldots+n_r)/(n_{l+1}+\ldots+n_r)=n_1+\ldots +n_l$ 

Then $w:=(0,\ldots,n_{j_0},\ldots,0,\ldots,n_{i_0}\ldots 0)\in
\tau_{00}$,  where the
$i_0$-th coordinate is $n_{j_0}$ and the $j_0$-th coordinate is $n_{i_0}$.
Thus  $v-w\in\tau_{0k}$ and $v$ is not $k$-minimal.
Consequently, all minimal $k$-vectors satisfy 
$\max\{x_i(v)\mid 1\leq i\leq r\}< 
(n_1+\ldots +n_l)(n_{l+1}+\ldots +n_r)+|k|$. \qed

As a corollary we get
\begin{lemma} \label{le: decomposition2} Let $K[\widetilde{X}_{\sigma}]=\bigoplus 
K[\widetilde{X}_{\sigma}]^k$ be the decomposition
according to weights with respect to the $\Gamma_\sigma=K^*$-action.
Then $K[\widetilde{X}_{\sigma}]^k$ is a finitely
generated $K[\widetilde{X}_{\sigma}]^0$-module.
\qed
\end{lemma}

\begin{lemma} \label{le: tild3} Let $\sigma$  be a
dependent cone in $\Sigma$.
Denote by ${j}_{\pi_{\sigma_-}}:
{X}_{\sigma_{-}}\to {X}_{\sigma_{-}}/\Gamma_\sigma$  and \\ 
$
\widetilde{j}_{\pi_{\sigma_-}}:
\widetilde{X}_{\sigma_{-}}\to \widetilde {X}_{\sigma_{-}}/\Gamma_{\sigma}$
the quotient morphisms.
Then \begin{enumerate}
\item $\widetilde{X}_{\sigma_{-}}/\Gamma_\sigma\simeq 
\widetilde{X}_{\pi_{\sigma_-}(\sigma_{-})}$.

\item $\tilde{j}_{\pi_{\sigma_-}}: \widetilde{X}_{\sigma_{-}}\to 
\widetilde{X}_{\pi_{\sigma_-}({\sigma_{-}})}$ is $G_\sigma$-equivariant.

\item For any $v\in \pi(\sigma_{-})$, ${j}_{\pi}^*
(\val(v,{X}_{\pi_{\sigma_-}({\sigma_{-}})} ))
=\val(\pi_{\sigma_-}^{-1}(v),{X}_{\sigma_{-}})$.

\item If
$\val(v,\widetilde{X}_{\pi_{\sigma_-}({\sigma_{-}})} )$
is $G^\sigma$ invariant then
 ${\tilde{j}_{\pi_{\sigma_-}}}^*
(\val(v,\widetilde{X}_{\pi_{\sigma_-}({\sigma_{-}})} ))
=\val(\pi_{\sigma_-}^{-1}(v),\widetilde{X}_{\sigma_{-}})$
is $G^\sigma$-invariant.

\end{enumerate}
\end{lemma}

\noindent{\bf Proof.}
(1)  For $\tau\in\sigma_{-}$ set
$\widetilde{X}_{\tau}:={X}_{(\tau,N_\sigma)}\times_{{X}_{\sigma}}
\widetilde{X}_{\sigma}$,\,\,\,\\ 
$\widetilde{X}_{\pi_{\sigma_-}({\tau},N_\sigma)}
:={X}_{\pi_{\sigma_-}({\tau})}
\times_{{X}_{\pi({\sigma)}}}
\widetilde{X}_{{\pi({\sigma})}}$.

Then $K[\widetilde{X}_{\tau}]=
K[{X}_{(\tau,N_{\sigma})}]
\otimes_{K[\widetilde{X}_{\sigma}]}K[\widetilde{X}_{\sigma}]$.
 
The elements of this ring are finite sums
$\sum\chi_if_i$, where $f_i\in
K[\widetilde{X}_{\sigma}]$ and $\chi_i\in \tau^\vee$
is  a character.Set

$R_1:=K[\widetilde{X}_{\tau}/\Gamma_\sigma]=
K[\widetilde{X}_{\tau}^{\Gamma_\sigma}]=
(K[{X}_{\tau}]
\otimes_{K[{X}_{\sigma}]}
K[\widetilde{X}_{\sigma}])^{\Gamma_\sigma}$, 

$R_2:=K[\widetilde{X}_{\pi_{\sigma_-}({\tau})}]=
K[{X}_{\pi_{\sigma_-}({\tau})}]\otimes
_{K[{X}_{\pi({\sigma)}}]}
K[\widetilde{X}_{\pi({\sigma})}]\subset R_1$.

$R_2$
consists of elements   
$\sum\chi_if_i$, where $\chi_i\in\tau^\vee$ and 
$f_i\in K[\widetilde{X}_{\sigma}]$ have weight $0$.
$R_1$
consists of elements   
$\sum\chi_if_i$, where $\chi_i\in\tau^\vee$ and 
$f_i\in K[\widetilde{X}_{\sigma}]$ have 
opposite weights. By Lemma \ref{le: decomposition2} there
is a finite
number of characters $\chi_{ij}\in
\sigma^{\vee}$ such that
$f_i=\sum\chi_{ij}f_{ij}$, where $f_{ij}$ have weight zero.
Finally,
$\sum\chi_if_i=\sum_i\chi_i(\sum_j\chi_{ij}f_{ij})=
\sum_i\sum_j(\chi_i\chi_{ij})f_{ij}\in R_2$.

(2) Follows from (1). 

(3) By Lemma \ref{le: projection} the projection $\pi:=\pi_{\sigma_{-}}$ maps
cones of $\sigma_{-}$ isomorphically onto cones of
$\pi(\sigma_-)$. 
  
 By \cite{KKMS}, Ch. I Th. 9, the sheaf of 
ideals $\cI:=\cI_{\val(\pi(v)),d}$ on
$X_{\pi_{\sigma_{-}}(\sigma_-)}$ corresponds to the strictly
convex piecevise linear function $\ord(\cI)$. By definition
and Lemma \ref{le: center}, $\ord(\cI)$ equals $0$ on 
all cones not containing $v$.
Let $\tau\in \sigma_-$ be a cone containing $v$. 
 By the proof of Lemma \ref{le: blow-up valuation}, if $d$ is
sufficiently divisible then for any face
$\delta$ of $\tau$ such that $\delta$ does not contain $v$, 
there is an integral vector $m_{\pi(\delta),\pi(\tau)}\in \pi(\tau)^\vee$ for which$ (m_{\pi(\delta),\pi(\tau)},\pi(v))=d$ , $m_{\pi(\delta),\pi(\tau)}$ is
$0$ on $\pi(\delta)$ and
$\ord(\cI)$ equals $m_{\pi(\delta),\pi(\tau)}$ on $\pi(\delta+\langle
v\rangle)$. 

In particular $\cI$ is generated on each cone $X_{\pi(\tau)}$ by
$m_{\pi(\delta),\pi(\tau)}$, where of $\delta$ is a face of
$\tau$  that does not contain
$v$. Let $m_{\delta}=m_{\pi(\delta)}\circ\pi\in \tau^\vee$ be the
induced functional on $M_\tau$. 
 Thus the ideal $\pi^*(\cI)$ is generated on each cone $X_{\tau}$, where
$\tau\in\sigma_-$ contains $v$, by
$m_{\delta,\tau}$, where $\delta$ doesn not contain $v$.
Then $(m_{\delta,\tau},v)=d$ , $m_{\delta}$ is
$0$ on $\delta$ and
$\ord(\pi^*(\cI))$ equals $m_{\delta,\tau}$ on $\delta+\langle
v\rangle$. This gives by the proof of Lemma \ref{le: blow-up valuation},
$\pi^*(\cI_{\val(\pi(v)),d})=\cI_{\val(v),d}$.

(4) Follows from (2) and (3). \qed

\begin{lemma} \label{le: inclu} 
Let $\tau\leq\sigma$ be  faces of $\Sigma$. Then either the inclusion 
$\imath^{\sigma \Gamma}_\tau$ maps $\tau^\Gamma$ to a face of 
$\sigma^\Gamma$, or $\tau$ is independent, $\sigma$ is dependent, and 
$\tau^\Gamma$ is mapped isomorphically onto the face of $\pi(\sigma_{-})$ (or
 $\pi(\sigma_{+})$).  
\end{lemma}

\noindent{\bf Proof.} If both faces $\tau$ and $\sigma$ are independent then
it follows that $\pi_\sigma$ and $\pi_\tau$ are linear isomorphisms. Then 
$\imath^{\sigma \Gamma}_\tau(\tau^\Gamma)$ is a face of $\sigma^\Gamma$.
If $\tau$ and $\sigma$ are both dependent then the kernel
of the projection 
$\pi_\sigma$ is contained in $\lin(\tau)$. Then $\pi(\tau)$
is a face of $\pi(\sigma)$. 
 If $\tau$ is an independent face of a dependent cone $\sigma$ then $\tau$ is a 
face of $\sigma_-$ or $\sigma_+$. Both fans consist of independent faces of 
$\sigma$ and project onto the subdivisions $\pi(\sigma_{-})$ and
 $\pi(\sigma_{+})$ of $\pi(\sigma)$.
\qed

\begin{definition}
 A  vector $v\in {\rm
int}({\sigma^\Gamma})$, where $\sigma\in \Sigma$,  
is $\Sigma$-{\it stable} if for any $\tau\geq
\sigma$ the corresponding
valuation ${\rm val}(\imath^{\tau \Gamma}_{\sigma}(v))$ on 
$\widetilde{X}_{\tau^\Gamma}$
is $G^\tau$-invariant. A vector $v\in \sigma^\Gamma$ is $\Sigma$-{\it
stable} if there is a $\Sigma$-stable vector $v_0\in {\rm
int}({\sigma^\Gamma_{0}})$, where $\sigma_0\leq \sigma$, for which
$v=\imath^{\sigma\Gamma}_{\sigma_0}(v_0)$.

\end{definition}

\begin{proposition} \label{pr: G-stab} Let $\sigma$ be a
semicone in $\Sigma$.
A vector $v\in\inte(\sigma)$ is 
$\Sigma$-stable if $\pi_\sigma(v)\in \sigma^\Gamma$
 is $\Sigma$-stable. 

\end{proposition}

\noindent{\bf Proof.}
 A  stable vector $v\in 
{\rm int}(\sigma)$ determines a $G^\tau$-invariant
valuation on any $\widetilde{X}_{\tau}$ for $\tau\geq
\sigma$. Consequently, it
determines a $G^\tau$-invariant
valuation on any $\widetilde{X}_{\tau}/\Gamma_{\tau}$ and
finally $\pi_{\sigma}(v)$ is $\Sigma$-stable.

Now
let $v\in {\rm int}(\sigma)$ be an integral vector 
 such that $\pi_{\sigma}(v)\in 
\inte(\sigma^\Gamma)$ is stable. Then  
$\val{(\pi_{\sigma}(v))}$ is $G^\tau$-invariant on 
$\widetilde{X}_{\tau}/\Gamma_{\tau}$. 
We have to prove that $v$ is stable, or equivalently, that for any $\tau\geq \sigma$,
${\rm val}(v)$ is $G^\tau$-invariant on
$\widetilde{X}_{\tau}$. Consider two cases:  

(1)  $\tau\in \Ind(\Sigma)$. 
The morphism $\tilde{j}_{\pi_{\tau}}: \widetilde{X}_{\tau}\to 
\widetilde{X}_{\tau}/\Gamma_\tau=\widetilde{X}_{\tau^\Gamma}$
is $G^\tau$-equivariant.
The valuation
  $\tilde{j}_{\pi_{\tau}*}(\val{(v)})=\val{(\pi_{\tau}(v))}$ is 
$G^\tau$-invariant on 
$\widetilde{X}_{\tau}/\Gamma_\tau$. Thus for
any $g\in G^\tau$, $g_*\tilde{j}_{\pi_\tau *}(\val{(v)})=
\tilde{j}_{{\pi_{\tau}}*}g_*
(\val{(v)})=\tilde{j}_{{\pi_{\sigma}}*}(\val{(v)})$.
This means that $\val(v)$ and $g_*(\val(v))$ define the same
functional on the lattice of characters $M^\Gamma_\tau$ and
consequently $M_\tau$ (since $M^{\bf
Q}_\tau=(M^\Gamma_\tau)^{\bf Q}$ ). This shows by Lemma \ref{le: monomial} that
$g_*(\val(v))\geq \val(v)$ and finally by Lemma \ref{le:
monomial2} we conclude  $g_*(\val(v)) = \val(v)$.

(2) $\tau\in \Dep(\Sigma)$. By Lemma \ref{le: tild2}(2), 
$\widetilde{ i_{/K^*}}:\widetilde X_{\pi(\tau_{-})} \to
\widetilde{X}_{\pi(\tau)}$ is  a 
$G^\tau$-equivariant proper birational morphism.
Then the valuation ${\rm val}(v,\widetilde{X}_{\pi(\tau_-)})=
(\widetilde{i_{/K^*}}_*^{-1}({\rm
val}(v,\widetilde{X}_{\pi(\tau)}))$ is $G^\tau$-invariant.
By Lemma \ref{le: tild3}(4), 
$\widetilde{j}_{\pi_{\tau_{-}}}^*({\rm
val}(\pi(v),\widetilde{X}_{\pi(\tau_-)}))=  {\rm
val}(\pi_{\tau_{-}}^{-1}\pi(v),\widetilde{X}_{\tau_-})$ is 
$G^\tau$-equivariant.
We get  $\pi_{\tau_{-}}^{-1}(\pi(v))\in {\rm Inv}(\tau)$.
Analogously $\pi_{\tau_{+}}^{-1}(\pi(v))\in
{\rm Inv}(\tau)$.  By Lemma \ref{le: convex} this  gives  
$v\in \langle\pi_{\tau_{-}}^{-1}(\pi(v)),
\pi_{\tau_{+}}^{-1}(\pi(v))\rangle\subset
{\rm Inv}(\tau)$.

 \qed

\begin{proposition} \label{pr: G-stab 2} Let $\sigma \in
\Ind(\Sigma)$ be an 
independent face of $\Sigma$. Then 

\begin{enumerate}

\item All vectors in $ \overline{\rm
par}(\sigma^\Gamma)\cap \inte(\sigma^\Gamma)$ 
are $\Sigma$-stable.

\item All vectors in ${\rm
par}(\sigma^\Gamma)$ are $\Sigma$-stable. 

\end{enumerate}
\end{proposition}

\noindent {\bf Proof} 
(1) Let $v\in
\overline{\rm par}(\sigma^\Gamma)\cap {\rm
int}(\sigma^\Gamma)$. We have to prove that it defines a
$G^\tau$-invariant valuation on any
$\widetilde{X}_{\pi_\tau(\tau)}$ for $\tau\geq \sigma$.

By Lemma \ref{le: inclu} either  $\sigma^\Gamma$ is a face of $\tau^\Gamma$ or 
$\tau\in \Dep(\Sigma)$ and $\sigma^\Gamma$ is a face of $\pi_\tau(\tau_-)$ (or 
$\pi_\tau(\tau_+)$).
In the first case $\sigma^\Gamma$ is a is a
$G^\tau$-invariant face of
$\tau^\Gamma$. By Lemma
\ref{le: parr}(2), $\val(v)$ defines a $G^\tau$-invariant
valuation on 
$\widetilde{X}_{\tau^\Gamma}$. In the second case $\sigma$ is
a $G^\tau$-invariant face of $\tau_-$ (or $\tau_+$). Consequently, 
the closure of the
orbit $O_{\sigma^\Gamma}$ in $\widetilde{X}_{\pi_\tau(\tau_-)}$
is $G_\tau$-invariant.
Let $\widetilde{i_{/K^*}}:
\widetilde{X}_{\pi_\tau(\tau_{-})} \longrightarrow
\widetilde{X}_{\pi_\tau(\tau)}$ be 
the $G_\tau$-equivariant morphism defined by the
subdivision $\pi_\tau(\tau_{-})$ of
$\pi_\tau(\tau)$. Then by Lemma
\ref{le: parr}(2), $\overline{\rm
par}(\sigma^\Gamma)\cap \inte(\sigma^\Gamma)$ is contained 
$\Inv(\pi_\tau(\tau))$.

(2) Let $v\in
{\rm par}(\pi_\sigma(\sigma))$. Then $v\in \overline{\rm par}
(\pi_\sigma(\tau))\cap {\rm
int}(\tau)$ for some $\Gamma$-indecomposable
face $\tau$ of $\sigma$. We apply (1).

 \begin{lemma} \label{le: centers2} 
 Let
$\sigma=\langle v_1,\ldots,v_k\rangle$ be a
circuit in $\Sigma$ with the unique relation
$\sum\alpha_iw_i=0,$
where all $\alpha_i\neq 0$ and $w_i=\prim(\pi_\sigma(v_i))$. 
 Then
${\rm Ctr}_-(\sigma):=\sum_{\alpha_i<0}w_i$ and 
${\rm Ctr}_+(\sigma):=\sum_{\alpha_i>0}w_i$
are $\Sigma$-stable.
\end{lemma}

\noindent{\bf Proof.} We have to show that
for any $\tau\geq \sigma$ the valuation $\val({\rm
Ctr}_-(\sigma),\widetilde{X}_{\tau^\Gamma})$ is
$G_\tau$-invariant on $\widetilde{X}_{\tau^\Gamma}$.
By Lemmas \ref{le: weights2} and \ref{le: tild2}(1),
the face $\sigma^-=\langle
v_i\rangle_{\alpha_i>0}\in \sigma_{-}\subset \tau_{-}$ 
corresponds to the $G^\tau$-invariant closed subscheme of
$\widetilde{X}_{\tau_-}\subset \widetilde{X}_{\tau}$.  
Therefore the orbit closure 
 scheme $\overline{O_{\pi_\sigma(\sigma^-)}}\subset 
\widetilde{X}_{\pi_\tau(\tau_-)}$ is $G^\tau$-invariant. 
By Lemma \ref{le: parr}(2), applied to the subdivision
${\pi_\tau(\tau_-)}$ of ${\pi_\tau(\tau)}$,
 $v:={\rm Ctr}_-(\sigma)$ defines a $G^\tau$-invariant valuation on $
\widetilde{X}_{\pi_\tau(\tau)}$. 
\qed

\begin{lemma} \label{le: centers} Let $\Delta$ be a
canonical subdivision of  $\Sigma$.
\begin{enumerate}
 \item If 
$v\in {\rm par}(\pi(\delta))$, where $\delta \in
\Delta^\sigma$, is an independent cone,  then
${\rm Mid}(v,\delta)$ is $\Sigma$-stable.

\item Let 
$\delta\in\Delta^\sigma$  be a
circuit.   
 Then
${\rm Mid}({\rm Ctr}_{+}(\delta),\delta)$,
${\rm Mid}({\rm Ctr}_{-}(\delta),\delta)$
are $\Sigma$-stable.
\end{enumerate}
\end{lemma}

\noindent{\bf Proof}
(1)  There is a $\Gamma$-indecomposable cone $\delta'\preceq\delta$
such that $v\in {\rm par}(\pi(\delta'))$. Then $\delta'\in \Delta^\sigma_\stab$.
We apply Lemma \ref{le: centers2} and Proposition \ref{pr: G-stab} to $\delta'\in
\Sigma'=\Delta_\stab$ to deduce  that 
${\rm Mid}(v,\delta)={\rm Mid}(v,\delta')$ is
$\Sigma'$-stable. By Proposition \ref{pr: canonical orientation}(3), it is
$\Sigma$-stable. 

(2) The circuit $\delta\in \Delta^\sigma$ is
$\Gamma$-indecomposable. Thus  $\delta\in \Delta_\stab$.
We apply Lemma \ref{le: centers2} and Proposition \ref{pr: G-stab} to $\delta\in
\Sigma'=\Delta_\stab$ to infer that ${\rm Mid}({\rm Ctr}_{+}(\delta),\delta)$,
${\rm Mid}({\rm Ctr}_{-}(\delta),\delta)$
are $\Sigma'$-stable. Consequently, by Proposition \ref{pr:
canonical orientation}(3) they are
$\Sigma$-stable. \qed

\subsection{The $\pi$-desingularization Lemma of Morelli}

\begin{definition} An independent cone $\delta \in
\Delta^\sigma$ is $\pi$-{\it nonsingular} if
 $\pi_{\sigma}(\delta)$ is regular. A subdivision $\Delta$ is   
$\pi$-{\it nonsingular} if all independent cones in $\Delta^\sigma$, 
where $\sigma\in\Sigma$, are $\pi$-{\it nonsingular}.
\end{definition}

\begin{definition}(Morelli \cite{Morelli1}, \cite{Morelli2},
\cite{Abramovich-Matsuki-Rashid})\label{de: Morelli} 
Let $\pi:N^{{\bf Q}+}:=N^{\bf Q}\oplus {\bf Q}\to N^{\bf Q}$ be the
natural projection and $v=(\{0\}\times 1)\in N^{\bf Q}$. 
A fan $\Sigma$ in $N^{{\bf Q}+}$ is
called a {\it polyhedral cobordism} or simply a cobordism if
\begin{enumerate}
\item For any cone $\sigma\in\Sigma$ the image $\pi(\sigma)$ is
strictly convex (contains no line).
\item The sets of cones 
 $$\Sigma_+: =\{\sigma \in\Sigma\mid \mbox{there  exists
} \quad \epsilon>0, \mbox{such that } \sigma+\epsilon\cdot v\not\in
|\Sigma| \}$$ \noindent and
$$\Sigma_-=\{\sigma \in\Sigma\mid \quad\mbox{there exists} \quad 
\epsilon>0, \mbox{such that}\sigma+\epsilon\cdot
v\not\in |\Sigma|\}$$ are
subfans of $\Sigma$ and
$\pi(\Sigma_-):=\{\pi(\tau)\mid\tau\in \Sigma_-\}$ and 
$\pi(\Sigma_+):=\{\pi(\tau)\mid\tau\in \Sigma_+\}$
are fans in $N^{\bf Q}$.

\end{enumerate}
\end{definition}

\begin{lemma}(Morelli \cite{Morelli1}, \cite{Morelli2},
\cite{Abramovich-Matsuki-Rashid}) \label{le: Morelli} 
 Let $\Sigma$ be a simplicial cobordism in $N^{{\bf Q}+}$. Then there exists a
simplicial cobordism $\Delta$ obtained from $\Sigma$ by a sequence of star
subdivisions such that $\Delta$ is $\pi$-nonsingular.
Moreover, the sequence can be taken so that any independent
and already $\pi$-nonsingular face of $\Sigma$ remains
unaffected during the process. Moreover
all  centers of star subdivisions are of
the form (1), (2) from Lemma \ref{le: centers}.$\qed$
\end{lemma}

A simple corollary of the above lemma is
\noindent \begin{lemma} \label{le: pi-lemma}
Let $\Sigma$
be a simplicial strictly convex $K^*$-semicomplex. Then there exists a
canonical subdivision  $\Delta$ of 
$\Sigma$ which is a sequence of star
subdivisions such that $\Delta$ is $\pi$-nonsingular.
Moreover, the sequence can be taken so that any independent
and already $\pi$-nonsingular face of $\Sigma$ remains
unaffected during the process. Moreover all centers are of
the form (1), (2) from Lemma \ref{le: centers} and
therefore are $\Sigma$-stable.
\end{lemma}

\begin{lemma} \label{le: normal} (\cite{Wlodarczyk1},Lemma
10, \cite{Morelli1}) Let $w_1,\ldots,w_{k+1}$ be integral vectors in
${\bf Z}^{k}\subset {\bf Q}^{k}$ 
which are not contained in a proper vector subspace of ${\bf Q}^{k}$.
Then $$\sum_{i=1}^n(-1)^i\det(w_1,\ldots,\check{w_i},\ldots,w_{k+1})\cdot w_i=0 $$ is the
unique (up to proportionality) linear relation between $w_1,\ldots,w_{k+1}$. 
\end{lemma}
\noindent{\bf Proof.}
 Let $v:=
\sum_{i=1}^n (-1)^i\det(w_1,\ldots,\check{w_i},\ldots,w_{k+1})\cdot w_i$.
Then  for any $i<j$, 
\[\begin{array}{rc}
&\det(w_1,\ldots,\check{w_i},\ldots,\check{w_j},\ldots,w_{k+1},v)=
\det(w_1,\ldots,\check{w_i},\ldots,w_{k+1})\cdot 
\det(w_1,\ldots,\check{w_i},\ldots,\check{w_j},\ldots,w_{k+1},w_i)+
\\
&\det(w_1,\ldots,\check{w_j},\ldots,w_{k+1})\cdot 
\det(w_1,\ldots,\check{w_i},\ldots,\check{w_j},\ldots,w_{k+1},w_j)=\\
&(-1)^i(-1)^{k-i}\det(w_1,\ldots,\check{w_j},\ldots,w_{k+1})
\cdot\det(w_1,\ldots,\check{w_i},\ldots,w_{k+1})+\\
&(-1)^j(-1)^{k-j+1}\det(w_1,\ldots,\check{w_i},\ldots,w_{k+1})
\cdot\det(w_1,\ldots,\check{w_j},\ldots,w_{k+1})=0
\end{array}\]
Therefore $v\in \bigcap_{i,j}
\lin\{w_1,\ldots,\check{w_i},\ldots,\check{w_j},\ldots,w_{k+1}\}=\{0\}$.\qed

\noindent \begin{lemma} \label{le: centers3} If $\Sigma$ is a simplicial
cobordism which is $\pi$-nonsingular and $\sigma\in \Sigma$
is a circuit then \\${\rm Mid}({\rm Ctr}_{+}(\sigma),\sigma)=
{\rm Mid}({\rm Ctr}_{-}(\sigma),\sigma)$ and 
${\rm Mid}({\rm Ctr}_{+}(\sigma),\sigma)\cdot\Sigma$ is $\pi$-nonsingular.
\end{lemma}

\noindent{\bf Proof.}
 Let  $\sigma=\langle v_1,\ldots,v_k\rangle$ be a
curcuit and $w_i:=\prim(\pi(\langle v_i \rangle))$..
\\Then by Lemma \ref{le: normal}, in the relation the unique relation 
$\sum\alpha_iw_i=0$ all
$\alpha_i=\det[w_1,\ldots,\check{w}_i,\ldots,
w_k]=\pm 1$
and projections of ''new independent faces''are regular
since they are obtained by
regular star subdivisions applied to projections of "old"
independent faces. \qed

\noindent{\bf Proof of Lemma \ref{le: pi-lemma}.} 
First note that 
the local projections $\pi_{\sigma}: N_\sigma\to
N^\Gamma_\sigma$ glue together 
and commute with subdivisions by Lemma \ref{le: projections2} and therefore the notion of
$\pi$-nonsingularity is coherent.
Let $\sigma_1,\ldots,\sigma_k$ denote all the faces of $\Sigma$.
 We shall construct by induction on $i$ a canonical
$\Gamma$-subdivision $\Delta_i$ of $\Sigma$ such that the
fan ${\Delta_i^{\sigma_j}}$ is $\pi$-nonsingular for all
$j\leq i$.  Suppose $\Delta_i$ is
already constructed. By Lemma \ref{le: Morelli}
there exist centers $v_1,\ldots,v_k$ of the appropriate type for the cobordism
$\Delta^{\sigma_{i+1}}_i$ in $N_{\sigma_{i+1}}$ such that
$\langle v_k \rangle\cdot\ldots\cdot\langle v_1
\rangle\cdot\Delta^{\sigma_{i+1}}_i$ is
$\pi$-nonsingular. It follows from Lemma
 \ref{le: centers}  that all centers
in the $\pi$-desingularization process are $\Sigma$-stable, hence by
 Proposition \ref{pr: blow-ups}, $\Delta_{i+1}:
=\langle v_k \rangle\cdot\ldots\cdot\langle v_1
\rangle\cdot\Delta_{i}$ is canonical. By construction 
$\Delta^{\sigma_{i+1}}_{i+1}$ is $\pi$-nonsingular. Note that 
$\Delta^{\sigma_{j}}_{i+1}$ for $j<i+1$ is obtained by star
subdivisions of the $\pi$-nonsingular 
cobordisms $\Delta^{\sigma_j}_{i}$  at
centers of type $\Mid(\Ctr_\pm(\delta),\delta)$ only and hence by Lemma \ref{le:
centers3} it remains $\pi$-nonsingular.  \qed

\bigskip

For completeness we give a simple proof of Lemma
\ref{le: Morelli}. The original
proof of Morelli used different centers of star subdivisions but finally it 
was modified by the author
(\cite{Morelli1},\cite{Morelli2}). 
The complete published proof of the lemma is given in
\cite{Abramovich-Matsuki-Rashid} (\cite{Morelli2}). 
 The present proof is based upon the algorithm developed
in \cite{Wlodarczyk1}. All algorithms in their  final versions use the same centers of subdivisions and are
closely related. 

\noindent{\bf Proof of Lemma \ref{le: Morelli}.} 

\begin{definition} Let  $\delta=\langle v_1,\ldots,v_k\rangle$ be a
dependent cone and $w_i:=\prim(\pi(\langle v_i \rangle))$. 
Then we shall call the relation  $\sum_{i=1}^k r_iw_i=0$ {\it
normal} if it is positively normalized and 
$|r_i|=|\det(w_1,\ldots,\check{w_i},\ldots,w_{k})|$ for $i=1,\ldots,k$. 
An independent cone $\sigma$ is called an {\it $n$-cone}
if $|\det(\sigma)|=n$. A dependent cone $\sigma$ 
is called an {\it $n$-cone} if one of its independent faces 
is an $n$-cone and the others are $m$-cones, where $m\leq n$. 

\end{definition}

\begin{lemma} \label{le: normal2} Let $\Sigma$ be a
simplicial coborism and $\delta=\langle v_1,\ldots,v_k\rangle$ be a
maximal dependent cone in $\Sigma$. Set
$w_i:=\prim(\pi(\langle v_i \rangle))$. Let  
$$\sum_i r_iw_i=0 \eqno(0),$$ \noindent where
$|r_i|=|\det(w_1,\ldots,\check{w_i},\ldots,w_{k})|$, 
be its normal relation (up to sign).

\noindent (1)  Let $v=\Mid(\Ctr_+(\delta),\delta)\in\inte
\langle v_i\mid r_i\neq 0\rangle$. Let $m_w\geq 1$ be the integer
such that the vector
$w=\frac{1}{m_w}(w_1+\ldots+w_k)$ is primitive. Then the maximal
dependent cones in $\langle v \rangle\cdot\delta$ are of the
form $\delta_{i}=\langle
v_1,\ldots,\check{v}_{i},\ldots,v_k,v\rangle$.

\noindent 1a.  Let $r_{i_0}>0$. Then for the maximal dependent cone
$\delta_{i_0}=\langle
v_1,\ldots,\check{v_{i_0}},\ldots,v_k,v\rangle$ in 
$\langle v \rangle\cdot\delta$, the  normal relation
is given (up to sign) by $$\sum_{r_i> 0, i\neq i_0}\frac{r_i-r_{i_0}}{m_w}w_i+
\sum_{r_i< 0} \frac{r_i}{m_w}w_i+r_{i_0}w=0. \eqno(1a)$$

\noindent 1b. Let $r_{i_0}<0$. Then for the maximal dependent cone
$\delta_{i_0}=\langle
v_1,\ldots,\check{v}_{i_0},\ldots,v_k,v\rangle$ in $\langle
v \rangle\cdot\delta$, the  normal relation
is given (up to sign) by $$\sum_{r_i< 0}-\frac{r_{i_0}}{m_w}w_i+
r_{i_0}w=0. \eqno(1b)$$

(2) Let $\sigma=\langle v_i\mid i\in I\rangle$ be a
codefinite face of
$\delta$. For simplicity assume that $r_i\geq 0$ for $i\in I$.
  Let
$w=\sum_{i\in I}\alpha_iw_i\in
\para(\pi(\sigma))\cap\inte(\pi(\sigma))$ and $v=\Mid(w,\sigma)\in\inte
(\sigma)$. Then the maximal dependent cones in 
 $\langle v \rangle\cdot\delta$
are of the form $\delta_{i_0}=
\langle v_1,\ldots,\check{v}_{i_0},\ldots,v_k,v\rangle$,
where $i_0\in I$.

2a. Let $i_0\in I$ and $r_{i_0}>0$. Then for the maximal
dependent cone $\delta_{i_0}=
\langle v_1,\ldots,\check{v}_{i_0},\ldots,v_k,v\rangle$ in
$\langle v  \rangle\cdot\delta$, 
 the  normal relation
is given (up to sign) by 
$$\sum_{i\in I\setminus \{i_0\},r_i>0}
(\alpha_{i_0}r_{i}-
\alpha_{i}r_{i_0})w_{i}+\sum_{i\not\in I, r_i>0}
\alpha_{i_0}r_{i}w_i+\sum_{
i\in I,r_i=0} -\alpha_{i}r_{i_0}w_{i}+ 
\sum_{r_i<0}\alpha_{i_0} r_iw_{i}+r_{i_0}w=0. \eqno(2a)$$ 

2b. Let  $i_0\in I$ and $r_{i_0}=0$. For  the maximal dependent cone
$\delta_{i_0}=\langle
v_1,\ldots,\check{v}_{i_0},\ldots,v_k,v\rangle$,
  the  normal relation
is given (up to sign) by
$$\sum_{i\in I\setminus\{i_0\}} \,\,\alpha_{i_0}r_iw_i+0w_{i_0}=0. \eqno(2b)$$

\end{lemma}
{\bf Proof.} 
It is straightforward to see that the above equalities hold. We only
need to show that the relations considered are normal. For that
it suffices to show that one of the coefficients is equal
(up to sign) to the corresponding determinant 

 Comparing the coefficients of $w$ in the above
relations with the normal
relations from Lemma \ref{le: normal} we get

1a. $\det(w_1,\ldots,\check{w}_{i_0},\ldots,w_k)=
r_{i_0}$.

1b. The coefficient of $w$ is equal to
$\det(w_1,\ldots,\check{w}_{i_0},\ldots,w_k))=r_{i_0}$, 

2a. The coefficient of $w$ is equal to 
$\det(w_1,\ldots,\check{w}_{i_0},\ldots,w_k)=r_{i_0}$. 

2b. The coefficient of $w_i$, where $r_i>0$, is equal to \\ 
$\det(w_1,\ldots,
\check{w}_{i_0},\ldots,\check{w_{i}},\cdots,w_k,w)=\alpha_{i_0}
\det(w_1,\ldots,
\check{w}_{i_0},\ldots,\check{w_{i}},\cdots,w_k,w_{i_0})=\\
(-1)^{k-i_0}\alpha_{i_0}
\det(w_1,\ldots,\check{w}_{i},\cdots,w_k,w)=
(-1)^{k-i_0} \alpha_{i_0}r_i$.

\qed

 Let $\delta=\langle v_1,\ldots,v_k\rangle$ be a dependent cone. Let 
$\sum\,\, r_iw_i=0$ be its normal relation.  

We shall consider the following kinds of dependent cones.

Type $I$: the cones which are {\it pointing down}: $\Card\{r_i\mid
r_i<0\}=1$. 

and the cones {\it pointing up}: $\Card\{r_i\mid r_i>0\}=1$.

Type $I(n,n)$: the cones of type $I$ for which 
$\max\{|r_i|\mid r_i>0\}=\max\{|r_i|\mid r_i<0\}=n$.

Type $I(n,*)$: the cones pointing down for which
$\max\{|r_i|\mid r_i>0\}<n$ and $\max\{|r_i|\mid r_i<0\}=n$
and the cones pointing up for which
$\max\{|r_i|\mid r_i<0\}<n$ and $\max\{|r_i|\mid r_i>0\}=n$. 

Type $I(*,n)$: the cones pointing down for which
$\max\{|r_i|\mid r_i<0\}<n$ and $\max\{|r_i|\mid r_i>0\}=n$,
and the cones pointing up for which
$\max\{|r_i|\mid r_i>0\}<n$ and $\max\{|r_i|\mid r_i<0\}=n$.

\bigskip
Type $II$: the cones which are neither pointing down nor
pointing up: $ \Card\{r_i\mid
r_i>0\}>1$ and $ \Card\{r_i\mid
r_i<0\}>1$.

Type $II(n,n)$: the cones of type $II$ for which 
$\max\{|r_i|\mid r_i>0\}=\max\{|r_i|\mid r_i<0\}=n$.

Type $II(n,*)$: the cones of type $II$ for which
$\max\{|r_i|\mid r_i\neq 0\}=n$ but 
$\max\{|r_i|\mid r_i>0\}<n$ or $\max\{|r_i|\mid r_i<0\}<n$.

\bigskip
Type $III$: the cones which are pointing down and pointing up
at the same time:\\ $\Card\{r_i\mid
r_i>0\}=\Card\{r_i\mid
r_i<0\}=1$.

Type $III(n,n)$: the cones of type $III$ for which 
$\max\{|r_i|\mid r_i>0\}=\max\{|r_i|\mid r_i<0\}=n$.

\bigskip
Types $I(*,*)$, $II(*,*)$, $III(*,*)$ the cones of type $I$, $II$,
$III$, respectively, for which $\max\{|r_i|\mid r_i\neq 0\}<n$. 

\bigskip
 Denote
by $n$  the maximum of the determinants of the projections
of independent cones in the simplicial cobordism $\Sigma$.
The $\pi$-desingularization algorithm consists of
eliminating all dependent $n$-cones according to the above classification.

All normal relations below are considered up to sign.

\bigskip
{\bf Step 1}. {\it Eliminating all cones of type $II(n,n)$} 

Let $\delta$ be a maximal dependent cone of type $II(n,n)$
in $\Sigma$ with normal
relation (0) (see Lemma \ref{le: normal2}). By Lemma \ref{le:
normal2} after the star
subdivision at $\Mid(\Ctr_+,\delta)$ all new dependent 
$\delta'=\delta_{i_0}$ cones have normal relations of type
(1a) if $r_{i_0}> 0$ or (1b) if  $r_{i_0}> 0$.
If $r_{i_0}> 0$ in the relation (1a) the absolute values of
the coefficients  satisfy $\frac{1}{m_w}|r_i-r_{i_0}|<n$, the absolute values of
the other coefficients  can be  $n$ only for some $r_i<0$
and $r_{i_0}$.  If $r_{i_0}<n$ then  there is no
positive coefficient in (1a) which is equal to $n$. 
 Hence $\delta'$ can be of type
$I(n,*)$, $I(*,n)$, $I(*,*)$, $II(n,*)$ or $II(*,n)$. If $r_{i_0}=n$ then there is only
one positive coefficient (of $w$) in (1a). The other
coefficients are
$\frac{r_i-r_{i_0}}{m_w}=\frac{r_i-n}{m_w} \leq 0$ and $r_i<0$. Thus
$\delta'$ is of the form $I(n,*)$ or $I(n,n)$.

If $\delta'=\delta_{i_0}$ for $r_{i_0}<0$ 
then the normal relation is of the form (1b). If
$|r_{i_0}|=n$ then $\delta'$
is of the form $I(n,n)$. If $|r_{i_0}|<n$  then $\delta'$
is of the form  $I(*,*)$.

By taking the star subdivision we eliminate dependent cones
in  $\Star(\delta_+,\Sigma)$, with the normal relation
proportional to (0). All the dependent $n$-cones  
we create are either of type $I$ or $II(n,*)$.

\bigskip
{\bf Step 2}. {\it Eliminating all cones of type $I(n,n)$}.
Let $\delta$ be a maximal dependent cone of type $I(n,n)$.
Without loss of generality we may assume that $\delta$ is
pointing down (one negative ray). Then $\delta^+$ is of
codimension $1$ and 
$\det(\pi(\delta^+))=n$.   

Find
$v_+\in \para(\pi(\delta^+))$. Then $v_+\in \inte(\pi(\sigma_v))$, where 
$\sigma_v$ is a face of $\delta^+$. Thus $\sigma_v$ is
independent.

{\bf Step 2a}. {\it Making $\sigma_v$ codefinite with respect to all
dependent cones in $\Star(\sigma_v,\Sigma)$}.

By applying certain star subdivisions we eliminate all
cones for which $\sigma_v$ is not codefinite.
Let $\delta'\in\Star(\sigma_v,\Sigma)$ be any dependent cone 
for which $\sigma_v$ is not codefinite. Note that $\delta'$
is not of type $III$ since in this case the circuit consists
of one positive and one negative ray and the noncodefinite face $\sigma_v$
containing positive and negative rays would contain a
circuit, which is impossible for an independent cone.
Therefore $\delta'$ is either of type $I$ or $II$. 

Suppose first that $\delta'=\langle
v_1,\ldots,v_{k}\rangle$ is of type $I(n,n)$  
or
$I(*,n)$. Without loss of generality we may assume that
there is only one coefficient $r_{i_1}>0$ (equal to $n$) 
in the normal relation
(0) from Lemma \ref{le: normal2}(1). 

Since $\sigma_v$ is not codefinite
with respect to $\delta'$, we have $v_{i_1}\in \sigma_v$. Note
that $\Ctr_+(\delta')=\prim(\pi(\langle v_{i_1}\rangle))$ and  
for $v:=\Mid(\Ctr_+(\delta'),\delta')$ we have 
$\prim(\pi(\langle v_{i_1}\rangle))=\prim(\pi(\langle
v_{i_1}\rangle))$. Let $w_i:=\prim(\pi(\langle
v_{i_1}\rangle))$ and $w:=\prim(\pi(\langle v\rangle))$.
Then $w=w_{i_1}$.
After the star subdivision of $\delta'$ at 
$\langle v\rangle$  by Lemma \ref{le: normal2}(1a),
we create a dependent cone $\delta'_{i_1}$ in $\langle
v\rangle\cdot\delta'$ 
 with the 
normal relation obtained from the relation (0) by replacing $w$ with
$w_{i_1}$ . Since $\sigma_v$ contains only negative rays of
$\delta'_{i_1}$ it is codefinite.  We also create  
dependent cones $\delta'_{i_0}$, where
$i_0\neq i_1$, of type $III$
($III(n,n)$ or $III(*,*)$) with the normal relation
$r_{i_0}w-r_{i_0}w_{i_1}=0$ (case 1b). The face $\sigma_v$ is
codefinite with respect to $\delta'_{i_0}$.
Thus after blow-up we create one cone of the type $I(n,n)$
and cones of type $III(n,n)$ or $III(*,*)$.
Since  at the same
time we eliminate the cone $\delta'$, the number of maximal cones of type $I(n,n)$ 
remains the same. After the star subdivision  
$\sigma_v$ is 
codefinite with respect to all new dependent cones.

Suppose $\delta'$ is of type $I(n,*)$, $I(*,*)$ $II(n,*)$
or $II(*,*)$. Without loss of generality
we may 
assume that in
the normal relation (0) for $\delta'$ (see Lemma \ref{le: normal2}(1)):

 $\max\{r_i\mid r_i>0\}\geq \max\{|r_i|\mid r_i<0\}$. 

In particular $\max\{|r_i|\mid r_i<0\}< n$.

We shall assign to such a cone $\delta'$ for which
$\sigma_v$ is not codefinite, the invariant  
$\inv(\delta'):=(\max\{r_i\mid r_i>0\},\max\{|r_i|\mid
r_i>0\}$), ordered lexicogragraphically.

By  Lemma \ref{le: normal2}(1)), after the star subdivision  at
$\Mid(\Ctr_+(\delta'),\delta')$ all new dependent cones
$\delta''=\delta'_{i_0}$  have normal relations of type
(1a) if $r_{i_0}> 0$ or (1b) if  $r_{i_0}< 0$.

If $r_{i_0}=\max\{r_i\mid r_i>0\}$ then in the relation
(1a) there is only one positive coefficient of $w$ which is
equal to $r_{i_0}$. All other coefficients:
$\frac{r_i-r_{i_0}}{m_w}$ and $r_i<0$ are nonpositive with
the absolute values $<n$. Thus $\delta'_{i_0}$ is of the
type $I(n,*)$ or $I(*,*)$. The face $\sigma_v$ is
codefinite with respect to $\delta'_{i_0}$.

If $0<r_{i_0}<\max\{r_i\mid r_i>0\}$ then 
the absolute values of
the coefficients $\frac{r_i-r_{i_0}}{m_w}$, $r_{i_0}$ and
$r_i\leq 0$, are $<n$. The cone $\delta'_{i_0}$ is of type
$I(*,*)$, $II(*,*)$, $III(*,*)$ and 
$\inv(\delta'_{i_0})<\inv(\delta')$.

If $r_{i_0}< 0$ then the absolute values of
all coefficients  in the relation
(1b) are equal to $|r_{i_0}|<n$ or $|\frac{r_{i_0}}{m_w}|<n$ and
$\delta'_{i_0}$  is a dependent cone of type $I(*,*)$  for which
$\sigma_v$ is not codefinite.

We repeat the procedure for all cones $\delta''$ for which
$\sigma_v$ is not codefinite.
The procedure
terminates since for all new dependent cones for which $\sigma_v$ is
not codefinite the invarinat $\inv(\delta')$ drops.

{\bf Step 2b}. {\it Eliminating $\delta$}. 
Apply the star subdivision at $\langle v\rangle$ to the
resulting cobordism.  We eliminate all cones in
$\Star(\sigma_v)$. In particular $\delta\in\Star(\sigma_v)$ will be eliminated.      
 By Step 2(a) $\sigma_v$
is codefinite with respect to all dependent cones in
$\Star(\sigma_v)$. All such cones $\delta'$ can be of the
type $I(n,n)$, $I(*,n)$ $I(n,*)$, $I(*,*)$, $II(n,*)$, $II(*,*)$, $III(n,n)$, $III(*,*)$.  

By Lemma \ref{le: normal2}(2) all new
dependent cones after the star subdivision 
have  normal relations of the form (2a)  and (2b). 

Since all the coefficients in normal relations of the form
(2b) are $<n$, the corresponding dependent cones are of
the form $I(*,*)$, $II(*,*)$, $III(*,*)$.

In the normal relation of the form (2a) for $\delta'_{i_0}$, the absolute 
values 
of the coefficients of $w_i$, for $i\neq i_1$, are $<n$. The
coefficient $r_{i_0}$ of $w$ is the only one which can be
equal to $n$. Therefore all new dependent cones we create
are of the type $I(n,*)$, $I(*,*)$, $II(n,*)$, $II(*,*)$, $III(*,*)$.

\bigskip
{\bf Step 3}. {\it Eliminating all cones of type  $I(*,n)$ and $II(n,*)$}.

Let $\delta=\langle v_1,\ldots,v_{k}\rangle$ be a cone of
the type $I(*,n)$ or $II(n,*)$.

We may assume that  the  positive coefficients
$r_{i}$ in the normal relation are  
$\leq n$ and the absolute values
of all negative coefficients are $< n$.  Then $\delta^+$
contains at least two positive rays.

By Lemma \ref{le: normal2}(1) after the star subdivision at
$\Mid(\Ctr_+(\delta),\delta)$ we create new dependent
$n$-cones $\delta_{i_0}$ only for $r_{i_0}=n$ in case
(1a). Then by the same reasoning as in Step1, 
$\delta_{i_0}$ is of type $I(n,*)$.

\bigskip
{\bf Step 4}. {\it Eliminating all cones of type  $III(n,n)$}.

Let $\delta$ be a cone of type $III(n,n)$. 
Then $\pi(\delta^+)=\pi(\delta^-)$.
Find 
$w\in \para(\pi(\delta^+))$ and let $v:=\Mid (w, \delta^+)$.
 Then $v$ is in the relative interior  of
  a face $\sigma_v$ of $\delta^+$.

{\bf Step 4a}. {\it Making $\sigma_v$ codefinite with respect to all
dependent cones in $\Star(\sigma_v,\Sigma)$}.

Let $\delta'\in\Star(\sigma_v,\Sigma)$ be any dependent cone 
for which $\sigma_v$ is not codefinite. As in Step 2,
$\delta'$ can be only of type $I$ or $II$. Then $\delta'$ is
of type $I(n,*)$, $I(*,*)$ or $II(*,*)$. 

Suppose $\delta'$ is of type  $I(n,*)$. 
Without loss of generality we may assume that in the normal
relation for $\delta'$
there is only one positive coefficient $r_{i_1}$ equal to $n$  and
$|r_{i}|<n$ for  all  negative coefficients $r_{i}$.
 
Since $\sigma_v$ is not codefinite
with respect to $\delta'$,  we have $v_{i_1}\in \sigma_v$. 
Note that $\Ctr_+(\delta')=\langle v_{i_1}\rangle$ and  
for $v:=\Mid(\Ctr_+(\delta'),\delta')$ we have $w=w_{i_1}$.
After the star subdivision at 
$\langle v'\rangle$  by Lemma \ref{le: normal2}(1)
we create one dependent cone $\delta'_{i_1}$ 
 with the 
normal relation obtained from the normal relation for
$\delta'$ by replacing $v$ with
$v_{i_1}$, and some dependent cones $\delta'_{i_0}$, for
$i_0\neq i_1$, which are of type
 $III(*,*)$ with the normal relation
$r_{i_0}w=r_{i_0}w_{i_1}$. 
After the star subdivision $\sigma_v$ is not
codefinite with respect to any new dependent cones.

Suppose $\delta'$ is of type  $I(*,*)$ 
or $II(*,*)$. Without loss of generality
we may 
assume that in
the normal relation (0) for $\delta'$ (see Lemma \ref{le: normal2}(1)):

 $\max\{r_i\mid r_i>0\}\geq \max\{|r_i|\mid r_i<0\}$. 

In particular $\max\{|r_i|\mid r_i<0\}< n$.

We shall assign to such a cone $\delta'$ for which
$\sigma_v$ is not codefinite, the invariant  
$\inv(\delta'):=(\max\{r_i\mid r_i>0\},\max\{|r_i|\mid
r_i>0\}$), ordered lexicogragraphically.

By  Lemma \ref{le: normal2}(1)), after the star subdivision  at
$\Mid(\Ctr_+(\delta'),\delta')$ all new dependent cones
$\delta''=\delta'_{i_0}$  have normal relations of type
(1a) if $r_{i_0}> 0$ or (1b) if  $r_{i_0}< 0$.

If $r_{i_0}=\max\{r_i\mid r_i>0\}$ then in the relation
(1a) there is only one positive coefficient of $w$ which is
equal to $r_{i_0}<n$. All other coefficients:
$\frac{r_i-r_{i_0}}{m_w}$ and $r_i<0$ are nonpositive with
the absolute values $<n$. Thus $\delta'_{i_0}$ is of the
type  $I(*,*)$. The face $\sigma_v$ is
codefinite with respect to $\delta'_{i_0}$.

If $0<r_{i_0}<\max\{r_i\mid r_i>0\}$ then 
the absolute values of
the coefficients $\frac{r_i-r_{i_0}}{m_w}$, $r_{i_0}$ and
$r_i\leq 0$, are $<n$. The cone $\delta'_{i_0}$ is of type
$I(*,*)$, $II(*,*)$, $III(*,*)$ and 
$\inv(\delta'_{i_0})<\inv(\delta')$.

If $r_{i_0}< 0$ then the absolute values of
all coefficients  in the relation
(1b) are equal to $|r_{i_0}|<n$ or $|\frac{r_{i_0}}{m_w}|<n$ and
$\delta'_{i_0}$  is a dependent cone of type $I(*,*)$  for which
$\sigma_v$ is not codefinite.

We repeat the procedure for all cones $\delta''$ for which
$\sigma_v$ is not codefinite.
The procedure
terminates since for all new dependent cones for which $\sigma_v$ is
not codefinite the invarinat $\inv(\delta')$ drops.

{\bf Step 4b}. {\it Eliminating $\delta$}. 
Apply the star subdivision at $\langle v\rangle$ to the
resulting cobordism.  We eliminate all cones in
$\Star(\sigma_v)$. In particular $\delta\in\Star(\sigma_v)$ will be eliminated.      
 By Step 4(a) $\sigma_v$
is codefinite with respect to all dependent cones in
$\Star(\sigma_v)$. All such cones $\delta'$ can be of the
type $I(n,*)$, $I(*,*)$, $II(*,*)$, $III(n,n)$, $III(*,*)$.  

By Lemma \ref{le: normal2}(2) all new
dependent cones after the star subdivision 
have  normal relations of the form (2a)  and (2b). 

Since all the coefficients in normal relations of the form
(2b) are $<n$, the corresponding dependent cones we create are of
the form $I(*,*)$, $II(*,*)$, $III(*,*)$.

If $\delta'$ is of type $I(*,*)$, $II(*,*)$, $III(*,*)$ then
all the coefficients in normal relations of the form
(2a) are $<n$ and the corresponding dependent cones we
create are of
the form $I(*,*)$, $II(*,*)$, $III(*,*)$.

If $\delta'$ is of type $I(n,*)$ then we can assume that
$\sigma_v$ is a face of ${\delta'}^+$. If ${\delta'}$ is
pointing down then all positive
coefficients $r_i$ in the normal relation for $\delta'$ are
$<n$. In the normal relation (2a) for $\delta'_{i_0}$,
where $r_{i_0}>0$, all coefficients are $<n$. If ${\delta'}$ is
pointing up then there is only one  positive
coefficient $r_{i_1}$. In the normal relation (2a) for $\delta'_{i_1}$,
the first two sums disappear and there is only one 
positive coefficient of $w$ equal to $r_{i_0}\leq n$. By
applying the star subdivision at $\langle v \rangle$ to
$\delta'$ we create one maximal dependnent cone of type
$I(n,*)$. At the same time we eliminate the cone $\delta'$ of type $I(n,*)$.

If $\delta'$ is of type $III(n,n)$ (for instance if
$\delta'=\delta$) then  there is one positive coefficient $r_{i_1}=n$
and one negative coefficient $r_{i_2}=-n$ in the normal
relation for $\delta'$. 

After the star subdivision we create only one maximal 
dependent cone $\delta'_{i_0}$ for 
$i_0=i_1$ having the normal relation of the form (2a). 
In the  relation
(2a) the first two sums disappear and there is only one 
positive coefficient of $w$ equal to $r_{i_0}\leq n$. The
absolute values of
the other negative coefficients are $< n$.
Thus we create one maximal dependent  
cone  of type $I(n,*)$.

\bigskip
{\bf Step 5}. {\it Eliminating all dependent cones of type $I(n,*)$
and all n-cones which are not faces of dependent cones}.

Let $\delta$ be a cone of type $I(n,*)$ or 
an $n$-cone which is not a face of a dependent cone.  

In the first case without loss
of generality we may assume that 
 in the normal relation there is one negative
coefficient equal to $-n$ and all positive coefficients are
$<n$ ($\delta$ is pointing down). 

Then $\pi(\delta^+)=\pi(\delta)$.
Find $w\in \para(\pi(\delta^+))$   and let $v=\Mid
(w,(\delta_+))$.
Then $v\in \inte (\sigma_v)$, where 
$\sigma_v$ is a face of $\delta_+$.

{\bf Step 5a}. {\it Making $\sigma_v$ codefinite with respect to all
dependent cones in $\Star(\sigma_v,\Sigma)$}.
We repeat  word by word the procedure in Step 4a.

{\bf Step 5b}. {\it Eliminating $\delta$}. Apply the star subdivision at $\langle v\rangle$ to the
resulting cobordism. We eliminate all cones in
$\Star(\sigma_v)$. 
 By Step 5(a) $\sigma_v$
is codefinite with respect to all dependent cones in
$\Star(\sigma_v)$. All such cones $\delta'$ can be of the
type $I(n,*)$, $I(*,*)$, $II(*,*)$, $III(*,*)$. We repeat
word by word the procedure in Step 5a excluding the last case $III(n,n)$.  
By applying the star subdivision at $\langle v\rangle$ we
decrease the number of maximal cones of type $I(n,*)$ (or
independent $n$-cones).  All new dependent
cones we create are of type $I(n,*)$, $I(*,*)$, $II(*,*)$.

\qed

\section{Birational cobordisms}
\subsection{Definition of a birational cobordism}

\noindent \begin{definition} (\cW): 
Let $X_1 $ and $X_2$ be two birationally equivalent
 normal varieties.
By a {\it birational cobordism} or simply a {\it cobordism} 
$B:=B(X_1,X_2)$ between them we understand a
normal variety $B$ with an algebraic action of $K^*$ such that the
sets
\[\begin{array}{rccc}
&B_-:=\{x \in B\mid  \, \lim_{{\bf t}\to 0} \, {\bf t}x \,\, \mbox{does
not exist}\}& \, \mbox{and}& \\
&B_+:=\{x \in B\mid \, \lim_{{\bf t}\to\infty} \, {\bf t}x \,\, 
\mbox{does not exist}\} && 
\end{array}\]
\noindent are nonempty and open and there exist 
 geometric quotients
$B_-/ K^*$ and $B_+\// K^*$ such that $B_+/ K^*\simeq X_1$ and $B_-\//
K^*\simeq X_2$ and the birational map $X_1\mathrel{{-}\,{\rightarrow}}
X_2$ is given by the above
isomorphisms and the open embeddings 
$B_+ \cap B_-\// K^*$ into
$B_+\//K^*$ and $B_-\//K^*$ respectively.
\end{definition}

\bigskip
\subsection{Collapsibility}
Let $X$ be a variety with an action of $K^*$. Let
$F\subset X$ be a set consisting of fixed points. Then we define  
$$F^+(X)=F^+=\{x\in X\mid 
\, \lim_{{\bf t}\to 0} {\bf t}x \in F\},\,  F^-(X)=F^-=\{x\in X\mid 
\, \lim_{{\bf t}\to \infty} {\bf t}x \in F\}.$$

\begin{definition} (\cW).\label{de: order} Let $X$ be a
cobordism or any variety with a $K^*$-action.
\begin{enumerate}
\item We say that a connected component $F$ of the fixed point
set  is an {\it  immediate predecessor} of a component
$F'$ if there exists a nonfixed point $x$ such that
$\lim_{{\bf t}\to 0} {\bf t}x\in F$ and
$\lim_{{\bf t}\to
\infty} {\bf t}x\in F'$. 
\item
We say that $F$  {\it precedes}  $F'$ and
write $F<F'$ if there exists a sequence of connected fixed
point set components $F_0=F ,F_1,\ldots,F_l=F'$ such that
$F_{i-1}$ is an immediate predecessor of $F_i$ 
(see \cite{BB-S} ,  Def. 1. 1).
\item
We call a cobordism (or a variety with $K^*$-action) {\it collapsible} (see
also Morelli \cite{Morelli1}) if the relation $<$ on its set of
connected components of the fixed point set  
is an order. (Here an
order is just required to be transitive.)
\end{enumerate}
\end{definition}

\begin{lemma}(\cW). A projective  cobordism is
collapsible. \qed 
\end{lemma}

\begin{lemma}\label{le: collapsible} Let $B$ be a collapsible
variety and $I$ be a $K^*$-invariant sheaf of ideals.
Then the blow-up $B':=\bl_I(B)$ is also collapsible.
\end{lemma}

\noindent {\bf Proof.} Set $\pi:B'\to B$. 
For any connected fixed point component $F'$ on $B'$ 
by abuse of notation denote by $\pi(F')$ the fixed point
component containing the image $\pi(F')$. Then for any
$F'_0$ and $F'_1$ on $B'$ the relation $F'_0<F'_1$ implies
that either $\pi(F'_0)<\pi(F'_1)$ or $\pi(F'_0)=\pi(F'_1)$.
Consider the latter case. Let $x\in F'$ and $u_1,\ldots,u_k$
be semiinvariant generators. Then at least one of the
functions $\pi^*(u_{i})$, say $\pi^*(u_{i_0})$, 
generates $\pi^{-1}(I)\cdot{\cO}_{B'}$ at $x$. We let ${\rm
ord}(F)$ denote  the weight of 
$\pi^*(u_{i_0})$ with respect to the $K^*$-action. This definition does not depend on the
choice of the semiinvariant generator. For any two generators
$\pi^*(u_{i_0})$ and $\pi^*(u_{i_1})$ their quotient is a
semiinvariant function which is invertible at the fixed point
$x\in B'$, which implies that it is invariant function. This
definition is locally constant, which means that it does not
depend on the choice of $x\in X$. Moreover the functions
$\pi^{*}(u_1),\ldots,\pi^*(u_k)$ are sections generating
$\cO_{\pi^{-1}(U)}(-D)$, where $D$ is the exceptional
divisor of $\pi$. This gives a $K^*$-equivariant morphism 
$\phi:\pi^{-1}(U)\to {\bf P}^k$
into a projective
space  with semiinvariant
coordinates $x_1,\ldots x_k$ having the same weights as $u_1,\ldots,u_k$, say
$a_1,\ldots,a_k$. We can assume that $a_1\leq\ldots\leq a_k$. Then  
the function ${\rm ord}(F)=\min\{a_i \mid x_i(\phi(F))\neq 0\}$ coincides with the one
 induced on ${\bf P}^k$, and  $F < F'$ implies
${\rm ord}(F)<{\rm ord}(F')$. Define $F\prec F'$ if
either $\pi(F)<\pi(F')$, or $\pi(F)=\pi(F')$ and 
${\rm ord}(F)\leq{\rm ord}(F')$. Since $F<F'$ implies
$F\prec F'$ it follows that $\prec$ is an order on the fixed
point components on $B'$. \qed

\subsection{Decomposition of a birational cobordism}

\begin{definition} (\cW). Let $B$ be a collapsible 
cobordism and $F_0$ be a minimal
component. By an {\it elementary collapse with respect to $F_0$}
we mean the cobordism $B^{F_0}:=B\setminus F_0^-$ (see
section \ref{se: pi-stable}).
By an {\it elementary cobordism with respect to
$F_0$} we mean the cobordism $B_{F_0}:=B\setminus \bigcup_
{F\neq F_0} F^+$.
\end{definition}

\begin{proposition} (\cW). Let $F_0$ be a minimal component
 of the fixed point
set in a collapsible cobordism $B$. Then the
elementary collapse $B^{F_0}$ with respect to $F_0$ is again a
collapsible cobordism, in particular it satisfies:
\begin{enumerate}
\item $B^{F_0}_+=B_+$ is an open subset of $B$,

\item $F_0^-$ is a closed subset of $B$ and equivalently $B^{F_0}$
is an
 open subset of $B$.

\item $B^{F_0}_-$ is an open subset of $B^{F_0}$  and
$B^{F_0}_-=B^{F_0}\setminus \bigcup_{F\neq F_0} F^+$.

\item The elementary cobordism $B_{F_0}$ is  an open subset
of $B$ such that
$$\displaylines{ (B_{F_0})_-=B_{F_0} \setminus F_0^+=B_-\cr
(B_{F_0})_+=B_{F_0} \setminus F_0^- =B^{F_0}_-.\cr} $$

\item There exist good and respectively geometric quotients  $B_{F_0}//{K^*}$
and $B^{F_0}_-/{K^*}$ and moreover the natural embeddings $ i_- : B_-\subset
B_{F_0}$ and $i_+ : B^{F_0}_-\subset B_{F_0}$ induce  proper morphisms 
 $ i_{-/K^*} : B^{F_0}_-/{K^*}\rightarrow B_{F_0}//{K^*}$
 and  $ i_{+/K^*} : B_-/{K^*}\rightarrow B_{F_0}//{K^*}$.  \qed  
\end{enumerate}
\end{proposition}
We can extend the order $<$ from Definition \ref{de: order}
to a total order. 

As a corollary from the above we obtain:
\begin{proposition} \label{pr: factorization} Let $F_1<\ldots<F_k$ denote
the connected fixed point components of a cobordism $B$.

Set $B_i:=B\setminus (\bigcup_{F_j<F_i}F_j^-\cup \bigcup_{F_j>F_i}F_j^+)$.

 Then \begin{enumerate}
\item
$(B_1)_-=B_-$, $(B_k)_+=B_+$. 

\item $(B_{i+1})_-=(B_i)_+$.

\item There is a factorization
$$B_-/K^*=(B_1)_-/K^*{{-}\to} (B_1)_+/K^*= (B_2)_-/K^*{{-}\to}
\ldots (B_{k-1})_+/K^*= (B_k)_-/K^* {{-}\to}
(B_k)_+/K^*=B_+/K^*.$$ \qed
\end{enumerate} 
\end{proposition}

\subsection{Existence of a smooth birational cobordism}\label{se: construction}
\begin{proposition} (see also \cW)\label{pr: construction} Let $\phi: Y\to X$ be a birational projective
morphism between smooth complete varieties over an algebraically
closed field of characteristic $0$. Let $U\subset X, Y$ be
an open subset , where $\phi$ is an isomorphism. 
\begin{enumerate} 
\item There exists a smooth
collapsible cobordism $B=B(X,Y)$ between $X$ and $Y$ such that
 $B\supset U\times K^*$

\item If $D_X:=X\setminus U$ and $D_Y:=Y\setminus U$ are divisors
with normal crossings then there exists 
a cobordism $\tilde{B}$ between $\tilde{X}$ and
$\tilde{Y}$ such that 
\begin{itemize}
\item $\tilde{X}$ and $\tilde{Y}$ are obtained from $X$ and $Y$ by a sequence of
blow-ups at centers which have normal crossings with
components of the total transforms of $D_X$ and $D_Y$ respectively.
\item $U\times K^* \subset \tilde{B} $ and
$\tilde{B}\setminus (U\times K^*)$ is a
divisor with normal crossings.
\end{itemize}
\end{enumerate} 
\end{proposition}

\noindent {\bf Proof.} (1) If $X$ and $Y$ are projective the
construction of $B$ is given in Proposition 2 of \cW; in
general we use the construction of D.Abramovich (Remark after
Proposition 2 in \cW). $B$ is an open subset of $\overline{B}$
where $\overline{B}$ is
obtained from $X\times {\bf P}^1$ by a sequence of blow-ups at
ideals, hence by a single blow-up of an ideal. By Lemma
\ref{le: collapsible}, $\overline{B}$ is collapsible. Consequently
$B$ is a  collapsible cobordism. 

(2) The sets $B_-$ and $B_+$ are locally
trivial $K^*$-bundles with projections $\pi_-:B_-\to B_-/K^*\simeq
X$ and $\pi_+:B_+\to B_+/K^*\simeq
Y$. Let $Z:=B\setminus (U\times K^*)$. Then $Z\cap B_-$ and
$Z\cap B_+$ are normal crossing divisors and $\pi_-(Z\cap
B_-)=D_X$ and $\pi_+(Z\cap
B_+)=D_Y$.
Let $f:\tilde{B}\to B$ be a canonical principalization of
$I_Z$ (see Hironaka \cite{Hironaka1},
Villamayor \cite{Villamayor} and Bierstone-Milman \cite{Bierstone-Milman}). Let $f_+:f^{-1}(B_+)\to B_+$
(resp.$f_-:f^{-1}(B_-)\to B_-$) be 
its restriction. By functoriality  
$f_+$ (resp. $f_-$)  is  a canonical
principalization of $B_+$ (resp. $B_-$) which commutes
with a canonical principalization $\tilde{Y}$ of $I_{D_Y}$ on $Y$ 
(resp. $\tilde{X}$ of $I_{D_X}$ on $D_X$)
In particular  all centers are
$K^*$-invariant of the form $\pi_+^{-1}(C)$ (resp.
$\pi_-^{-1}(C)$), where $C$ have
normal crossings with components of the total transform of $D_Y$ (resp. $D_X$).
\qed

\section{Weak factorization theorem}
\subsection{Construction of a $\pi$-regular cobordism}

\begin{definition} \begin{enumerate} 
\item We say that a cobordism $B$ is 
{\it $\pi$-regular} if for any $x$ in $B$ which is not a fixed point of
the $K^*$-action in $B$ there exists an
affine invariant
neighborhood $U$ without fixed points such that $U/K^*$ is
smooth.
\item We call a cobordism $B$  {\it $K^*$-stratified} if it
is a $K^*$-stratified toroidal variety.  
\end{enumerate}
\end{definition}

\begin{proposition} \label{pr: pi-lemma2} Let $(B,S)$ be a smooth stratified
cobordism
between smooth varieties $X$ and $X'$.  There
exists a $K^*$-toroidal modification $\tilde{B}$ of $B$ obtained
as a sequence of blow-ups of stable valuations such that
$\tilde{B}$ is  a $\pi$-regular
$K^*$-stratified cobordism between $X$ and $X'$. Moreover if $U\subset X$
is an open affine invariant fixed point free subset such
that $U/K^*$ is smooth then all centers of blow-ups are
disjoint from $U$. 
\end{proposition}

\noindent {\bf Proof.} By Lemma \ref{le: asemicomplex} we can associate with $(B,S)$ a
$K^*$-semicomplex $\Sigma$. Since for any fixed point component $F$
in $B$ the sets $F^+$ and $F^-$ are not dense it follows
that $\Sigma_S$ is strictly convex.  By Lemma \ref{le: pi-lemma} 
there is a $K^*$-canonical subdivision $\Delta$ of
$\Sigma$ which is $\pi$-nonsingular and does not affect
any $\pi$-nonsingular cones in $\Sigma$. By
Theorem \ref{th: correspondence}, $\Delta$ corresponds
to a $K^*$-stratified cobordism between $X$ and $Y$. \qed

\subsection{Factorization determined by a $\pi$-regular cobordism}

\begin{proposition} \label{pr: regular factorization}
 Let $(B,R)$ be  a $\pi$-regular collapsible 
 $K^*$-stratified cobordism between smooth stratified toroidal varieties
$(X,S_X)=(B_-/K^*,(R\cap B_-)/K^*)$, $(Y,S_Y)=(B_+/K^*,(R\cap B_+)/K^*)$.
Let $B_i$ be elementary cobordisms as in Proposition
\ref{pr: factorization}.
Then there exists a sequence of smooth stratified toroidal varieties
$(X_i,S_i):=((B_i)_-/K^*, (R\cap (B_i)_-)/K^*)$
 and  a factorization
\[\begin{array}{rcccccccccccccccc}
&&Y_0&&&&Y_1&&&&&&&Y_{n-1}&&&\\
&g_0\swarrow&& \searrow f_0&& g_1\swarrow&& \searrow
f_1&&&&&g_{n-1}\swarrow&& \searrow f_{n-1}&&\\
X=X_0&& {{-}{\to}}&& X_1&&
  {{-}{\to}}&&X_2&&\ldots& X_{n-1}&&
{{-}{\to}}& X_n&=Y&,
\end{array}\]

\noindent where each $X_i$ is a smooth variety and $f_i:Y_i\to X_{i+1}$
and $g_i:Y_i\to X_i$ for $i=0,\ldots,n-1$ are blow-ups at
smooth centers which have normal crossings with the closures of
strata in $S_{i}$ (resp. $S_{i+1}$). 

\end{proposition}
\noindent{\bf Proof.} Follows from the following Lemma 
\ref{le: pi-lemma2} and Proposition \ref{pr: factorization}.
\qed

\begin{lemma} \label{le: pi-lemma2} Let $(B,R)$ be a
$\pi$-regular elementary
$K^*$-stratified cobordism. Then $X=B_-/K^*$ and $Y=B_+/K^*$ are smooth
stratified toroidal varieties with stratifications $S_X$ and
$S_Y$ induced
by $R\cap B_+$ and $R\cap B_-$ and there exists a smooth
variety $Z$, obtained by a blow-up  
at a smooth center $C_X$ (resp. $C_Y$), 
having only normal crossings with closures of strata in
$S_X$ (resp. $S_X$).
\end{lemma}

\noindent{\bf Proof.} 
For any $X$ let $\widetilde{X}$ denote its normalization.

 The fixed 
point component is equal to the closure
$\overline{\strat_B(\sigma)}$ of the maximal stratum
, corresponding to a circuit
$\sigma$.

Consider the 
commutative diagram 

\[\begin{array}{rccccc}
&B_+/K^* &&&&B_-/K^*\\
&&\searrow&& \swarrow&\\  
&&&B//K^*&&\\
\end{array}\]

This diagram can be completed to
\[\begin{array}{rccccccc}
&&&Z:=({B_+/K^* 
\times_{B//K^*} B_-/K^*})\widetilde{}&&\\
&&&\psi\downarrow&&&&\\
&&&B_+/K^* \times_{B//K^*} B_-/K^*&&\\
&&\swarrow&& \searrow&&&\\
&B_+/K^* &&&&B_-/K^*&& (*)\\
&&\searrow&& \swarrow&&&\\  
&&&B//K^*&&&&\\

\end{array}\]
\noindent  where the morphism $\psi$ is given by the normalization.

For any $x\in \overline{\strat(\sigma)}$ find $\tau\geq\sigma$ such that
$x\in\strat(\tau)$. 
Let $\phi_{\tau}: U_{\tau}\to X_{\tau}$ be the
relevant chart. 

The above diagram defines locally a diagram
\[\begin{array}{rccccccc}
&&&((U_{\tau})_+/K^*\times_{U_{\tau}/K^*}
(U_{\tau})-/K^*)\widetilde{} &&&&\\
&&&\downarrow&&&&\\
&&&(U_{\tau})_+/K^*\times_{U_{\tau}/K^*}{(U_{\tau})_-/K^*}&&&&\\
&&\swarrow&& \searrow&&&\\
&(U_{\tau})_+/K^*&&&&(U_{\tau})_-/K^*&&\\
&&\searrow&& \swarrow&&&  \\
&&&U_{\tau}//K^*&&&&\\
\end{array}\]
This diagram is a pull-back via a smooth morphism
$\phi_{\tau/K^*}:U_{\tau}//K^*\to X_{\tau}/K^*$
of the  diagram of toric varieties 

\[\begin{array}{rccccc}
&&&X_{\Sigma} :=((X_{\tau})_-/K^*
\times_{X_{\tau}/K^*}(X_{\tau})_+/K^*)\widetilde{}
 &&\\
&&&\downarrow&&\\
&&&(X_{\tau})_-/K^*\times_{X_{\tau}/K^*} (X_{\tau})_+/K^*&&\\
&&\swarrow&& \searrow&\\
&X_{\Sigma_1}:=(X_{\tau})_-/K^*&&&&X_{\Sigma_2}:=(X_{\tau})_+/K^*\\
 &&\searrow&& \swarrow&\\  
&&&X_{\tau^\Gamma}=X_{\tau}//K^*&&\\
\end{array}\]

\noindent It follows from the universal property of the  fiber product that 
$X_\Sigma$ is a
normal toric variety whose fan consists of the cones $
\{\tau_1\cap\tau_2 \mid   \tau_1\in \Sigma_1, \tau_2\in \Sigma_2\}.$

The cone $\tau$ is $\pi$-nonsingular and by Lemma \ref{le: normal} we can write the
only dependence relation  of
$\tau^\Gamma=\langle v_1,\ldots,v_k,w_1,\ldots,w_m,q_1,\ldots,q_l \rangle$, where
$\langle v_1,\ldots,v_k,w_1,\ldots,w_m \rangle$  is a curcuit,  
as follows:
$$e:=v_1+\ldots+v_k=w_1+\ldots+w_m.$$ 
Finally, $\Sigma=\langle e \rangle\cdot\Sigma_1$ and 
$\Sigma=\langle e\rangle \cdot \Sigma_2$
are regular star subdivisions.

We have shown that the morphisms $Z\to X$ and 
$Z\to Y$ defined by the diagram (*) are
blow-ups at smooth centers. These centers have normal
crossings with closures of strata since locally the centers and
the closures of strata are determined by the closures of orbits on a
smooth toric  variety.
 \qed

\subsection{Nagata's factorization}
\begin{lemma}(Nagata)\label{le: Nagata} 
Let $X\supset U \subset Y$ be two complete varieties over
an algebraically closed field, containing an open subset $U$.
Then there exists  a variety $Z$ which is  simultaneously a
 blow-up of a coherent sheaf of ideals $\cI_X$ on $X$ and 
a blow-up of a coherent sheaf of ideals $\cI_Y$ on $Y$.
Moreover the supports of $\cI_X$ and $\cI_Y$ are disjoint from $U$
\end{lemma}

 The above lemma is a refinement of the
original lemma of Nagata
\begin{definition}(Nagata \cite{Nagata}) \label{de: Na}
Let $X_i\supset U $, for $i=1,\ldots,n$ be  complete varieties over
an algebraically closed field, containing an open subset $U$.
By the {\it join} $X_1*\ldots *X_k$ we mean   
the closure of the diagonal $\Delta(U) \subset
X_1\times\ldots\times X_k$.
\end{definition}

\begin{lemma}(Nagata \cite{Nagata}) \label{le: Na}
Let $X\supset U \subset Y$ be two complete varieties over
an algebraically closed field, containing an open subset $U$.
Then there exist coherent sheaves of ideals $\cI_1,\ldots,\cI_k$ on
$X$ with supports disjoint from $U$ and blow-ups 
$X_1=\bl_{\cI_1}X,\ldots,X_k=\bl_{\cI_k}X$ 
such that the  join 
$X_1*\ldots *X_k$
 dominates $Y$.\qed
\end{lemma}

\noindent {\bf Proof of Lemma \ref{le: Nagata}}
 Let $X_i$ be as in the assertion of Lemma \ref{le: Na}. 
Note that $X_1\times\ldots\times X_k$ is projective over
$X\times\ldots \times X$. This implies that $X_1*\ldots *X_k$ is
projective over $\Delta(X)\simeq X$. Thus there exists a
sheaf of ideals $\cJ$ on $X$ such that $X_1*\ldots *X_k\to X$ is
the blow-up at $\cJ$.
 By Nagata and the above find  a blow-up $X'$ of $X$ at $\cJ$
such that $X'$ dominates $Y$. Denote the blow-up by $\phi:X'\to Y$. 
 Analogously find  a blow-up $Z$ of $Y$ at $\cI_Y$
such that $Z$ dominates $X'$. 
 Thus $Z\to X'$ is a blow-up at $\phi^{-1}(\cI_Y)$ and $Z\to
X$ is a composite of blow-ups and consequently a single blow-up 
at some sheaf of ideals $\cI_X$. \qed

\subsection{Proof of the Weak Factorization Theorem}\label{se: factorization}

\bigskip

\noindent {\bf Step 1} By Nagata's Lemma \ref{le: Nagata}
we find a $Z'$ obtained from $X$ and $Y$ respectively
by blow-ups at
centers disjoint from $U$. If $X$ and $Y$ are projective we
can take $Z'$ to be the graph of $\phi: X{-\to} Y$. 
Let $Z''$ be a
resolution of singularities of $Z'$ and $Z$ be a canonical
principalization of the ideal of the set $Z''\setminus U$
(this is needed for part (2) of the theorem only). Then
$Z\setminus U$  
 is a divisor with normal crossings. 
It suffices to prove the theorem for the projective
morphisms $Z\to X$ (or $Z\to Y$). 

\noindent {\bf Step 2} 
By Proposition \ref{pr: construction} there exists
a smooth collapsible cobordism $B=B(Z,X)$. If $D_X:=X\setminus U$ and $D_Z:=Z\setminus U$
are divisors with normal crossings then there exists a
cobordism $\tilde{B}$ between $\tilde{X}$ and $\tilde{Y}$
as in the proposition such that $\tilde{B}\setminus
(U\times K^*)$ is a normal crossings divisor. For
simplicity denote $\tilde{X}, \tilde{Z},\tilde{B}$ by $X
,Z, B$ respectively.

\noindent {\bf Step 3} 
By Lemma
\ref{le: existence2} there exists a $K^*$-invariant stratification $S$
on $B$  such that $(B,S)$ is a $K^*$-stratified
toroidal variety. If $B\setminus
(U\times K^*)$ is a normal crossings divisor
 then the
components of this divisor are the closures of some strata
in $S$.

\noindent {\bf Step 4} 
By Proposition 
\ref{pr: pi-lemma2} we can find a $\pi$-regular collapsible $K^*$-stratified
cobordism $(\overline{B},\overline{R})$ between $X$ and $Z$
by applying a sequence of blow-ups at stable valuations to $(B,S)$
and
using the combinatorial algorithm. The subsets
$B_+=\overline{B}_+$ and $B_-=\overline{B}_-$ are unaffected.

\noindent {\bf Step 5} By Proposition \ref{pr: regular
factorization} the $\pi$-regular
cobordism $(\overline{B},\overline{R})$ determines a
factorization of the map $Z\to X$ into a sequence of blow-ups and blow-downs.
All centers have only normal crossings with closures of the
induced strata on intermediate varieties.  
\qed

\section{Comparison with another proof. Torification
and stable support}\label{se: comparison}

\subsection{Comparison with another proof}

In \cite{AKMW} another proof of the Weak
Factorization Theorem is given , where the theorem is stated
an proved in a more general version. In particular the
assumption of the base field to be algebraically closed was
removed. The theorem was also formulated and proved for
bimeromorphic maps of complex compact analytic manifolds.  
The proof in \cite{AKMW}  uses the idea of torific ideals due to
Abramovich and De Jong.
 The idea of the proof is
to consider an ideal on a smooth cobordism whose blow-up
determines a structure of toroidal embedding compatible
with the $K^*$-action. Such an ideal, called {\it ''torific''},
and the the blow-up procedure , called {\it
torification} were introduced  by D.Abramovich and
J.de Jong. The torific ideal is defined in an invariant
neighborhood $U$ of
the fixed point component of the smooth variety with  $K^*$-action.
Let $u_1,\ldots,u_k$ denote the semiinvariant parameters with weights
$a_1,\ldots,a_k$. We define the torific ideal to be 
$I=I_{a_1}\cdot\ldots\cdot I_{a_k}$, where $I_a$ denotes the
ideal generated by all semiinvariant functions of weight $a$.
 The blowing-up of the torific ideal
induces a  structure of toroidal embedding, locally on each elementary piece
of cobordism. Torification allows one to pass to the category of toroidal
embeddings, where the factorization problem has already been
solved using combinatorial algorithms. However the torific
ideal can be constructed only locally on elementary
cobordisms. Hence after the torification is done on
each elementary cobordism some glueing procedure should be
applied. Using canonical principalizations of torific ideals
and canonical resolution of singularities of intermediate varieties
gives a factorization on an elementary cobordism between
canonical resolution of singularities.
Patching these factorizations together one
gets a decomposition of the given birational map. 

\subsection{Torification and stable support}
The notion of torification can be understood and  generalized
on the ground of the theory  of stratified toroidal varieties.

\noindent
\begin{definition} Let $(X,S)$ be an oriented stratified toroidal
variety (with or without group action). By a {\it stable
ideal} $\cI$ we mean a coherent sheaf of ideals $\cI$ satisfying the following
conditions :
\begin{enumerate}
\item The support of $\cI$ is a union of strata in $S$.

\item For any points $x_1$ and $x_2$ in a stratum $s$ and any
isomorphism $\phi:\widehat{X}_{x_1}\to \widehat{X}_{x_2}$
preserving strata and orientation (respectively equivariant
with respect to the $\Gamma_\sigma$-action),
$\phi^*(\cI)\cdot
\widehat{\cO}_{x_2}=\cI\cdot\widehat{\cO}_{x_1}$.
\end{enumerate}
\end{definition}

\begin{proposition} Let $\cI$ be a stable ideal on an oriented stratified toroidal
variety $(X,S)$.
Then for any $\sigma\in \Sigma$ there is a $G^\sigma$-invariant toric
ideal $I_\sigma$ on ${X}_{\sigma}$ such that 

\begin{enumerate}
 \item The induced ideal
on $\widetilde{X}_{\sigma}$ is $G^\sigma$-invariant. 

\item For any
chart $\phi_\sigma:U\to X_{\sigma}$, $\cI_U=\phi_\sigma^*(I_\sigma)$.

\item  For any $\tau\leq \sigma$, $I_\tau\cdot
\cO_{x,X_{\sigma}}$ is the restriction of $I_{\sigma}$
to the open subset $X_{(\tau, N_\sigma)}\subset X_{\sigma}$.

\end{enumerate}
\end{proposition}

One can rephrase the above proposition in combinatorial
language, using the correspondence between toric ideals $\cI$
and $\ord(\cI)$-functions in \cite{KKMS} and Lemma \ref{le: vertices}.

\begin{proposition} Let $(X,S)$ be an oriented stratified toroidal
variety with an associated oriented semicomplex $\Sigma$.
The stable complete coherent sheaves of ideals $\cI$ on $(X,S)$ are in $1-1$
correspondence with collections of functions $\ord(I_\sigma)$,
$\sigma\in \Sigma$, 
satisfying:
\begin{enumerate}
\item
$\ord(I_\sigma):\sigma\to \QQ$ is a convex piecewise
linear function $\ord(I_\sigma)({N}_\sigma)\to {\bf Z}$. 

\item Denote by $\Delta(I_\sigma)$ the subdivision
of $\sigma$ into the maximal cones where $\ord(I_\sigma)$ is
linear. Then $\Ver(\Delta(I_\sigma))\subset\stab(\sigma)\cup
\Ver(\overline{\sigma})$, and $\ord(I_\sigma)(v)=0$ for any $v\in
\Ver(\Delta(I_\sigma))\setminus \stab(\sigma)$. 

\item For any $\tau\leq \sigma$, the restriction of $\ord(I_{\sigma})$ to
$\tau$ is equal to $\ord(I_{\tau})$. \qed

\end{enumerate}
\end{proposition}
It follows that the functions $\ord(I_\sigma)$ patch together and
give a function $\ord(\cI)$ on the totality of vectors in
faces of $\Sigma$.
\begin{lemma}
Denote by $\Delta(\cI)$ the subdivision of $\Sigma$ obtained
by glueing the subdivisions $\Delta(I_\sigma)$. Then 
$\Delta(\cI)$ is a canonical subdivision of $\Sigma$
correspondinding to the blow-up of $(X,S)$ at $\cI$. \qed
 \end{lemma}

\begin{definition} Let $(X,S)$ be an oriented stratified toroidal
variety (with or without group action). Then an ideal
$\cI$ is called {\it torific} if it is stable and the
normalization of its  blow-up is a
canonical toroidal morphism $(Y,R)\to (X,S)$ such that $(Y,R)$ is a
 toroidal embedding. Such a blow-up is called a {\it
torification}. 
\end{definition}

\begin{lemma} 
Let $(X,S)$ be an oriented stratified toroidal variety with
an
associated semicomplex $\Sigma$. Then a stable ideal $\cI$ is
torific iff  any $\sigma \in \Delta(\cI)$ whose 
relative interior intersects the stable support $\Stab(\Sigma)$,
is contained in $\Stab(\Sigma)$. 
\end{lemma}
\noindent{\bf Proof.}  
If $\inte(\sigma)\cap \Stab(\Sigma) \neq\emptyset$ then
$\sigma\in \Delta(\cI)_{\stab}$. By Theorem 
\ref{th: correspondence} and Lemma \ref{le: associated
semicomplex}, $\Delta(\cI)_{\stab}$ is a complex corresponding
to the toroidal embedding. Hence all vectors from $\sigma$
are stable. \qed

Torification is not always possible.
However we can  perform torification locally in the following sense.

\begin{proposition}\label{pr: torification} Let $(X,S)$ be an
oriented stratified toroidal variety with an associated
oriented semicomplex $\Sigma$. 
Then for any $\sigma\in \Sigma$ there is a $G^\sigma$-invariant
ideal  $$I_\sigma:=(G^\sigma\cdot u_1)\cdot\ldots\cdot (G^\sigma\cdot
u_k),$$ where $u_1,\ldots,u_k$ are toric paramaters on $\widetilde{X}_{\sigma}$ 
such that $\Delta(I_\sigma)$ contains the face $\Inv(\sigma)$. 
In particular $$\Inv(\sigma)=
{\rm conv}((\Ver(\Delta(I_\sigma))\setminus
\Ver(\sigma)) \cup \Ver^\semic(\sigma)),$$
\noindent where $\Ver^\semic(\sigma)$ denotes the set of
$1$-dimensional faces in the semicone $\sigma$.
\end{proposition}

\noindent{\bf Proof.} Denote by  $v_\nu$
 the vector corresponding to a toric valuation $\nu$.
 Note that a toric valuation $\nu$ is $G^\sigma$-invariant iff
$\nu(g_*(u_i))=\nu(u_i))$ for any $i=1,\ldots,k$ and $g\in
G^\sigma$. It follows from Lemma \ref{le: monomial2} that the latter
statement is equivalent to the following one:
$\nu(g_*(u_i))\geq\nu(u_i)$
for any $g\in G^\sigma$. This can be reformulated as follows: $\nu$
is $G^\sigma$-invariant iff 
$\ord{(I_\sigma)}(v_\nu)= \nu(u_1\cdot\ldots\cdot u_k)$.
Finally, $\Inv(\sigma)$
is a set  where $\ord{(I_\sigma)}$ is linear and equal to
the functional of $u_1\cdot\ldots\cdot u_k$. Hence $\Inv(\sigma)$ occurs as a cone in
$\Delta(I_\sigma)$. All vertices of this cone belong to 
$\Inv(\sigma)$
and hence they have to be in the form as in the
proposition. \qed

\begin{remark}The above proposition generalizes the construction of
torific ideal due to Abramovich and de Jong (\cite{Abramovich-de-Jong}). 
\end{remark} 

\subsection{Computing the stable support}

Proposition \ref{pr: torification} provides a method of
computing the stable support of semicomplexes.

\begin{example}\label{ex: isolated4} Let $(X,S)$ be a toroidal variety with an isolated
singularity $x_1x_2=x_3x_4$ as in Example \ref{ex:
isolated}. The associated semicomplex $\Sigma$
consists of  the cone over a square and the apex of this cone. Then $G^\sigma\cdot
x_i=(x_1,x_2,x_3,x_4)$ is the ideal $m_p$ of the singular point $p$. Hence the
torific ideal $I_\sigma$, where $s=\{p\}$, is equal to $m_p^4$.
$\Delta(I_\sigma)$ is the star
subdivision at the ray $\varrho$ over the centre of the square. By
Proposition \ref{pr: torification}, $$\Stab(\Sigma)=
\Inv(\sigma) =\varrho.$$ 
\end{example}    

\begin{example}\label{ex: Hironaka2} Let $(X,S)$ be a
smooth stratified toroidal variety of dimension $3$
with the stratification determined by two curves $l_1$ and
$l_2$ meeting transversally at some point $p$. The
associated semicomplex $\Sigma$ consists of the regular
$3$-dimensional cone $\sigma=\langle w_1,w_2,w_3 \rangle $ and its two $2$-dimensional
faces $\tau_1=\langle w_1,w_3 \rangle $ and $\tau_2=
\langle w_2,w_3 \rangle $. Set $v_1:=w_1+w_3$,
$v_2:=w_2+w_3$, $v:=w_1+w_2+w_3$. 
Write $l_1: x_1=x_3=0$, $l_2: x_2=x_3=0$. Then $(G^\sigma\cdot
x_1)=I_{l_1}=(x_1,x_3)$, $(G^\sigma\cdot
x_2)=(x_2,x_3)$,
$(G^\sigma\cdot x_3)=I_{l_1\cup l_2}=(x_3,x_1x_2)$. Then
$\Delta(I_\sigma)$ is determined by two star subdivisions:
$\langle v_1 \rangle\cdot\sigma$, $\langle v_2 \rangle\cdot\sigma$,
and the subdivision $\Sigma_0$ of $\sigma$ determined by 
$\ord((x_3,x_1x_2)):=\min(w^*_3, w_1^*+w_2^*)$. The latter
consists of two cones separated by the $2$-dimensional
face $\langle v_1, v_2 \rangle$.
Then $\Inv(\sigma)=\langle v_1, v_2, v \rangle$.    
Also $\Inv(\tau_i)= \langle v_i \rangle$. Finally, 
$$\Stab(\Sigma)=\Inv(\sigma)=\langle v_1, v_2, v \rangle.$$
\end{example}

The above considerations can be generalized to any regular
semicomplex. We shall skip the details of the proof.  
\begin{proposition}
Let $\Sigma$ be a regular semicomplex. For any cone
$\sigma=\langle v_1,\ldots, v_{k_\sigma} \rangle$ set \\ $v_\sigma:=
v_1+\ldots+v_{k_\sigma}$.
Then $$\stab(\sigma)=\bigoplus_{\tau\leq \sigma}  {\bf Q}\cdot v_{\tau}.$$
\qed
\end{proposition}
This formula can be extended to any oriented simplicial semicomplex.
\begin{conjecture}
Let $\Sigma$ be an oriented simplicial semicomplex. For any 
$\sigma\in \Sigma$ let $V_\sigma$ denote the set of all the minimal
internal vectors of $\sigma$ and $W_\sigma:=\bigcup_{\tau\leq \sigma}
V_{\tau}$. Then
$$\stab(\sigma)= \bigoplus_{v\in W_\sigma}  {\bf Q}\cdot v.$$
\end{conjecture}

\end{document}